\documentclass[a4paper,12pt]{amsbook}

\usepackage{amsmath,amssymb,amsthm,graphicx,stmaryrd}
\usepackage{verbatim,enumitem}

\usepackage[refpage]{nomencl}

\makenomenclature

\newcommand{\BB}{\mathbb{B}} 
\newcommand{\CC}{\mathbb{C}}
\newcommand{\EE}{\mathbb{E}} 
\newcommand{\FF}{\mathbb{F}}

\newcommand{\HH}{\mathbb{H}}
\newcommand{\KK}{\mathbb{K}}
\newcommand{\LL}{\mathbb{L}}
\newcommand{\NN}{\mathbb{N}}
\newcommand{\PP}{\mathbb{P}} 
\newcommand{\RR}{\mathbb{R}} 
\newcommand{\RRp}{\RR_+}
\renewcommand{\SS}{\mathbb{S}}
\newcommand{\TT}{\mathbb{T}} 
\newcommand{\VV}{\mathbb{V}} 
\newcommand{\WW}{\mathbb{W}} 
\newcommand{\XX}{\mathbb{X}}
\newcommand{\YY}{\mathbb{Y}}
\newcommand{\ZZ}{\mathbb{Z}}

\newcommand{\cA}{\mathcal{A}}
\newcommand{\cB}{\mathcal{B}} 
\newcommand{\cI}{\mathcal{I}}
\newcommand{\cC}{\mathcal{C}}
\newcommand{\cD}{\mathcal{D}} 
\newcommand{\cF}{\mathcal{F}}
\newcommand{\cE}{\mathcal{E}}
\newcommand{\cG}{\mathcal{G}} 
\newcommand{\cK}{\mathcal{K}}
\newcommand{\cL}{\mathcal{L}} 
\newcommand{\cH}{\mathcal{H}}
\newcommand{\cN}{\mathcal{N}}
\newcommand{\cO}{\mathcal{O}}
\newcommand{\cT}{\mathcal{T}}
 
\newcommand{\cR}{\mathcal{R}} 
\newcommand{\cQ}{\mathcal{Q}}
\newcommand{\oD}{D}
\newcommand{\tc}{\mathrm{c}}
\newcommand{\oper}{\nu}
\newcommand{\wt}{w}
\newcommand{\Qab}{\cQ^{A,B}}
\newcommand{\Qabxy}{\cQ^{A\sd xy,B\sd xy}}

\newcommand{\fr}{\phi_{\rho}}
\newcommand{\frr}{\phi_{\rho'}}

\renewcommand{\a}{\alpha}
\newcommand{\D}{\Delta} 
\renewcommand{\d}{\delta} 
\newcommand{\G}{\Gamma}
\newcommand{\Gh}{\G}
\newcommand{\g}{\gamma} 
\renewcommand{\L}{\Lambda} 
\renewcommand{\i}{\iota} 
\newcommand{\LLam}{\mathbf{\Theta}} 
\renewcommand{\l}{\lambda}
\renewcommand{\b}{\beta} 
\renewcommand{\k}{\kappa} 
\newcommand{\Om}{\Omega}
\newcommand{\om}{\omega} 
\renewcommand{\S}{\Sigma} 
\newcommand{\s}{\sigma}
  
\newcommand{\boundary}{\hat\partial}
\newcommand{\eps}{\varepsilon}
\newcommand{\dual}{\mathrm{d}} 
\newcommand{\rs}{\mathrm{s}}
\renewcommand{\max}{\mathrm{max}}

\newcommand{\el}{\langle} 
\newcommand{\er}{\rangle}
\newcommand{\tr}{\mathrm{tr}}
\newcommand{\ev}{\mathrm{ev}}
\newcommand{\odd}{\mathrm{odd}}
\newcommand{\rf}{\mathrm{f}}
\newcommand{\sd}{\triangle}
\newcommand{\dlr}{{\sc dlr}}
\newcommand{\fk}{{\sc fk}}
\newcommand{\fkg}{{\sc fkg}}
\newcommand{\gks}{{\sc gks}}
\newcommand{\ghs}{{\sc ghs}}
\newcommand{\sle}{{\sc sle}}
\newcommand{\lra}{\leftrightarrow}
\newcommand{\nlra}{\nleftrightarrow}

\newcommand{\bc}{\rho_{\mathrm{c}}}
\newcommand{\bs}{\rho_{\mathrm{s}}}
\renewcommand{\b}{\beta}
\newcommand{\oo}{\infty}
\newcommand{\ghost}{\Gamma}
\newcommand{\qq}{\quad\quad}
\newcommand{\rc}{random-cluster}
\newcommand{\bra}[1]{\langle#1|}
\newcommand{\ket}[1]{|#1\rangle}
\newcommand{\Si}{\Sigma}
\newcommand{\sm}{\setminus}
\newcommand{\resp}{respectively}
\newcommand{\even}{\mathrm{even}}
\renewcommand{\odd}{\mathrm{odd}}
\renewcommand{\o}{\mathrm{o}}
\newcommand{\wtilde}{\widetilde}
\newcommand{\what}{\widehat}
\newcommand{\pd}{\partial}

\newcommand{\es}{\varnothing}
\newcommand{\se}{\subseteq}
\newcommand{\ul}{\underline}
\newcommand{\od}{d}
\newcommand{\ol}{\overline}
\newcommand{\lrao}[1]{\overset{#1}{\lra}}
\newcommand{\bigmid}{\,\big|\,}
\newcommand{\Bigmid}{\,\Big|\,}
\newcommand{\wired}{\mathrm{w}} 
\newcommand{\free}{\mathrm{f}} 
\newcommand{\swired}{\mathrm{sw}} 
\newcommand{\sfree}{\mathrm{sf}} 
\newcommand{\crit}{\mathrm{c}} 

\newcommand{\one}{\hbox{\rm 1\kern-.27em I}}

\newtheoremstyle{slthm}
     {}
     {\baselineskip}
     {\slshape}
     {\parindent}
     {\scshape}
     {.}
     { }
     {}

\theoremstyle{slthm}
\newtheorem{definition}{Definition}[chapter]
\newtheorem{theorem}[definition]{Theorem}
\newtheorem{proposition}[definition]{Proposition}
\newtheorem{lemma}[definition]{Lemma}
\newtheorem{corollary}[definition]{Corollary}
\newtheorem{remark}[definition]{Remark}

\newtheorem{assumption}[definition]{Assumption}
\newtheorem{example}[definition]{Example}

\newcounter{mycount}
\newenvironment{romlist}{\begin{list}{\rm(\roman{mycount})}%
   {\usecounter{mycount}\labelwidth=1cm\itemsep 0pt}}{\end{list}}

\newenvironment{letlist}{\begin{list}{(\alph{mycount})}%
   {\usecounter{mycount}\labelwidth=1cm\itemsep 0pt}}{\end{list}}

\numberwithin{section}{chapter}
\numberwithin{subsection}{section}
\numberwithin{equation}{section}
\numberwithin{definition}{section}
\numberwithin{figure}{chapter}

\begin{document}

\title{Graphical representations \\of Ising and Potts models\\
\vspace{1cm}
\large{Stochastic geometry of the quantum Ising model 
and the space--time Potts model}}

\author{\vspace{0.5cm}
\includegraphics{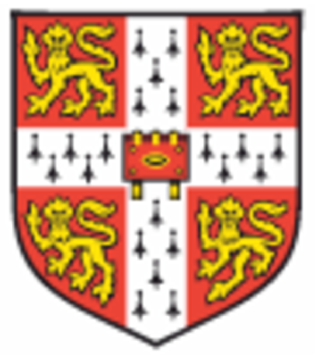}\\
\vspace{1.5cm}
Jakob Erik Bj\"ornberg\\
\vspace{0.8cm}
\large{Gonville \& Caius College and Statistical Laboratory\\
University of Cambridge}\\
\vspace{1cm}
\large{This dissertation is submitted for the degree of\\
\textit{Doctor of Philosophy}\\
June 2009}
}

\maketitle

\frontmatter

\chapter*{Preface}

This dissertation is the result of my own work and includes nothing
which is the outcome of work done in collaboration except where
specifically indicated in the text.

I would like to thank my PhD supervisor Geoffrey Grimmett.
Chapter~\ref{qim_ch} and Section~\ref{sec_1d} were done in
collaboration with him.  We have agreed that $65\%$
of the work is mine.  This work has appeared in a  journal
as a joint publication~\cite{bjogr2}.  
I am the sole author of the remaining material.
Section~\ref{starlike_sec} has been published in a journal~\cite{bjo0}.  

I would also like to thank the following.
Anders Bj\"orner and the
Royal Institute of Technology (KTH) in Stockholm, Sweden,
made this work possible through extremely generous support and funding.
The House of Knights (Riddarhuset) in Stockholm, Sweden, 
has supported me very generously throughout my studies.
I have received further generous support from the
Engineering and Physical Sciences Research Council under a Doctoral
Training Award to the University of Cambridge.
The final writing of this thesis took place during a very 
stimulating stay at the Mittag-Leffler Institute for Research in 
Mathematics, Djursholm, Sweden, during the spring of 2009.

\newpage

\chapter*{Summary}

Statistical physics seeks to explain macroscopic properties of
matter in terms of microscopic interactions.  Of particular interest
is the phenomenon of phase transition:  the sudden changes in macroscopic
properties as external conditions are varied.  
Two models in particular are of great interest to mathematicians,
namely the Ising model of a magnet and the percolation model of a 
porous solid.  These models in turn are part of the unifying framework of
the random-cluster representation, a model for random graphs
which was first studied
by Fortuin and Kasteleyn in the 1970's.  The random-cluster representation
has proved extremely useful in proving important facts about the
Ising model and similar models.

In this work we study the corresponding graphical framework for
two related models.  The first model is the transverse field quantum
Ising model, an extension of the original Ising model
which was introduced by Lieb, Schultz and Mattis in the 1960's.  
The second model
is the space--time percolation process, which is closely related
to the contact model for the spread of disease.  In 
Chapter~\ref{st_ch} we define the appropriate `space--time'
random-cluster model and explore a range of useful probabilistic
techniques for studying it.  The space--time Potts model emerges
as a natural generalization of the quantum Ising model.  The
basic properties of the phase transitions in these models are 
treated in this chapter, such as the fact that there is at most one
unbounded \fk-cluster, and the resulting  lower bound on the
critical value in~$\ZZ$.

In Chapter~\ref{qim_ch} we develop an alternative graphical
representation of the quantum Ising model, called the random-parity
representation.  This representation is based on the random-current
representation of the classical Ising model, and allows us to study
in much greater detail the phase transition and critical behaviour.
A major aim of this chapter is to prove sharpness of the phase
transition in the quantum Ising 
model---a central issue in the theory---and to establish bounds
on some critical exponents.  We address these issues by using the
random-parity representation to establish certain differential inequalities,
integration of which give the results.

In Chapter~\ref{appl_ch} we explore some consequences 
and possible extensions of
the results established in Chapters~\ref{st_ch} and~\ref{qim_ch}.
For example,  we determine the critical point for the quantum
Ising model in $\ZZ$ and in `star-like' geometries.


\tableofcontents

\newpage
\printnomenclature[2cm]

\mainmatter

\chapter{Introduction and background}

Many physical and mathematical systems undergo a \emph{phase transition},
of which some of the following examples may be familiar to the reader:
water boils at $100^\circ$C and freezes at $0^\circ$C;  
Erd\H{o}s-R\'enyi random graphs
produce a `giant component' if and only if the edge-probability
$p>1/n$;  and magnetic materials exhibit `spontaneous magnetization'
at temperatures below the Curie point.  In physical terminology, these
phenomena may be unified by saying that there is an `order parameter'
$M$ (density, size of largest component, magnetization) which behaves 
non-analytically on the parameters of the system at certain points.
In the words of Alan Sokal:  ``at a phase transition $M$ may be 
discontinuous, or continuous but not differentiable, or 16 times
differentiable but not 17 times''---any behaviour of this sort
qualifies as a phase transition.

\begin{figure}[hbt]
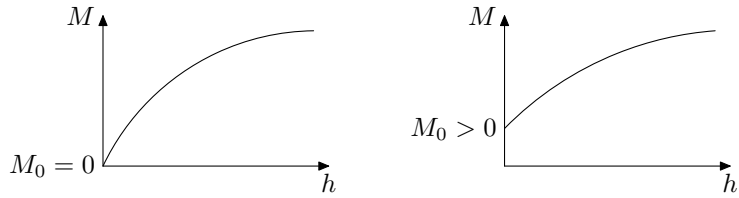

\centering
\includegraphics{thesis.31}
\qquad\includegraphics{thesis.32}
\caption[Magnetization at different temperatures]
{Magnetization $M$ when $T>T_\crit$ (left)
and when $T<T_\crit$ (right).  The residual magnetization $M_0$
is zero at high temperature and positive at low temperature.}
\label{mag_fig}   
\end{figure}

Since it is the example closest to the topic of this work, let us 
look  at the case of spontaneous magnetization.  For the moment
we will stay on an entirely intuitive level of description.
If one takes a piece of iron and places it in a magnetic field,
one of two things will happen.  When the strength of the external
field is decreased to nought, the iron piece may retain magnetization,
or it may not.  Experiments confirm that there is a critical value
$T_\crit$ of the  temperature $T$ such that:  if $T<T_\crit$
there is a residual (`spontaneous') magnetization, 
and if $T>T_\crit$ there is not.
See Figure~\ref{mag_fig}.
Thus the order parameter $M_0(T)$ 
(residual magnetization) is non-analytic at $T=T_\crit$
(and it turns out that the phase transition is of the
`continuous but not differentiable' variety, see
Theorem~\ref{crit_val_cor}).  Can we account for this behaviour
in terms of the `microscopic' properties of the material, that is
in terms of individual atoms and their interactions?

Considerable ingenuity has, since the 1920's and earlier, gone in to
devising mathematical models that strike a
good balance between three desirable properties:
physical relevance, mathematical (or computational) tractability, and
`interesting' critical behaviour.  A whole arsenal of mathematical tools,
rigorous as well as non-rigorous, have been developed to study such models.
One of the most exciting aspects of the mathematical theory of phase 
transition is the abundance of amazing conjectures originating
in the physics literature;
attempts by mathematicians to `catch up' with the physicists and
rigorously prove some of these conjectures have led to the 
development of many beautiful mathematical theories.  As an example
of this one can hardly at this time fail to mention the theory
of {\sle} which has finally established some long-standing
conjectures in two-dimensional 
models~\cite{smirnov01,smirnov_ising}.

This work is concerned with  the representation
of physical models using stochastic geometry, in particular what 
are called percolation-, \fk-, and random-current
representations.  A major focus of this work is on the 
\emph{quantum Ising model} of a magnet (described below);  
on the way to studying
this model we will also study `space--time' random-cluster
(or \fk) and Potts models.  Although a lot of attention has
been paid to the graphical representation of classical Ising-like
models, this is less true for quantum models, hence the current work.
Our methods are rigorous, and mainly utilize the mathematical
theory of probability.  Although graphical methods may give less 
far-reaching results 
than the `exact' methods favoured by mathematical physicists, they are
also more robust to changes in geometry:  towards the end of this work
we will see some examples of results on high-dimensional, and
`complex one-dimensional', models where exact methods cannot be used.

\section{Classical models}

\subsection{The Ising model}

The best-known, and most studied, model in statistical physics
is arguably the Ising model of a magnet, given as follows.  One
represents the magnetic material at hand by a finite graph 
$L=(V,E)$ where the vertices $V$ represent individual particles
(or atoms) and an edge is placed between particles that
interact (`neighbours').  A `state' is an assignment of the numbers
$+1$ and $-1$ to the vertices of $L$;  these numbers are usually
called `spins'.  The set $\{-1,+1\}^V$
of such states is  denoted $\S$, and an element of $\S$ is
denoted $\s$.  The model has two parameters, namely the
temperature $T\geq 0$  and the external magnetic
field $h\geq0$.  The probability of seeing a particular configuration
$\s$ is then proportional to the number
\begin{equation}\label{a1}
\exp\Big(\b\sum_{e=xy\in E}\s_x\s_y+\b h\sum_{x\in V}\s_x\Big).
\end{equation}
Here $\b=(k_{\mathrm{B}}T)^{-1}>0$ is the `inverse temperature',
where $k_{\mathrm{B}}$ is a constant called the `Boltzmann constant'.
Intuitively, the number~\eqref{a1} is bigger if more spins agree, since
$\s_x\s_y$ equals $+1$ if $\s_x=\s_y$ and $-1$ otherwise;  similarly
it is bigger if more spins `align with the external field' in that
$\s_x=+1$.  In particular, the spins at different sites are not
in general statistically independent, 
and the structure of this dependence
is subtly influenced by the geometry of the graph $L$.
This is what makes the model interesting.

The Ising model was introduced around 1925 (not originally by but \emph{to}
Ising by his thesis advisor Lenz) as a candidate for
a model that exhibits a phase transition~\cite{ising1925}.  
It turns out that 
the magnetization $M$, which is by definition the expected value 
of the spin at some given vertex, behaves (in the limit as the graph
$L$ approaches an infinite graph $\LL$) 
non-analytically on the parameters $\b,h$ at
a certain point $(\b=\b_\crit,h=0)$ in the $(\b,h)$-plane.

The Ising model is therefore the second-simplest physical model with an
interesting phase transition;  the simplest such model is the following.
Let $\LL=(\VV,\EE)$ be an infinite, but countable, graph.
(The main example to bear in mind is the lattice $\ZZ^d$ with
nearest-neighbour edges.)  Let $p\in[0,1]$ be given, and examine
each edge in turn, keeping it with probability $p$ and deleting
it with probability $1-p$, these choices being independent for
different edges.  The resulting subgraph of $\LL$ is typically
denoted $\om$, and the set of such subgraphs is denoted $\Om$.
The graph $\om$ will typically not be connected, but will break into
a number of connected components.  Is
one of these components infinite?  The model possesses a
phase transition in the sense that the probability that there
exists an infinite component jumps from 0 to 1 at a
critical value $p_\crit$ of~$p$.

This model is called \emph{percolation}.  It
was introduced by Broadbent and Hammersley in 1957
as a model for a porous material immersed in a 
fluid~\cite{broadbent_hammersley}.  Each edge in
$\EE$ is then thought of as a small hole which may be open
(if the corresponding edge is present in $\om$) or closed
to the passage of fluid.
The existence of an infinite component corresponds to the fluid
being able to penetrate from the surface to the `bulk' of the material.
Even though we are dealing here
with a countable set of independent random variables, the theory
of percolation is a genuine departure from the traditional
theory of sequences of independent variables, 
again since geometry plays such a vital role.

\subsection{The random-cluster model}

At first sight, the Ising- and percolation models seem unrelated,
but they have a common generalization.  On a finite graph $L=(V,E)$,
the percolation configuration $\om$ has probability
\begin{equation}\label{a2}
p^{|\om|}(1-p)^{|E\sm \om|},
\end{equation}
where $|\cdot|$ denotes the number of elements in a finite set, and 
we have identified the subgraph $\om$ with its edge-set.  A natural
way to generalize~\eqref{a2} is to consider absolutely continuous 
measures, and it turns out that the distributions defined by
\begin{equation}\label{a3}
\phi(\om):=p^{|\om|}(1-p)^{|E\sm \om|}\frac{q^{k(\om)}}{Z}
\end{equation}
are particularly interesting.  Here $q>0$ is an additional
parameter, $k(\om)$ is the number of connected components in $\om$,
and $Z$ is a normalizing constant.  The `cluster-weighting factor'
$q^{k(\om)}$ has the effect of skewing the distribution in favour
of few large components (if $q<1$) or many small components 
(if $q>1$), respectively.  This new model is called the
random-cluster model, and it contains percolation
as the special case $q=1$.  By considering limits as
$L\uparrow\LL$, one may see that the random-cluster models
(with $q\geq 1$) also have a phase transition in the same sense
as the percolation model, with associated critical probability
$p_\crit=p_\crit(q)$.

There is also a natural way to
generalize the Ising model.  This is easiest to describe
when $h=0$, which we assume henceforth.  
The relative weights~\eqref{a1} depend 
(up to a multiplicative constant) only on the number of 
adjacent vertices with equal spin, so the same model is obtained
by using the weights
\begin{equation}\label{a4}
\exp\Big(2\b\sum_{e=xy\in E}\d_{\s_x,\s_y}\Big),
\end{equation}
where $\d_{a,b}$ is 1 if $a=b$ and 0 otherwise.  
(Note that $\d_{\s_x,\s_y}=(\s_x\s_y+1)/2$.)
In this formulation
it is natural to consider the more general model when the spins $\s_x$
can take not only two, but $q=2,3,\dotsc$ different values,
that is each $\s_x\in\{1,\dotsc,q\}$.
Write $\pi$ for the corresponding distribution on spin
configurations;
the resulting model is called the $q$-state Potts model.
It turns out that the $q$-state Potts models is closely related
to the random-cluster model, one manifestation of this being
the following.  
(See~\cite{fortuin_kasteleyn_1972}, 
or~\cite[Chapter~1]{grimmett_rcm} for a modern proof.)

\begin{theorem}\label{thm1}
If $q\geq2$ is an integer and $p=1-e^{-2\b}$ then for all $x,y\in V$
\[
\pi(\s_x=\s_y)-\frac{1}{q}=
\Big(1-\frac{1}{q}\Big)\phi(x\lra y)
\]
\end{theorem}
Here $\pi(\s_x=\s_y)$ denotes the probability that, in the Potts
model, the spin at $x$ takes the same value as the spin at $y$.
Similarly, $\phi(x\lra y)$ is the probability that, in the 
random-cluster model, $x$ and $y$ lie in the same component of $\om$.
Since the right-hand-side concerns a typical graph-theoretic property
(connectivity), the random-cluster model 
is called a `graphical representation' of the Potts model.
The close relationship between the random-cluster and Potts models
was unveiled by Fortuin and Kasteleyn during the 1960's and 1970's
in a series of papers including~\cite{fortuin_kasteleyn_1972}.
The random-cluster model is therefore sometimes called the
`\fk-representation'.  In other words, Theorem~\ref{thm1} says that
the correlation between distant spins in the Potts model
is translated to the existence of paths between the sites in the
random-cluster model.
Using this and related facts one can deduce
many important things about the phase transition of the Potts model
by studying  the random-cluster model.
This can be extremely useful since the random-cluster formulation
allows geometric arguments that are  not present in the 
Potts model.  Numerous examples of this may be found 
in~\cite{grimmett_rcm};  very recently, in~\cite{smirnov_ising},
the `loop' version of the
random-cluster model was also used to prove conformal invariance for the
two-dimensional Ising model, a major breakthrough
in the theory of the Ising model.

\subsection{Random-current representation}
\label{intro_rcr}

For the  Ising model there exists also another
graphical representation, distinct from the random-cluster
model.  This is called the 
`random-\emph{current} representation' and was developed in a sequence
of papers in the late 1980's \cite{aiz82,abf,af}, building on
ideas in~\cite{ghs}.
These papers answered many questions for the Ising model 
on $\LL=\ZZ^d$ with $d\geq 2$ that are 
still to this day unanswered for general Potts models.  Cast
in the language of the $q=2$ random-cluster model, these
questions include the following [answers in square brackets].
\begin{itemize}
\item If $p<p_\crit$, is the expected size of a component 
finite or infinite? [Finite.]
\item If $p<p_\crit$, do the connection
probabilities $\phi(x\lra y)$ go to zero exponentially fast 
as $|x-y|\rightarrow\oo$? [Yes.]
\item At $p=p_\crit$, does $\phi(x\lra y)$ go to zero 
exponentially fast as $|x-y|\rightarrow\oo$?  [No.]
\end{itemize}
In fact, even more detailed information could be obtained,
especially in the case $d\geq4$, giving at least partial
answer to the question
\begin{itemize}
\item How does the magnetization $M=M(\b,h)$ behave as the critical 
point $(\b_\crit,0)$ is approached?
\end{itemize}
It is one of the main objectives of this work to develop a random-current
representation for the \emph{quantum} Ising model (introduced in the
next section), and answer the above
questions also for that model.  

Here is a very brief sketch of the random-current representation of
the classical Ising model.
Of particular importance is the normalizing constant or
`partition function' that makes~\eqref{a1} a probability distribution,
namely
\begin{equation}\label{a5}
\sum_{\s\in\S}\exp\Big(\b\sum_{e=xy\in E}\s_x\s_y\Big)
\end{equation}
(we assume that $h=0$ for simplicity).  
We rewrite~\eqref{a5} using the following steps.
Factorize the exponential
in~\eqref{a5} as a product over $e=xy\in E$, and then expand each
factor as a Taylor series in the variable $\b\s_x\s_y$.  By interchanging
sums and products we then obtain a weighted sum over vectors
$\ul m$ indexed by $E$ of a quantity which (by $\pm$ symmetry)
is zero if a certain condition on $\ul m$ fails
to be satisfied, and a positive
constant otherwise.  The condition on $\ul m$ is that: for each $x\in V$
the sum over all edges $e$ adjacent to $x$ of $m_e$ is a multiple
of 2.  

Once we have rewritten the partition function in this way, we may
interpret the weights on $\ul m$ as probabilities.
It follows that the partition function is (up to a multiplicative 
constant) equal to the probability that the random graph $G_{\ul m}$
with each edge $e$ replaced by $m_e$ parallel edges is 
\emph{even} in that each vertex has even total degree.  
Similarly, other quantities of interest may be expressed
in terms of the probability that only a given set of vertices fail to have
even degree in $G_{\ul m}$;  for example, the correlation between
$\s_x$ and $\s_y$ for $x,y\in V$ is expressed in terms of
the probability that only $x$ and $y$ fail to have even degree.
By elementary graph theory, the latter event implies the
existence of a path from $x$ to $y$ in $G_{\ul m}$.  
By studying connectivity
in the above random graphs with restricted degrees
one obtains surprisingly detailed
information about the Ising model.  Much more will be said about this
method in Chapter~\ref{qim_ch}, see for example the Switching Lemma
(Theorem~\ref{sl}) and its applications in Section~\ref{sw_appl_sec}.

\section{Quantum models and space--time models}

There is a version of the Ising model formulated to meet the 
requirements of quantum theory, introduced in~\cite{lieb}.
We will only be concerned with the \emph{transverse field}
quantum Ising model.
Its definition and physical motivation bear a certain level
of complexity which it is beyond the scope of this work to justify
in an all but very cursory manner.  One is given, as before, a finite
graph $L=(V,E)$, and one is interested
in the properties of certain matrices (or `operators') acting on
the Hilbert space $\cH=\bigotimes_{v\in V}\CC^2$%
\nomenclature[H]{$\cH$}{Hilbert space $\bigotimes_{v\in V}\CC^2$}.
The set $\S=\{-1,+1\}^V$ may now be identified with a basis
for $\cH$, defined by letting each factor $\CC$ in the tensor
product have basis consisting of the two vectors 
$\ket{+}:=\big(\begin{smallmatrix} 1 \\ 0\end{smallmatrix}\big)$ 
and $\ket{-}:=\big(\begin{smallmatrix} 0 \\ 1\end{smallmatrix}\big)$.
We write $\ket{\s}=\bigotimes_{v\in V}\ket{\s_v}$%
\nomenclature[>]{$\ket{\s}$}{Basis vector in $\cH$}
 for these basis vectors.
In addition to
the inverse temperature $\b>0$, one is given parameters $\l,\d>0$,
interpreted as spin-coupling and transverse field intensities,
respectively.  The latter specify the \emph{Hamiltonian}  
\begin{equation}\label{qi_ham_eq}
H=-\tfrac{1}{2}\l\sum_{e=uv\in E}\s_u^{(3)}\s_v^{(3)}-
\d\sum_{v\in V}\s_v^{(1)},
\end{equation}
where the `Pauli spin-$\frac12$ matrices' are given as 
\begin{equation}
\s^{(3)}=
\Bigg(\begin{matrix} 
1 & 0 \\
0 & -1
\end{matrix}\Bigg),
\qquad
\s^{(1)}=
\Bigg(\begin{matrix}
0 & 1 \\
1 & 0
\end{matrix}\Bigg),
\end{equation}%
\nomenclature[s]{$\s^{(1)},\,\s^{(3)}$}{Pauli matrices}%
and $\s^{(i)}_v$ acts on the copy of $\CC^2$ in $\cH$ indexed 
by $v\in V$.  Intuitively, the matrices $\s^{(1)}$ and $\s^{(3)}$
govern spins in `directions' $1$ and $3$ respectively (there is 
another matrix $\s^{(2)}$ which does not feature in this model).
The external field is called `transverse' since it acts in a different
`direction' to the internal interactions.  When $\d=0$ this model
therefore reduces to
the (zero-field) classical Ising model (this will be obvious
from the space--time formulation below).

The basic operator of interest is $e^{-\b H}$, which is thus 
a (Hermitian) matrix acting on $\cH$;  
one usually normalizes it and studies
instead the matrix $e^{-\b H}/\tr(e^{-\b H})$.  
Here the \emph{trace} of the Hermitian matrix $A$ is
defined as
\begin{equation*}
\tr(A) = \sum_{\s\in\Si} \bra{\s}A\ket{\s},
\end{equation*}%
\nomenclature[>]{$\ket{\pm}$}{Basis of $\CC^2$}%
\nomenclature[t]{$\tr(\cdot)$}{Trace}%
where $\bra{\s}$%
\nomenclature[<]{$\bra{\cdot}$}{Conjugate transpose}
 is the adjoint, or conjugate transpose, of 
the column vector $\ket{\s}$,
and we are using the usual matrix product.
An eigenvector of $e^{-\b H}/\tr(e^{-\b H})$ may be
thought of as a `state' of the system, and is now a `mixture'
(linear combination) of classical states in $\S$;  the corresponding
eigenvalue (which is real since the matrix is Hermitian) 
is related to the `energy level' of the state.

In this work we will not be working directly with
this formulation of the quantum Ising 
model, but a (more probabilistic) `space--time' formulation, which
we describe briefly now.  It is by now standard that many properties of
interest in the transverse field 
quantum Ising model may be studied by means of
a `path integral' representation, which maps the model onto
a type of classical Ising model on the continuous space
$V\times[0,\b]$.  (To be precise, the endpoints of the interval
$[0,\b]$ must be identified for this mapping to hold.)
This was first used in~\cite{ginibre69}, 
but see also for
example~\cite{akn,aizenman_nacht,campanino_klein_perez,chayes_ioffe_curie-weiss,GOS,nachtergaele93} 
and the recent surveys to be found in \cite{G-pgs,ioffe_geom}.
Precise definitions will be given in Chapter~\ref{st_ch},
but in essence we must consider piecewise constant
\emph{functions} $\s:V\times[0,\b]\rightarrow\{-1,+1\}$,
which are random and have a distribution reminiscent of~\eqref{a1}.
The resulting model is called the `space--time Ising model'.
As for the classical case, it is straightforward to generalize
this to a space--time \emph{Potts} model with $q\geq 2$
possible spin values, and also to give a graphical representation
of these models
in terms of a space--time random-cluster model.
Although the partial continuity of the underlying geometry
poses several technical difficulties,
the corresponding theory is very similar to the classical
random-cluster theory.  The most important basic properties of the
models are developed in detail in Chapter~\ref{st_ch}.  
On taking limits as $L$ and/or $\b$ become infinite, one may
speak of the existence of unbounded connected components, and
one finds (when $\b=\oo$) 
that there is a critical dependence on the ratio
$\rho=\l/\d$ of the probability of seeing such a component.
One may also develop, as we do in Chapter~\ref{qim_ch}, a type of
random-current representation of the space--time Ising model
which allows us to deduce many facts about the critical behaviour
of the \emph{quantum} Ising model.

Other models of space--time type have been around for a long time
in the probability literature.  Of these the most relevant for us
is the \emph{contact process} (more precisely, its graphical
representation), see for example~\cite{liggett85,liggett99} 
and references therein.  In the contact process, one imagines
individuals placed on the vertices of a graph, such as $\ZZ^2$.
Initially, some of these individuals may be infected with a contagious
disease.  As time passes, the individuals themselves stay fixed but
the disease may spread:  individuals may be infected by their
neighbours, or by a `spontaneous' infection.
Infected individuals may recover spontaneously.
Infections and recoveries are governed by
Poisson processes, and depending on the ratio of infection rate
to recovery rate the infection may or may not persist indefinitely.
The contact model may be regarded as the $q=1$ or `independent'
case of the space--time random-cluster model (one difference
is that we in the space--time model regard time as `undirected').
Thus one may get to general space--time random-cluster models
in a manner reminiscent of the classical case, by skewing
the distribution by an appropriate `cluster weighting factor'.
This approach will be treated in detail in Section~\ref{basics_sec}.

\section{Outline}

A brief outline of the present work follows.  In Chapter~\ref{st_ch},
the space--time random-cluster and Potts models are defined.  As for the 
classical theory, one of the most important tools is 
\emph{stochastic comparison},
or the ability to compare the probabilities of certain events under
measures with different parameters.  A number of results of this type
are presented in Section~\ref{stoch_ineq_sect}.  We then consider the
issue of defining random-cluster and Potts measures on infinite graphs,
and of their phase transitions.
We etablish the existence of weak limits of Potts and random-cluster
measures as $L\uparrow\LL$, and introduce the central question of
when there is a unique such limit.  It turns out that this question is
closely related to the question if there can be an unbounded
connected component;  this helps us to define a critical 
value $\rho_\crit(q)$.
In general not a lot can be said about the precise value of 
$\rho_\crit(q)$, but in the case when $\LL=\ZZ$ there are additional
geometric (duality) arguments that can be used to show that
$\rho_\crit(q)\geq q$.

Chapter~\ref{qim_ch} deals exclusively with the quantum Ising model
in its space--time formulation.  We develop the 
`random parity representation', which is the space--time analog of the
random-current representation, and the tools associated with it,
most notably the switching lemma.  This representation allows us to
represent truncated correlation functions in terms of single geometric
events.  Since truncated correlations are closely related to the
derivatives of the magnetization $M$, we can use this 
to prove a number of inequalities between the different partial derivatives
of $M$, along the lines of~\cite{abf}.  Integrating these
differential inequalities gives the information
on the critical behaviour that was referred to in Section~\ref{intro_rcr},
namely the sharpness of the phase transition, bounds on critical
exponents, and the vanishing of the mass gap.
Chapter~\ref{qim_ch} (as well as Section~\ref{sec_1d})
is joint work with Geoffrey Grimmett, and appears in the
article \emph{The phase transition of the quantum Ising model
is sharp}~\cite{bjogr2}, published by 
the Journal of Statistical Physics.

Finally, in Chapter~\ref{appl_ch}, we combine the results of 
Chapter~\ref{qim_ch} with the results of Chapter~\ref{st_ch}
in some concrete cases.  Using  duality
arguments  we prove that the critical
ratio $\rho_\crit(2)=2$ in the case $\LL=\ZZ$.  We then develop some 
further geometric arguments for the random-cluster representation to
deduce that the critical ratio is the same as for $\ZZ$ on a much
larger class of `$\ZZ$-like' graphs.  These arguments
(Section~\ref{starlike_sec}) appear in the article
\emph{Critical value of the quantum Ising model on star-like
graphs}~\cite{bjo0}, 
published in the Journal of Statistical
Physics.  We conclude by describing some future
directions for research in this area.

\chapter{Space--time models: \\
random-cluster, Ising, and Potts}
\label{st_ch}

\begin{quote}
{\it Summary.}  We provide basic definitions and facts
pertaining to the space--time random-cluster and -Potts models.
Stochastic inequalities, a major tool in the theory, are proved
carefully, and the notion of phase transition
is defined.  We also introduce the notion of graphical duality.
\end{quote}

\section{Definitions and basic facts}\label{basics_sec}

The space--time models we consider
live on the product of a graph with the real line.  To
define space--time random-cluster and Potts models we first work on
bounded subsets of this product space, and then pass to a limit.
The continuity of $\RR$ makes the definitions of boundaries
and boundary conditions more delicate than in the discrete case.

\subsection{Regions and their boundaries}

Let $\LL=(\VV,\EE)$%
\nomenclature[L]{$\LL$}{Infinite graph}%
\nomenclature[V]{$\VV$}{Vertex set of $\LL$}%
\nomenclature[E]{$\EE$}{Edge set of $\LL$}
 be a countably infinite, 
connected, undirected graph, which is \emph{locally finite} in that
each vertex has finite degree.  Here $\VV$ is the vertex set 
and $\EE$ the edge set.  For simplicity we assume that $\LL$
does not have multiple edges or loops.  An edge of $\LL$
with endpoints $u$, $v$ is denoted by $uv$. 
We write $u \sim v$ if $uv \in \EE$.
The main example to bear in mind is when $\LL=\ZZ^d$ is the $d$-dimensional
lattice, with  edges between points that differ
by one in exactly one coordinate.

\begin{figure}[hbt]
\centering
\includegraphics{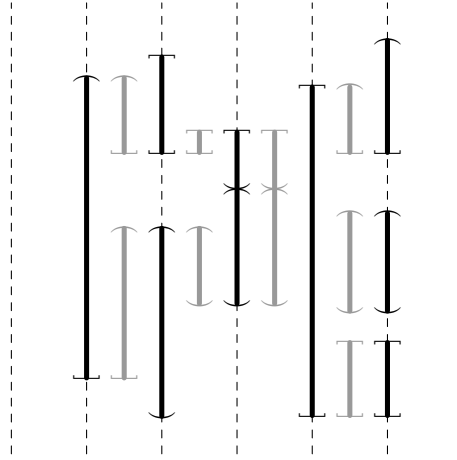}
\caption[A region]
{A region $\L=(K,F)$ as a subset of $\LLam$  
when $\LL=\ZZ$.  Here $\KK$ is drawn dashed, $K$ is drawn
bold black, and $F$ is drawn bold grey.
An endpoint of an interval in $K$ (respectively, $F$) is drawn as a square
bracket if it is included in $K$ (respectively, $F$) 
or as a rounded bracket if it is not.}
\label{region_fig}   
\end{figure}

Let
\begin{gather}
\KK:=\bigcup_{v\in\VV} (v\times \RR),\qq
\FF:=\bigcup_{e\in\EE} (e\times \RR),
\label{o8}\\
\LLam:=(\KK,\FF).
\end{gather}%
\nomenclature[K]{$\KK$}{The product $\VV\times\RR$}%
\nomenclature[F]{$\FF$}{The product $\EE\times\RR$}%
\nomenclature[T]{$\LLam$}{The pair $(\KK,\FF)$}%
Let $L=(V,E)$%
\nomenclature[L]{$L$}{Finite subgraph of $\LL$}%
\nomenclature[V]{$V$}{Vertex set of $L$}%
\nomenclature[E]{$E$}{Edge set of $L$}
 be a finite 
connected subgraph of $\LL$.
In the case when $\LL=\ZZ^d$, the main example for
$L$ is the `box' $[-n,n]^d$.
For each $v\in V$, let $K_v$ be a finite union of
(disjoint) bounded intervals in $\RR$.  
No assumption is made whether 
the constituent intervals are open, closed, or half-open.  For 
$e=uv\in E$ let $F_e:=K_u\cap K_v\se\RR$.  Let
\begin{equation}
K:=\bigcup_{v\in V}(v\times K_v),\quad
F:=\bigcup_{e\in E}(e\times F_e).
\end{equation}%
\nomenclature[K]{$K$}{Subset of $\KK$}%
\nomenclature[F]{$F$}{Subset of $\FF$}%
We define a \emph{region} to be a pair
\begin{equation}\label{def-newL}
\L=(K,F)
\end{equation}%
\nomenclature[L]{$\L$}{Region}%
for $L$, $K$ and $F$ defined as above.  We will often think of
$\L$ as a subset of $\LLam$ in the natural way, see 
Figure~\ref{region_fig}.
Since a region $\L=(K,F)$ is completely determined
by the set $K$, we will sometimes
abuse notation by writing $x\in\L$ when we mean $x\in K$, and
think of subsets of $K$ (respectively, $\KK$) as subsets 
of $\L$ (respectively, $\LLam$).

An important type of a region is a \emph{simple region},
defined as follows.  For $L$  as above, let $\b>0$ and let 
$K$ and $F$ be given by letting each $K_v=[-\b/2,\b/2]$.  Thus
\begin{gather}
K=K(L,\b):=\bigcup_{v\in V} (v\times [-\b/2,\b/2]),\\
F=F(L,\b):=\bigcup_{e\in E} (e\times [-\b/2,\b/2]),\\
\L=\L(L,\b):=(K,F).
\label{def-oldL}
\end{gather}
Note that in a simple region, the intervals constituting $K$ are
all closed.  (Later, in the quantum Ising model of 
Chapter~\ref{qim_ch}, the parameter $\b$ will be interpreted as 
the `inverse temperature'.)

Introduce an additional point $\G$%
\nomenclature[G]{$\G$}{Ghost site}
 external to $\LLam$, 
to be interpreted as a 
`ghost-site' or `point at infinity';  the use of $\G$ will be
explained below, when the space--time random-cluster and Potts
models are defined.  Write $\LLam^\G=\LLam\cup\{\G\}$,
$\KK^\G=\KK\cup\{\G\}$, and similarly for other notation.

\begin{figure}[hbt]
\centering
\includegraphics{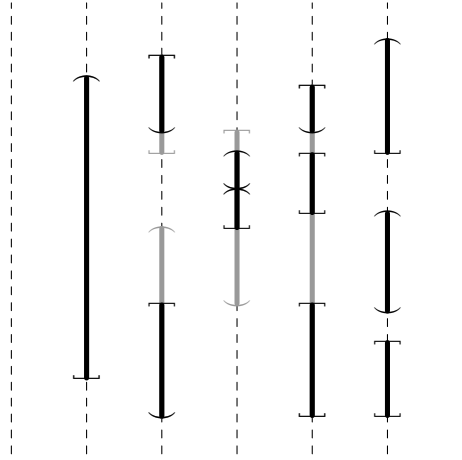}
\caption[Boundary of a region]
{The (inner) boundary $\boundary\L$ 
of the region $\L$ of Figure~\ref{region_fig}
is marked black, and $K\setminus\boundary\L$ is marked grey.  An
endpoint of an interval in $\boundary\L$ is drawn as a square bracket
if it lies in $\boundary\L$ and as a round bracket otherwise.}
\label{boundary_fig}   
\end{figure}

We will require two distinct notions of boundary for
regions $\L$.  For $I\subseteq\RR$ we denote the closure
and interior of $I$ by $\ol I$ and $I^\circ$, respectively.
For $\L$ a region as in~\eqref{def-newL}, define the
\emph{closure} to be the region 
$\ol\L=(\ol K,\ol F)$ given by
\begin{equation}
\ol K:=\bigcup_{v\in V}(v\times \ol K_v),\quad
\ol F:=\bigcup_{e\in E}(e\times \ol F_e);
\end{equation}%
\nomenclature[L]{$\ol \L$}{Closure of the region $\L$}%
similarly define the \emph{interior} of $\L$ to be the region
$\L^\circ=(K^\circ,F^\circ)$ given by
\begin{equation}
K^\circ:=\bigcup_{v\in V}(v\times K_v^\circ),\quad
F^\circ:=\bigcup_{e\in E}(e\times F_e^\circ).
\end{equation}%
\nomenclature[L]{$\L^\circ$}{Interior of the region $\L$}%
Define the \emph{outer boundary}%
\nomenclature[d]{$\partial\L$}{Outer boundary}
 $\partial \L$ of $\L$ to be the union
of $\ol K\setminus K^\circ$ with the set of points
$(u,t)\in K$ such that $u\sim v$ for some $v\in\VV$ such that
$(v,t)\not\in K$.  Define the \emph{inner boundary}
$\boundary\L$  of $\L$ by
$\boundary\L:=(\partial\L)\cap K$%
\nomenclature[d]{$\boundary\L$}{(Inner) boundary}.
The inner boundary of $\L$ will often simply be called 
the \emph{boundary} of $\L$.
Note that if $x$ is an endpoint
of a closed interval in $K_v$, then $x\in\partial\L$ if and only if
$x\in\boundary\L$, but if $x$ is an endpoint
of an open interval in $K_v$, then $x\in\partial\L$ but
$x\not\in\boundary\L$.  In particular, if $\L$ is a simple region
then $\partial\L=\boundary\L$.
  A word of caution:  this
terminology is nonstandard, in that for example the interior and the
boundary of a region, as defined above, need not be disjoint.
See Figure~\ref{boundary_fig}.  
We define 
$\partial\L^\G=\partial\L\cup\{\G\}$ and
$\boundary\L^\G=\boundary\L\cup\{\G\}$.

A subset $S$ of $\KK$ will be called \emph{open} if it equals a union
of the form 
\[
\bigcup_{v\in\VV}(v\times U_v),
\]
where each $U_v\se\RR$ is an open set.  Similarly for subsets of $\FF$.
The $\s$-algebra generated by this topology on $\KK$ (respectively, on $\FF$)
will be denoted $\cB(\KK)$%
\nomenclature[B]{$\cB(\KK)$}{Borel $\s$-algebra}
 (respectively, $\cB(\FF)$) and will be
referred to as the Borel $\s$-algebra.

Occasionally, especially in Chapter~\ref{qim_ch}, we will in place
of $\LLam$ be using the finite $\b$ space
$\LLam_\b=(\KK_\b,\FF_\b)$ given by
\begin{equation}\label{def-bL}
\KK_\b:=\bigcup_{v\in\VV} (v\times [-\b/2,\b/2]),\qq
\FF_\b:=\bigcup_{e\in\EE} (e\times [-\b/2,\b/2]).
\end{equation}%
\nomenclature[T]{$\LLam_\b$}{Finite-$\b$ space}%
This is because in the quantum Ising model $\b$ is thought of as
`inverse temperature', and then both $\b<\oo$ (positive temperature)
and $\b=\oo$ (ground state) are interesting.  

In what follows, proofs will often, for simplicity, be given 
for simple regions only;  proofs
for general regions will in these cases be straightforward adaptations.
We will frequently be using integrals of the forms
\begin{equation}
\int_K f(x)\,dx\qquad\text{and}\qquad\int_F g(e)\, de.
\end{equation}
These are to be interpreted, respectively, as
\begin{equation}\label{KF_int}
\sum_{v\in V}\int_{K_v} f(v,t)\,dt,
\qq 
\sum_{e\in E} \int_{F_e} g(e,t) \,dt.
\end{equation}
If $A$ is an event, we will write $\one_A$ or $\one\{A\}$ for the 
indicator function of~$A$.

\subsection{The space--time percolation model}

Write $\RRp=[0,\oo)$ and let
$\l:\FF\rightarrow\RRp$,%
\nomenclature[l]{$\l$}{Intensity of $B$}
 $\d:\KK\rightarrow\RRp$,%
\nomenclature[d]{$\d$}{Intensity of $D$}
 and 
$\g:\KK\rightarrow\RRp$%
\nomenclature[g]{$\g$}{Intensity of $G$}
 be bounded functions.
We assume throughout that
$\l,\d,\g$ are all Borel-measurable.
We retain the notation $\l$, $\d$, $\g$
for the restrictions of these functions to $\L$, given in \eqref{def-newL}.
Let $\Om$%
\nomenclature[O]{$\Om$}{Percolation configuration space}
 denote the set of triples $\om=(B,D,G)$%
\nomenclature[o]{$\om$}{Percolation configuration}
 of countable subsets
$B\se \FF$,%
\nomenclature[B]{$B$}{Process of bridges}
 $D,G\se \KK$;%
\nomenclature[D]{$D$}{Process of deaths}%
\nomenclature[G]{$G$}{Process of ghost-bonds}
  these triples will often be called \emph{configurations}.
Let $\mu_\l$,%
\nomenclature[m]{$\mu_\l$}{Law of $B$}
 $\mu_\d$,%
\nomenclature[m]{$\mu_\d$}{Law of $D$}
 $\mu_\g$%
\nomenclature[m]{$\mu_\g$}{Law of $G$}
 be the probability measures associated with independent Poisson processes
on $\KK$ and $\FF$ as appropriate, with
respective intensities $\l$, $\d$, $\g$. 
Let $\mu$%
\nomenclature[m]{$\mu$}{Law of space--time percolation}
 denote the probability measure $\mu_\l\times\mu_\d\times\mu_\g$
on $\Om$.
Note that, with $\mu$-probability 1, each of the countable sets 
$B,D,G$ contains no accumulation points;  
we call such a set \emph{locally finite}.
We will sometimes write $B(\om),D(\om),G(\om)$ for clarity.
\begin{remark}\label{rem-as} 
For simplicity of notation we will
frequently overlook events of probability zero, and will thus assume for
example that $\Om$ contains only triples
$(B,D,G)$ of locally finite  sets, such that
no two points in $B\cup D\cup G$ have the same $\RR$-coordinates.  
\end{remark}
For the purpose of defining a metric and a $\s$-algebra on $\Om$,
it is convenient to identify each $\om\in\Om$ with a collection
of step functions.  To be definite, we then regard each 
$\om\cap(v\times\RR)$ and each $\om\cap(e\times\RR)$ as an
\emph{increasing, right-continuous} step function, which equals 0 at 
$(v,0)$ or $(e,0)$ respectively.
There is a metric on the space of right-continuous
step functions on $\RR$, called the Skorokhod metric, which may
be extended in a straightforward manner to a metric on
$\Om$.  Details may be found in  Appendix~\ref{skor_app},
alternatively 
see~\cite{bezuidenhout_grimmett}, and~\cite[Chapter~3]{ethier_kurtz}
or~\cite[Appendix~1]{lindvall}.
We let $\cF$%
\nomenclature[F]{$\cF$}{Skorokhod $\s$-algebra on $\Om$}
 denote the $\s$-algebra on 
$\Om$ generated by the Skorokhod metric.  Note that 
the metric space $\Om$ is 
\emph{Polish}, that is to say separable (it contains a countable
dense subset) and complete (Cauchy sequences converge).

\begin{figure}[hbt]
\centering
\includegraphics{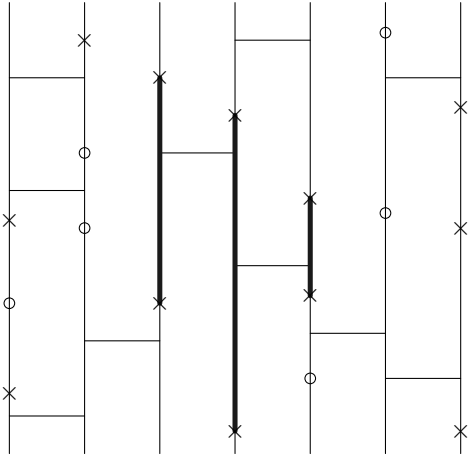}
\caption[Space--time percolation configuration]
{Part of a configuration $\om$ when $\LL=\ZZ$.  Deaths are marked as
crosses and bridges as horizontal line segments;  the positions of
ghost-bonds are marked as small circles.  One of the connected components 
of $\om$ is drawn bold.}\label{sample_fig}  
\end{figure}

However, in the context of percolation,
here is how we usually want to
think about elements of $\Om$.   Recall the
`ghost site' or `point at infinity' $\G$.
Elements of $D$ are thought of as `deaths', or missing points; 
elements of $B$
as `bridges' or line segments between points $(u,t)$ and $(v,t)$, 
$uv\in\EE$; and elements of $G$ as `bridges to $\ghost$'.  
See Figure~\ref{sample_fig} for an
illustration of this.  Elements of $B$ will sometimes
be referred to as \emph{lattice bonds} and elements of $G$ as 
\emph{ghost bonds}.  A lattice bond $(uv,t)$ is said to 
have \emph{endpoints} $(u,t)$ and $(v,t)$;  a ghost bond at
$(v,t)$ is said to have endpoints $(v,t)$ and~$\G$.

For two points $x,y\in\KK$ we say that
there is a \emph{path}, or an \emph{open path}, in $\om$ between $x$ and $y$
if there is a sequence $(x_1,y_1),\dotsc,(x_n,y_n)$ of pairs of elements
of $\KK$ satisfying the following:
\begin{itemize}
\item Each pair $(x_i,y_i)$ consists either of the two endpoints
of a single lattice bond (that is, element of $B$) or of the endpoints
in $\KK$ of two distinct ghost bonds (that is, elements of $G$),
\item Writing $y_0=x$ and $x_{n+1}=y$, we have that
for all $0\leq i\leq n$, there is a $v_i\in\VV$ such that
$y_{i},x_{i+1}\in (v_i\times\RR)$,
\item For each $0\leq i\leq n$, the (closed) interval 
in $v_i\times\RR$ with endpoints
$y_i$ and $x_{i+1}$ contains no elements of $D$.
\end{itemize}
In words, there is a path between $x$ and $y$ if
$y$ can be reached from $x$ by traversing bridges and ghost-bonds,
as well as subintervals of $\KK$ which do not contain elements of $D$.
For example,
in Figure~\ref{sample_fig} there is an open path between any two points
on the line segments that are drawn bold.
By convention, there is always an open path from $x$ to itself.
We say that there is a path between $x\in\KK$
and $\G$ if there is a $y\in G$ such that there is a path 
between $x$ and $y$.
Sometimes we say that $x,y\in\KK^\G$ are
\emph{connected} if there is an open path between them.  
Intuitively, elements of $D$ break connections on vertical lines, and 
elements of $B$ create connections between neighbouring lines.
The use of $\G$, and the process $G$, is to provide a 
`direct link to $\infty$';  two points that are joined to $\G$ are
automatically joined to eachother.

We write $\{x\lra y\}$ for
the event that there is an open path between $x$ and $y$.
We say that
two subsets $A_1,A_2\se\KK$ are connected, and write
$A_1\lra A_2$, if there exist $x\in A_1$
and $y\in A_2$ such that $x\lra y$.
For a region $\L$, we say that there is
an open path between $x,y$ \emph{inside $\L$} if $y$ can be 
reached from $x$ by traversing
death-free line segments, bridges, and ghost-bonds 
that all lie in $\L$.  Open paths \emph{outside} $\L$ are defined similarly.

\begin{definition}
With the above interpretation, the measure $\mu$ on $(\Om,\cF)$ is called 
the space--time percolation measure
on $\LLam$ with parameters $\l,\d,\g$.
\end{definition}

The measure $\mu$ coincides with the law of the graphical representation
of a contact process
with spontaneous infections, 
see~\cite{aizenman_jung,bezuidenhout_grimmett}.  
In this work, however, we regard
`time' as undirected, and thus think of $\om$ as a geometric object
rather than as a process evolving in time.  

\subsection{Boundary conditions}

Any $\om\in\Om$ breaks into \emph{components}, where a component is
by definition the maximal subset
of $\KK^\G$ which can be reached 
from a given point in $\KK^\G$ by traversing
open paths.  See Figure~\ref{sample_fig}.
One may imagine $\KK$ as
a collection of infinitely long strings, which are cut at deaths, tied
together at bridges, and also tied to $\G$ at ghost-bonds.  
The components are
the pieces of string that `hang together'.
The random-cluster measure, which is defined in the next subsection,
is obtained by `skewing' the percolation measure $\mu$ in favour
of either many small, or a few big, components.  Since the total
number of components in a typical $\om$ is infinite, we must first,
in order to give an analytic definition,
restrict our attention to the number of components which intersect
a fixed region $\L$.  We consider a number of different rules
for counting those components which intersect the boundary of $\L$.
Later we will be interested in limits as the region $\L$ grows,
and whether or not these `boundary conditions' have an effect on 
the limit.

\begin{figure}[hbt]
\centering
\includegraphics{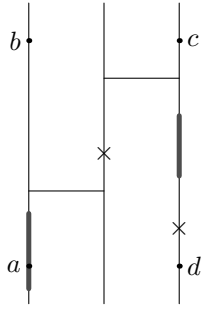}
\caption[Connections with boundary condition]
{Connectivities with respect to the boundary condition
$b=\{P_1\}$, where $P_1\setminus\{\G\}$ is the subset drawn bold.  
The following connectivities hold:
$a\lra b$, $a\lra c$, $a\not\lra d$.
(This picture does not specify which endpoints of the 
subintervals of $P_1$ lie in $P_1$.)}\label{bc_fig}  
\end{figure}

Let $\L=(K,F)$ be a region.  We define a 
\emph{random-cluster boundary condition} $b$%
\nomenclature[b]{$b$}{Boundary condition}
 to be a 
finite nonempty collection
$b=\{P_1,\dotsc,P_m\}$, where the $P_i$ are disjoint, 
nonempty subsets of $\boundary\L^\G$, such that each 
$P_i\setminus\{\G\}$ is a finite union of intervals.
(These intervals may be open, closed, or half-open, and may
consist of a single point.)
We require that $\G$ lies in one of the $P_i$, and by convention
we will assume that $\G\in P_1$.  Note that the union of
the $P_i$ will in general be a proper subset of $\boundary\L^\G$.
For $x,y\in\L^\G$ we say that $x\lra y$ \emph{with respect to $b$}
if there is a sequence $x_{1},\dotsc,x_{l}$ (with $0\leq l\leq m$)
such that 
\begin{itemize}
\item Each $x_{j}\in P_{i_j}$ for some $0\leq i_j\leq m$;
\item There are open paths inside $\L$ from $x$ to $x_1$
and from $x_l$ to $y$;
\item For each $j=1,\dotsc,l-1$ there is some point 
$y_{j}\in P_{i_j}$ such that there is a path inside $\L$ from $y_j$ to
$x_{j+1}$. 
\end{itemize}
See Figure~\ref{bc_fig} for an example.

When $\L$ and $b$ are fixed and $x,y\in\L^\G$, we will typically 
without mention
use the symbol $x\lra y$ to mean that there is a path between $x$
and $y$ in $\L$ with respect to $b$.
Intuitively, each $P_i$ is thought of as \emph{wired together};
as soon as you reach one point $x_j\in P_{i_j}$ 
you automatically reach all
other points $y_j\in P_{i_j}$.
It is important in the definition that each
$P_i$ is a subset of the \emph{inner} boundary $\boundary\L^\G$ and
not $\partial\L^\G$.

Here are some important examples of random-cluster boundary 
conditions.
\begin{itemize}
\item If $b=\{\boundary\L^\G\}$ then the entire boundary 
$\boundary\L$ is
wired together;  we call this the \emph{wired} boundary condition
and denote it by $b=\wired$;%
\nomenclature[w]{$\wired$}{Wired boundary condition}
\item If $b=\{\{\G\}\}$ then $x\lra y$ with respect to $b$
if and only if there is an open path between $x,y$ inside $\L$;
we call this the \emph{free} boundary condition, and denote it
by $b=\free$.%
\nomenclature[f]{$\free$}{Free boundary condition}
\item Given any $\tau\in\Om$, the boundary condition $b=\tau$
is by definition obtained by letting the $P_i$ consist of
those points in $\boundary\L^\G$ which are connected by open paths of 
$\tau$ \emph{outside} $\L$.
\item
We may also impose a number of \emph{periodic} boundary conditions
on simple regions.  One
may then regard $[-\b/2,\b/2]$ as a circle by identifying its endpoints,
and/or in the case $L=[-n,n]^d$ identify the latter with the 
torus $(\ZZ/[-n,n])^d$.  Notation for periodic boundary conditions
will be introduced when necessary.  Periodic boundary conditions will
be particularly important in the study of the quantum Ising model
in Chapter~\ref{qim_ch}.
\end{itemize}

For each boundary condition $b$ on $\L$, define the function 
$k^b_\L:\Om\rightarrow\{1,2,\dotsc,\oo\}$%
\nomenclature[k]{$k^b_\L$}{Number of connected components}
 to count the number
of components of $\om$ in $\L$, 
counted with respect to the boundary condition $b$.
There is a natural partial order on boundary conditions given by:
$b'\geq b$ if $k^{b'}_\L(\om)\leq k^{b}_\L(\om)$ for all $\om\in\Om$.
\label{bc_po}
For example, for any boundary condition $b$ we have 
$k^\wired_\L\leq k^b_\L\leq k^\free_\L$ and hence
$\wired\geq b\geq\free$.  
(Alternatively, $b'\geq b$ if $b$ is a refinement of $b'$.
Note that for $b=\tau\in\Om$,
this partial order agrees with the natural partial order on $\Om$,
defined in Section~\ref{stoch_ineq_sect}.)

\subsection{The space--time random-cluster model}\label{strcm_subsec}

For $q>0$ and $b$ a boundary condition, define the 
(random-cluster) \emph{partition functions}
\begin{equation}
Z^b_\L=Z^b_\L(\l,\d,\g,q):=\int_\Om q^{k^b_\L(\om)}\:d\mu(\om).
\end{equation}%
\nomenclature[Z]{$Z^b_\L$}{Random-cluster model partition function}%
It is not hard to see that each $Z^b_\L<\infty$.
\begin{definition}\label{rc_def}
We define the \emph{finite-volume random-cluster measure}
$\phi^b_\L=\phi^b_{\L;q,\l,\d,\g}$%
\nomenclature[f]{$\phi$}{Random-cluster measure}
 on $\L$
to be the probability measure on
$(\Om,\cF)$ given by
\[
\frac{d\phi^b_\L}{d\mu}(\om):=\frac{q^{k^b_\L(\om)}}{Z^b_\L}.
\]
\end{definition}
Thus, for any bounded, $\cF$-measurable $f:\Om\rightarrow\RR$ we have that
\begin{equation}
\phi^b_\L(f)=\frac{1}{Z^b_\L}\int_\Om f(\om)q^{k^b_\L(\om)}
\, d\mu(\om).
\end{equation}

We say that an event $A\in\cF$ is \emph{defined} on a pair $(S,T)$ 
of subsets $S\subseteq\KK$ and $T\se\FF$ 
if whenever $\om\in A$, and $\om'\in\Om$ is such that
$B(\om)\cap T=B(\om')\cap T$, $D(\om)\cap S=D(\om')\cap S$
and $G(\om)\cap S=G(\om')\cap S$, then also $\om'\in A$.  Let $\cF_{(S,T)}\subseteq\cF$ be the $\s$-algebra of
events defined on $(S,T)$.  For $\L=(K,F)$ a region we write $\cF_{\L}$%
\nomenclature[F]{$\cF_\L$}{Restricted $\s$-algebra}
 for $\cF_{(K,F)}$;  we abbreviate $\cF_{(S,\varnothing)}$ and
$\cF_{(\varnothing,T)}$ by $\cF_S$ and $\cF_T$, respectively.
Let $\cT_{(S,T)}=\cF_{(\KK\sm S,\FF\sm T)}$%
\nomenclature[T]{$\cT_\L$}{Events defined outside $\L$}
 denote the
$\s$-algebra of  events defined \emph{outside} $S$ and $T$.    
We call $A\in\cF$ a \emph{local}
event if there is a region $\L$ such that $A\in\cF_\L$
(this is sometimes also called a \emph{finite-volume} event
or a \emph{cylinder} event).

Note that the version of $d\phi^b_\L/d\mu$ given in Definition~\ref{rc_def}
is $\cF_\L$-measurable;  thus we may either regard $\phi^b_\L$ as a measure
on the full space $(\Om,\cF)$, or, by restricting consideration to events
in $\cF_\L$, as a measure on $(\Om,\cF_\L)$.

For $\D=(K,F)$ a region and $\om,\tau\in\Om$, let 
\[
\begin{split}
B_\D(\om,\tau)&=(B(\om)\cap F)\cup (B(\tau)\cap (\FF\sm F)),\\
D_\D(\om,\tau)&=(D(\om)\cap K)\cup (D(\tau)\cap (\KK\sm K)),\\
G_\D(\om,\tau)&=(G(\om)\cap K)\cup (G(\tau)\cap (\KK\sm K)).
\end{split}
\]
We write
\[
(\om,\tau)_\D=(B_\D(\om,\tau),D_\D(\om,\tau),G_\D(\om,\tau))
\]
for the configuration that agrees with $\om$ in $\D$ and with
$\tau$ outside $\D$.
The following result is a very useful `spatial Markov' property
of random-cluster measures;  it is sometimes referred to as
the \dlr-, or Gibbs-, property.
The proof follows standard arguments and
may be found in Appendix~\ref{pfs_app}.

\begin{proposition}\label{cond_meas_rcm}
Let $\L\subseteq\D$ be regions, $\tau\in\Om$, and
$A\in\cF$.  Then 
\[
  \phi^\tau_\D(A\mid\cT_\L)(\om)=\phi_\L^{(\om,\tau)_\D}(A),
  \qquad \phi^\tau_\D\mbox{-a.s.}
\]
\end{proposition}

Analogous results hold for $b\in\{\free,\wired\}$.
The following is an immediate consequence of 
Proposition~\ref{cond_meas_rcm}. 

\begin{corollary}[Deletion-contraction property]\label{del_contr}
Let $\L\subseteq\D$ be regions such that
$\boundary\L\cap\boundary\D=\varnothing$, and let $b$ be a boundary
condition on $\D$.
Let $\cC$ be the event that all components 
inside $\L$ which intersect $\boundary\L$ are connected in $\D\sm\L$; 
let $\cD$ be the
event that \emph{none} of these components are connected in $\D\sm\L$.
Then 
\[
\phi^b_\D(\cdot\mid \cC)=\phi_\L^\wired(\cdot)\quad\text{and}\quad
\phi^b_\D(\cdot\mid \cD)=\phi_\L^\free(\cdot).
\]
\end{corollary}

\subsection{The space--time Potts model}\label{stp_sec}

The classical random-cluster model 
is closely related to the Potts model of statistical mechanics.
Similarly
there is a natural `space--time Potts model' which may be coupled
with the space--time random-cluster model.  
A realization of the space--time Potts measure is a piecewise constant
`colouring' of $\KK^\G$.  As for the random-cluster model, we will
be interested in specifying different boundary conditions, and these
will not only tell us which parts of the boundary are `tied together',
but may also specify the precise colour on certain parts of the boundary.

Let us fix a region $\L$ and $q\geq 2$ an \emph{integer}.
Let $\cN=\cN_q$%
\nomenclature[N]{$\cN$}{Potts model configuration space} 
 be the set of functions $\nu:\KK^\G\rightarrow\{1,\dotsc,q\}$%
\nomenclature[n]{$\nu$}{Potts configuration}
 which  have the property that
their restriction to any $v\times\RR$ is piecewise
constant and right-continuous. 
Let $\cG$%
\nomenclature[G]{$\cG$}{$\s$-algebra for the Potts model}
 be the $\s$-algebra on $\cN$ generated by all the functions
$\nu\mapsto(\nu(x_1),\dotsc,\nu(x_N))\in\RR^N$ as $N$ ranges through
the integers and $x_1,\dotsc,x_N$ range through $\KK^\G$
(this coincides with the $\s$-algebra generated by the Skorokhod
metric, see Appendix~\ref{skor_app} 
and~\cite[Proposition~3.7.1]{ethier_kurtz}).
For $S\subseteq\KK$ define the $\s$-algebra 
$\cG_S\subseteq\cG$ of events defined on $S^\G$.
Although we canonically let $\nu\in\cN$ be right-continuous, 
we will usually identify such $\nu$ which agree off sets of Lebesgue
measure zero, compare Remark~\ref{rem-as}.  Thus we will
without further mention allow $\nu$ to be any 
piecewise constant function with values
in $\{1,\dotsc,q\}$, and we will frequently even allow $\nu$ to be
undefined on a set of  measure zero.  
We call elements of $\cN$ `spin configurations' and will
usually write $\nu_x$ for $\nu(x)$.

Let $b=\{P_1,\dotsc,P_m\}$ be any random-cluster boundary condition and
let $\a:\{1,\dotsc,m\}\rightarrow\{0,1,\dotsc,q\}$.%
\nomenclature[a]{$\a$}{Part of Potts boundary condition}
 We call the pair $(b,\a)$ a \emph{Potts boundary condition}. 
We assume that
$\G\in P_1$, and write $\a_\G$ for $\a(1)$;  we also require that
$\a_\G\neq 0$.
Let $D\se K$ be a finite set, and let $\cN^{b,\a}_\L(D)$%
\nomenclature[N]{$\cN(D)$}{Potts configurations permitted by $D$}
 be the set of $\nu\in\cN$ with the following properties.
\begin{itemize}
\item For each $v\in V$ and each interval $I\se K_v$ such that
$I\cap D=\varnothing$, $\nu$ is constant on $I$,
\item if $i\in\{1,\dotsc,m\}$ is such that $\a(i)\neq 0$ then
$\nu_x=\a(i)$ for all $x\in P_i$,
\item if $i\in\{1,\dotsc,m\}$ is such that $\a(i)=0$ and $x,y\in P_i$
then $\nu_x=\nu_y$,
\item if $x\not\in\L$ then $\nu_x=\a_\G$.
\end{itemize}
Intuitively, the boundary condition $b$ specifies which parts of the 
boundary are forced to have the same spin, and the function $\a$ specifies
the \emph{value} of the spin on some parts of the boundary;  
$\a(i)=0$ is taken to mean that the  value on $P_i$ is not
specified.  (The value of $\a$ at $\G$ is special, in that it takes on
the role of an external field, see~\eqref{st_potts_def_eq}.)

Let $\l:\FF\rightarrow\RR$, $\g:\KK\rightarrow\RR$ and 
$\d:\KK\rightarrow\RR_+$ be bounded and Borel-measurable;  note that
$\l$ and $\g$ are allowed to take negative values. 
For $a,b\in\RR$, let $\d_{a,b}=\one_{\{a=b\}}$, and for $\nu\in\cN$
and $e=xy\in\EE$, let $\d_\nu(e)=\d_{\nu_x,\nu_y}$.
Let $\pi_\L^{b,\a}$%
\nomenclature[p]{$\pi$}{Potts measure}
 denote the probability measure on $(\cN,\cG)$ defined
by, for each bounded and $\cG$-measurable $f:\cN\rightarrow\RR$, 
letting  $\pi^{b,\a}_\L(f(\nu))$ be a constant multiple of
\begin{equation}\label{st_potts_def_eq}
\int d\mu_\d(D)\,\sum_{\nu\in\cN^{b,\a}_\L(D)} f(\nu)
\exp\Big(\int_F\l(e)\d_\nu(e)de+
\int_K\g(x) \d_{\nu_x,\a_\G}dx\Big)
\end{equation}
(with constant determined by the requirement that $\pi_\L^{b,\a}$
be a probability measure).
The integrals 
in~\eqref{st_potts_def_eq} are to be interpreted as in~\eqref{KF_int}.
\begin{definition}\label{potts_def}
The probability
measure $\pi^{b,\a}_\L=\pi^{b,\a}_{\L;q,\l,\g,\d}$ on $(\cN,\cG)$ 
defined by~\eqref{st_potts_def_eq} is called 
the space--time Potts measure with $q$ states on~$\L$.
\end{definition}

Note that, as with $\phi^b_\L$, we may regard 
$\pi^{b,\a}_{\L}$ as a measure on $(\cN,\cG_\L)$.
Here is a word of motivation for \eqref{st_potts_def_eq}
in the case $b=\free$ and $\a_\G=q$;  similar constructions
hold for other $b,\a$.
See Figure~\ref{bridge_part_fig} in Section~\ref{ssec-rpr},
and also~\cite{GOS}.  The set $(v\times K_v)\sm\oD$ is
a union of maximal death-free intervals $v\times J_v^k$,
where $k=1,2,\dotsc,n$ and $n=n(v,\oD)$ is the number of such intervals.  
We write $V(\oD)$ for the collection of all such intervals
as $v$ ranges over $V$, together with the ghost-vertex
$\Gh$, to which we assign spin $\nu_\Gh=q$.  The set $\cN^{\free,\a}_\L(\oD)$ 
may be identified with $\{1,\dotsc,q\}^{V(\oD)}$, and
we may think of $V(D)$ as the set of vertices of a graph with edges
given as follows. 
An edge is placed between $\Gh$ and each $\bar v\in V(D)$. 
For $\bar u, \bar v\in V(D)$,
with $\bar u=u\times I_1$ and $\bar v=v\times I_2$ say, we place an edge 
between $\bar u$ and $\bar v$ if and only if: (i) $uv$ is an edge of $L$, and
(ii) $I_1\cap I_2\neq\es$.  Under the space--time Potts measure
\emph{conditioned on $D$}, a spin-configuration
$\nu\in \cN^{\free,\a}_\L(D)$ on this graph
receives a (classical) Potts weight
\begin{equation}
\exp\left\{\sum_{\bar u\bar v}J_{\bar u\bar v}\d_{\nu}(\bar u\bar v)
+ \sum_{\bar v} h_{\bar v}\d_{\nu_{\bar v},q}\right\},
\end{equation}
where $\nu_{\bar v}$ denotes the common value of $\nu$ along $\bar v$, 
and where
$$
J_{\bar u\bar v}=\int_{I_1\cap I_2}\l(uv,t)\,dt\qquad\text{and}
\qquad h_{\bar v}=\int_{\bar v}\g(x)\,dx.
$$
This observation will be pursued further for the
Ising model in Section~\ref{ssec-rpr}.

The space--time Potts measure may, for special boundary conditions,
be coupled to the space--time random-cluster measure, as follows.
For $\a$ of the form $(\a_\G,0,\dotsc,0)$, we call $(b,\a)$ a
\emph{simple} Potts boundary condition.  Thus, under a simple boundary
condition, the only spin value which is specified in advance 
is that of $\G$.
Let $\om=(B,D,G)\in\Om$ be sampled from $\phi^b_\L$ and   
write $\cN^{b,\a}_\L(\om)$ 
for the set of $\nu\in\cN$ such that (i) $\nu_x=\a_\G$ for $x\not\in\L$, and
(ii) if $x,y\in\L$ and $x\lra y$ in $\om$ under the boundary 
condition $b$ in $\L$ then $\nu_x=\nu_y$.  In particular, since
$\G\not\in\L$ we have that $\nu_\G=\a_\G$.
Note that each $\cN^{b,\a}_\L(\om)$ is a finite set.
With $\om$ given, we sample $\nu\in\cN^{b,\a}_\L(\om)$ as follows.
Set $\nu_\G:=\a_\G$ and set $\nu_x=\a_\G$ for all $x\not\in\L^\G$;  then
choose the spins of the other components of $\om$ in $\L$ uniformly 
and independently at random.  The resulting pair $(\om,\nu)$ has 
a distribution $\PP^{b,\a}_\L$ on $(\Om,\cF)\times(\cN,\cG)$%
\nomenclature[P]{$\PP$}{Edwards--Sokal coupling}
 given by
\begin{equation}\label{es_def}
\begin{split}
\PP^{b,\a}_\L(f(\om,\nu))&= \int_\Om d\phi^b_\L(\om) 
\,\frac{1}{q^{k^b_\L(\om)-1}}\sum_{\nu\in\cN^{b,\a}_\L(\om)} f(\om,\nu)\\
&\propto \int_\Om d\mu(\om) \,\sum_{\nu\in\cN^{b,\a}_\L(\om)} f(\om,\nu),
\end{split}
\end{equation}
for all bounded $f:\Om\times\cN\rightarrow\RR$, measurable in the product
$\s$-algebra $\cF\times\cG$.  We call
the measure $\PP^{b,\a}_\L$ of~\eqref{es_def}  the Edwards--Sokal
measure.  This definition is completely analogous to a coupling 
in the discrete model, which was  was found in~\cite{edwards_sokal}.
Usually we take  $\a_\G=q$ and in this case 
we will often suppress reference to $\a$,
writing for example $\cN^b_\L(\om)$ and similarly for other
notation.

The marginal of $\PP^{b,\a}_\L$ on $(\cN,\cG)$ is computed
as follows.  Assume that 
$f(\om,\nu)\equiv f(\nu)$ depends only on $\nu$,
and let $D\subseteq K$ be a finite set.
For $\nu\in\cN^{b,\a}_\L(D),$ let  
$\{\nu\sim\om\}$ be the event that $\om$ has no
open paths \emph{inside} $\L$ that violate the condition that
$\nu$ be constant on the components of $\om$.
We may rewrite~\eqref{es_def} as 
\begin{equation}
\PP^{b,\a}_\L(f(\nu))
\propto \int d\mu_\d(D)\int d(\mu_\l\times\mu_\g)(B,G)
\sum_{\nu\in\cN^{b,\a}_\L(D)} f(\nu)\one\{\nu\sim\om\}.
\end{equation}
With $D$ fixed, the probability under $\mu_\l\times\mu_\g$ of 
the event $\{\nu\sim\om\}$ is
\begin{equation}\label{i5}
\exp\Big(-\int_F\l(e)(1-\d_\nu(e))de
-\int_K\g(x)(1-\d_{\nu_x,\a_\G})dx\Big).
\end{equation}
Taking out a constant, it follows that $\PP^{b,\a}_\L(f(\nu))$
is proportional to
\begin{align}
\int d\mu_\d(D)\,\sum_{\nu\in\cN^{b,\a}_\L(D)} f(\nu)
\exp\Big(\int_F\l(e)\d_\nu(e)de+
\int_K\g(x) \d_{\nu_x,\a_\G}dx\Big).
\end{align}
Comparing this with~\eqref{st_potts_def_eq}, and noting that both
equations define probability measures, it follows that
$\PP^{b,\a}_\L(f(\nu))=\pi^{b,\a}_\L(f)$.

We may ask for a description of how to obtain an $\om$ with law
$\phi^b_\L$ from a $\nu$ with law $\pi^{b,\a}_\L$.
In analogy with the discrete case this is as follows:
\begin{quote}
Given $\nu\sim\pi^{b,\a}_\L(\cdot)$, place a death wherever 
$\nu$ changes spin in $\L$, and also place additional deaths elsewhere
in $\L$ at rate $\d$;  place bridges between
intervals in $\L$ of the same spin at rate $\l$;  
and place ghost-bonds in intervals in $\L$
of spin $\a$ at rate $\g$.  The outcome $\om$ has law 
$\phi^b_\L(\cdot)$.
\end{quote}
It follows that we have the following correspondence between 
$\phi=\phi^b_\L$ and $\pi=\pi^{b,\a}_{\L,q}$ when $(b,\a)$ is  simple.
The result is completely 
analogous to the corresponding result for the discrete Potts model
(Theorem~\ref{thm1}), and 
the proof is included only for completeness.  
\begin{proposition}\label{corr_conn_prop}
Let $x,y\in\L^\G$.  Then
\[
\pi(\nu_x=\nu_y)=\Big(1-\frac{1}{q}\Big)\phi(x\leftrightarrow y)+\frac{1}{q}.
\]
\end{proposition}
\begin{proof}
Writing $\PP$ for the Edwards--Sokal coupling, we have that
\begin{align*}
q\pi(\nu_x=\nu_y)-1&=\PP(q\cdot\PP(\nu_x=\nu_y\mid\om)-1)\\
&=\PP\Big(q\big(\one\{x\lra y\mbox{ in }\om\}+
\frac{1}{q}\one\{x\not\lra y\mbox{ in }\om\}\big)-1\Big)\\
&=\PP((q-1)\cdot\one\{x\leftrightarrow y\mbox{ in }\om\})\\
&=(q-1)\phi(x\leftrightarrow y).
\end{align*}
\end{proof}

The case $q=2$ merits special attention.  In this case it is customary to 
replace the states $\nu_x=1,2$ by $-1,+1$ respectively, and we
thus define $\s_x=2\nu_x-3$.  For $\a$ taking values in $\{0,-1,+1\}$,
we let $\S,\S^{b,\a}_\L(\om),\S^{b,\a}_\L(D)$%
\nomenclature[S]{$\S$}{Ising configuration space}%
\nomenclature[S]{$\S(D)$}{Ising configurations permitted by $D$}
 denote the
images of $\cN,\cN^{b,\a}_\L(\om),\cN^{b,\a}_\L(D)$ respectively under the map 
$\nu\mapsto\s$.  Reference to $\a$ may be suppressed if $(b,\a)$ is 
simple and $\a_\G=+1$.   

We have that
\begin{equation}\label{sigma_nu}
\one\{\s_x=\s_y\}=\frac{1}{2}(\s_x\s_y+1),\qquad 
\one\{\s_x=\a_\G\}=\frac{1}{2}(\a_\G\s_x+1).
\end{equation}
Consequently, $\pi^{b,\a}_{\L;q=2}(f(\s))$ is proportional to
\begin{equation}\label{ising_def_eq}
\int d\mu_\d(D)\,\sum_{\s\in\S^{b,\a}_\L(D)} f(\s)
\exp\Big(\frac{1}{2}\int_F\l(e)\s_e\,de+
\frac{1}{2}\int_K\g(x)\a_\G\s_x \,dx\Big),
\end{equation}
where we have written $\s_e$ for $\s_x\s_y$ when $e=xy$.  In this formulation,
we call the measure of~\eqref{ising_def_eq} 
the \emph{Ising measure}.  Expected values
with respect to this measure will typically be written 
$\el\cdot\er^{b,\a}_\L$;%
\nomenclature[<]{$\el\cdot\er$}{Expectation under Ising measure}
  thus for
example Proposition~\ref{corr_conn_prop} says that when $q=2$ 
and $(b,\a)$ is simple, then
\begin{equation}\label{ising_corr_conn}
\el\s_x\s_y\er^{b,\a}_\L=\phi^b_\L(x\lra y).
\end{equation}

For later reference, we make a note here of the constants of proportionality
in the above definitions.  Let 
\begin{equation}
Z^b_{\mathrm{RC}}=Z^b_{\mathrm{RC}}(q)=\int_{\Om} q^{k^b_\L(\om)}\,d\mu(\om)
\end{equation}
denote the partition function of the random-cluster model, and 
\begin{equation}
Z_{\mathrm{Potts}}^{b,\a}(q)=\int d\mu_\d(D)\,\sum_{\nu\in\cN^{b,\a}_\L(D)} 
\exp\Big(\int_F \d_\nu(e)\l(e)\,de+\int_K \d_{\nu_x,\a_\G}\g(x)\, dx\Big)
\end{equation}
that of the $q$-state Potts model.  Also, let 
\begin{equation}\label{ising_pf}
Z^{b,\a}_{\mathrm{Ising}}=\int d\mu_\d(D)\,\sum_{\s\in\S^{b,\a}_\L(D)} 
\exp\Big(\frac{1}{2}\int_F\l(e)\s_e\,de+
\frac{1}{2}\int_K\g(x)\a_\G\s_x\, dx\Big)
\end{equation}
be the partition function of the Ising model.  By keeping track of the 
constants in the above calculations we obtain the following result,
which for simplicity is stated only for $\a_\G=q$.
\begin{proposition}\label{pfs}  Let $b$ be a random-cluster boundary
condition.  Then
\begin{align}
Z^{b}_{\mathrm{Potts}}(q)&
=\frac{1}{q}Z^b_{\mathrm{RC}}(q)\cdot
\exp\Big(\int_F\l(e)\,de+\int_K\g(x)\,dx\Big)\\
Z^{b}_{\mathrm{Ising}}&=Z^{b}_{\mathrm{Potts}}(2)\cdot
\exp\Big(-\frac{1}{2}\int_F\l(e)\,de-\frac{1}{2}\int_K\g(x)\,dx\Big)\\
&=\frac{1}{2}Z^b_{\mathrm{RC}}(2)\cdot
\exp\Big(\frac{1}{2}\int_F\l(e)\,de+\frac{1}{2}\int_K\g(x)\,dx\Big).\nonumber
\end{align}
\end{proposition}

It is easy to check, by a direct computation, that the Potts model
behaves in a similar manner to the random-cluster model upon 
conditioning on the value of $\nu$ in part of a region, i.e. that
analogs of Proposition~\ref{cond_meas_rcm} and 
Corollary~\ref{del_contr} hold.  We will not state these results explicitly
in full generality, but will record here the following special case
for later reference.
\begin{lemma}\label{potts_cond}
Let $\L\subseteq\D$ denote two regions, and consider the
boundary condition $(\wired,\a)$.
Then for all $\cG_\L$-measurable $f$ we have that
\[
\pi_\L^{\wired,\a}(f(\nu))=
\pi_\D^{\wired,\a}(f(\nu)\mid\s\equiv \a_\G\mbox{ on }\D\sm\L).
\]
\end{lemma}

\section{Stochastic comparison}\label{stoch_ineq_sect}

The ability to compare the probabilities of events under a range
of different measures is extremely important in the theory of 
random-cluster measures.  In this section we develop in detail the
basis for such a methodology in the space--time setting.  We also 
prove versions of the {\gks-} and {\fkg} inequalities suitable
for the space--time Potts and Ising measures, respectively.

Let $\L$ be a  region.
Let the pair $(E,\cE)$ denote one of $(\Om,\cF)$, $(\Om,\cF_\L)$,
$(\S,\cG)$ and $(\S,\cG_\L)$.  Thus $E$, equipped with the
Skorokhod metric, is a Polish metric space.
Given a partial
order $\geq$ on $E$, a measurable function $f:E\rightarrow\RR$
is called \emph{increasing} if for all $\om,\xi\in E$ such
that $\om\geq\xi$ we have $f(\om)\geq f(\xi)$.  An event
$A\in\cE$ is increasing if the indicator function $\one_A$ is.
We assume that the set $\{(\om,\xi)\in E^2:\om\geq\xi\}$ is closed
in the product topology;  this will hold automatically in our
applications.

Let $\psi_1,\psi_2$ be two probability measures
on $(E,\cE)$.
\begin{definition}\label{stoch_dom_def} 
We say that $\psi_1$
\emph{stochastically dominates} $\psi_2$, and we write $\psi_1\geq\psi_2$,
if $\psi_1(f)\geq\psi_2(f)$ for all bounded, increasing local functions $f$.
\end{definition}
By a standard approximation argument using the
monotone convergence theorem, $\psi_1\geq\psi_2$ 
holds if for all increasing local
\emph{events} $A$ we have $\psi_1(A)\geq\psi_2(A)$.

The following general result lies at the heart of stochastic
comparison and will be used repeatedly.  It goes back 
to~\cite{strassen};
see also~\cite[Theorem~IV.2.4]{lindvall} and~\cite[Theorem~4.6]{ghm}.
\begin{theorem}[Strassen]\label{strassen_thm}
Let $\psi_1,\psi_2$ be probability measures on $(E,\cE)$.
The following statements are equivalent.
\begin{enumerate}
\item $\psi_1\geq\psi_2$;
\item For all \emph{continuous} bounded increasing 
local functions
$f:E\rightarrow\RR$ we have $\psi_1(f)\geq\psi_2(f)$;
\item There exists a probability
measure $P$ on $(E^2,\cE^2)$ such that
\[
P(\{(\om_1,\om_2):\om_1\geq\om_2\})=1.
\]
\end{enumerate}
\end{theorem}
Note that the equivalence of (1) and (3) extends to \emph{countable}
sequences $\psi_1,\psi_2,\psi_3,\dotsc$;  
see~\cite[Theorem~IV.6.1]{lindvall}.

\begin{definition}\label{pos_assoc_def}
A measure $\psi$ is on $(E,\cE)$ is called
\emph{positively associated} if for all local
increasing events $A,B$ we have that 
$\psi(A\cap B)\geq\psi(A)\psi(B)$.
\end{definition}
The inequality $\psi(A\cap B)\geq\psi(A)\psi(B)$ for local increasing
events is sometimes referred to as the \fkg-inequality as the systematic
study of such inequalities was initiated by Fortuin, Kasteleyn
and Ginibre~\cite{fkg}.

\subsection{Stochastic inequalities for the random-cluster
model}\label{rcm_si}

The results in
this section are applications, and slight modifications, of stochastic
comparison results for point processes that appear in~\cite{preston75}
and~\cite{georgii_kuneth}.  See also~\cite[Theorem~10.4]{ghm}.
Some of the results, 
such as positive association in the space--time random-cluster
model, have been stated before, sometimes with additional assumptions;
see for example~\cite{akn,aizenman_nacht,bezuidenhout_grimmett}. 
We do not believe detailed proofs for space--time models
have appeared before.  The results
presented are satisfyingly similar to those for the discrete case,
compare~\cite[Chapter~3]{grimmett_rcm} and~\cite{grimmett_gks}.

We will follow the method of~\cite{preston75}
rather than the later (and more general)~\cite{georgii_kuneth}.
This is because the former method avoids discretization
and is closer to the standard 
approach of~\cite{holley74} (also~\cite[Chapter~2]{grimmett_rcm}) 
for the classical random-cluster model.  The method makes use
of coupled Markov chains on $\Om$ (specifically, jump-processes,
see~\cite[Chapter~X]{feller71_vol2}).

For $\om\in\Om$, write $B(\om),D(\om),G(\om)$ for the sets
of bridges, deaths and ghost-bonds in $\om$, respectively.
We define a partial order on $\Om$ by saying that
$\om\geq\xi$ if
$B(\om)\supseteq B(\xi)$, $D(\om)\subseteq D(\xi)$ and
$G(\om)\supseteq G(\xi)$.

We will in this section only consider measures on $\cF_\L$,
that is we take $(E,\cE)=(\Om,\cF_\L)$.  We
will regard $B,G,D$ as subsets of $K$ and $F$ as appropriate.
The symbol $x$ will be used to denote a generic point of
$\L\equiv K\cup F$, interpreted either as a bridge, a ghost-bond, 
or a death, as specified.  More formally, we may regard $x$
as an element of 
$F\cup(K\times\{\mathrm{d}\})\cup(K\times\{\mathrm{g}\})$, where
the labels $\mathrm{d},\mathrm{g}$ 
allow us to distinguish between deaths and 
ghost-bonds, respectively.
We let $X=(X_t:t\geq 0)$ be a continuous-time
stochastic process with state space $\Om$, defined as follows.
If $X_t=(B,G,D)$, there are 6 possible transitions.  The
process can either jump to one of
\begin{equation}\label{iml1}
(B\cup\{x\},G,D),\quad\mbox{or }
(B,G\cup\{x\},D),\quad\mbox{or }
(B,G,D\cup\{x\}),
\end{equation}
where $x\in\L$;
the corresponding move is called a \emph{birth} at $x$.
Alternatively, in the case where $x\in B$, the process can jump to
\[
(B\sm\{x\},G,D),
\]
and similarly for $x\in G$ or $x\in D$; the corresponding
move is called a \emph{demise} at $x$.
If $\om=(B,G,D)\in\Om$, we will often abuse notation and write
$\om^x$ for the configuration~\eqref{iml1}
with a point at $x$ added, making it clear
from the context whether $x$ is a bridge, ghost-bond, or death.
Similarly, if $x\in B\cup G\cup D$, we will write $\om_x$ for
the configuration with the bridge, ghost-bond or death
at $x$ removed.

The transitions described above happen at the following rates.
Let $\cL$ denote the Borel $\s$-algebra on 
$\L\equiv F\cup(K\times\{\mathrm{d}\})\cup(K\times\{\mathrm{g}\})$, 
and let $\cB:\Om\times\cL\rightarrow\RR$ be a given function,
such that for each $\om\in\Om$, $\cB(\om;\cdot)$ is a finite
measure on $(\L,\cL)$.   Also let 
$\cD:\Om\times\L\rightarrow\RR$ be such that for all $\om\in\Om$
we have that $\cD(\om;x)$ is a non-negative measurable function of $x$.
If for some $t\geq 0$ we have that $X_t=\om$,
then there is a birth in the measurable set $H\subseteq\L$ before
time $t+s$ with probability $\cB(\om; H)s+o(s)$.
Alternatively, there is a demise at the point $x\in\om$
before time $t+s$ with probability $\cD(\om_x;x)s+o(s)$.

We may give an equivalent `jump-hold' description of the chain, 
as follows.  Let
\begin{equation}\label{car1}
\cA(\om):=\cB(\om;\L)+\sum_{x\in\om}\cD(\om_x;x).
\end{equation}
For $A\in\cF_\L$ let
\begin{equation}\label{car2}
\cK(\om,A):=\frac{1}{\cA(\om)}\Big(\cB(\om;\{x\in\L:\om^x\in A\})
  +\sum_{\substack{x\in\om\\\om_x\in A}}\cD(\om_x;x)\Big).
\end{equation}
Then given that $X_t=\om$, the holding
time until the next transition has the exponential distribution 
with parameter
$\cA(\om)$;  once the process jumps it goes to some state $\xi\in A$ with
probability $\cK(\om,A)$.
Existence and basic properties of such Markov chains are discussed
in~\cite{preston75}.

We will aim to construct such chains $X$ which are in detailed
balance with a given probability measure $\psi$ on $(\Om,\cF_\L)$.
It will be necessary to make some assumptions on $\psi$, and
these will be stated when appropriate.  For now the main assumption
we make is the following.  Let $\k=\mu_{1,1,1,}$ 
denote the probability measure
on $(\Om,\cF_\L)$ given by letting $B,G,D$ all be independent Poisson
processes of constant intensity $1$.
\begin{assumption}\label{MC_ass_1}
The probability measure $\psi$ is absolutely continuous
with respect to $\k$; there exists a version of the density
\[
f=\frac{d\psi}{d\k}
\]
which has the property that for all $\om\in\Om$ and $x\in\L$,
if $f(\om)=0$ then $f(\om^x)=0$.
\end{assumption}
\begin{example}
The space--time percolation measures 
(restricted to $\L$)
satisfy Assumption~\ref{MC_ass_1}, 
because by standard properties of Poisson processes, 
if  $\mu=\mu_{\l,\d,\g}$ then a version of the density is given by
\begin{equation}
\frac{d\mu}{d\k}(\om)\propto \prod_{x\in B}\l(x)
\prod_{y\in D}\d(y)\prod_{z\in G}\g(z).
\end{equation}
Moreover, the random-cluster measure
$\phi^b_\L=\phi^b_{\L;q,\l,\d,\g}$ also satisfies
Assumption~\ref{MC_ass_1}, having density 
\begin{equation}
\frac{d\phi^b_\L}{d\k}(\om)=
\frac{d\phi^b_\L}{d\mu}(\om)\frac{d\mu}{d\k}(\om) 
\propto q^{k^b_\L(\om)}\prod_{x\in B}\l(x)
\prod_{y\in D}\d(y)\prod_{z\in G}\g(z)
\end{equation}
against $\k$.  
\end{example}

\begin{definition}
The \emph{Papangelou intensity} of $\psi$ is the function
$\i:\Om\times\L\rightarrow\RR$ given by
\begin{equation}
\i(\om,x)=\frac{f(\om^x)}{f(\om)}
\end{equation}
(where we take $0/0$ to be 0).
\end{definition}

The following construction will not itself be used, but
serves as a helpful illustration.  To construct a birth-and-death chain 
which has equilibrium distribution $\psi$ we
would simply take $\cD\equiv 1$ and $\cB(\om;dx)=\i(\om,x)dx$.
(Here $dx$ denotes Lebesgue measure on 
$F\cup(K\times\{\mathrm{d}\})\cup(K\times\{\mathrm{g}\})$.)
The corresponding chain $X$ is in
detailed balance with $\psi$, since 
$d\psi(\om_x)\cdot\cB(\om_x; dx)=d\k(\om_x)f(\om^x)dx=d\psi(\om^x)\cdot 1$.
In light of this
one may may think of $\i(\om,x)$ as the intensity with which
the chain $X$, in equilibrium with $\psi$, attracts a birth at~$x$.

\begin{example}
For the random-cluster measure $\phi^b_\L$, 
\begin{equation}\label{rcm_intensity}
\i(\om,x)=
q^{k^b_\L(\om^x)-k^b_\L(\om)}\cdot
\left\{\begin{array}{ll}
\l(x), & \mbox{for $x$ a bridge} \\
\d(x), & \mbox{for $x$ a death} \\
\g(x), & \mbox{for $x$ a ghost-bond}.
\end{array}\right.
\end{equation}
\end{example}

In the rest of this section we let 
$\psi,\psi_1,\psi_2$ be three probability measures satisfying
Assumption~\ref{MC_ass_1}, and let $f,f_1,f_2$ and $\i,\i_1,\i_2$
denote their density functions against $\k$ and their
Papangelou intensities, respectively.

\begin{definition}\label{preston_lattice}
We say that the pair $(\psi_1,\psi_2)$ 
satisfies the \emph{lattice condition}
if the following hold whenever $\om\geq\xi$:
\begin{enumerate}
\item $\i_1(\om,x)\geq \i_2(\xi,x)$
whenever $x$ is a bridge or ghost-bond such that $\xi^x\not\leq\om$;
\item $\i_2(\xi,x)\geq \i_1(\om,x)$
whenever $x$ is a death such that $\xi\not\leq\om^x$.
\end{enumerate}
We say that $\psi$ has the \emph{lattice property} if
the following hold whenever $\om\geq\xi$:
\begin{enumerate}[resume]
\item $\i(\om,x)\geq \i(\xi,x)$ whenever
$x$ is a bridge or ghost-bond such that $\xi^x\not\leq\om$;
\item $\i(\xi,x)\geq \i(\om,x)$ whenever
$x$ is a death such that $\xi\not\leq\om^x$.
\end{enumerate}
\end{definition}

(We use the term `lattice' in the above definition 
in the same sense as~\cite{fkg}; `lattice' is the name for 
any partially ordered set in which
any two elements have a least upper bound and greatest lower bound.)

The next result states that `well-behaved'
measures $\psi_1,\psi_2$ which satisfy the lattice condition
are stochastically ordered, in that $\psi_1\geq\psi_2$.
Intuitively, the lattice condition implies that a chain
with equilibrium distribution $\psi_1$ acquires bridges and 
ghost-bonds faster than, but deaths slower than, the 
chain corresponding to $\psi_2$.  Similarly, we will see that 
measures with the lattice property are positively associated;
a similar intuition holds in this case.

\begin{theorem}\label{preston_thm_1}
Suppose $\psi_1,\psi_2$ satisfy the lattice condition, and 
that the Papangelou intensities $\i_1,\i_2$ are 
bounded.  Then $\psi_1\geq\psi_2$.
\end{theorem}
\begin{theorem}\label{FKG}
Suppose $\psi$ has the lattice property, and that 
$\i$ is bounded.  Then $\psi$ is positively associated.
\end{theorem}
\begin{proof}[Sketch proof of Theorem~\ref{preston_thm_1}]
This essentially follows from \cite{preston75}, the main difference 
being that our order on $\Om$ is different, in that 
`deaths count negative'.  The method of~\cite{preston75}
is to couple two jump-processes $X$ and $Y$, which
have the respective equilibrium distributions $\psi_1$
and $\psi_2$.  One may define a jump process on the product space
$\Om\times\Om$ in the same way as described in~\eqref{car1}
and~\eqref{car2};  here is the specific instance we require.

Let $T:=\{(\om,\xi)\in\Om^2:\om\geq\xi\}$, and for $a,b\in\RR$
write $a\vee b$ and $a\wedge b$ for the maximum and minimum
of $a$ and $b$, respectively.
We let $Z=(X,Y)$ be the birth-and-death process on $T$
started at $(\varnothing,\varnothing)$ and given by
the $\cA$ and $\cK$ defined below.  First,
\begin{multline}
\cA(\om,\xi):=\int_\L (\i_1(\om,x)\vee \i_2(\xi,x))\:dx+\\
+(|B(\om)|\vee|B(\xi)|)+(|D(\om)|\vee|D(\xi)|)+(|G(\om)|\vee|G(\xi)|).
\end{multline}
Write $\om\cap\xi$ for the element
$(B(\om)\cap B(\xi),D(\om)\cap D(\xi),G(\om)\cap G(\xi))$
of $\Om$;  similarly let
$\om\setminus\xi=(B(\om)\sm B(\xi),D(\om)\sm D(\xi),G(\om)\sm G(\xi))$.
For $A\subseteq T$ measurable in the product topology, let
\begin{equation}
\cK(\om,\xi;A):=\frac{1}{\cA(\om,\xi)}
\big(\cK_{\mathrm{b}}(\om,\xi;A)+\cK_{\mathrm{d}}(\om,\xi;A)\big)
\end{equation}
where
\begin{multline}
\cK_{\mathrm{d}}(\om,\xi;A):=|\{x\in\om\cap\xi:(\om_x,\xi_x)\in A\}|+\\
+|\{x\in\om\sm\xi:(\om_x,\xi)\in A\}|
+|\{x\in\xi\sm\om:(\om,\xi_x)\in A\}|
\end{multline}
and
\begin{multline}
\cK_{\mathrm{b}}(\om,\xi;A):=
\int_\L \one_A(\om^x,\xi^x)(\i_1(\om,x)\wedge \i_2(\xi,x))\:dx+\\
+\int_\L \one_A(\om^x,\xi)[\i_1(\om,x)-(\i_1(\om,x)\wedge \i_2(\xi,x))]\:dx+\\
+\int_\L \one_A(\om,\xi^x)[\i_2(\xi,x)-(\i_1(\om,x)\wedge \i_2(\xi,x))]\:dx.
\end{multline}
Thanks to the lattice condition, $Z$ is indeed a 
process on $T$.  In other words, if
$\om\geq\xi$ then $\cK(\om,\xi;T)=1$.  It is also not hard to see that $X$ 
and $Y$ are birth-and-death processes on $\Om$ with transition intensities
$\cB_1,\cD_1$ and $\cB_2,\cD_2$ respectively, where
$\cD_k\equiv 1$ and $\cB_k(\om;dx)=\i_k(\om,x)dx$, for $k=1,2$.

Define, for $n\geq 0$ and $k\in\{1,2\}$,
\begin{equation}
\cB^{(n)}_k=\sup_{|\om|=n}\cB_k(\om;\L),
\end{equation}
where $|\om|$ is the total number of bridges, ghost-bonds and 
deaths in $\om$.
The boundedness of $\i_1,\i_2$ ensures that the following 
properties, which appear as conditions in~\cite{preston75}, hold.
First, the expectation
\begin{equation}\label{iml2}
\k(\cB_k(\cdot;\L))<\oo,
\end{equation}
and second,
\begin{equation}\label{iml3}
\sum_{n=1}^\infty\frac{\cB_k^{(0)}\dotsb\cB_k^{(n-1)}}{n!}
<\infty.
\end{equation}
Theorems 7.1 and 8.1 of~\cite{preston75} therefore combine to give
that the chain $Z$ has a unique invariant distribution $P$ such that
$Z_t\Rightarrow P$, and such that $P(F\times\Om)=\psi_1(F)$ and 
$P(\Om\times F)=\psi_2(F)$.  Since $P(T)=1$, the result follows:
if $A\in\cF_\L$ is increasing then
\begin{equation}
\psi_1(A)=P(\om\in A,\,\om\geq\xi)\geq
P(\xi\in A,\,\om\geq\xi)=\psi_2(A).
\end{equation}
\end{proof}

\begin{remark}
The two technical properties~\eqref{iml2} and~\eqref{iml3}
are not strictly necessary for the main results of~\cite{preston75},
as shown in~\cite{georgii_kuneth},
but they do seem necessary for the proof method in~\cite{preston75}.  
See~\cite[Remark~1.6]{georgii_kuneth}.
\end{remark}

Theorem~\ref{FKG} follows from Theorem~\ref{preston_thm_1} 
using the following standard argument~\cite{holley74}.

\begin{proof}[Proof of Theorem~\ref{FKG}]
Let $g,h$ be two bounded, increasing and $\cF_\L$-measurable 
functions.  By adding
constants, if necessary, we may assume that
$g,h$ are strictly positive.  Let $\psi_2=\psi$ 
and let $\psi_1$ be given by
\begin{equation}
f_1(\om)=\frac{d\psi_1}{d\k}(\om):=\frac{h(\om)f(\om)}{\psi(h)}.
\end{equation}
We have that
\begin{equation}
\i_1(\om,x)=\frac{h(\om^x)f(\om^x)}{h(\om)f(\om)},
\qquad\i_2(\xi,x)=\frac{f(\xi^x)}{f(\xi)}.
\end{equation}
Clearly $\i_1,\i_2$ are uniformly bounded;
we check that $\psi_1,\psi_2$ satisfy the lattice condition.
Let $\om\geq\xi$.  If $x$ is a bridge or a ghost-bond then
$h(\om^x)/h(\om)\geq1$, so by the lattice property of $\psi$
we have that $\i_1(\om,x)\geq\i_2(\xi,x)$.
Similarly, if $x$ is a death then $h(\xi^x)/h(\xi)\leq1$
so $\i_1(\om,x)\leq\i_2(\xi,x)$, as required.

We thus have that
\begin{equation}
\psi(gh)=\psi(h)\psi_1(g)\geq\psi(h)\psi_2(g)=
\psi(h)\psi(g).
\end{equation}
\end{proof}

For the next result we let $\l,\d,\g,\l',\d',\g'$ be non-negative,
bounded and Borel-measurable, and write $\l'\geq\l$ if $\l'$
is pointwise no less that $\l$ (and similarly for other functions).
For $a\in\RR$, write $a\l$ or $\l a$ for the function
$x\mapsto a\cdot\l(x)$ (and similarly for other functions).
Recall also the ordering of boundary conditions defined in
Section~\ref{basics_sec} (page~\pageref{bc_po}).

\begin{theorem}\label{correlation}
If $q\geq1$ and  $0<q'\leq q$ then for any boundary condition
$b$ we have that
\begin{equation*}
\begin{split}
  \phi^b_{\L;q,\l,\d,\g}&\leq \phi^b_{\L;q',\l',\d',\g'},\qquad
  \mbox{if } \l'\geq\l,\:\d'\leq\d \mbox{ and } \g'\geq\g\\
  \phi^b_{\L;q,\l,\d,\g}&\geq \phi^b_{\L;q',\l',\d',\g'},\qquad
  \mbox{if } \l'\leq\l q'/q,\:\d'\geq\d q/q', \mbox{ and }
  \g'\leq\g q'/q.
\end{split}
\end{equation*}
Moreover, if $b'\geq b$ are two boundary conditions,  then
\[
\phi^{b'}_{\L;q,\l,\d,\g}\geq \phi^b_{\L;q,\l,\d,\g}. 
\]
\end{theorem}
\begin{corollary}\label{comp_perc}
Let $b$ be any boundary condition.  If $q\geq 1$ then
\begin{equation}
  \phi^b_{\L;q,\l,\d,\g}\leq \mu_{\l,\d,\g}\quad\mbox{and}\quad
  \phi^b_{\L;q,\l,\d,\g}\geq \mu_{\l/q,q\d,\g/q}
\end{equation}
and if $0<q<1$ then
\begin{equation}
  \phi^b_{\L;q,\l,\d,\g}\geq \mu_{\l,\d,\g}\quad\mbox{and}\quad
  \phi^b_{\L;q,\l,\d,\g}\leq \mu_{\l/q,q\d,\g/q}.
\end{equation}
\end{corollary}

\begin{proof}[Proof of Theorem~\ref{correlation}]
We prove the first inequality;  the rest are similar.  
The proof (given Theorem~\ref{preston_thm_1}) is completely
analogous to the one for the discrete random-cluster model,
see~\cite[Theorem~3.21]{grimmett_rcm}.
Recall the formula~\eqref{rcm_intensity} for $\i(\cdot,\cdot)$ in the
random-cluster case.  Let
$\psi_1=\phi^b_{\L;q',\l',\d',\g'}$ and 
$\psi_2=\phi^b_{\L;q,\l,\d,\g}$.  
Clearly $\i_1,\i_2\leq qr$ for all $\om,x$, where $r$
is an upper bound for all of $\l,\d,\g,\l',\d',\g'$.
Let us check the lattice conditions
of Definition~\ref{preston_lattice}.  Let
$\om\leq\xi$ and let $x$ be a bridge such that $\xi^x\not\leq\om$.  Then 
$\i_1(\om,x)=\l'(x)(q')^{k^b_\L(\om^x)-k^b_\L(\om)}$ and
$\i_2(\xi,x)=\l(x)q^{k^b_\L(\xi^x)-k^b_\L(\xi)}$. 
Note that the powers of $q,q'$ both take values in $\{0,-1\}$. 
Since $\l'\geq\l$ and $q'\leq q$, we are done if we
show that $k^b_\L(\om^x)-k^b_\L(\om)\geq k^b_\L(\xi^x)-k^b_\L(\xi)$. 
The left-hand-side is $-1$ if and only if $x$ ties together two different
components of $\om$.  But if it does, then certainly it does the same to $\xi$
since $\xi\leq\om$;  so then also the right-hand-side is $-1$, as required. 
It follows that $\i_1(\om,x)\geq \i_2(\xi,x)$.  The cases when $x$ is a 
death or a ghost-bond are similar.  
\end{proof}

\begin{theorem}[Positive association]\label{rc_fkg}
Let $q\geq1$.  The random-cluster measure $\phi^b_{\L;q,\l,\d,\g}$
is positively associated.
\end{theorem}
Presumably positive association fails when $q<1$, as it does
in the discrete random-cluster model.
\begin{proof}
We only have to verify that $\phi^b_{\L;q,\l,\d,\g}$ has
the lattice property.  Since $q\geq1$ this follows from the fact that
$k^b_\L(\om^x)-k^b_\L(\om)\geq k^b_\L(\xi^x)-k^b_\L(\xi)$ if
$\om\geq\xi$ and $x$ is a bridge or ghost-bond, and the other way
around if $x$ is a death, as in the proof of Theorem~\ref{correlation}.
\end{proof}
The next result is a step towards the `finite energy property'
of Lemma~\ref{fin_en};  it provides upper and lower bounds on the 
probabilities of seeing or not seeing any bridges, deaths
or ghost-bonds in small regions.  These bounds are useful because
they are uniform in $\L$.  For the statement of the result, we let
$q>0$, let $\L=(K,F)$ be a region and $I\subseteq K$ and 
$J\subseteq F$ intervals.  Define
\begin{equation}
\ol\l=\sup_{x\in J}\l(x),\qquad \ul\l=\inf_{x\in J}\l(x)
\end{equation}
and similarly for $\ol\d,\ul\d,\ol\g,\ul\g$ with $J$ replaced by $I$.
Write
\begin{align}
\eta_\l&=\min\{e^{-\ol\l|J|},e^{-\ol\l|J|/q}\},&
\qquad \eta^\l&=\max\{e^{-\ul\l|J|},e^{-\ul\l|J|/q}\},
\nonumber\\
\eta_\d&=\min\{e^{-\ol\d|I|},e^{-q\ol\d|I|}\},&
\qquad \eta^\d&=\max\{e^{-\ul\d|I|},e^{-q\ul\d|I|}\},
\nonumber\\
\eta_\g&=\min\{e^{-\ol\g|I|},e^{-\ol\g|I|/q}\},&
\qquad \eta^\g&=\max\{e^{-\ul\g|I|},e^{-\ul\g|I|/q}\}.
\nonumber
\end{align}
These are to be interpreted as six distinct quantities.

\begin{proposition}\label{fin_en1}
For any boundary condition $b$ we have that
\begin{align}
\eta_\l
&\leq\phi^b_{\L;q,\l,\d,\g}(|B\cap J|=0\mid\cF_{\L\sm J})
\leq \eta^\l\nonumber\\
\eta_\d
&\leq\phi^b_{\L;q,\l,\d,\g}(|D\cap I|=0\mid\cF_{\L\sm I})
\leq \eta^\d \nonumber\\
\eta_\g
&\leq\phi^b_{\L;q,\l,\d,\g}(|G\cap I|=0\mid\cF_{\L\sm I})
\leq \eta^\g \nonumber
\end{align}
\end{proposition}
\begin{proof}
Follows from Proposition~\ref{cond_meas_rcm}
and Corollary~\ref{comp_perc}.
\end{proof}

\begin{remark}
It is convenient, but presumably not optimal, 
to deduce finite energy from stochastic ordering
as we have done here.  For discrete models it is straightforward to
prove the analog of Proposition~\ref{fin_en1} without using
stochastic domination, see~\cite[Theorem~3.7]{grimmett_rcm}.
\end{remark}

\subsection{The FKG-inequality for the Ising model}
\label{ising_fkg_sec}

There is a natural partial order on the set $\S^{b,\a}_\L$ of space--time
Ising configurations, given by:
$\s\geq\tau$ if $\s_x\geq\tau_x$ for all $x\in K$.  In 
Section~\ref{ising_uniq_sec} we will require a \fkg-inequality
for the Ising model, and we prove such a result in this section.  It will
be important to have a result that is valid for all boundary conditions
$(b,\a)$ of Ising type, and when the function $\g$ is allowed to
take negative values.  
The result will be proved by expressing the space--time 
Ising measure as a weak limit of discrete Ising measures, for which
the \fkg-inequality is known.  The same approach was used for the 
space--time percolation model in~\cite{bezuidenhout_grimmett}.
We let $\l,\d$ denote non-negative functions,
as before, and we let $b=\{P_1,\dotsc,P_m\}$ and $\a$ be
fixed.

\begin{figure}[hbt]
\centering
\includegraphics{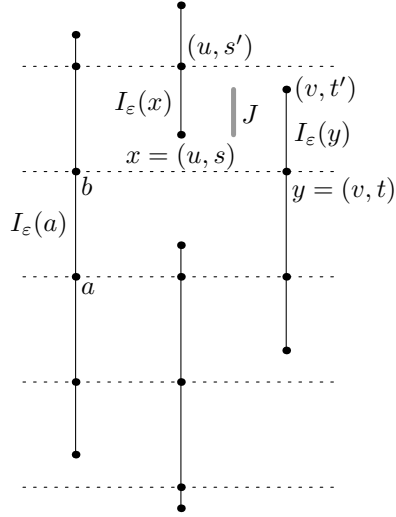}
\caption[Discretization]
{Discretized Ising model.  $K$ is drawn as solid vertical lines, and is 
the union of four closed, disjoint intervals.
Dotted lines indicate the levels $k\eps$ for $k\in\ZZ$.
Elements of $K^\eps$ are drawn as black dots.
The interval $J=uv\times\{[s,s')\cap[t,t')\}$, which appears in the
integral in~\eqref{pp2}, is drawn grey.
In this illustration $b=\free$.}\label{discr_fig}  
\end{figure}

Recall that $K$ consists of a collection of
disjoint intervals $I^v_i$.
Write $\cE$ for the set of endpoints $x$ of the $I^v_i$ 
for which $x\in K$.
Similarly, each $P_i\setminus\{\G\}$ is a finite union of disjoint
intervals;  write $\cB$ for the set of endpoints $y$ of these intervals
for which $y\in K$.
For $\eps>0$, let
\begin{equation}
K^\eps=\cE\cup\cB\cup\{(v,\eps k)\in K: k\in\ZZ\}.
\end{equation}
Let $\S^\eps$ denote the set of vectors $\s'\in\{-1,+1\}^{K^\eps\cup\{\G\}}$
that respect the boundary condition $(b,\a)$;  that is,
(i) if $x,y\in K^\eps\cup\{\G\}$ are such that $x,y\in P_i$ for 
some $i$, then $\s'_x=\s'_y$, and (ii) if in addition $\a(i)\neq 0$
then $\s'_x=\a(i)$.  For each $x=(v,t)\in K^\eps$, let $t'>t$ be maximal
such that the interval $I_\eps(x):=v\times[t,t')$ lies in $K$
but contains no other element
of $K^\eps$;  if no such $t'$ exists let $I_\eps(x):=\{x\}$.
See Figure~\ref{discr_fig}.

We now define the appropriate coupling constants for the discretized model.
Let $x,y\in K^\eps$, $x\neq y$.  First suppose $I_\eps(x)$ and $I_\eps(y)$ 
share an endpoint, which we may assume to be the right endpoint of
$I_\eps(x)$.  Then define
\begin{equation}\label{pp1}
p^\eps_{xy}=1-\int_{I_\eps(x)}\d(z)\, dz.
\end{equation}
Next, suppose $x=(u,s)$ and $y=(v,t)$ are such that $uv\in E$,
and such that $I_\eps(x)=\{u\}\times[s,s')$ and 
$I_\eps(y)=\{v\}\times[t,t')$ satisfy $[s,s')\cap[t,t')\neq\es$.
Then let $J=uv\times\{[s,s')\cap[t,t')\}$ and define
\begin{equation}\label{pp2}
p^\eps_{xy}=\int_J\l(e)\, de.
\end{equation}
For all other $x,y\in K^\eps$ we let $p^\eps_{xy}=0$.  Finally, 
for all $x\in K^\eps$ define
\begin{equation}\label{pp3}
p^\eps_{x\G}=\int_{I_\eps(x)}\g(z)\, dz.
\end{equation}
Note that $p^\eps_{x\G}$ can be negative.  

Let $J^\eps_{xy}$ and $h^\eps_x$ ($x,y\in K^\eps$) be defined by
\begin{equation}\label{js}
1-p^\eps_{xy}=e^{-2J^\eps_{xy}},\qquad
1-p^\eps_{x\G}=e^{-2h^\eps_x}.
\end{equation}
Let $\pi'_\eps$ be the Ising measure on $\S^\eps$ with these
coupling constants, that is 
\begin{equation}\label{mrk2}
\pi'_\eps(\s')=\frac{1}{Z^\eps}\exp\Big(
\frac{1}{2}\sum_{x,y\in K^\eps}J^\eps_{xy}\s'_x\s'_y+
\sum_{x\in K^\eps}h^\eps_x\s'_x\a_\G\Big),
\end{equation}
where $Z^\eps$ is the appropriate normalizing constant.
In the special case when $\g\geq 0$ and $(b,\a)$ is simple, all the
$p^\eps_{xy}$ and $p^\eps_{x\G}$ lie in $[0,1]$
for $\eps$ sufficiently small, and $\pi'_\eps$ is 
coupled via the standard Edwards--Sokal 
measure~\cite[Theorem~1.10]{grimmett_rcm} to the $q=2$ 
random-cluster measure with these edge-probabilities.

There is a natural way to map each element $\s'\in\S^\eps$ to an
element $\s$ of $\S^\free_\L$, namely by letting $\s$ take the value
$\s'_x$ throughout $I_\eps(x)$.  Let $\pi_\eps$ denote the law of 
$\s$ under this mapping.  By a direct computation 
using~\eqref{mrk2} (for example by splitting off the factor
corresponding to `vertical' interactions in the sum over $x,y$)
one may see that 
\begin{equation}\label{mrk3}
\pi_\eps\Rightarrow\el\cdot\er^{b,\a}_\L\quad\mbox{as }\eps\downarrow 0,
\end{equation}
where $\el\cdot\er^{b,\a}_\L$ is the space--time Ising measure defined 
at~\eqref{ising_def_eq}.

For $S\in\cG_\L$ an event, we write $\partial S$ for the boundary
of $S$ in the Skorokhod metric.  We say that $S$ is a 
\emph{continuity set} if $\el \one_{\partial S}\er^{b,\a}_\L=0$.
By standard facts about weak convergence, \eqref{mrk3} implies
that $\pi_\eps(S)\rightarrow\el \one_S\er^{b,\a}_\L$ for
each continuity set $S$.
Note that $\partial(S\cap T)\se\partial S\cup \partial T$,
so if $S,T\in\cG_\L$ are continuity sets then so is $S\cap T$.

\begin{lemma}\label{ising_fkg}
Let $S,T\in\cG_\L$ be increasing continuity sets.  Then
\[
\el\one_{S\cap T}\er^{b,\a}_\L\geq
\el\one_{S}\er^{b,\a}_\L \el\one_{T}\er^{b,\a}_\L.
\]
\end{lemma}
\begin{proof}
By the standard \fkg-inequality for the classical Ising model, we have
for each $\eps>0$ that
\[
\pi_\eps(S\cap T)\geq\pi_\eps(S)\pi_\eps(T).
\]
The result follows from~\eqref{mrk3}.
\end{proof}

In the next result, we write $\el\cdot\er_\g$ for the 
space--time Ising measure $\el\cdot\er^{b,\a}_\L$ with ghost-field
$\g$.
\begin{lemma}\label{ising_mon}
Let $S$ be an increasing continuity set, and let $\g_1\geq\g_2$
pointwise.  Then $\el \one_S\er_{\g_1}\geq\el \one_S\er_{\g_2}$.
\end{lemma}
\begin{proof}
Follows from~\eqref{mrk3} and the fact that $\pi'_\eps$ is
increasing in~$\g$.
\end{proof}

\begin{example}\label{cty_ex}
Here is an example of a continuity set.  Let 
$R\se K$ be a finite union of intervals, some of which may
consist of a single point.  Let $a\in\{-1,+1\}$.  Then the event
\[
S=\{\s\in\S:\s_x=a\mbox{ for all }x\in R\}
\]
is a continuity set, since $\s\in\partial S$ only if
$\s$ changes value exactly on an endpoint of one of the intervals
constituting $R$.  
\end{example}

The assumption above that $S,T$ be continuity sets is
an artefact of the proof method and can presumably be removed.
It should be possible to establish versions of 
Theorems~\ref{preston_thm_1} and~\ref{FKG} also for Ising spins, using
a Markov chain approach.  The auxiliary process $D$
complicates this.  The author would like to thank Jeffrey Steif
for pointing out an error in an earlier version of
this subsection.

\subsection{Correlation inequalities for the Potts model}
\label{gks_sec}

A cornerstone in the study of the classical Ising model 
is provided by the so-called \gks- or Griffiths' inequalities
(see~\cite{griffiths67_I,griffiths67_II,kelly_sherman}) 
which state that certain
covariances are non-negative.
Recently, in~\cite{ganikhodjaev_razak} and~\cite{grimmett_gks},
it was demonstrated that these inequalities follow
from the \fkg-inequality for the
random-cluster representation, using
an argument that also extends to the Potts
models.  In this section we adapt 
the methods of~\cite{grimmett_gks} to the space--time
setting.

Let $q\geq 2$ be fixed, $\L$ a fixed region, and
$b$ a fixed random-cluster boundary condition.
We let $\a$ be such that $(b,\a)$ is a simple boundary
condition with $\a_\G=q$.
It is important to note that the proofs in this section are only
valid for this choice of $\a$.  Therefore, some of the results here 
are less general than what we require for detailed study of the
space--time Ising model, and we will then resort to the results of the 
previous subsection.

Let $\pi,\phi$ denote the Potts- and random-cluster measures
with the given parameters, respectively.
We will be using the complex variables
\begin{equation}
\s_x=\exp\Big(\frac{2\pi i\nu_x}{q}\Big),
\end{equation}%
\nomenclature[s]{$\s$}{Ising- or Potts spin}%
where $i=\sqrt{-1}$.
Note that when $q=2$ this agrees with the previous definition
on page~\pageref{sigma_nu}.
(In~\cite{grimmett_gks} many alternative possibilities
for $\s$ are explored;  similar results hold
at the same level of generality here,
but we refrain from
treating this added generality for simplicity of presentation.)

Define for $A\subseteq K$ a finite set
\begin{equation}
\s_A:=\prod_{x\in A}\s_x.
\end{equation}
More generally, if $\ul r=(r_x:x\in A)$ is a vector of integers
indexed by $A$, define
\begin{equation}
\s^{\ul r}_A:=\prod_{x\in A}\s^{r_x}_x.
\end{equation}
Thus $\s_A\equiv\s^{\ul 1}_A$ where $\ul 1$ is a constant
vector of $1$'s.  The set $B$ in the following should not be
confused with the bridge-set $B=B(\om)$.

\begin{lemma}[\gks\ inequalities]\label{gks_lem}
Let $A,B\subseteq K$ be finite sets, not necessarily disjoint,
and let $\ul r=(r_x:x\in A)$ and $\ul s=(s_y:y\in B)$.  Then
\begin{equation}\label{gks_1_eq}
\pi(\s_A^{\ul r})\geq 0
\end{equation}
and
\begin{equation}\label{gks_2_eq}
\pi(\s_A^{\ul r};\s_B^{\ul s}) 
:= \pi(\s_A^{\ul r}\s_B^{\ul s}) - 
\pi(\s_A^{\ul r})\pi(\s_B^{\ul s}) \geq 0.
\end{equation}
In particular, $\pi(\s_A)\geq 0$ and $\pi(\s_A;\s_B) \geq 0$.
\end{lemma}

A result similar to Lemma~\ref{gks_lem} holds for $A,B\subseteq\ol K$,
but then care must be taken to define $\s_x$ appropriately for points
$x\in\partial\L$ that do not lie in $\L$.  
For example, if $x=(v,t)$ is an isolated 
point in $\KK\setminus K$ then the corresponding result holds if we replace
$\s_x$ by one of $\s_{x+}$ or $\s_{x-}$, where
$\s_{x+}=\lim_{\eps\downarrow0}\s_{(v,t+\eps)}$
and $\s_{x-}=\lim_{\eps\downarrow0}\s_{(v,t-\eps)}$
(these limits exist almost surely but are in general different
for such $x$).

For $\om\in\Om$ let $k=k^b_\L(\om)$, and 
let $C_1(\om),\dotsc,C_{k}(\om)$
denote the components of $\om$ in $\L$, defined according 
to the boundary condition $b$.  We assume that $\G\in C_k(\om)$,
and thus $C_1(\om),\dotsc,C_{k-1}(\om)$ are the 
`$\G$-free' components of $\om$.
Lemma~\ref{gks_lem}  will follow from 
Theorems~\ref{correlation} and~\ref{rc_fkg}
using the following representation.
\begin{lemma}\label{gks_fkg_lem}
Let $\ul r=(r_x:x\in A)$ and write $r_j=\sum_{x\in A\cap C_j}r_x$
(for $j=1,\dotsc,k-1$).  Then
\[
\pi(\s^{\ul r}_A)=\phi(r_j\equiv 0 
\mbox{ (mod $q$), for } j=1,\dotsc,k-1).
\]
\end{lemma}
Note that the event on the right-hand-side is increasing; also
note that if $r_x=1$ for all $x$ then $r_j=|A\cap C_j|$.
\begin{proof}
Let $U_1,U_2,\dotsc$ be independent random variables with
the uniform distribution on $\{e^{2\pi im/q}:m=1,\dotsc,q\}$, and
let $\PP$ denote the Edwards--Sokal coupling~\eqref{es_def}
of $\pi$ and $\phi$.  We have that
\begin{equation}\label{hej2}
\PP(\s_A^{\ul r}\mid\om)=E\Big(1\cdot\prod_{j=1}^{k-1}U_j^{r_j}\Big)
=\prod_{j=1}^{k-1}E(U_j^{r_j}),
\end{equation}
where $E$ denotes expectation over the $U_j$
(recall that $\nu_\G=q$, so $\s_\G=1$).  Since $U_j$ is uniform
we have that
\begin{equation}
E(U_j^r)=\frac{1}{q}\sum_{m=1}^q \big(e^{2\pi i m/q}\big)^r
=\left\{
\begin{array}{ll}
1, & \mbox{if $r\equiv 0$ (mod $q$)}, \\
0, & \mbox{otherwise}.
\end{array}\right.  
\end{equation}
The result follows on taking the expectation of~\eqref{hej2}.
\end{proof}

\begin{proof}[Proof of Lemma \ref{gks_lem}]
It is immediate from Lemma~\ref{gks_fkg_lem} that 
$\pi(\s^{\ul r}_A)\geq0$, which is~\eqref{gks_1_eq}.
For~\eqref{gks_2_eq} we note that
$\s_A^{\ul r}\s_B^{\ul s}=\s_{A\cup B}^{\ul t}$, where
$\ul t$ is the vector indexed by $A\cup B$ given by
$t_x=r_x+s_x$ if $x\in A\cap B$,
$t_x=r_x$ if $x\in A\setminus B$, and
$t_x=s_x$ if $x\in B\setminus A$.
Thus, with the obvious abbreviations,
\begin{align*}
\pi(\s^{\ul r}_A\s^{\ul s}_B)&=\phi(t_j\equiv 0\;\forall j)\\
&\geq\phi(r_j\equiv 0\;\forall j\mbox{ and } 
s_j\equiv 0\;\forall j)\\
&\geq \phi(r_j\equiv 0\;\forall j)\phi(s_j\equiv 0\;\forall j)\\
&=\pi(\s^{\ul r}_A)\pi(\s^{\ul s}_B),
\end{align*}
where the second inequality follows from positive association
of $\phi$, Theorem~\ref{rc_fkg}.
\end{proof}

In the Ising model, the covariance~\eqref{gks_2_eq} is related to
the derivative of $\el\s_A\er$ with respect to 
the coupling strengths;  thus it follows from~\eqref{gks_2_eq}
that $\el\s_A\er$ is increasing in these quantities.  Here is the
corresponding result for the Potts model.

\begin{figure}[hbt]
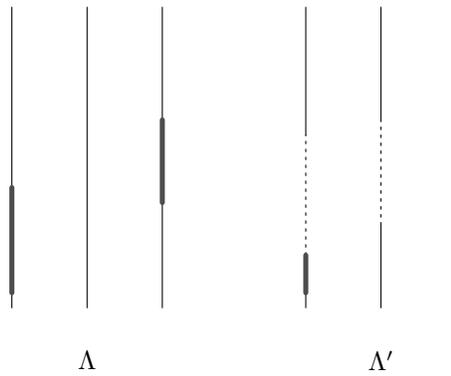

\centering
\includegraphics{thesis.36}
\qquad\qquad\includegraphics{thesis.37}
\caption[Subracting part of a region]
{Left:  a region $\L$ with the boundary condition
$b=\{P_1\}$, where $P_1\setminus\{\G\}$ is drawn bold.  
Right:  the corresponding region
$\L'$ when the set $R$, drawn dashed, has been removed;  the 
boundary condition is $b'=\{P'_1\}$ where $P'_1=P_1\setminus R$
and $P_1'\setminus\{\G\}$ is drawn bold.  In this picture we have not
specified which endpoints of $R$ belong to $R$.}
\label{R_fig}  
\end{figure}

Let $A\subseteq K$ be a finite set, and let $R\subseteq K$ 
be a finite union of positive length 
intervals whose interiors are disjoint from $A$.
We write $\L'$ for the region corresponding to
$K'=K\sm R$.  If $b=(P_1,\dotsc,P_m)$
we define the boundary condition $b'=(P'_1,\dotsc,P'_m)$, where
$P'_i=P_i\setminus R$.  Thus $b'$ agrees with $b$ on $\boundary\L$,
but is `free' on $\boundary\L'\sm\boundary\L$.  
See Figure~\ref{R_fig}.
Similar results hold
for other $b'$.

\begin{lemma}\label{cor_mon_lem}
The average $\pi^b_\L(\s^{\ul r}_A)$ is increasing in
$\l$ and $\g$ and decreasing in $\d$.  Moreover,  
\begin{equation}\label{cor_mon_eq}
\pi^{b'}_{\L'}(\s^{\ul r}_A)\leq \pi^b_\L(\s^{\ul r}_A).
\end{equation}
\end{lemma}
We interpret $\pi^{b'}_{\L'}(\s^{\ul r}_A)$ as $0$ when $A$ intersects
the interior of~$R$.

\begin{proof}
The claim about monotonicity in $\l,\g,\d$ follows from
the stochastic ordering of random-cluster measures,
Theorem~\ref{correlation}, and the representation
in Lemma~\ref{gks_fkg_lem}.
Let us prove~\eqref{cor_mon_eq}.
It suffices to consider the case when
$R=I$ is a single interval.  First note that
\begin{equation}
\pi^b_\L(\s^{\ul r}_A)=\phi^b_\L(T)\geq\tilde\phi^b_\L(T),
\end{equation}
where $T$ is the event on the right-hand-side of Lemma~\ref{gks_fkg_lem},
and $\tilde\phi^b_\L$ is the measure $\phi^b_\L$ with $\g$ set to
zero on $I$, and $\l(e)$ set to zero whenever $e\not\in F'$.
Hence, using also Corollary~\ref{del_contr},
\begin{equation}\label{car3}
\pi^{b'}_{\L'}(\s^{\ul r}_A)=\phi^{b'}_{\L'}(T)=
\tilde\phi^{b}_\L(T\mid D\cap I\neq\varnothing)\leq
\frac{\tilde\phi^b_\L(T)}{1-e^{-\d(I)}}\leq
\frac{\pi^b_\L(\s^{\ul r}_A)}{1-e^{-\d(I)}},
\end{equation}
where
\[
\d(I)=\int_I\d(x)dx.
\]
The left-hand-side of~\eqref{car3} does not depend on the value
of $\d$ on $I$, so we may let $\d\rightarrow\infty$ on $I$
to deduce the result.
\end{proof}

\begin{example}
Here is a consequence of Lemma~\ref{gks_lem} when $\ul r$ is
not constant.  Let $x,y\in K$, and write $\tau_{xy}=\s_x\s_y^{-1}$.
Then $\tau_{xy}$ is a $q$th root of unity, and it follows that
\begin{equation}
\one\{\nu_x=\nu_y\}=\one\{\s_x=\s_y\}=\frac{1}{q}\sum_{r=0}^{q-1}\tau_{xy}^r.
\end{equation}
So if $z,w\in K$ too then
\begin{equation}\label{tjo4}
\begin{split}
\pi^b_\L(\nu_x=\nu_y,\,\nu_z=\nu_w)&=
\frac{1}{q^2}\sum_{r,s=0}^{q-1}\pi^b_\L(\tau_{xy}^r\tau_{zw}^s)\\
&\geq\frac{1}{q^2}\sum_{r,s=0}^{q-1}
\pi^b_\L(\tau_{xy}^r)\pi^b_\L(\tau_{zw}^s)\\
&=\pi^b_\L(\nu_x=\nu_y)\pi^b_\L(\nu_z=\nu_w).
\end{split}
\end{equation}
This inequality does not quite follow  from
the correlation/connection property of 
Proposition~\ref{corr_conn_prop} when $q>2$.  In the case 
when $\g=0$ it follows straight away from the 
Edwards--Sokal coupling, without using stochastic
domination properties of the random-cluster model;
see~\cite[Corollary~6.5]{ghm}.
\end{example}

\section{Infinite-volume random-cluster measures}
\label{inf_rc_sec}

In this section we define 
random-cluster measures on the \emph{unbounded} spaces $\LLam,\LLam_\b$
of~\eqref{def-newL} and~\eqref{def-bL}, 
for which  Definition~\ref{rc_def} cannot make sense
(since $k$ will be infinite).  
One standard approach in statistical physics
is to study the class of measures which satisfy a conditioning
property similar to that of Proposition~\ref{cond_meas_rcm}
for all bounded regions;  the first task is then
to show that this class is nonempty.  The 
book~\cite{georgii88} is dedicated to this approach for classical
models.  We will instead  follow the
route of proving weak convergence as the bounded
regions $\L$ grow.  In doing so we follow standard methods
(see~\cite[Chapter~4]{grimmett_rcm}), adapted to the current setting.
See also~\cite{aizenman_nacht} for results of this type.

Central to the topic of infinite-volume measures is the question
when there is a unique such measure.  There
may in general be multiple such measures, obtainable by passing to
the limit using different boundary conditions.
Non-uniqueness of infinite-volume measures is intimately related
to the concept of phase transition described in the Introduction.
Intuitively, if there is not a unique limiting measure this
means that the boundary conditions have an `infinite range'
effect, and that the system does not know what state to favour,
indicating a transition from one preferred state to another.

\subsection{Weak limits}\label{rc_wl_sec}

We fix $q\geq1$ and non-negative bounded measurable 
functions $\l,\d,\g$.
Let $L_n$ be a sequence of subgraphs of $\LL$ and $\b_n$ a sequence of
positive numbers.  Writing $\L_n$ for the simple 
region given by $L_n$ and
$\b_n$ as in~\eqref{def-oldL}, we say
that $\L_n\uparrow\LLam$ if $L_n\uparrow\LL$ and $\b_n\rightarrow\oo$.
We assume throughout that $L_n$ and $\b_n$ are strictly increasing.
Versions of the results in this section are valid also when $\b<\oo$ is kept
fixed as $L_n\uparrow\LL$ so that $\L_n\uparrow\LLam_\b$ given
in~\eqref{def-bL}.  
We will only supply proofs in the $\b_n\rightarrow\oo$ case
as the $\b<\oo$ case is similar.

Recall that a sequence $\psi_n$ of probability measures on $(\Om,\cF)$
is \emph{tight} if for each $\eps>0$ there is a compact set $A_\eps$
such that $\psi_n(A_\eps)\geq1-\eps$ for all $n$.  Here compactness
refers, of course, to the Skorokhod topology outlined in
Section~\ref{basics_sec} and defined in detail in  
Appendix~\ref{skor_app}.

Let $\phi^b_n:=\phi^b_{\L_n}$.  The proof of the following result
is given in Appendix~\ref{skor_app}.

\begin{lemma}\label{tight_lem}
For any sequence of boundary conditions $b_n$ on $\L_n$,
the sequence of measures $\{\phi^{b_n}_n:n\geq 1\}$ is tight.
\end{lemma}

For $x=(e,t)\in\FF$
with $t\geq0$ (respectively $t<0$), 
let $V_x(\om)$%
\nomenclature[V]{$V_x(\om)$}{Element count of $B$, $G$ or $D$}
 denote the number of
elements of the set $B\cap(\{e\}\times[0,t])$
(respectively $B\cap(\{e\}\times(-t,0])$).  Similarly, for 
$x\in\KK\times\{\mathrm{d}\}$ and 
$x\in\KK\times\{\mathrm{g}\}$, define $V_x$ to count the number
of deaths and ghost-bonds between $x$ and the origin, respectively.
An event of the form
\[
R=\{\om\in\Om: V_{x_1}(\om)\in A_1,\dotsc,
V_{x_m}(\om)\in A_m\}\in\cF
\]
for $m\geq 1$ and the $A_i\subseteq\ZZ$
is called a \emph{finite-dimensional cylinder event}.
For $z=(z_1,\dotsc,z_m)$ and $z'=(z'_1,\dotsc,z'_m)$
elements of $\ZZ^m$, we write $z'\geq z$ if 
$z'_i\geq z_i$ for all $i=1,\dotsc, m$;  we write $z'> z$
if $z'\geq z$ and $z'\neq z$.

\begin{theorem}\label{conv_lem}
Let $b\in\{\free,\wired\}$ and $q\geq1$.  
The sequence of measures $\phi^b_n$ converges
weakly to a probability measure.  The limit measure does not depend on
the choice of sequence $\L_n\uparrow\LLam$.
\end{theorem}
The limiting measure in Theorem~\ref{conv_lem} will be denoted
$\phi^b$, or $\phi^{b,\b}_{q,\l,\d,\g}$ if the parameters need to 
be emphasized;  here $\b\in(0,\oo]$.
\begin{proof}
Consider the case $b=\wired$.
Let $\L$ be a simple region and $f:\Om\rightarrow\RR$ an increasing,
$\cF_\L$-measurable function.  Let $n$ be large enough so that 
$\L_n\supseteq\L$ and let $\cC$ be the event that all 
components inside $\L_n$
which intersect $\boundary\L_n$ are connected in $\L_{n+1}$.  Then
by Corollary~\ref{del_contr} and the
\fkg-property we have that
\begin{equation}\label{wcs1}
\phi^\wired_n(f)=\phi^\wired_{n+1}(f\mid \cC)\geq \phi^\wired_{n+1}(f),
\end{equation}
which is to say that $\phi^\wired_n\geq\phi^\wired_{n+1}$.
At this point we could appeal to Corollary~IV.6.4 of~\cite{lindvall}, 
which proves that a sequence of probability measure which is tight
and stochastically ordered as in~\eqref{wcs1} 
necessarily converges weakly.  However, we shall later need to know
that the finite dimensional distributions converge, so we prove this 
now;  it then follows from tightness and standard properties of
the Skorokhod topology that the sequence converges weakly.

Let 
$x_1,\dotsc,x_k\in F\cup(K\times\{\mathrm{g}\})$ and let
$x_{k+1},\dotsc,x_m\in K\times\{\mathrm{d}\}$.
For $z=(z_1,\dotsc,z_m)\in\ZZ^m$, write 
\[
\tilde z=(z_1,\dotsc,z_k,-z_{k+1},\dotsc,-z_m).
\]
Let $V=V(\om)=(V_{x_1}(\om),\dotsc,V_{x_m}(\om))$ and for $A\se\ZZ^m$
consider the finite-dimensional cylinder event 
$R=\{V\in A\}$.  We have that
\begin{equation}\label{wcs2}
\begin{split}
\phi^\wired_n(R)&=
\sum_{z\in A}\phi^\wired_n(V=z)
=\sum_{z\in A}\phi^\wired_n(\tilde V=\tilde z)\\
&=\sum_{z\in A}[\phi^\wired_n(\tilde V\geq\tilde z)
-\phi^\wired_n(\tilde V>\tilde z)].
\end{split}
\end{equation}
The events $\{\tilde V\geq\tilde z\}$ and $\{\tilde V>\tilde z\}$
are both increasing, so by~\eqref{wcs1} the limits
\[
\ol\phi(\tilde V\geq\tilde z)=\lim_{n\rightarrow\oo}
\phi^\wired_n(\tilde V\geq\tilde z)\quad\text{and}\quad
\ol\phi(\tilde V>\tilde z)=\lim_{n\rightarrow\oo}
\phi^\wired_n(\tilde V>\tilde z)
\]
exist.  Define $\ol\phi$ by
\[
\ol\phi(R):=\sum_{z\in A}[\ol\phi(\tilde V\geq\tilde z)
-\ol\phi(\tilde V>\tilde z)].
\]
Then, by the bounded convergence theorem, $\ol\phi$ defines
a probability measure on the algebra of 
finite-dimensional cylinder events in $\cF_\L$.  Thus $\ol\phi$ 
extends to 
a unique probability measure $\phi^\wired$ on $\cF_\L$
(see~\cite[Theorem~3.1]{billingsley_probmeas}).
Since
$\phi^\wired_n(R)\rightarrow\phi^\wired(R)$ for all
finite-dimensional cylinder events in $\cF_\L$ and since the 
sequence $(\phi^\wired_n:n\geq1)$ is tight, it follows that
$\phi^\wired_n\Rightarrow\phi^\wired$ on $(\Om,\cF_\L)$.
Since $\L$ was arbitrary and the $\cF_\L$ generate
$\cF$, the convergence for $b=\wired$ follows.  

For the independence of the choice of sequence
$\L_n$, let also $\D_n\uparrow\LLam$.  
Let $m$ be an integer,
and choose $l=l(m)$ and $n=n(m)$ so that 
$\L_l\subseteq\D_m\subseteq\L_n$.    We have that
\[
\phi^\wired_{\L_l}\geq\phi^\wired_{\D_m}\geq\phi^\wired_{\L_n},
\]
so letting $m\rightarrow\infty$ tells
us that the limits are equal (see Remark~\ref{po_rk}).

The arguments for $b=\free$ are similar.
\end{proof}

\begin{remark}\label{po_rk}
If $\psi_1,\psi_2$ are two probability 
measures on $(\Om,\cF)$ such 
that both $\psi_1\geq\psi_2$ and $\psi_2\geq\psi_1$ then
$\psi_1=\psi_2$.  To see this, note that for $R$ any 
finite-dimensional cylinder event,
we may as in~\eqref{wcs2} write
\[
\psi_j(R)=\sum_{z\in A} [\psi_j(\tilde V\geq\tilde z)
-\psi_j(\tilde V>\tilde z)],\quad j=1,2.
\]
It follows that $\psi_1(R)=\psi_2(R)$ for all such $R$, and hence that 
$\psi_1=\psi_2$ (see Appendix~\ref{skor_app}).
\end{remark}

For any sequence $b_n$ of boundary conditions, \emph{if} the sequence of
measures $(\phi^{b_n}_n:n\geq1)$ has a weak limit $\phi$, then
$\phi^\free\leq\phi\leq\phi^\wired$;  
this follows from the second part of Theorem~\ref{strassen_thm}.
Hence there is a unique random-cluster measure if and only if
$\phi^\free=\phi^\wired$.  It turns out that the set of real triples 
$(\l,\d,\g)$ such that there is \emph{not} a unique 
random-cluster measure has Lebesgue measure zero, 
see Theorem~\ref{rc_uniq_thm}.

\subsection{Basic properties}

Some further properties of the measures $\phi^b$, for
$b\in\{\free,\wired\}$, follow, all being straightforward adaptations
of standard results, as summarized in~\cite[Section~4.3]{grimmett_rcm}.
First, recall the upper and lower bounds on the probabilities
of seeing no bridges, deaths or ghost-bonds in small regions
which is provided by Proposition~\ref{fin_en1}, as well as
the notation introduced there.  
\begin{lemma}[Finite energy property]\label{fin_en}
Let $q\geq1$ and let $I\subseteq\KK$ and 
$J\subseteq\FF$ be bounded intervals.  
Then for $b\in\{\free,\wired\}$ we have that
\begin{align}
\eta_\l
&\leq\phi^b(|B\cap J|=0\mid\cT_{J})
\leq \eta^\l\nonumber\\
\eta_\d
&\leq\phi^b(|D\cap I|=0\mid\cT_{I})
\leq\eta^\d \nonumber\\
\eta_\g
&\leq\phi^b(|G\cap I|=0\mid\cT_{I})
\leq \eta^\g \nonumber
\end{align}
\end{lemma}

The same result holds for any weak limit of random-cluster measures
with $q>0$;  we assume that $q\geq1$ and $b\in\{\free,\wired\}$
only because then we know that the
measures $\phi^b_\L$ converge.  

\begin{proof}
Recall the notation $V_x(\om)$ introduced 
before Theorem~\ref{conv_lem}, and note that
the event $\{|B\cap J|=0\}$ is a finite-dimensional cylinder event.
For $J\subseteq\FF$ as in the statement, let 
$x_1,x_2,\dotsc$ be an enumeration of the points in
$(\KK\times\{\mathrm{d}\})\cup(\KK\times\{\mathrm{g}\})\cup(\FF\sm J)$
with rational $\RR$-coordinate.  We have that
$\cT_J=\s(V_{x_1},V_{x_2},\dotsc)$
so by the martingale convergence theorem
\[
\phi^b(|B\cap J|=0\mid\cT_{J})=
\lim_{n\rightarrow\oo}\phi^b(|B\cap J|=0\mid V_{x_1},\dotsc,V_{x_n}).
\]
For $\ul z\in\ZZ^n$, let $A_{\ul z}=\{(V_{x_1},\dotsc,V_{x_n})=\ul z\}$.
Then
\[
\begin{split}
\phi^b(|B\cap J|=0\mid \cF_n)&=
\sum_{\ul z\in \ZZ^n}\frac{\phi^b(A_{\ul z},\{|B\cap J|=0\})}
{\phi^b(A_{\ul z})}\one_{A_{\ul z}}\\
&=\lim_{\D}\sum_{\ul z\in\ZZ^n}
\frac{\phi_\D^b(A_{\ul z},\{|B\cap J|=0\})}
{\phi_\D^b(A_{\ul z})}\one_{A_{\ul z}}\\
&=\lim_{\D}\phi^b_\D(|B\cap J|=0\mid \cF_n).
\end{split}
\]
The result now follows from Proposition~\ref{fin_en1}.
A similar argument holds for $\{|D\cap I|=0\}$ and $\{|G\cap I|=0\}$.
\end{proof}

Define an \emph{automorphism} on $\LLam$ to be a bijection 
$T:\LLam\rightarrow\LLam$ of the form $T=(\a,g):(x,t)\mapsto(\a(x),g(t))$ 
where $\a:\VV\rightarrow\VV$ is an automorphism of the 
graph $\LL$, and $g:\RR\rightarrow\RR$ is a continuous bijection.
Thus $\a$ has the property that $\a(x)\a(y)\in\EE$ if and only if
$xy\in\EE$.  
For $T$ an automorphism and $\om=(B,D,G)\in\Om$, let
$T(\om)=(T(B),T(D),T(G))$.  For $f:\Om\rightarrow\RR$ measurable,
let $(f\circ T)(\om)=f(T(\om))$, and for $\phi$ a measure on 
$(\Om,\cF)$ define $\phi\circ T(f)=\phi(T(f))$.
\begin{lemma}\label{invar_lem}
Let $b\in\{\free,\wired\}$ 
and let $T$ be an automorphism of $\LLam$ such that
$\l=\l\circ T$, $\g=\g\circ T$ and $\d=\d\circ T$.  Then $\phi^b$
is invariant under $T$, that is $\phi^b=\phi^b\circ T$.
\end{lemma}
\begin{proof}
Let $f$ be a measurable function.  Under the given assumptions,
we have that for any region $\L$,
\begin{equation*}
\phi^b_{\L}(f\circ T)=\int f(T(\om))\,d\phi^b_{\L}(\om)
=\int f(\om)\,d\phi^b_{T^{-1}(\L)}(\om)=\phi^b_{T^{-1}(\L)}(f).
\end{equation*}
The result now follows from Theorem~\ref{conv_lem}.
\end{proof}

\begin{proposition}\label{tail_triv}
The tail $\s$-algebra $\cT$ is trivial under the measures 
$\phi^\free$ and $\phi^\wired$, in that $\phi^b(A)\in\{0,1\}$
for all $A\in\cT$.
\end{proposition}
\begin{proof}
Let $\L\subseteq\D$ be two regions.  We treat the case when $b=\free$,
the case $b=\wired$ follows similarly on reversing several of the
inequalities below.  Let $A\in\cF_\L$ be an
increasing finite-dimensional
cylinder event, and let $B\in\cF_{\D\sm\L}\subseteq\cT_\L$
be an arbitrary finite-dimensional
cylinder event.  We may assume without loss of generality
that $\phi^\free_\D(B)>0$.  By the conditioning
property Proposition~\ref{cond_meas_rcm} and the stochastic ordering of
Theorem~\ref{correlation}, we have that
\begin{equation}
\phi^\free_\D(A\cap B)=\phi^\free_\D(A\mid B)\phi^\free_\D(B)
\geq \phi^\free_\L(A)\phi^\free_\D(B).
\end{equation}
Let $\cR$ denote the set of finite-dimensional
cylinder events in $\cT_\L$.  
Letting $\D\uparrow\LLam$ implies that
\begin{equation}\label{hej3}
\phi^\free(A\cap B)\geq \phi^\free_\L(A)\phi^\free(B)
\end{equation}
for all $B\in\cR$ and all increasing 
finite-dimensional cylinder events $A\in\cF_\L$.  
The set $\cR$ is an algebra, so for fixed
$A$ the difference between the left and right sides of~\eqref{hej3}
extends to a finite measure 
$\psi$ on $\cT_\L$, and by the uniqueness of this extension
it follows that 
$0\leq\psi(B)=\phi^\free(A\cap B)- \phi^\free_\L(A)\phi^\free(B)$ 
for all $B\in\cT_\L\se\cT$.  Thus we may let 
$\L\uparrow\LLam$  to deduce that 
\begin{equation}\label{hej4}
\phi^\free(A\cap B)\geq \phi^\free(A)\phi^\free(B)
\end{equation}
for all increasing finite-dimensional cylinder
events $A\in\cF_\L$ and all $B\in\cT$.
However,~\eqref{hej4} also holds with $B$ replaced by its complement 
$B^c$;  since 
\[
\phi^\free(A\cap B)+\phi^\free(A\cap B^c)=
\phi^\free(A)\phi^\free(B)+\phi^\free(A)\phi^\free(B^c)
\]
it follows that 
\begin{equation}\label{hej5}
\phi^\free(A\cap B)= \phi^\free(A)\phi^\free(B)
\end{equation}
for all increasing finite-dimensional cylinder events
$A\in\cF_\L$ and all $B\in\cT$.  For fixed $B$, the left and right sides 
of~\eqref{hej5} are finite measures which agree on all increasing
events $A\in\cF_\L$.  Using the reasoning of Remark~\ref{po_rk},
it follows that~\eqref{hej5}
holds for all $A\in\cF_\L$, and hence also
for all $A\in\cF$.  Setting $A=B\in\cT$ gives the result.
\end{proof}

In the case when $\LL=\ZZ^d$ and $\l,\d,\g$ are constant, define the 
automorphisms $T_x$, for $x\in\ZZ^d$, by
\[
T_x(y,t)=(y+x,t).
\]
The $T_x$ are called \emph{translations}.  An event $A\in\cF$ is
called $T_x$-invariant if $A=T_x^{-1}A$.  The following ergodicity 
result is a
standard consequence of Proposition~\ref{tail_triv}, 
see for example~\cite[Proposition~14.9]{georgii88}
(here $0$ denotes the element $(0,\dotsc,0)$ of $\ZZ^d$).
\begin{lemma}\label{ergod_lem}
Let $x\in\ZZ^d\setminus\{0\}$ and $b\in\{\free,\wired\}$.  If
$A\in\cF$ is $T_x$-invariant then $\phi^b(A)\in\{0,1\}$.
\end{lemma}

\subsection{Phase transition}

In the random-cluster model, the probability that there is an
unbounded connected component serves as `order parameter':  depending on
the values of the parameters $\l,\d,\g$ this probability may be
zero or positive.  We show in this section that one may define
a critical point for this probability, and then establish some very
basic facts about the phase transition.
We assume throughout this section that $\g=0$, 
that $q\geq1$,  that $\l\geq0$, $\d>0$ are constant, 
and that $\LL=\ZZ^d$ for some $d\geq1$.  Some 
of the results hold for more general $\LL$, but we will not pursue
this here.  The boundary condition $b$ will denote either
$\free$ or  $\wired$  throughout.

Let $\{0\lra\oo\}$ denote the event that the origin lies in an unbounded
component.  Define for $0<\b\leq\oo$,
\begin{equation}
\theta^{b,\b}(\l,\d,q):=\phi^{b,\b}_{q,\l,\d}(0\leftrightarrow\oo).
\end{equation}%
\nomenclature[t]{$\theta$}{Percolation probability}%
When $\b=\oo$ a simple rescaling argument implies that
$\theta^{b,\oo}(\l,\d,q)$ depends on $\l,\d$ through the ratio
$\rho=\l/\d$ only.  Hence we will often in what follows set
$\d=1$ and $\l=\rho$, and define for $0<\b\leq\oo$
\begin{equation}
\theta^{b,\b}(\rho)=
\theta^{b,\b}(\rho,q):=\phi^{b,\b}_{q,\rho,1}(0\leftrightarrow\oo).
\end{equation}
By the stochastic monotonicity of Theorem~\ref{correlation}, 
and a small argument justifying its application
to the event $\{0\lra\oo\}$, the quantity $\theta^b(\rho)$
is increasing in $\rho$.  
\begin{definition}\label{crit_def}
For $b\in\{\free,\wired\}$ and $0<\b\leq\oo$ 
we define the \emph{critical value}
\[
\rho^{b,\b}_\crit(q):=\sup\{\rho\geq 0:\theta^{b,\b}(q,\rho)=0\}.
\]
\end{definition}%
\nomenclature[r]{$\rho_\crit(q)$}{Percolation threshold}%
In what follows we will usually suppress reference to $\b$.
We will see in Section~\ref{press_sec} that 
$\rho^\free(q)=\rho^\wired(q)$ for all $q\geq1$.  Therefore we
will write $\rho_\crit(q)$ for their common value.
We write $\phi^{b}_\rho$ for $\phi^{b,\b}_{q,\rho,1}$.

One may adapt standard methods (see~\cite[Theorem~5.5]{grimmett_rcm}) 
to prove the following:
\begin{theorem}\label{crit_nontriv}
Unless $d=1$ and $\b<\infty$ we have that 
\[
0<\rho_c(q)<\infty.
\]
\end{theorem}
(If $d=1$ and $\b<\oo$ then a standard zero-one argument,
involving comparison to percolation and
the second Borel--Cantelli lemma,
implies that $\rho_\crit=0$.)

Fix $\rho>0$ and for $\om\in\Om$ let 
$N=N(\om)$ denote the number of distinct unbounded components
in $\om$.  By Lemma~\ref{ergod_lem}, using for example the translation
$T:(x,t)\mapsto(x+1,t)$, we have that $N$ is almost surely constant under
the measures $\phi^b_\rho(\cdot)$, $b\in\{\free,\wired\}$.

\begin{theorem}\label{inf_clust_uniq}
The number $N$ of unbounded components is either $0$ or $1$ 
almost surely under~$\phi^b_\rho$.
\end{theorem}
\begin{proof}
We follow the strategy of~\cite{burton_keane},
and as previously we provide details only in the $\b=\oo$
case.  We first show that
$N\in\{0,1,\oo\}$ almost surely.  Suppose to the contrary that
there exists $2\leq m<\oo$ such that $N=m$ almost surely.  Then we may
choose (deterministic) $n,\b$ sufficiently large that the 
corresponding simple region $\L_n=\L_n(\b)$, 
regarded as a subset of $\LLam$,
has the property that $\phi^b_\rho(A)>0$, where $A$ is the event
that the $m$ distinct unbounded components all meet
$\partial\L_n$.
Let $C$ be the event that all points in $\partial\L_n$ are connected
inside $\L_n$.  By the finite energy property, Lemma~\ref{fin_en},
we have that $\phi^b_\rho(C\mid A)>0$, and hence
$\phi^b_\rho(C\cap A)>0$.  But on $\{C\cap A\}$ we have $N=1$, a
contradiction.  Thus $N\in\{0,1,\oo\}$.

Now suppose that $N=\oo$ almost surely.  Let $\b=2n$, and
for $v\in\VV$ and $r\in\ZZ$  let
\begin{equation}
I_{v,r}=\{v\}\times[r,r+1]\subseteq \KK.
\end{equation}
We call $I_{v,r}$ a \emph{trifurcation} if (i) it is contained
in exactly one unbounded component, and 
(ii) if one removes all bridges incident on $I_{v,r}$ and places a 
least one death in $I_{v,r}$, then the unbounded 
component containing it breaks into three distinct unbounded 
components.   See Figure~\ref{trif_fig}.
\begin{figure}[hbt]
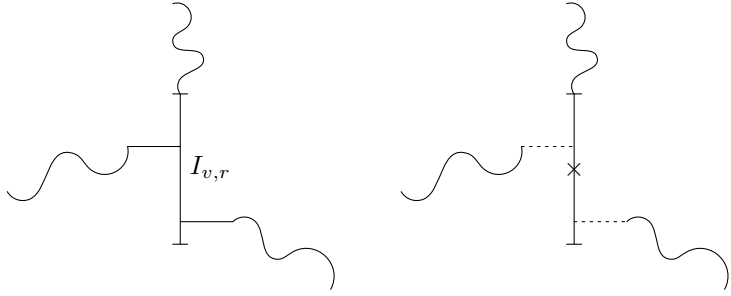

\centering
\includegraphics{thesis.3}\qquad
\includegraphics{thesis.4}
\caption[Trifurcation interval]
{A trifurcation interval (left);  upon removing all 
incident bridges and placing a death in the interval, the
unbounded cluster breaks in three (right).}\label{trif_fig}   
\end{figure}

We claim that
\begin{equation}\label{pos1}
\phi^b_\rho(I_{0,0}\mbox{ is a trifurcation})>0.
\end{equation}
To see this let $n$ be large enough so that $\partial\L_n$ meets
three distinct unbounded components with positive probability.
Conditional on $\cT_{\L_n}$, the finite energy property 
Lemma~\ref{fin_en} allows us to modify the configuration inside
$\L_n$ so that, with positive probability, $I_{0,0}$ is a trifurcation.

We note from translation invariance, Lemma~\ref{invar_lem}, that
the number $T_n$ of trifurcations in $\L_n$ satisfies
\begin{equation}\label{pos2}
\begin{split}
\phi^b_\rho(T_n)&=\sum_{\substack{v\in[-n,n]^d\\r=-n,\dotsc,n-1}} 
\phi^b_\rho(I_{v,r}\mbox{ is a trifurcation})\\
&=2n(2n+1)^d\phi^b_\rho(I_{0,0}\mbox{ is a trifurcation}).
\end{split}
\end{equation}
Define the \emph{sides} of $\L_n$ to be the union of all intervals
$v\times[-n,n]$ where $v$ has at least one coordinate which is $\pm n$.
Topological considerations imply that $T_n$
is bounded from above by the total number of 
deaths on the sides of $\L_n$ plus twice the number of vertices
in $[-n,n]^d$.  
(Each trifurcation needs at least one unique point of exit
from $\L_n$).  Using the stochastic domination
in Corollary~\ref{comp_perc} or otherwise, it follows that
$\phi^b_\rho(T_n)\leq 2(2n+1)^d+\d\cdot 4dn(2n+1)^{d-1}$.
In view of~\eqref{pos1} and~\eqref{pos2} this is a contradiction.
See~\cite[Chapter~5]{br_perc} for more 
details on the topological aspects of this  argument.
\end{proof}

It follows from Theorem~\ref{inf_clust_uniq} that $N=0$ almost surely
under $\phi^b_\crit$ if $\rho<\rho_c$ and that
$N=1$ almost surely if $\rho>\rho_c$.  
It is crucial for the proof that $\LL=\ZZ^d$ is `amenable' in 
the sense that the boundary of $[-n,n]^d$ is an order of magnitude
smaller than the volume.  The result fails, for example, when $\LL$
is a tree, in which case $N=\oo$ may occur;  
see~\cite{pemantle92} for the corresponding phenomenon
in the contact process.

\subsection{Convergence of pressure}\label{press_sec}

In this section we adapt the well-known `convergence of 
pressure' argument to the space--time random-cluster model.
By relating the question of uniqueness of measures to that of 
the existence of certain derivatives, we are able to
deduce that there is a unique infinite-volume measure at
almost every $(\l,\d,\g)$, see Theorem~\ref{rc_uniq_thm} below.
Arguments of this type are
`folklore' in statistical physics, and appear in many places
such as~\cite{ellis85:LD,georgii88,israel79}.
We follow closely the corresponding method
for the discrete random-cluster model given 
in~\cite[Chapter~4]{grimmett_rcm}.  

Let $\l,\d,\g>0$ be constants.  
We will for simplicity of presentation 
be treating only the case when $\g>0$ and $q\geq 1$, though similar
arguments hold when $\g=0$ and when $0<q<1$. 
The partition function
\begin{equation}
Z^b_\L(\l,\d,\g,q)=\int_{\Om} q^{k^b_\L(\om)}
\,d\mu_{\l,\d,\g}(\om)
\end{equation}
is now a function $\RR_+^4\rightarrow\RR$.  In this section we will
study the related \emph{pressure functions}
\begin{equation}\label{pr_eq}
P^b_\L(\l,\d,\g,q)=\frac{1}{|\L|}\log Z^b_\L(\l,\d,\g,q).
\end{equation}
Here, and in what follows, we have abused notation by writing
$|\L|$ for the (one-dimensional) Lebesgue measure $|K|$ of $K$, 
where $\L=(K,F)$.   We will be considering
limits of $P^b_\L$ as the region $\L$ grows.  To be concrete we
will be considering regions of the form 
\begin{equation}\label{pr_reg}
\L=\L_{\ul n,\b}\equiv\{1,\dotsc,n_1\}\times\dotsb\{1,\dotsc,n_d\}
\times[0,\b]
\end{equation}
and limits when $\L\uparrow\LLam$, that is to say all
$n_1,\dotsc,n_d,\b\rightarrow\oo$ (simultaneously).
Strictly speaking such regions do not tend to $\LLam$,
but the $P^b_\L$ are not affected by translating $\L$.
It will be clear from the arguments that  
one may deal in the same way
with limits as $\L\uparrow\LLam_\b$ with $\b<\oo$ fixed.
When $\ul n$ and $\b$ need to be emphasized we will
write $\L_{\ul n,\b}=(K_{\ul n,\b},F_{\ul n,\b})$.

Here is a simple observation about $Z^b_\L$.  Writing
\begin{equation}\label{seb3}
r=\log\l,\quad s=\log \d, \quad t=\log \g,\quad
u=\log q,
\end{equation}
and 
\begin{equation}
D_\L=D\cap K, \quad G_\L=G\cap K, \quad B_\L=B\cap F,
\end{equation}
we have that 
\begin{equation}
\begin{split}
Z^b_\L(r,s,t,u)&\equiv Z^b_\L(\l,\d,\g,q)\\
&=\int_\Om d\mu_{1,1,1}(\om)
\exp\big(r|B_\L|+s|D_\L|+t|G_\L|+u k^b_\L\big).
\end{split}
\end{equation}
(Where $\mu_{1,1,1}$ is the percolation measure where
$B,D,G$ all have rate 1.)
This follows from basic properties of Poisson processes.  
It will sometimes be more convenient to work with
$Z^b_\L(r,s,t,u)$ in this form.  We will also write 
$P^b_\L(r,s,t,u)$ for the pressure~\eqref{pr_eq} using these
parameters~\eqref{seb3}.

Let
$\ul h=(h_1,\dotsc,h_4)$ be a unit vector in $\RR^4$, and let $y\in\RR$.
It follows from a simple computation that the function
$f(y)=P^b_\L((r,s,t,u)+y h)$ has non-negative second derivative.
Indeed, $f''(y)$ is the variance under the appropriate 
random-cluster measure of the quantity 
\[
h_1|B_\L|+h_2|D_\L|+h_3|G_\L|+h_4 k^b_\L.
\]
Since variances are non-negative, have proved
\begin{lemma}\label{press_convex_lem}
Each $P^b_\L(r,s,t,u)$ is a convex function $\RR^4\rightarrow\RR$.
\end{lemma}

Our first objective in this section is the following result.
\begin{theorem}\label{press_converge_thm}
The limit 
\[
P(r,s,t,u)=\lim_{\L\uparrow\LLam} P^b_\L(r,s,t,u)
\]
exists for all $r,s,t,u\in\RR$ and all sequences $\L\uparrow\LLam$
of the form~\eqref{pr_reg},
and is independent of the boundary condition~$b$.
\end{theorem}

The function $P$ is usually called the 
\emph{specific Gibbs free energy}, or \emph{free energy} for short.
It follows that $P$ is a convex function $\RR^4\rightarrow\RR$,
and hence that the set $\cD$ of points in $\RR^4$ at which one or more
partial derivative of $P$ fails to exist has 
zero Lebesgue measure.
We will return to this observation after the proof of
Theorem~\ref{press_converge_thm}.

\begin{proof}[Proof of Theorem~\ref{press_converge_thm}]
We first prove convergence of $P^\free_\L$ 
with free boundary, and then deduce the result for general $b$.  
For each $i=1,\dotsc,d$ let $0<m_i\leq n_i$ and also let $0<\a<\b$.
Write $|\ul m|=m_1\dotsb m_d$.
We may regard the region $\L_{\ul m,\a}$ as a subset of 
$\L_{\ul n,\b}$.  Write $T^{\ul n,\b}_{\ul m,\a}$ for 
the set of points in $F_{\ul n,\b}\setminus F_{\ul m,\a}$ 
adjacent to at least one point in $K_{\ul m,\a}$.  We have that
\begin{equation}
k^\free_{\L_{\ul n,\b}}
\left\{\begin{array}{l}
\leq k^\free_{\L_{\ul m,\a}}+
k^\free_{\L_{\ul n,\b}\setminus{\L_{\ul m,\a}}} \\
\geq k^\free_{\L_{\ul m,\a}}+
k^\free_{\L_{\ul n,\b}\setminus{\L_{\ul m,\a}}} 
-|\ul m|-|B\cap T^{\ul n,\b}_{\ul m,\a}| -1.
\end{array}\right.
\end{equation}
The lower bound follows because the number of `extra' components created by 
`cutting out' $\L_{\ul m,\a}$ from $\L_{\ul n,\b}$ is bounded by
the number of intervals constituting $K_{\ul m,\a}$, plus the number of 
bridges that are cut, plus 1 (for the component of $\G$).  The upper
bound is similar but simpler.  Thus
\begin{equation}\label{seb1}
\begin{split}
\log Z^\free_{\L_{\ul n,\b}}&=
\log\mu_{\l,\d,\g}(q^{k^\free_{\L_{\ul n,\b}}})\\
&\left\{\begin{array}{l}
\leq \log Z^\free_{\L_{\ul m,\a}}
+\log Z^\free_{\L_{\ul n,\b}\setminus \L_{\ul m,\a}}\\
\geq  \log Z^\free_{\L_{\ul m,\a}}
+\log Z^\free_{\L_{\ul n,\b}\setminus \L_{\ul m,\a}}-\\
\quad-(\log q) |\ul m|-
\l(1-1/q)\a d|\ul m|\sum_{i=1}^d\frac{1}{m_i}-\log q.
\end{array}\right.
\end{split}
\end{equation}
We have used the fact that
\[
|T^{\ul n,\b}_{\ul m,\a}|\leq 
\a d|\ul m|\sum_{i=1}^d\frac{1}{m_i}.
\]

There are 
$\prod_{i=1}^d\lfloor n_i/m_i \rfloor\cdot \lfloor\b/\a\rfloor$
`copies' of $\L_{\ul m,\a}$ in $\L_{\ul n,\b}$, each being
a translation of $\L_{\ul m,\a}$ by a vector
\[
\ul{l}\in\{(b_1m_1,\dotsc,b_dm_d,c\a):
b_i=1,\dotsc,\lfloor n_i/m_i \rfloor,\:
c=1,\dotsc,\lfloor \b/\a \rfloor \}.
\]
Write
\begin{equation}
\L=\Big(\bigcup_{\ul{l}} (\L_{\ul{m},\a}+\ul{l})\Big)\cup \L';
\end{equation}
this union is disjoint up to a set of measure zero.
Let $\L'=(K',F')$.  Repeating the argument
leading up to~\eqref{seb1} once for each `copy' of $\L_{\ul m,\b}$
we deduce that $Z^\free_{\L_{\ul n,\b}}$
is bounded above by
\begin{equation}\label{abc1}
\Big(\prod_{i=1}^d\lfloor n_i/m_i \rfloor\cdot \lfloor\b/\a\rfloor
\Big)\log Z^\free_{\L_{\ul m,\a}}+\log Z^\free_{\L'}
\end{equation}
and below by the same quantity~\eqref{abc1} minus
\begin{equation}
\prod_{i=1}^d\lfloor n_i/m_i \rfloor\cdot \lfloor\b/\a\rfloor
\Big((\log q) |\ul m|+
\l(1-1/q)\a d|\ul m|\sum_{i=1}^d\frac{1}{m_i}+\log q\Big).
\end{equation}
We will prove shortly that 
\begin{equation}\label{seb2}
\lim_{n_i,\b\rightarrow\oo} \frac{1}{|\L_{\ul n,\b}|}
\log Z^\free_{\L'}=0;
\end{equation}
once this is done it follows on
dividing by $|\L_{\ul n,\b}|=\b\cdot|\ul n|$ and letting all
$n_i,\b\rightarrow\oo$ that
\begin{equation}
\begin{split}
\frac{1}{|\L_{\ul m,\a}|}\log Z^\free_{\L_{\ul m,\a}}&\leq
\liminf_{n_i,\b\rightarrow\oo} P^\free_{\L_{\ul n,\b}}
\leq \limsup_{n_i,\b\rightarrow\oo} P^\free_{\L_{\ul n,\b}}\\
&\leq \frac{1}{|\L_{\ul m,\a}|}\log Z^\free_{\L_{\ul m,\a}}
+\frac{1}{\a}\log q+\\
&\qquad+\l(1-1/q)d\sum_{i=1}^d\frac{1}{m_i}
+\frac{1}{|\L_{\ul m,\a}|}\log q,
\end{split}
\end{equation}
and hence that $\lim_\L P^\free_\L$ exists and is finite.

Let us prove the claim~\eqref{seb2}.  The set $K_{\L'}$ consists
of a number of disjoint intervals, of which
\[
\prod_{i=1}^d m_i\lfloor n_i/m_i\rfloor
\]
have length $\b-\a\lfloor\b/\a\rfloor$, and 
\[
\prod_{i=1}^d n_i-
\prod_{i=1}^d m_i\lfloor n_i/m_i\rfloor
\]
have length $\b$.  The number $k^\free_{\L'}$ of components
is bounded above by the sum over all such intervals $L$ of
$|D\cap L|+2$ (we have added $1$ for the component of $\G$).  Hence
\begin{equation}
\begin{split}
0\leq \log Z^\free_{\L'}&=\mu_{\l,\d,\g}(q^{k^\free_{\L'}})\\
&\leq \Big(\prod_{i=1}^d m_i\lfloor n_i/m_i\rfloor\Big)\cdot
(q-1)\d (\b-\a\lfloor \b/\a\rfloor)+\\
&\qquad+\Big(\prod_{i=1}^d n_i-
\prod_{i=1}^d m_i\lfloor n_i/m_i\rfloor\Big)\cdot
(q-1)\d\b+2\log q.
\end{split}
\end{equation}
Equation~\eqref{seb2} follows.

Finally, we must prove convergence with arbitrary boundary
condition.  It is clear that for any boundary condition $b$
we have
\[
k^\wired_\L\leq k^b_\L\leq k^\free_\L.
\]
On the other hand 
\[
k^\wired_\L\geq k^\free_\L-2|\ul n|-|D\cap \partial\L|-1.
\]
The result follows.
\end{proof}

We now switch parameters to $r,s,t,u$, given in~\eqref{seb3}.
For fixed $u$ (i.e. fixed $q$) let $\cD_u=\cD_q$ be the set of
points $(r,s,t)\in\RR^3$ at which at least one of the partial
derivatives 
\[
\frac{\partial P}{\partial r}, \quad
\frac{\partial P}{\partial s}, \quad
\frac{\partial P}{\partial t}
\]
fails to exist.  Since $P$ is convex, $\cD_q$ has zero (three-dimensional)
Lebesgue measure.  By general properties of convex functions,
the partial derivatives 
\[
\frac{\partial P^b_\L}{\partial r}, \quad
\frac{\partial P^b_\L}{\partial s}, \quad
\frac{\partial P^b_\L}{\partial t}
\]
converge to the corresponding derivatives of $P$ whenever
$(r,s,t)\not\in\cD_q$, for any $b$.  Now observe that
\begin{equation}
\begin{split}
\frac{\partial P^\free_\L}{\partial r}&=
\frac{1}{|\L|}\phi^\free_\L(|B_\L|)\leq
\frac{1}{|\L|}\phi^\free(|B_\L|)\\
&\leq\frac{1}{|\L|}\phi^\wired(|B_\L|)\leq
\frac{1}{|\L|}\phi^\wired_\L(|B_\L|)=
\frac{\partial P^\wired_\L}{\partial r},
\end{split}
\end{equation}
so if $(r,s,t)\not\in\cD_q$ then
\begin{equation}\label{seb4}
\lim_{\L\uparrow\LLam}\frac{1}{|\L|}\phi^\free(|B_\L|)=
\lim_{\L\uparrow\LLam}\frac{1}{|\L|}\phi^\wired(|B_\L|)=
\frac{\partial P}{\partial r}.
\end{equation}
Recall from Lemma~\ref{invar_lem} that $\phi^\free$
and $\phi^\wired$ are both invariant under translations.  
The set $B$ is a point process on $\FF$, which
is therefore stationary under both $\phi^\free$
and $\phi^\wired$, and hence has constant intensities under
these measures.  Said
another way, the \emph{mean measures} $m^\free,m^\wired$
on $(\FF,\cB(\FF))$, given respectively by
\[
m^\free(F):=\phi^\free(|B\cap F|),\quad\text{and}\quad
m^\wired(F):=\phi^\wired(|B\cap F|)
\]
are translation invariant measures.  It is therefore a general fact
that there are constants $c^\free_{\mathrm{b}}$ and
$c^\wired_{\mathrm{b}}$ such that for all regions $\L=(K,F)$,
\[
m^\free(F)=\phi^\free(|B_\L|)=c^\free_{\mathrm{b}}|F|,\quad\text{and}\quad
m^\wired(F)=\phi^\wired(|B_\L|)=c^\wired_{\mathrm{b}}|F|,
\]
where $|\cdot|$ denotes Lebesgue measure.  Similarly, there are
constants $c^\free_{\mathrm{d}}$, $c^\wired_{\mathrm{d}}$,
$c^\free_{\mathrm{g}}$ and $c^\wired_{\mathrm{g}}$ such that
\[
\phi^\free(|D_\L|)=c^\free_{\mathrm{d}}|K|,\quad\text{and}\quad
\phi^\wired(|D_\L|)=c^\wired_{\mathrm{d}}|K|,
\]
and
\[
\phi^\free(|G_\L|)=c^\free_{\mathrm{g}}|K|,\quad\text{and}\quad
\phi^\wired(|G_\L|)=c^\wired_{\mathrm{g}}|K|,
\]
for all regions $\L=(K,F)$.

Note that 
\[
\lim_{n_i,\b\rightarrow\oo} \frac{|F_{\ul n,\b}|}{|K_{\ul n,\b}|}=d.
\]
It follows from~\eqref{seb4}, and similar calculations
for $D$ and $G$, that
\begin{equation}
c^\free_{\mathrm{b}}=c^\wired_{\mathrm{b}},\quad
c^\free_{\mathrm{d}}=c^\wired_{\mathrm{d}},\quad\text{and}
\quad c^\free_{\mathrm{g}}=c^\wired_{\mathrm{g}}
\quad\text{whenever } (r,s,t)\not\in\cD_q.
\end{equation}
Recall the condition given at the end of
Section~\ref{rc_wl_sec} for the uniqueness of the 
infinite-volume random-cluster measures,
namely that $\phi^\free=\phi^\wired$.  We will use the facts listed
above to prove

\begin{theorem}\label{rc_uniq_thm}
There is a unique random-cluster measure, in that
$\phi^\free=\phi^\wired$, whenever $(r,s,t)\not\in\cD_q$.
\end{theorem}
The corresponding results holds when $\g\geq0$ is fixed, in that
$\phi^\free=\phi^\wired$ except on a set of points
$(r,s)$ of zero (two-dimensional) Lebesgue measure.
For also $\d>0$ fixed, the corresponding set of $\l$
where uniqueness fails is countable, again by general properties of
convex functions.  Presumably this latter
set consists of a single point, namely the point corresponding
to $\rho=\rho_\crit$, but this has not been proved even for the
discrete models.

\begin{proof}
Since $\phi^\wired\geq\phi^\free$, there is by Theorem~\ref{strassen_thm} 
a coupling $\PP$
of the two measures such that 
\[
\PP(\{(\om^\wired,\om^\free)\in\Om^2:\om^\wired\geq\om^\free\})=1,
\]
and such that $\om^\wired$ and $\om^\free$ have marginal distributions
$\phi^\wired$ and $\phi^\free$ under $\PP$, respectively.  
Write $B^b$, $b\in\{\free,\wired\}$ for the bridges of $\om^b$, and
similarly for deaths and ghost-bonds.
Let $A\in\cF_\L$ be an increasing event.  Then
\begin{equation}
\begin{split}
0\leq \phi^\wired(A)-\phi^\free(A)&\leq
\PP(\om^\wired\in A,\om^\wired\neq\om^\free\text{ in }\L)\\ &\leq
\PP(|B^\wired_\L\setminus B^\free_\L|+|D^\free_\L\setminus D^\wired_\L|
+|G^\wired_\L\setminus G^\free_\L|)\\
&=\phi^\wired(|B_\L|)-\phi^\free(|B_\L|)+
\phi^\free(|D_\L|)-\phi^\wired(|D_\L|)+\\
&\qquad +\phi^\wired(|G_\L|)-\phi^\free(|G_\L|)\\
&=|\L|(c^\wired_{\mathrm{b}}-c^\free_{\mathrm{b}}+
c^\free_{\mathrm{d}}-c^\wired_{\mathrm{d}}+
c^\wired_{\mathrm{g}}-c^\free_{\mathrm{g}})\\
&=0,
\end{split}
\end{equation}
so $\phi^\wired=\phi^\free$ as required.
\end{proof}

Here is a consequence when $\g=0$.  Recall that we set
$\l=\rho$ and $\d=1$.  Suppose $0<\rho<\rho'$ are given.
We may pick $\l_1=\rho_1$ so that $\rho<\rho_1<\rho'$ 
and so that there is a unique 
infinite-volume measure with parameters $\l_1=\rho_1,\d_1=1$ 
and $\g=0$.   Hence
\begin{equation}\label{seb5}
\phi^\wired_\rho\leq\phi^\wired_{\rho_1}=
\phi^\free_{\rho_1}\leq\phi^\free_{\rho'}.
\end{equation}
It follows that the critical values $\rho^\free_\crit(q)$
and $\rho^\wired_\crit(q)$ of Definition~\ref{crit_def}
are equal for all $q\geq1$.

\section{Duality in $\ZZ\times\RR$}\label{duality_sec}

In this section we let $\LL=\ZZ$.
Thanks to the notion of planar duality for graphs, much more is 
known about the discrete random-cluster model in two dimensions
than in general dimension.  In particular, the critical value
for $q=1,2$ and $q\geq 25.72$ has been calculated in two dimensions,
see~\cite{abf,kesten_half,laanait_etal_91,laanait_etal_86}.
In the space-time setting, the $d=1$ model occupies the 
two-dimensional space $\ZZ\times\RR$, so we may adapt duality arguments
to this case;  that is the objective of this section.  Such arguments
have been applied when $q=1$ to prove that $\rho_\crit(1)=1$, 
see~\cite{bezuidenhout_grimmett}.  We will see in Chapter~\ref{appl_ch}
that $\rho_\crit(2)=2$, and Theorem~\ref{zhang_thm} in
the present section is a first step towards this result.

Throughout this section we assume that $\g=0$, and hence
suppress reference to both $\g$ and $G$.  
We also assume that $q\geq1$ and 
that $\l,\d$ are positive constants.  In light 
of Theorem~\ref{crit_nontriv} 
we may disregard the $\b<\oo$ case, hence we deal in
this section only with the $\b\rightarrow\oo$ case.
We think of $\LLam\equiv\ZZ\times\RR$ as embedded in $\RR^2$
in the natural way.

\begin{figure}[hbt]
\centering
\includegraphics{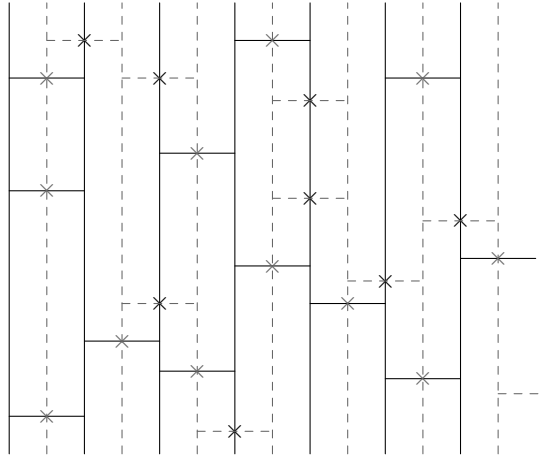}
\caption[Planar duality]
{An illustration of duality.
The primal configuration $\om$ is drawn solid black, the
dual $\om_\dual$ dashed grey.}\label{duality_fig}   
\end{figure}

We write $\LL_\dual$ for $\ZZ+1/2$;  of course $\LL$ and $\LL_\dual$
are isomorphic graphs.  With any $\om=(B,D)\in\Om$ we associate
the `dual' configuration $\om_\dual:=(D,B)$%
\nomenclature[o]{$\om_\dual$}{Dual configuration}
 regarded as a configuration
in $\LLam_\dual=\LL_\dual\times\RR$.  Thus each bridge in $\om$
corresponds to a death in $\om_\dual$, and each death in $\om$
corresponds to a bridge in $\om_\dual$.
This correspondence is illustrated
in Figure~\ref{duality_fig}.  We identify $\om_\dual=(D,B)$ with the
element $(D-1/2,B-1/2)$ of $\Om$.  Under this identification we may
for any measurable $f:\Om\rightarrow\RR$ define $f_\dual:\Om\rightarrow\RR$
by $f_\dual(\om)=f(\om_\dual)$.

In the case when $q=1$ it is clear that for any measurable function
$f:\Om\rightarrow\RR$, we have the relation
$\mu_{\l,\d}(f_\dual)=\mu_{\d,\l}(f)$, since the roles of $\l$ and
$\d$ are swapped under the duality transformation.  We will see that
a similar result holds when $q>1$.
\begin{definition}
Let $\psi_1,\psi_2$ be probability measures on $(\Om,\cF)$.  
We say that $\psi_2$ is dual to $\psi_1$ if for all 
measurable $f:\Om\rightarrow\RR$ we have that 
\begin{equation}\label{dual_meas_cond}
\psi_1(f_\dual)=\psi_2(f).
\end{equation}
\end{definition}
Thus the dual of $\mu_{\l,\d}$ is $\mu_{\d,\l}$.
Clearly it is enough to check~\eqref{dual_meas_cond} on some determining
class of functions, such as the local functions.

It will be convenient in what follows to denote the free and wired
random-cluster measures on a  region $\L$ by 
$\phi^0_{\L;q,\l,\d}$%
\nomenclature[f]{$\phi^0$}{Free random-cluster measure}
 and $\phi^1_{\L;q,\l,\d}$%
\nomenclature[f]{$\phi^1$}{Wired random-cluster measure}
 respectively, instead of $\phi^\free_{\L;q,\l,\d}$ and 
$\phi^\wired_{\L;q,\l,\d}$.  The following result is stated in terms
of infinite-volume measures, but from the proof we see that an
analogous result holds also in finite volume.
\begin{theorem}\label{duality_thm}
Let $b\in\{0,1\}$.  The dual of the measure $\phi^b_{q,\l,\d}$
is $\phi^{1-b}_{q,q\d,\l/q}$.
\end{theorem}
\begin{proof}
Fix $\b>0$ and $q\geq 1$;  later we will let
$\b\rightarrow\infty$.    
We write $[m,n]$ for the graph $L\subseteq\LL$ induced by the set 
$\{m,m+1,\dotsc,n\}\subseteq\ZZ$ and $\L_{m,n}=(K_{m,n},F_{m,n})$ 
for the  corresponding simple region.  We write $\phi^b_{m,n;\l,\d}$ for the
random-cluster measure on the region $\L_{m,n}$, with similar adjustments
to other notation.

In what follows it will be useful to
restrict attention to the bridges and deaths of $\om\in\Om$ 
that fall in $\L_{m,n}$ only.  It is then most natural to consider
only those (dual) bridges and deaths of $\om_\dual$ that fall in
$\L_{m,n-1}+1/2$.  In line with this we define
\begin{equation}
B_{m,n}(\om):=B(\om)\cap F_{m,n},\qquad
D_{m,n}(\om):=D(\om)\cap K_{m,n};
\end{equation}
and for the dual
\begin{equation}
B_{m,n-1}(\om_\dual):=D(\om)\cap K_{m+1,n-1},\qquad
D_{m,n-1}(\om_\dual):=B(\om)\cap F_{m,n}.
\end{equation}

The first step is to establish an analog of the Euler
equation for planar graphs.  We claim that
\begin{equation}\label{tjo1}
k^1_{m,n}(\om)-k^0_{m,n-1}(\om_\dual)+
|B_{m,n}(\om)|-|D_{m,n}(\om)|=1-n+m.
\end{equation}
(A similar result was obtained in~\cite[Lemma~3.3]{aizenman_nacht}.)
This is best proved inductively by successively adding elements to
the sets $B_{m,n}(\om)$ and $D_{m,n}(\om)$.  If both sets are empty,
the claim follows on inspection.  For each bridge you add to 
$B_{m,n}(\om)$, either $k^1_{m,n}(\om)$ decreases by one 
or $k^0_{m,n-1}(\om_\dual)$ increases by one, but never both.  
Similarly when you add deaths
to $D_{m,n}(\om)$,  either $k^1_{m,n}(\om)$ increases by one 
or $k^0_{m,n-1}(\om_\dual)$ decreases by one for each death, but never both.  
That establishes~\eqref{tjo1}.

Let $\mu_{m,n;\l,\d}$ denote the percolation measure restricted to 
$\L_{m,n}$.  For $f:\Om\rightarrow\RR$ any $\cF_{\L_{m,n-1}}$-measurable,
bounded and continuous function, we have, using~\eqref{tjo1}, that
\begin{equation}\label{tjo2}
\begin{split}
\phi^1_{m,n;\l,\d}(f_\dual)&\propto
\int d\mu_{m,n;\l,\d}(\om) q^{k^1_{m,n}(\om)}f(\om_\dual)\\
&\propto \int d\mu_{m,n;\l,\d}(\om) q^{k^0_{m,n-1}(\om_\dual)}
q^{|D_{m,n}(\om)|}q^{-|B_{m,n}(\om)|}f(\om_\dual)\\
&\propto \int d\mu_{m,n-1;\d,\l}(\om_\dual) q^{k^0_{m,n-1}(\om_\dual)}
q^{|B_{m,n-1}(\om_\dual)|}q^{-|D_{m,n-1}(\om_\dual)|}f(\om_\dual)\\
&\propto \int d\mu_{m,n-1;q\d,\l/q}(\om_\dual) 
q^{k^0_{m,n-1}(\om_\dual)}f(\om_\dual)\\
&\propto\phi^0_{m,n-1;q\d,\l/q}(f).
\end{split}
\end{equation}
We have used the fact that 
\begin{equation}
\frac{d\mu_{m,n-1;q\d,\l/q}}{d\mu_{m,n-1;\d,\l}}(\om)
\propto q^{|B_{m,n-1}(\om)|}q^{-|D_{m,n-1}(\om)|},
\end{equation}
a simple statement about Poisson processes.

Since both sides of~\eqref{tjo2} are probability measures,
it follows that
\begin{equation}\label{tjo3}
\phi^1_{m,n;\l,\d}(f_\dual)=\phi^1_{m,n-1;q\d,\l/q}(f).
\end{equation}
Letting $m,n,\b\rightarrow\infty$ in~\eqref{tjo3} and using
Theorem~\ref{conv_lem}, the result follows.
\end{proof}

Note that if $\l/\d=\rho$ then the corresponding ratio for
the dual measure is $q\d/(\l/q)=q^2/\rho$.  We therefore say
that the space--time random-cluster 
model is \emph{self-dual} if $\rho=q$.  
This self-duality was referred to 
in~\cite[Proposition~3.4]{aizenman_nacht}.

\subsection{A lower bound on $\rho_c$ when $d=1$}

In this section we adapt Zhang's famous and versatile argument 
(published in~\cite[Chapter~6]{grimmett_rcm}) to the space-time setting.
See~\cite{bezuidenhout_grimmett} for the special case of this
argument when $q=1$.

\begin{theorem}\label{zhang_thm}
If $d=1$  and $\rho=q$ then $\theta^\free(\rho,q)=0$;  hence the critical
ratio $\rho_\crit\geq q$. 
\end{theorem}
\begin{proof}
Assume  for a contradiction that with $\rho=q$ we have that
$\theta^\free(\rho,q)>0$.  Then by Theorem~\ref{inf_clust_uniq} 
there is almost surely a unique unbounded component in $\om$
under $\phi^\free$.  It follows from self-duality and the fact that
$\theta^\wired\geq\theta^\free$ that there is almost surely also 
a unique unbounded 
component in $\om_\dual$.  Let $D_n=\{(x,t)\in\RR^2:|x+1/2|+|t|\leq n\}$ be 
the `lozenge', and $D_n^\dual=\{(x,t)\in\RR^2:|x|+|t|\leq n\}$ its 
`dual', as in Figure~\ref{zhang_fig}. 
\begin{figure}[hbt]
\centering
\includegraphics{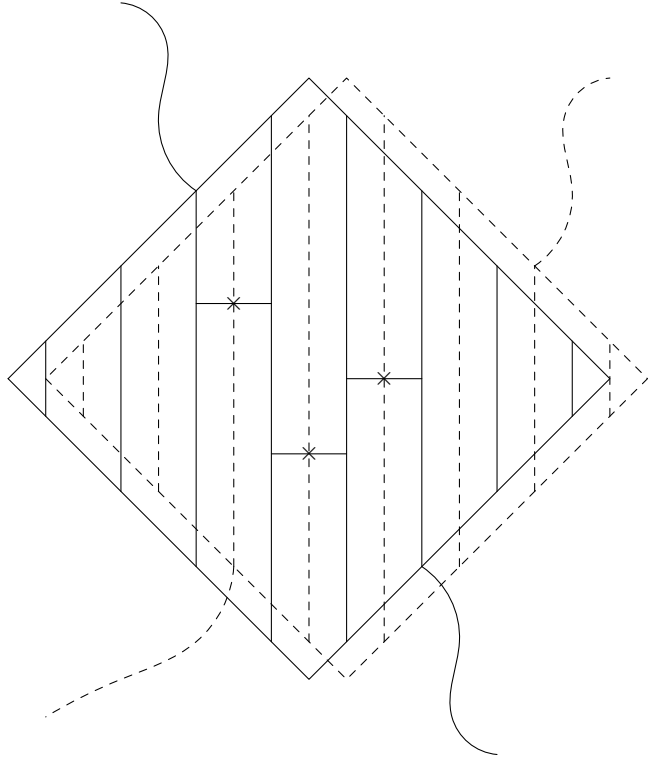}
\caption[Zhang's argument]
{On the event $A_2\cap A_4\cap A^\dual_1\cap A^\dual_3$
either the unbounded primal cluster breaks into 2 parts, or the dual 
one does.}\label{zhang_fig}   
\end{figure}
Number the four sides of each of  $D_n$ and $D_n^\dual$ counterclockwise,
starting in each case with the north-east side.  For $i=1,\dotsc,4$ let $A_i$
be the event that the $i$th side of $D_n$ is attached to an unbounded path
of $\om$, which does not otherwise intersect $D_n$.  
Similarly let $A^\dual_i$ be the event that the $i$th side of the
dual $D_n^\dual$ is attached to an unbounded path of $\om_\dual$.  
Clearly $\phi^\free(\cup_{i=1}^4 A_i)\rightarrow1$ as $n\rightarrow\infty$.  
However, all the $A_i$ are
increasing, and by symmetry under reflection they
carry equal probability.  It follows from positive
association, Theorem~\ref{rc_fkg}, that
\begin{equation*}
\phi^\free(\cup_{i=1}^4 A_i)\leq 1-\phi^\free(A_2^c)^4=1-(1-\phi^\free(A_2))^4,
\end{equation*}
and hence $\phi^\free(A_2)\rightarrow1$ too.  Hence 
for $n$ large enough we have
that $\phi^\free(A_2)=\phi^\free(A_4)\geq 5/6$, so 
by positive association again
$\phi^\free(A_2\cap A_4)\geq (5/6)^2>5/8$ for $n$ large enough.  
In the same way
it follows that for large $n$ we have 
$\phi^\free(A^\dual_1\cap A^\dual_3)>5/8$.  But then 
\begin{equation*}
\phi^\free(A_2\cap A_4\cap A^\dual_1\cap A^\dual_3)
\geq\frac{10}{8}-1=\frac{1}{4}.
\end{equation*}
Now a glance at Figure~\ref{zhang_fig} should convince the reader that
this contradicts the uniqueness of the unbounded cluster, either in $\om$
or $\om_\dual$.  This contradiction shows that $\theta^\free(\rho,q)=0$ as
required.
\end{proof}
\begin{remark}
It is natural to suppose that the critical value equals the self-dual
value $\l/\d=q$.  For $q=1$ this is proved in~\cite{bezuidenhout_grimmett} 
and in~\cite{aizenman_jung};  for $q=2$ it is proved in 
Theorem~\ref{crit_val_cor} (see also~\cite{bjogr2}).
\end{remark}

\section{Infinite-volume Potts measures}\label{inf_potts_sec}

Using the convergence results in Section~\ref{inf_rc_sec}, we will
in this section construct infinite-volume weak limits  of 
Potts measures.  We will also provide more details about
uniqueness of infinite-volume measure in the space-time Ising
model, extending in that case the arguments of 
Section~\ref{press_sec}.
The results in this section will form the foundation
for our study of the quantum Ising model in Chapter~\ref{qim_ch}.

\subsection{Weak limits of Potts measures}

Let $q\geq 2$ be an integer, and let $\a_\G=q$;
we will suppress reference to the simple 
boundary condition $(b,\a)$ throughout this subsection.
Recall the two random-cluster measures
$\phi_\L^\wired$ and $\phi_\L^\free$ as well as their Potts counterparts
$\pi_\L^\wired$ and $\pi_\L^\free$, connected via the 
coupling~\eqref{es_def}.  
For simplicity we assume in this section that $\LL=\ZZ^d$ for 
some $d\geq 1$;  similar arguments are valid in greater generality,
but we do not pursue this here.  All regions in this
section will be simple, as in~\eqref{def-oldL}.
We let $\L_n=(K_n,F_n)$ denote a
strictly increasing sequence of simple 
regions, containing the origin and increasing to either
$\LLam$ or $\LLam_\b$.  Denote by
$\phi_n^\wired$, $\phi_n^\free$, $\pi_n^\wired$ and $\pi_n^\free$
the corresponding random-cluster and Potts measures.
Proofs will be given for the $\b=\infty$ case,  the case
$\b<\oo$ is similar.

Throughout this subsection we will be making use of the
concept of \emph{lattice components}:  given 
$\om=(B,D,G)$ the lattice components of $\om$ are the connected
components in $\KK$ 
of the configuration $(B,D,\varnothing)$.  We will think
of the points in $G$ as green points, and of any lattice component
containing an element of $G$ as  \emph{green}.
In this subsection we will only use the notation $x\lra y$ to mean that
$x,y$ lie in the same \emph{lattice} component.  We write
$C_x(\om)$ for the lattice component of $x$ in~$\om$.

The following convergence result is an adaptation of arguments
in~\cite{ACCN}, see also~\cite[Theorem~4.91]{grimmett_rcm}.
\begin{theorem}\label{potts_lim_thm}
The weak limits 
\begin{equation}
\pi^\free=\lim_{n\rightarrow\oo}\pi^\free_n
\qquad\mbox{and}\qquad
\pi^\wired=\lim_{n\rightarrow\oo}\pi^\wired_n
\end{equation}
exist and are independent of the manner in which
$\L_n\uparrow\LLam$.  Moreover, $\pi^\free$ and $\pi^\wired$ are given
as follows:
\begin{itemize}
\item Let $\om\sim\phi^\free$ and assign to each green component 
of $\om$ spin $q$, and assign to the remaining
components uniformly independent spins from $1,\dotsc,q$;  then the
resulting spin configuration has law $\pi^\free$.
\item Let $\om\sim\phi^\wired$ and assign to each unbounded component 
and each green component of $\om$ spin $q$, and assign to the remaining
components uniformly independent spins from $1,\dotsc,q$;  then the
resulting spin configuration has law~$\pi^\wired$.
\end{itemize}
\end{theorem}
\begin{proof}
We will make use of a certain total order
on $\KK=\ZZ^d\times\RR$.  The precise details are not important,
except that the ordering be such that every (topologically)
closed set contains an earliest point.  We define such an
ordering as follows.  We say that 
$x=(x_1,\dotsc,x_d,t)<(x_1',\dotsc,x_d',t')=x'$ if (a)
for $k\in\{1,\dotsc,d\}$ minimal with $x_kx'_k<0$ we have $x_k>0$;
or if (a) fails but (b) $tt'<0$ with $t>0$;
or if (a) and (b) fail but (c)  
$|x|<|x'|$ lexicographically, where $|x|=(|x_1|,\dotsc,|x_d|,|t|)$.

Slightly different arguments are required for the two
boundary conditions.  We give the argument only for free boundary.
It will be necessary to modify the probability space
$(\Om,\cF)$, as follows (we omit some details).  
For each $n\geq 1$ and each  $\om=(B,D,G)\in\Om$, let 
$\tilde\om_n=(\tilde B_n,\tilde D_n,\tilde G_n)$
be given by
\[
\tilde B_n=B\cap F_n, \quad
\tilde D_n=(D\cap K_n)\cup(\KK\sm K_n), \quad
\tilde G_n=G\cap K_n. 
\]
Thus, in $\tilde\om_n$, no two points in $\KK\setminus K_n$ are
connected.
Let $\tilde\Om=\Om\cup\{\tilde\om_n:\om\in\Om,n\geq 1\}$, and define 
connectivity in elements of $\tilde\Om$ in the obvious way.  Define the
functions $V_x$ as before Theorem~\ref{conv_lem};  if 
$x\in\KK\times\{\mathrm{d}\}$ then $V_x$ may now take the value $+\oo$.
Let $\tilde\cF$ denote the $\s$-algebra generated by the $V_x$'s.
(Alternatively, $\tilde\cF$ is the $\s$-algebra generated by the 
appropriate Skorokhod metric when the associated step functions
are allowed to take the values $\pm\oo$.)  Let $\tilde\phi^\free_n$
denote the law of $\tilde\om_n$ when $\om$ has law $\phi^\free_n$.
Note that the number of components of $\tilde\om_n$ equals
$k^\free_n(\om)$.  

Extending the partial order on $\Om$ to $\tilde\Om$ in the natural
way, we see that for each $n$ we have 
$\tilde\phi^\free_n\leq\tilde\phi^\free_{n+1}$.
(It is here that we need to use $\tilde\Om$, since the stochastic 
ordering  $\phi^\free_n\leq\phi^\free_{n+1}$ holds only on $\cF_n$, 
not on the  full $\s$-algebra $\cF$.)
Hence there exists by Strassen's Theorem~\ref{strassen_thm} 
a probability measure $P$ on $(\tilde\Om^\NN,\tilde\cF^\NN)$ such that 
in the sequence $(\tilde\om_1,\tilde\om_2,\dotsc)$ the $n$th
component has marginal distribution $\tilde\phi^\free_n$, and such that
$\tilde\om_n\leq\tilde\om_{n+1}$ for all $n$, with $P$-probability one.
The sequence $\tilde\om_n$ increases to a limiting configuration
$\tilde\om_\infty$, which has law $\phi^\free$.  We have
that $\phi^\free(\Om)=1$.

For each fixed (bounded) region $\D$,
if $n$ is large enough then $\tilde\om_n$ agrees with $\tilde\om_\infty$
throughout $\D$.  Let $\L$ be a fixed region, and let
$\D=\D(\tilde\om_\oo)\supset\L$ be large enough so that
the following hold:
\begin{enumerate}
\item Each bounded lattice-component of $\tilde\om_\infty$ which
 intersects $\L$ is entirely contained in $\D$;
\item Any two points $x,y\in\L$ which are connected in
 $\tilde\om_\oo$ are connected inside $\D$;
\item Any lattice-component of $\tilde\om_\oo$ which is green
 has a green point inside $\D$.
\end{enumerate}
It is (almost surely) possible to choose such a $\D$ because
only finitely many lattice components intersect $\L$.  
We choose $n=n(\tilde\om_\oo)$ large enough so that $\tilde\om_n,\tilde\om_{n+1},\dotsc$
all agree with $\tilde\om_\oo$ throughout~$\D$.  

\emph{Claim:}  for all $x,y\in\L$, we have that $x\lra y$ in $\tilde\om_n$
if and only if $x\lra y$ in $\tilde\om_\oo$.  To see this, first note
that $C_x(\tilde\om_n)\subseteq C_x(\tilde\om_\oo)$ since $\tilde\om_n\leq\tilde\om_\oo$,
proving one of the implications.  Suppose now that $x\lra y$ in $\tilde\om_\oo$.
Then by our choice of $\D$, there is a path from $x$ to $y$ inside
$\D$.  But $\tilde\om_\oo$ and $\tilde\om_n$ agree on $\D$, so it follows that also
$x\lra y$ in~$\tilde\om_n$.

Let $\tilde\om\in\tilde\Om$, and let $C$ be a lattice component
of $\tilde\om$.  The (topological) closure of $C$ contains an earliest
point in the order defined above.  Order the lattice components
$C_1(\tilde\om),C_2(\tilde\om),\dotsc$ according to the earliest point
in their closure;  this ordering is almost surely well-defined
under any of $\tilde\phi^\free_n,\tilde\phi^\free$.
Note that the claim above implies that this ordering agrees
for those lattice components of $\tilde\om_n$ and $\tilde\om_\oo$
which intersect~$\L$.

Let $S_1,S_2,\dotsc$
be independent and uniform on $\{1,\dotsc,q\}$, and
define for $x\in\LLam$,
\begin{equation}
\tau_x(\tilde\om)=\left\{
\begin{array}{ll}
q, & \mbox{if } C_x(\tilde\om)\mbox{ is green}, \\
S_i, & \mbox{otherwise, where } C_x(\tilde\om)=C_i.
\end{array}\right.
\end{equation}
Then $\tau(\tilde\om_\oo)$ has the law $\pi^\free$ described 
in the statement of the 
theorem, and $\tau(\tilde\om_n)$ has the law $\pi^\free_n$ on events
in $\cG_\L$.  Moreover, from the claim it follows that 
$\tau_x(\tilde\om_\oo)=\tau_x(\tilde\om_n)$ for any $x\in\L$.
Hence for all continuous, bounded $f$, measurable
with respect to $\cG_\L$, we have that
$f(\tau(\tilde\om_n))\rightarrow f(\tau(\tilde\om_\oo))$
almost surely.  It follows from the bounded convergence theorem
that 
\begin{equation}
\pi^\free_n(f)=E(f(\tau(\tilde\om_n)))
\rightarrow E(f(\tau(\tilde\om_\oo)))
=\pi^\free(f).
\end{equation}
Since such $f$ are convergence determining it follows that
$\pi^\free_n\Rightarrow\pi^\free$.
\end{proof}

\begin{remark}\label{inf_vol_rk}
From the representation given in Theorem~\ref{potts_lim_thm}
it follows that the correlation/connectivity relation
of Proposition~\ref{corr_conn_prop} holds also for infinite-volume
random-cluster and Potts measures.  In particular, when
$q=2$, it follows (using the obvious notation) that 
the analogue of~\eqref{ising_corr_conn} holds, namely
\[
\el\s_x\s_y\er^b=\phi^b(x\lra y),
\]
for $b\in\{\free,\wired\}$.
Note also that when $\g=0$ then, as in Proposition~\ref{corr_conn_prop},
we have for for $b\in\{\free,\wired\}$ that
\begin{equation}\label{corr_conn_inf}
\el\s_x\er^b=\phi^b(x\lra\oo).
\end{equation}
\end{remark}

\subsection{Uniqueness in the Ising model}\label{ising_uniq_sec}

We turn our attention  now to 
the space--time Ising model on $\LL=\ZZ^d$ with
constant $\l,\d,\g$.  In this section we continue our discussion,
started in Section~\ref{press_sec}, about uniqueness of 
infinite-volume measures.  More information can be obtained
in the case of the Ising model, partly thanks to the so-called
\ghs-inequality which allows us to show the absence of a phase transition
when $\g\neq0$.  In contrast, using only results obtained via the
random-cluster representation one can say next to nothing about
uniqueness when $\g\neq0$ since there is no useful way of combining
a $+1$ external field with a $-1$ `lattice-boundary'.  
The arguments in this section follow very closely
those for the classical Ising model, as developed 
in~\cite{lebowitz_martin-lof} and~\cite{preston_ghs} 
(see also~\cite[Chapters IV and V]{ellis85:LD}).  We provide full
details for completeness.  

As remarked earlier, the Ising model admits
more boundary conditions than the corresponding random-cluster model.
It will therefore seem
like some of the arguments presented below repeat what was said
in Section~\ref{press_sec}.  It should be noted, however, 
that the arguments
in this section can deal with all boundary conditions that
occur in the Ising model.  
It will be particularly useful to consider the $+$ and $-$ boundary
conditions, defined as follows.%
\nomenclature[<]{$\el\cdot\er^\pm$}{Ising measure with $\pm$ boundary condition}
  Let $b=\{P_1,P_2\}$
where $P_1=\{\G\}$ and $P_2=\boundary\L$.  We define the
\emph{$+$ boundary condition} by letting $\a_1=\a_2=+1$;  
when $\g\geq 0$ this equals the wired
random-cluster boundary condition with $\a_\G=+1$.  We define the
\emph{$-$ boundary condition} by letting $\a_1=+1$ and $\a_2=-1$.
The measure $\el\cdot\er^-_\L$ does not have a
satisfactory random-cluster representation when $\g>0$.
(See~\cite{chayes_machta_redner} for an
in-depth treatment of some difficulties associated with
the graphical representation of the Ising model in an arbitrary
external field.)  
In line with physical terminology we will sometimes in this 
section refer to the measures $\el\cdot\er^{b,\a}_\L$ as `states'.

For simplicity of notation we will in this 
section replace $\l$ and $\g$ by $2\l$ and $2\g$ throughout.
We will be writing $Z^{b,\a}_\L$ for the Ising
partition function~\eqref{ising_pf}, which therefore becomes
\begin{equation}
Z^{b,\a}_\L=\int d\mu_\d(D)\,\sum_{\s\in\S^{b,\a}_\L(D)} 
\exp\Big(\int_F\l(e)\s_ede+\int_K\g(x)\s_x dx\Big).
\end{equation}
We will similarly write $P^{b,\a}_\L=(\log Z^{b,\a}_\L)/|\L|$.  
Thanks to Proposition~\ref{pfs}, the $P^{b,\a}_\L$ thus
defined converge to a function $P$ which is a multiple of the 
original $P$ in Theorem~\ref{press_converge_thm}. Straightforward
modifications of the argument in Theorem~\ref{press_converge_thm} let
us deduce that this convergence holds for all boundary
conditions $b$ of Ising-type.

We assume throughout this section that $\L=\L_n\uparrow\LLam$
in such a way that
\begin{equation}\label{vanhove}
\frac{|K_n\sm K_{n-1}|}{|K_n|}\rightarrow 0, 
\qq \mbox{as }n\rightarrow\infty,
\end{equation}
where $\L_n=(K_n,F_n)$ and
$|\cdot|$ denotes Lebesgue measure.  As previously, straightforward
modifications of the argument are valid when $\b<\oo$ is fixed
and $\L\uparrow\LLam_\b$.  

Here are some general facts about convex functions;  some facts
like these were already used in Section~\ref{press_sec}.
See e.g.~\cite[Chapter~IV]{ellis85:LD} for proofs.
Recall that for a function
$f:\RR\rightarrow\RR$, the left and right derivatives of 
$f$ are given respectively by
\begin{equation}
\frac{\partial f}{\partial \g^+}:=
\lim_{h\downarrow 0}\frac{f(\g+h)-f(\l)}{h}
\quad\mbox{and}\quad
\frac{\partial f}{\partial \g^-}:=
\lim_{h\downarrow 0}\frac{f(\g-h)-f(\l)}{-h}
\end{equation}
provided these limits exist.

\begin{proposition}\label{conv_prop}
Let $I\subseteq\RR$ be an open interval and
$f:I\rightarrow\RR$
a convex function;  also let $f_n:I\rightarrow\RR$ be
a sequence of convex functions.  Then
\begin{itemize}
\item The left and right derivatives of $f$ exist throughout $I$;
the right derivative is right-continuous and the left derivative
is left-continuous.
\item The derivative $f'$ of $f$ exists at all but countably
many points in $I$.
\item If all the $f_n$ are differentiable and $f_n\rightarrow f$
pointwise then the derivatives $f_n'$ converge to $f'$ whenever 
the latter exists.
\item If the $f_n$ are uniformly bounded above and below
then there exists a sub-sequence $f_{n_k}$  and a 
(necessarily convex) function $f$ such that $f_{n_k}\rightarrow f$
pointwise.
\end{itemize}
\end{proposition}

We will usually
keep $\l,\d$ fixed and regard $P=P(\g)$ as a function of $\g$, and similarly
for other functions.  Note that $P$ is an even function of $\g$:
we have for all $\g>0$ that $P^+_\L(-\g)=P^-_\L(\g)$,
and since the limit $P$ is independent of boundary condition
it follows that $P(-\g)=P(\g)$.

Let
\begin{equation}\label{diff1}
\bar M^{b,\a}_\L:=\frac{\partial P^{b,\a}_\L}{\partial \g}=
\frac{1}{|\L|}\int_\L dx\el\s_x\er^{b,\a}_\L,
\end{equation}
where we abuse notation to write $x\in\L$ (respectively, $|\L|$) 
in place of the
more accurate $x\in K$ (respectively, $|K|$).  Also let
\begin{equation}
M^{b,\a}_\L:=\el\s_0\er^{b,\a}_\L.
\end{equation}%
\nomenclature[M]{$M^{b,\a}_\L$}{Finite-volume magnetization}%
 Note that~\eqref{diff1} together with the first 
\gks-inequality~\eqref{gks_1_eq} imply that $P^\wired_\L$, and hence
also $P$, is increasing for $\g>0$ (and hence decreasing for $\g<0$).
Moreover, we see that
\begin{equation}\label{diff2}
\frac{\partial^2 P^\wired_\L}{\partial \g^2}=
\frac{1}{|\L|}\int_\L\int_\L dxdy\el\s_x;\s_y\er^\wired_\L\geq0,
\end{equation}
from the second \gks-inequality~\eqref{gks_2_eq}.  Thus $P$
is convex in $\g$.  

\begin{lemma}\label{pm_states}
The states $\el\cdot\er^+_\L$ and $\el\cdot\er^-_\L$ converge weakly
as $\L\uparrow\LLam$.  The limiting states 
$\el\cdot\er^+$ and $\el\cdot\er^-$ are independent of the
way in which $\L\uparrow\LLam$ and are translation invariant.
\end{lemma}

\begin{remark}
The convergence result for $+$ boundary follows from 
Theorem~\ref{potts_lim_thm} and Remark~\ref{inf_vol_rk},
since when $q=2$ the measure $\pi^\wired_\L$ there is precisely
the state $\el\cdot\er_\L^+$.  However, the result for $-$
boundary does not follow from that result
since the random-cluster representation
as employed there does not admit the spin at $\G$ to be
different from that at $\partial\L$.  (One would have to
condition on the event that, in the random-cluster model, the
boundary is disconnected from $\G$, and then one loses desired
monotonicity properties.)  
\end{remark}

In the proof of Lemma~\ref{pm_states} we will be applying
the \fkg-inequality, Lemma~\ref{ising_fkg}.  
For each $x\in\KK$, let $\nu'_x=(\s_x+1)/2$%
\nomenclature[n]{$\nu'_x$}{$(\s_x+1)/2$} and for $A\se\KK$ finite, 
write
\begin{equation}\label{mrk1}
\nu'_A=\prod_{x\in A} \nu'_x.
\end{equation}
Note that $\nu'_A=\one_S$, where $S$ is the event that $\s_x=+1$
for all $x\in A$.  This is an increasing event, and a continuity
set by Example~\ref{cty_ex}.  Similarly, if $\L\se\D$ are regions
and $T$ is the event that $\s=+1$ on $\D\setminus \L$, then
$T$ is an increasing event and a continuity set, also by
Example~\ref{cty_ex}.

\begin{proof}[Proof of Lemma~\ref{pm_states}]
It is easy to check that the
variables $\nu'_A$, as $A$ ranges over the finite subsets of $\KK$,
form a convergence determining class.  
By Lemma~\ref{potts_cond} and Lemma~\ref{ising_fkg}
we therefore see that for any regions $\L\subseteq\D$ we
have that
\begin{equation}
\el\nu'_A\er^+_\L=\el\nu'_A\mid\s\equiv+1\mbox{ on } \D\sm\L\er^+_\D
\geq \el\nu'_A\er^+_\D
\end{equation}
and
\begin{equation}
\el\nu'_A\er^-_\L=\el\nu'_A\mid\s\equiv-1\mbox{ on } \D\sm\L\er^-_\D
\leq \el\nu'_A\er^-_\D.
\end{equation}
Hence $\el\nu'_A\er^+_\L$ and $\el\nu'_A\er^-_\L$ converge for all
finite $A\subseteq\KK$, as required.
\end{proof}
The proof of Lemma~\ref{pm_states} shows in particular that
\begin{equation}\label{mo1}
\el\s_0\er^+_\L\downarrow\el\s_0\er^+
\qquad\mbox{and}\qquad \el\s_0\er^-_\L\uparrow\el\s_0\er^-,
\end{equation}
and indeed that all the $\el\s_A\er^\pm_\L$ converge to the
corresponding $\el\s_A\er^\pm$.
Recall that by convexity, the left and right derivatives of $P$
exist at all $\g\in\RR$.
\begin{lemma}\label{press_der_lem1}
For all $\g\in\RR$ we have that
\begin{equation}
\frac{\partial P}{\partial \g^+}=\el\s_0\er^+
\qquad\mbox{and}\qquad
\frac{\partial P}{\partial \g^-}=\el\s_0\er^-.
\end{equation}
\end{lemma}
\begin{proof}
As a preliminary step we first show that 
$\bar M^\pm_\L$ has the same infinite-volume
limit as $M^\pm_\L$, that is to say
\begin{equation}\label{s1}
\lim_{\L\uparrow\LLam}\bar M^\pm_\L=\el\s_0\er^\pm.
\end{equation}
We prove this in the case of $+$ boundary, the case of $-$ boundary
being similar.  First note that
\begin{equation}
\bar M^+_\L=\frac{1}{|\L|}\int_\L dx\el\s_x\er^+_\L\geq
\frac{1}{|\L|}\int_\L dx\el\s_x\er^+=\el\s_0\er^+,
\end{equation}
by~\eqref{mo1} and translation invariance.  Thus 
$\liminf_\L\bar M^+_\L\geq\el\s_0\er^+$.  
Next let $\eps>0$ and let $\L$ be
large enough so that $\el\s_0\er^+_\L\leq\el\s_0\er^++\eps$.  If 
$x\in\KK$ and $\D$ is large enough that the translated region
$\L+x\subseteq\D$ then
\begin{equation}
\el\s_x\er^+_\D\leq\el\s_x\er^+_{\L+x}=\el\s_0\er^+_\L\leq\el\s_0\er^++\eps.
\end{equation}
Let $\D':=\{x\in\D:\L+x\in\D\}$.  Then
\begin{equation}
\begin{split}
\bar M^+_\D&=\frac{1}{|\D|}\int_\D dx\el\s_x\er^+_\D\leq
\frac{1}{|\D|}\Big(\int_{\D'} dx\el\s_x\er^+_\D+|\D\sm\D'|\Big)\\& \leq
\frac{1}{|\D|}\Big(|\D'|\big(\el\s_x\er^++\eps\big)+|\D\sm\D'|\Big).
\end{split}
\end{equation}
It therefore follows from the assumption~\eqref{vanhove} that
$\limsup_\L \bar M^+_\L\leq\el\s_0\er^++\eps$, which 
gives~\eqref{s1}.

Next we claim that $\el\s_0\er^+$ and $\el\s_0\er^-$ are right-
and left continuous in $\g$, respectively.  
First consider $+$ boundary.  Then for $\g'>\g$, we have
for any $\L$ from Lemma~\ref{ising_mon} that
$\el\s_0\er^+_{\L,\g'}\geq\el\s_0\er^+_{\L,\g}$.  Thus
\begin{equation}\label{mo3}
\begin{split}
\el\s_0\er^+_\g&\leq\liminf_{\g'\downarrow\g}\el\s_0\er^+_{\g'}
\leq\limsup_{\g'\downarrow\g}\el\s_0\er^+_{\g'}\\
&\leq\limsup_{\g'\downarrow\g}\el\s_0\er^+_{\L,\g'}=
\el\s_0\er^+_{\L,\g}\xrightarrow[\L\uparrow\LLam]{}\el\s_0\er^+_\g.
\end{split}
\end{equation}
(We have used the fact that $\el\s_0\er^+_\L$ is  
continuous in $\g$.)  A similar calculation holds for
$-$ boundary.  

Now, by convexity of $P$, the right derivative 
$\frac{\partial P}{\partial\g^+}$ is right-continuous, and also
$\lim_\L \bar M^\pm_\L=\frac{\partial P}{\partial\g}$ whenever the right side
exists.  But it exists for all but countably many $\g$, so given $\g$ there is
a sequence $\g_n\downarrow\g$ such that 
$\frac{\partial P}{\partial\g}(\g_n)=\el\s_0\er^+_{\g_n}$ for all $n$,
and similarly for $-$ boundary.  The result follows.
\end{proof}

We say that there is a \emph{unique state} at $\g$ (or at $\l,\d,\g$)
if for all finite $A\se\KK$, the limit 
$\el\s_A\er:=\lim_\L\el\s_A\er^{b,\a}_\L$ 
exists and is independent of the boundary condition $(b,\a)$. 
Note that, by linearity, it is equivalent to require that all the 
limits $\el\nu'_A\er:=\lim_\L\el\nu'_A\er^{b,\a}_\L$ exist and are independent
of the boundary condition.  Alternatively, there
is a unique state if and only if the measures $\el\cdot\er^{b,\a}_\L$ 
all converge weakly to the same limiting measure.

\begin{lemma}\label{uniq_lem}
There is a unique state at $\g\in\RR$ if and only if $P$ is differentiable
at $\g$.  There is a unique state at any $\g\neq0$.
\end{lemma}
\begin{proof}
We have that 
\begin{equation}
f_A:=\sum_{x\in A}\nu'_x-\nu'_A
\end{equation}
is increasing in $\s$.  By the \fkg-inequality, Lemma~\ref{ising_fkg},
we have that $\el f_A\er^+_\L\geq\el f_A\er^-_\L$.
It follows on letting $\L\uparrow\LLam$, and using translation invariance
as well as Lemma~\ref{press_der_lem1}, that
\begin{equation}
0\leq \el\nu'_A\er^+-\el\nu'_A\er^-\leq
\frac{1}{2}\sum_{x\in A}(\el\s_x\er^+-\el\s_x\er^-)
=\frac{|A|}{2}\Big(\frac{\partial P}{\partial\g^+}-
\frac{\partial P}{\partial\g^-}\Big),
\end{equation}
where $|A|$ is the number of elements in $A$.
Hence $\el\nu'_A\er^+=\el\nu'_A\er^-$ whenever 
$\frac{\partial P}{\partial\g}$ exists.  Since 
$\el\nu'_A\er^-\leq\el\nu'_A\er^{b,\a}\leq\el\nu'_A\er^+$ for all $(b,\a)$
(a consequence of Lemma~\ref{ising_fkg}), the 
first claim follows. 

The next part makes use of the facts about
convex functions stated above;  this part of the argument originates
in~\cite{preston_ghs}.  Let $\g>0$, and use the free
boundary condition.  We already know that $P$ and each $P^\free_\L$
is convex.   The \ghs-inequality, which is standard for the 
classical Ising model and proved for the current model 
in Lemma~\ref{ghs_lem}, implies that each $\bar M^\free_\L$ has
non\emph{positive} second derivative for $\g>0$, and hence that
each $\bar M^\free_\L$ is concave.  Moreover, each $\bar M^\free_\L$ lies
between $-1$ and $1$.  There therefore exists a sequence $\L_n$
of simple regions such that the 
sequence $\bar M^\free_{\L_n}$ converges pointwise 
to a limiting function which  we denote by $M^\free_\oo$.  
If $0<\g<\g'$ then by the fundamental theorem of calculus
and the bounded convergence theorem, we have that
\begin{equation}
\begin{split}
P(\g')-P(\g)&=\lim_{n\rightarrow\infty}
\big(P^\free_{\L_n}(\g')-P^\free_{\L_n}(\g)\big)\\
&=\lim_{n\rightarrow\infty}\int_{\g}^{\g'}\bar M^\free_{\L_n}(\g)\:d\g
=\int_{\g}^{\g'}M^\free_\oo(\g)\:d\g.
\end{split}
\end{equation}
The function $M^\free_\oo$ is concave, and hence continuous,
in $\g>0$.  It therefore
follows from the above that $P$ is in fact differentiable at each 
$\g>0$ (with derivative $M^\free_\oo$).  
The result follows since $P(-\g)=P(\g)$ for all $\g>0$.
\end{proof}

Whenever there is a unique infinite-volume state at $\g$, we will
denote it by $\el\cdot\er=\el\cdot\er_\g$.

\begin{lemma}\label{press_der_lem2}
For each $\g\neq0$ and each $(b,\a)$, we have that 
\begin{equation}
M:=\frac{\partial P}{\partial \g}=
\lim_{\L\uparrow\LLam} M^{b,\a}_\L=\lim_{\L\uparrow\LLam} \bar M^{b,\a}_\L.
\end{equation}
\end{lemma}
\begin{proof}
The proof of Lemma~\ref{uniq_lem} shows that
at each $\g\neq0$ the derivative of $P$ is $M^\free_\oo$.  
Since for all $(b,\a)$ and $\L$, the function $P^{b,\a}_\L(\g)$ is convex
and differentiable with
\begin{equation}
\frac{\partial P^{b,\a}_\L}{\partial\g}=\bar M^{b,\a}_\L
\end{equation}
it follows from the properties of convex functions that
$\bar M^{b,\a}_\L(\g)\rightarrow M(\g)$ at all $\g\neq0$.
That also $M^{b,\a}_\L\rightarrow M$ for $\g\neq0$ follows from the 
the fact that $M^-_\L\leq M^{b,\a}_\L\leq M^+_\L$
and the fact that $\lim M^\pm_\L=\lim \bar M^\pm_\L$
as we saw at~\eqref{s1}. 
\end{proof}

Lemma~\ref{press_der_lem2} implies in particular that
\begin{equation}
M=\lim_{\L\uparrow\LLam}\el\s_0\er^\pm_\L
\end{equation}
at all $\g\neq0$.  We know from Lemma~\ref{pm_states} that the
limits
\begin{equation}
M_\pm:=\lim_{\L\uparrow\LLam}\el\s_0\er^\pm_\L
\end{equation}%
\nomenclature[M]{$M_+$}{Spontaneous magnetization}%
exist also at $\g=0$.  By Lemma~\ref{press_der_lem1} there is a
unique state at $\g=0$ if and only if $M_+(0)=M_-(0)$.
We sometimes call $M_+(0)$ the \emph{spontaneous magnetization}.

Note that
for all $\L$ and all $\g>0$ we have $M^+_\L(-\g)=-M^-_\L(\g)$, so that
$\lim M^+_\L(-\g)=-M(\g)$.  Hence $M$ is an \emph{odd} function
of $\g\neq0$.   Note also that 
\[
M_+(0)=\lim_{\g\downarrow 0} M(\g).
\]
Indeed, rather more is true:  by repeating the argument 
at~\eqref{mo3} with $\s_A$ in place of $\s_0$, it follows that the
state $\el\cdot\er^+$ of Lemma~\ref{pm_states} may be written as the
weak limit
\begin{equation}\label{plus_ising}
\el\cdot\er^+_{\g=0}=\lim_{\g\downarrow 0} \el\cdot\er_\g
\end{equation}
where $\el\cdot\er_\g$ is the unique state at $\g>0$.
Thus we may summarize the results of this section as follows.
\begin{theorem}\label{ising_summary_thm}
There is a unique state at all $\g\neq0$ and there is a unique state at 
$\g=0$ if and only if
\begin{equation}
M_+(0)\equiv\lim_{\g\downarrow0}M(\g)=0.
\end{equation}
\end{theorem}
We now recall the remaining parameters $\l$, $\d$ and $\b$.
As previously,  we set $\d=1$, $\rho=\l/\d$, and write
\[
M^\b(\rho,\g)=M^\b(\rho,1,\g).
\]
It follows from Lemma~\ref{cor_mon_lem} that $M^\b_+(\rho,0)$
is an increasing function of $\rho$.  This motivates the
following definition.
\begin{definition}\label{ising_crit_val_def}
We define the \emph{critical value}
\[
\rho^\b_\crit:=\inf\{\rho>0:M^\b_+(\rho,0)>0\}.
\]
\end{definition}%
\nomenclature[r]{$\rho^\b_\crit$}{Critical value}
From Remark~\ref{inf_vol_rk} and~\eqref{plus_ising} it follows that
this  $\rho^\b_\crit$ coincides with the `percolation threshold'
$\rho_\crit(2)$ for the $q=2$ space--time random-cluster model
as defined in Definition~\ref{crit_def}.
More information about $\rho^\b_\crit$ 
and the behaviour of $M^\b$ and
related quantities near the critical point may be found
in Section~\ref{cons_sec}.

\chapter{The quantum Ising model:
random-parity representation and
sharpness of the phase transition}
\label{qim_ch}

\begin{quote}
{\it Summary.}  We develop a `random-parity' representation for the
space--time Ising model;  this is the space--time analog of the 
random-current representation.  The random-parity 
representation is then used to
derive a number of differential inequalities, from which  one can
deduce many important properties of the phase transition of the 
quantum Ising model, such as sharpness of the transition.
\end{quote}

\section{Classical and quantum Ising models}\label{sec-backg}


Recall from the Introduction that
the (transverse field) quantum Ising model on the finite graph $L$ 
is given by the Hamiltonian
\begin{equation}
H=-\tfrac{1}{2}\l\sum_{e=uv\in E}\s_u^{(3)}\s_v^{(3)}
-\d\sum_{v\in V}\s_v^{(1)},
\end{equation}
acting on the Hilbert space $\cH=\bigotimes_{v\in V}\CC^2$.  
We refer to that chapter for definitions of the notation used.
In the quantum Ising model the number
$\b>0$ is thought of as the `inverse temperature'.  We 
define the positive temperature states
\begin{equation}\label{qi_states_eq}
\oper_{L,\b}(Q)=\frac{1}{Z_L(\b)}\tr(e^{-\b H}Q),
\end{equation}
where $Z_L(\b)=\tr(e^{-\b H})$ and $Q$ is a suitable matrix.  The
\emph{ground state} is defined as the limit $\oper_L$ of $\oper_{L,\b}$ as
$\b\rightarrow\infty$.  If $(L_n: n \ge 1)$ is an increasing sequence of graphs
tending to the infinite
graph $\LL$, then we may also make use of the \emph{infinite-volume} limits
$$
\oper_{L,\b}=\lim_{n\rightarrow\infty}\oper_{L_n,\b},\qq
\oper_L=\lim_{n\rightarrow\infty}\oper_{L_n}.
$$
The existence of such limits is
discussed in~\cite{akn}, see also the related discussion of limits of
space--time Ising measures in Section~\ref{inf_potts_sec}.  

The quantum Ising model is intimately
related to the space--time Ising model,
one manifestation of this being the following.  Recall that 
if $\ket{\psi}$ denotes a vector then $\bra{\psi}$ denotes its 
conjugate transpose.  The state $\oper_{L,\b}$
of \eqref{qi_states_eq} gives rise to a probability measure $\mu$ on
$\{-1,+1\}^V$ by 
\begin{equation}
\mu(\s)=\frac{\el\s|e^{-\b H}|\s\er}{\tr(e^{-\b H})},\qq \s\in\{-1,+1\}^V.
\end{equation}
When $\g=0$, it turns out that $\mu$ is the
law of the vector $(\s_{(v,0)}: v\in V)$ under the space--time Ising measure
of \eqref{ising_def_eq} (with periodic boundary, see below).  
See \cite{akn} and the references therein.
It therefore makes sense to study the phase diagram of the quantum Ising model
via its representation in the space--time Ising model.
Note, however, that in our analysis it is crucial to work with $\g>0$,
and to take the limit $\g \downarrow 0$ later.  
The role played in the classical model by the external
field will in our analysis be played by the `ghost-field' $\g$ rather than
the `physical' transverse field $\d$.  (In fact, $\g$ corresponds to
a $\s^{(3)}$-field, see~\cite{crawford_ioffe}.)

In most of this chapter we will be working with periodic
boundary conditions in the $\RR$-direction.  That is to say, for simple
regions of the form~\eqref{def-oldL} we will identify the endpoints
of the the `time' interval $[-\b/2,\b/2]$, and think of this interval 
as the circle of circumference $\b$.  We will denote this circle by
$\SS=\SS_\b$%
\nomenclature[S]{$\SS_\b$}{Circle of circumference $\b$}
 and thus our simple regions will be of the form
$L\times\SS$ for some finite graph $L$.  
We shall generally (until Section~\ref{cons_sec}) keep $\b>0$ fixed,
and thus suppress reference to $\b$.  Similarly, we will generally
suppress reference to the boundary condition.
Thus we will write for instance $\S(D)$ for the set of spin
configurations permitted by $D$ (see the discussion 
before~\eqref{ising_def_eq}).

General regions of the form~\eqref{def-newL} will usually be thought of
as subsets of the simple region $L\times\SS$.  Thus,
for $v\in V$, we let $K_v\subseteq\SS$ be a finite union of disjoint
intervals, and we write $K_v=\bigcup_{i=1}^{m(v)} I^v_i$.%
\nomenclature[m]{$m(v)$}{Number of intervals constituting $K_v$}%
\nomenclature[I]{$I^v_i$}{Maximal subinterval of $K_v$}
 As before, no assumption is made  on whether the $I^v_i$ are
open, closed, or half-open.  
With the $K_v$ given, we define $F$ and $\L$ as in~\eqref{def-newL}.

For simplicity of notation we replace in this chapter the functions
$\l,\g$ in~\eqref{ising_pf} by $2\l,2\g$, respectively.  Thus
the space--time Ising measure
on a region $\L=(K,F)$ has partition function
\begin{equation}
Z'=\int d\mu_\d(\oD)\sum_{\s\in\S(\oD)}\exp\left\{\int_F\l(e)\s_e\, de+
\int_K\g(x)\s_x\, dx\right\},
\label{o12}
\end{equation}
where $\s_e=\s_{(u,t)}\s_{(v,t)}$ if $e=(uv,t)$.
See~\eqref{ising_def_eq}.  As previously, we write $\el f\er$
for the mean of a $\cG_\L$-measurable 
$f:\S\rightarrow\RR$ under this measure.  Thus
for example
\begin{equation}\label{st_Ising_eq}
\el\s_A\er=
\frac{1}{Z'}\int d\mu_\d(\oD)\sum_{\s\in\S(\oD)}
\s_A\exp\left\{\int_F\l(e)\s_e\, de+
\int_K\g(x)\s_x\, dx\right\}.
\end{equation}
Note that in this chapter we denote the partition function
by~$Z'$.%
\nomenclature[Z]{$Z'$}{Ising partition function}

It is essential for our method in this chapter that we work on general
regions of the form given in \eqref{def-newL}. The reason for this
is that, in the geometrical analysis of currents, we shall at times
remove from $K$  a random subset called the `backbone', and the ensuing
domain has the form of \eqref{def-newL}.
Note that considering this general class of regions
also allows us to revert to a `free'
rather than a `vertically periodic' boundary condition.  
That is, by setting $K_v =[-\b/2,\b/2)$ for all $v\in V$, rather than
$K_v = [-\b/2,\b/2]$, we effectively remove the restriction that the
`top' and `bottom' of each $v\times\SS$ have the same spin.  

Whenever we wish to
emphasize the roles of particular $K$, $\l$, $\d$, $\g$, we include 
them as subscripts.  For example, we may 
write $\el\s_A\er_K$ or $\el\s_A\er_{K,\g}$ or $Z'_\g$, and so on.

\subsection{Statement of the main results}

Let $0$ be a given point of $V\times \SS$.
We will be particularly concerned with the
\emph{magnetization}  and \emph{susceptibility}
of the space--time Ising model on $\L=L\times \SS$, given respectively by
\begin{align}
M=M_\L(\l,\d,\g) &:=\el\s_0\er,\label{def-mag}\\ 
\chi=\chi_\L(\l,\d,\g) &:=\frac{\partial M}{\partial \g}
=\int_\L\el\s_0;\s_x\er\,dx,
\label{def-susc}
\end{align}
where we recall that
the \emph{truncated} two-point function $\el\s_0;\s_x\er$ is given by
\begin{equation}
\el\s_A;\s_B\er := \el \s_A\s_B\er - \el\s_A\er\el\s_B\er.
\label{o20}
\end{equation}
Note that, for simplicity of notation, we will in most of this 
chapter keep $M$ and $\chi$ free from sub- and superscripts even
though they refer to finite-volume quantities.
Some basic properties of these quantities were discussed in
Section~\ref{ising_uniq_sec}.

Our main choice for $L$ is a box 
$[-n,n]^d$ in the $d$-dimensional cubic
lattice $\ZZ^d$ where $d \ge 1$, with a periodic boundary condition.
That is to say, apart from the usual nearest-neighbour bonds, we also
think of 
two vertices $u$, $v$ as joined by an edge whenever there exists
$i\in\{1,2,\dots,d\}$ 
such that $u$ and $v$ differ by exactly $2n$ in the $i$th 
coordinate.  
Subject to this boundary condition, $M$
and $\chi$ do not depend on the choice of origin $0$.
We shall
pass to the infinite-volume limit as $L\uparrow \ZZ^d$. The model is
over-parametrized, and we shall, as before, normally assume $\d=1$, 
and write $\rho=\l/\d$.   The
critical point $\bc=\bc^\b$ is given as in 
Definition~\ref{ising_crit_val_def} by
\begin{equation}
\bc^\b:=\inf\{\rho: M^\b_+(\rho)>0\},
\label{o1}
\end{equation}
where
\begin{equation}
M^\b_+(\rho):=\lim_{\g\downarrow 0}M^\b(\rho,\g),
\label{o2}
\end{equation}
is the magnetization in the limiting state $\el\cdot\er^\b_+$ as
$\g\downarrow 0$.
As in Theorem~\ref{crit_nontriv},
we have that:
\begin{equation}
\begin{split}
\text{if $d \ge 2$}&:\quad 0<\bc^\b<\oo\text{ for $\b\in(0,\oo]$},\\ 
\text{if $d=1$}&:\quad \bc^\b=\oo \text{ for $\b\in(0,\oo)$}, \ 0<\bc^\oo<\oo. 
\end{split}
\label{critvals}
\end{equation}

Complete statements of our main results are deferred until 
Section~\ref{cons_sec}, but here are two examples of what can be proved.

\begin{theorem}\label{ed_thm}
Let $u,v\in \ZZ^d$ where $d \ge 1$, and $s,t\in\RR$. For $\b\in(0,\oo]$:
\begin{romlist}
\item
if $0 < \rho < \bc^\b$, the two-point correlation function 
$\langle \s_{(u,s)}\s_{(v,t)}\rangle^\b_+$ of the space--time
Ising model decays exponentially to $0$ 
as $|u-v|+|s-t|\to\oo$,
\item
if $\rho\ge\bc^\b$,
 $\langle \s_{(u,s)}\s_{(v,t)}\rangle^\b_+ \ge M^\b_+(\rho)^2>0$.
\end{romlist}
\end{theorem}
Theorem~\ref{ed_thm} is what is called
`sharpness of the phase transition':  there is no intermediate 
regime in which correlations decay to zero slowly. 
(See for example~\cite{chayes_chayes86} 
and~\cite{ghm} for examples of systems
where this does occur).

\begin{theorem}\label{mf_thm} Let $\b\in(0,\oo]$. 
In the notation of Theorem \ref{ed_thm},
there exists $c = c(d)>0$ such that 
$$
M^\b_+(\rho) \geq c(\rho-\bc^\b)^{1/2}\qq\text{for }
\rho>\bc^\b.
$$ 
\end{theorem}

These and other facts will be stated and proved in Section \ref{cons_sec}.
Their implications for the infinite-volume quantum model will be elaborated 
around \eqref{o4}--\eqref{o5}. 

The approach used here is to prove a family of differential
inequalities for the finite-volume
magnetization $M(\rho,\g)$. This parallels the methods
established in \cite{ab,abf} for the analysis of the phase transitions in
percolation and Ising models on discrete lattices, 
and indeed our arguments are closely
related to those of \cite{abf}. Whereas new problems arise in
the current context and require treatment, certain aspects of the analysis
presented here are simpler that the corresponding steps of \cite{abf}. 
The application to the quantum model imposes a periodic boundary condition in
the $\b$ direction; some of our conclusions are valid
for the space--time Ising model with a free boundary condition.

The following is the principal differential inequality
we will derive.  (Our results are in fact valid in greater generality,
see the statement before Assumption \ref{periodic_assump}.) 

\begin{theorem}\label{main_pdi_thm}
Let $d\ge 1$, $\b<\oo$, and $L = [-n,n]^d$ with periodic boundary.  Then
\begin{equation}\label{ihp18}
M\leq \g\chi+M^3+2\l M^2\frac{\partial M}{\partial \l}
-2\d M^2\frac{\partial M}{\partial \d}.
\end{equation}
\end{theorem}

A similar inequality was derived in~\cite{abf} for the classical Ising
model, and our method of proof is closely related to that used there. 
Other such inequalities
have been proved for percolation in \cite{ab} 
(see also \cite{grimmett_perc}), and for the contact model
in \cite{aizenman_jung,bezuidenhout_grimmett}.
As observed in \cite{ab,abf}, the powers of $M$ on the right side 
of \eqref{ihp18} determine the bounds of Theorems \ref{ed_thm}(ii) 
and \ref{mf_thm}  on the critical exponents.
The cornerstone of our proof
is a `random-parity representation' of the space--time Ising model.

The analysis of the differential inequalities, 
following \cite{ab,abf}, reveals a
number of facts about the behaviour of the model.  In
particular, we will show the exponential
decay of the correlations $\el\s_0\s_x\er_+$ when $\rho<\bc^\b$ and $\g=0$, 
as asserted in Theorem \ref{ed_thm}, and in addition
certain bounds on two critical exponents of the model.  See
Section~\ref{cons_sec} for further details. 

We draw from \cite{akn,aizenman_nacht} in the following 
summary of the relationship between the 
phase transitions of the quantum and 
space--time Ising models. Let $u,v\in V$, and 
\[
\tau^\b_{L}(u,v) := \tr\bigl(\oper_{L,\b}(Q_{u,v})\bigr),\qq
Q_{u,v} = \s^{(3)}_u\s^{(3)}_v.
\]%
\nomenclature[t]{$\tau^\b$}{Two-point function}%
It is the case that
\begin{equation}
\tau^\b_{L}(u,v) = \el \s_A \er^\b_L
\label{o3}
\end{equation}
where $A=\{(u,0),(v,0)\}$, and the role of $\b$ is stressed
in the superscript. Let $\tau_L^\oo$ denote the limit
of $\tau^\b_{L}$ as $\b\to\oo$. For $\b\in(0,\oo]$, let $\tau^\b$
be the limit of $\tau^\b_L$ as $L\uparrow \ZZ^d$.
(The existence of this limit may depend on the choice 
of boundary condition on $L$,
and we return to this at the end of Section \ref{cons_sec}.)
By Theorem \ref{ed_thm},
\begin{equation}
\tau^\b(u,v) \le c'e^{-c|u-v|},
\label{o4}
\end{equation}
where $c'$, $c$ depend on $\rho$, 
and $c>0$ for $\rho<\bc^\b$ and $\b\in(0,\oo]$.
Here, $|u-v|$ denotes the $L^1$ distance from $u$ to $v$. 
The situation when $\rho=\bc^\b$ is
more obscure, but one has that
\begin{equation}
\limsup_{|v|\to\oo}\tau^\b(u,v) \le M^\b_+(\rho),
\label{o4a}
\end{equation}
so that $\tau^\b(u,v) \to 0$ whenever $M^\b_+(\rho)=0$. It 
is proved at Theorem \ref{crit_val_cor} that $\bc^\oo=2$ 
and $M^\oo_+(2)=0$ when $d=1$.

By the {\fkg} inequality, and the uniqueness
of infinite clusters in the space--time \rc\ model 
(see Theorem~\ref{inf_clust_uniq}),
\begin{equation}
\tau^\b(u,v) \ge M^\b_+(\rho-)^2 > 0,
\label{o5}
\end{equation}
when $\rho > \bc^\b$ and $\b\in(0,\oo]$, where $f(x-):= \lim_{y\uparrow x}f(y)$.
The proof is discussed at the end of Section \ref{cons_sec}.

The critical value $\bc^\b$ depends of course on the number of dimensions.
We shall in the next chapter use Theorem~\ref{ed_thm} and 
planar duality to show that $\bc^\oo=2$ when $d=1$,
and in addition that the transition is of second 
order in that $M^\oo_+(2)=0$. See Theorem~\ref{crit_val_cor}. 
The critical point has been calculated by other means in the quantum case,
but we believe that the current proof is valuable. Two applications
to the work of \cite{bjo0,GOS} are summarized in Section \ref{sec_1d}.

Here is a brief outline of the contents of this chapter.  
Formal definitions are presented
in Section \ref{sec-backg}. The random-parity 
representation of the quantum Ising model is 
described in Section \ref{rcr_sec}.
This representation may at first sight seem quite different from 
the random-current representation of the 
classical Ising model on a discrete lattice.
It requires more work to set up than does its discrete cousin, 
but once in place it works in a very similar, and sometimes simpler, manner.  
We then state and prove, in Section~\ref{ssec-switching}, the fundamental
`switching lemma'.  In
Section~\ref{sw_appl_sec} are presented a number of important consequences
of the switching lemma, including {\ghs} and
Simon--Lieb 
inequalities, as well as other useful inequalities and identities.
In Section~\ref{pf_sec}, we prove the 
somewhat more involved differential inequality
of Theorem \ref{main_pdi_thm}, 
which is similar to the main inequality of \cite{abf}. Our main results  
follow from Theorem \ref{main_pdi_thm} in conjunction with the results of
Section \ref{sw_appl_sec}.
Finally, in Sections \ref{cons_sec} and \ref{sec_1d}, we give rigorous 
formulations and proofs of our main results.

This chapter forms the contents of the
article~\cite{bjogr2}, which has been published in the 
Journal of Statistical Physics.   The quantum
mean-field, or Curie--Weiss, model has been studied using 
large-deviation techniques in \cite{chayes_ioffe_curie-weiss}, see also
\cite{grimmett_stp}.  
There is a very substantial overlap between the results 
reported here and those of the independent and contemporaneous 
article~\cite{crawford_ioffe}. 
The basic differential inequalities of 
Theorems~\ref{main_pdi_thm} and~\ref{three_ineq_lem}
appear in both places.  The proofs are in essence the same despite 
some superficial differences.  We are grateful to the authors 
of~\cite{crawford_ioffe} for 
explaining the relationship between the random-parity representation 
of Section~\ref{rcr_sec} and the random-current representation 
of~\cite[Section~2.2]{ioffe_geom}.
As pointed out in~\cite{crawford_ioffe}, 
the appendix of~\cite{chayes_ioffe_curie-weiss} contains a type of 
switching argument for the mean-field model. A principal difference between 
that argument and those of~\cite{crawford_ioffe,ioffe_geom}
and the current work is that it uses 
the classical switching lemma developed in~\cite{aiz82}, applied to a discretized version of the mean-field system.

\section{The random-parity representation}\label{rcr_sec}

The classical Ising model on a discrete 
graph $L$ is a `site model', in the sense that
configurations comprise spins assigned to the vertices (or `sites') of
$L$.  As described in the Introduction, the 
classical random-current representation maps this into a bond-model,
in which the 
sites no longer carry random values, but instead the 
\emph{edges} $e$ (or `bonds')  of
the graph are replaced by a random number $N_e$ of parallel edges. The
bond $e$ is called \emph{even} (\resp, \emph{odd}) if $N_e$ is
even (\resp, odd). The odd bonds may be arranged into paths and
cycles.  One cannot proceed
in the same way in the above space--time
Ising model.  

There are two 
possible alternative
approaches.  The first uses the fact that,
conditional on the set $D$ of deaths, $\L$ may be 
viewed as a discrete structure 
with finitely many components, to which the 
random-current representation of \cite{aiz82} may be applied.
This is explained in
detail around \eqref{rcr_step1_eq} below.
Another approach is to forget about `bonds', and
instead to concentrate on the parity configuration associated 
with a current-configuration, as follows.  

The circle $\SS$ may be viewed as a continuous limit of a ring of 
equally spaced points.  If we apply the random-current
representation to the discretized system, but only record whether a bond is
even or odd, the representation has a well-defined limit as a
partition of $\SS$ into even and odd sub-intervals.  In the limiting
picture, even and odd intervals carry different weights, and it is the
properties of these weights that render the representation useful.  This is
the essence of the main result in this section, 
Theorem~\ref{rcr_thm}. We will prove this result without recourse to
discretization.

We now define two additional random processes associated with the
space--time Ising measure on $\L$.  The first is a random colouring of $K$,
and the second is a random (finite) weighted graph.  These two objects
will be the main components of the random-parity representation.

\subsection{Colourings}\label{ssec-col}

Let $\ol K$ be the closure of $K$.
A set of \emph{sources} is a finite set $A\se\ol K$ such that:
each $a\in A$ is the endpoint of at most one maximal subinterval
$I^v_i$.  (This last condition is for simplicity later.)
Let $B\subseteq F$ and $G\subseteq K$ be finite sets.  
Let $S=A \cup G \cup V(B)$,%
\nomenclature[S]{$S$}{Switching points}
 where $V(B)$ is the set of
endpoints of bridges of $B$,
and call members of $S$ \emph{switching points}. As usual
we shall assume that $A$, $G$ and $V(B)$ are disjoint.

We shall define a colouring $\psi^A=\psi^A(B,G)$%
\nomenclature[p]{$\psi^A$}{Colouring}
 of $K\sm S$ 
using the two colours (or labels) `even' and `odd'.  
This colouring is constrained to be `valid', where
a valid colouring is defined to be a mapping 
$\psi:K\sm S \to\{\even, \odd\}$ such that: 
\begin{romlist}
\item the label is constant between 
two neighbouring switching points, that is, $\psi$ 
is constant on any sub-interval
of $K$ containing no members of $S$,
\item  the label always switches at each switching point, 
which is to say that, for $(u,t) \in S$, $\psi(u,t-) \ne \psi(u,t+)$, 
whenever these two values are defined,
\item for any pair $v$, $k$ such that $I^v_k\neq\SS$, in the limit as we move
along $v\times I^v_k$ towards either endpoint of $v\times I^v_k$, 
the colour converges to `odd' if that endpoint lies in $S$, and
to `even' otherwise. 
\end{romlist} 

If there exists $v\in V$ and
$1\leq k\leq m(v)$ such that $v\times \ol{I^v_k}$ contains an
\emph{odd} number of switching points, then conditions (i)--(iii) cannot be
satisfied;  in this case we set the colouring $\psi^A$ to
a default value denoted~$\#$.  

Suppose that (i)--(iii) \emph{can} be satisfied, and let 
\begin{equation*}
W=W(K):=\{v\in V: K_v=\SS\}. 
\end{equation*}%
\nomenclature[W]{$W$}{Vertices $v\in V$ such that $K_v=\SS$}%
If $W = \es$, then there exists a unique valid colouring, 
denoted $\psi^A$.  If $r=|W|\ge 1$, there are exactly $2^r$ valid
colourings, one for each of the two possible colours assignable to the sites
$(w,0)$, $w \in W$; in this case we let $\psi^A$ be chosen uniformly at random
from this set, independently of all other choices. 

\begin{figure}[tbp]
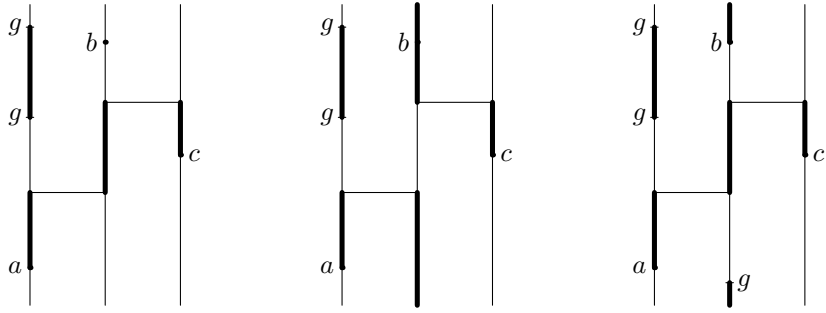

\includegraphics{thesis.6}
\hspace{1.3cm} 
\includegraphics{thesis.7} 
\hspace{1.3cm} 
\includegraphics{thesis.8} 
\caption[Colourings in the random-parity representation]
{Three examples of colourings for given 
$B\subseteq F$, $G\subseteq K$.
Points in $G$ are written $g$. Thick line
segments are `odd' and thin segments `even'.  In this illustration we have
taken $K_v=\SS$ for all $v$.  \emph{Left and middle}: two of the eight
possible colourings when the sources are $a$, $c$.  
\emph{Right}: one of the possible
colourings when the sources are $a$, $b$, $c$.}
\label{colouring_fig}
\end{figure}

We write $M_{B,G}$%
\nomenclature[M]{$M_{B,G}$}{Uniform measure on colourings}
 for the probability measure (or expectation when
appropriate) governing the randomization in the definition of $\psi^A$: 
$M_{B,G}$ is
the uniform (product) measure on the set of valid colourings, and it is a point
mass if and only if $W=\es$.  See Figure~\ref{colouring_fig}.  

Fix the set $A$ of sources.
For (almost every) pair $B$, $G$, one may construct as above
a (possibly random) colouring $\psi^A$. Conversely, it is easily seen that
the pair $B$, $G$ may (almost surely) be reconstructed from 
knowledge of the colouring $\psi^A$.
For given $A$, we may thus speak of a configuration as being either a pair $B$, $G$, or 
a colouring $\psi^A$. 
While $\psi^A(B,G)$ is a colouring of $K \sm S$ only, we shall
sometimes refer to it as a colouring of~$K$.

The next step is to assign weights $\pd\psi$ to colourings $\psi$. The
`failed' colouring $\#$ is assigned weight $\pd\# =0$.  
For every valid colouring $\psi$, let $\ev(\psi)$
\nomenclature[e]{$\ev(\psi)$}{Set of `even' points in $\psi$}%
 (\resp, $\odd(\psi)$%
\nomenclature[o]{$\odd(\psi)$}{Set of `odd' points in $\psi$}%
)
denote the subset of $K$ that is labelled even (\resp, odd), 
and let
\begin{equation}
\partial\psi :=\exp\bigl\{2\d(\ev(\psi))\bigr\},
\label{def-wt}
\end{equation}%
\nomenclature[d]{$\partial\psi$}{Weight of colouring $\psi$}%
where
$$
\d(U):=\int_{U}\d(x)\,dx, \qq U \subseteq K.
$$
Up to a multiplicative constant depending on $K$ and $\d$ only,
$\pd\psi$ equals the square of the 
probability that the odd part of $\psi$ is death-free.

\subsection{Random-parity representation}\label{ssec-rpr}
The expectation $E(\pd\psi^A)$ is taken over the sets $B$, $G$, and over
the randomization that takes place when $W \ne \es$, 
that is, $E$ denotes expectation with respect to the
measure $d\mu_\l(B) d\mu_\g(G) dM_{B,G}$. 
The notation has been chosen to harmonize with that
used in \cite{abf} in the discrete case:
the expectation $E(\partial\psi^A)$ will play the role of the probability
$P(\partial\underline n=A)$ of \cite{abf}.
The main result of this section now follows.

\begin{theorem}[Random-parity representation]\label{rcr_thm}
For any finite set $A \subseteq \ol K$ of sources,
\begin{equation}
\el\s_A\er=\frac{E(\partial\psi^A)}{E(\partial\psi^\es)}.
\label{ihp8}
\end{equation}
\end{theorem}

We introduce a second random object in advance of proving this.
Let $D$ be a finite subset of $K$.
The set $(v\times K_v)\sm\oD$ is
a union of maximal death-free intervals which we write
$v\times J^v_k$,%
\nomenclature[J]{$J^v_k$}{Subintervals of $K$ bounded by deaths}
 and
where $k=1,2,\dotsc,n$ and $n=n(v,\oD)$%
\nomenclature[n]{$n(v,D)$}{Number of death-free intervals in $K_v$}
 is the number of such intervals.  
We write $V(\oD)$%
\nomenclature[V]{$V(D)$}{Collection of maximal death-free intervals}
 for the collection of all such intervals.

For each $e=uv\in E$, and each $1\leq k\leq n(u)$ and $1\leq l\leq n(v)$, let 
\begin{equation}
J^e_{k,l}:=J^u_k\cap J^v_l,
\end{equation}%
\nomenclature[J]{$J^e_{k,l}$}{Element of $E(D)$}%
and 
\begin{equation}
E(D)=\bigl\{e\times J^e_{k,l}:e\in E,\ 1\leq k\leq n(u),\ 1\leq l\leq n(v),
\,J^e_{k,l}\neq\es\bigr\}.
\end{equation}%
\nomenclature[E]{$E(D)$}{Edge set of the graph $G(D)$}%
Up to a finite set of points, $E(D)$  forms a partition of the set
$F$ induced by the `deaths' in~$D$.  

\begin{figure}[tbp]
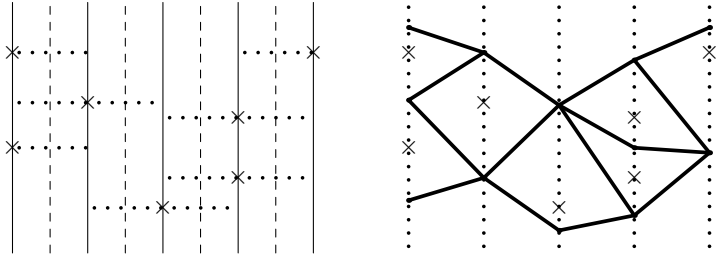

\includegraphics{thesis.9}
\qquad
\includegraphics{thesis.10} 
\caption[Applying the discrete random-current representation in
the quantum Ising model]
{\emph{Left}: The partition $E(\oD)$.  We have: $K_v=\SS$ for $v \in V$, 
the lines $v\times K_v$ are drawn as solid, the lines $e\times K_e$ 
as dashed, and elements of $\oD$ are marked as crosses.  The endpoints 
of the $e\times J^e_{k,l}$ are the points where the dotted lines meet the
dashed lines.
\emph{Right}:  
The graph $G(\oD)$.  In this illustration, the dotted lines are the
$v\times K_v$, and the solid lines are the edges of $G(\oD)$.}
\label{bridge_part_fig}
\end{figure}

The pair
\begin{equation}
G(\oD):=(V(\oD),E(\oD))
\end{equation}%
\nomenclature[G]{$G(D)$}{Discrete graph constructed from $D$}%
may be viewed as a graph, illustrated in Figure~\ref{bridge_part_fig}.
We will use the symbols $\bar v$ and $\bar e$ for
typical elements of $V(\oD)$ and $E(\oD)$, respectively.
There are natural weights on the edges and vertices of $G(\oD)$: for
$\bar e=e\times J^e_{k,l}\in E(\oD)$ and
$\bar v=v\times J^v_k\in V(\oD)$, let
\begin{equation}
J_{\bar e}:= \int_{J^e_{k,l}}\l(e,t)\,dt,
\qquad 
h_{\bar v}:= \int_{J^v_k}\g(v,t)\,dt.
\label{o10}
\end{equation}
Thus the weight of a vertex or edge is its measure, calculated according to
$\l$ or $\g$, respectively. By \eqref{o10},
\begin{equation}
\sum_{\bar e\in E(\oD)} J_{\bar e} + \sum_{\bar v\in V(\oD)}h_{\bar v}
= \int_{F}\l(e)\,de + \int_{K}\g(x)\,dx.
\label{o11}
\end{equation}

\begin{proof}[Proof of Theorem~\ref{rcr_thm}]
With $\L=(K,F)$ as in \eqref{def-newL}, we consider the partition function
$Z'=Z'_K$ given in \eqref{o12}.
For each $\bar v\in V(\oD)$,
$\bar e\in E(\oD)$, the spins $\s_v$ and $\s_e$ are constant for $x\in\bar v$
and $e\in\bar e$, respectively.  Denoting their common values by 
$\s_{\bar v}$ and $\s_{\bar e}$ respectively, the summation in \eqref{o12} equals
\begin{multline}
\sum_{\s\in\S(\oD)}\exp\left\{
\sum_{\bar e\in E(\oD)}\s_{\bar e}\int_{\bar e}\l(e)\, de+
\sum_{\bar v\in V(\oD)}\s_{\bar v}\int_{\bar v}\g(x)\, dx\right\}\\
=\sum_{\s\in\S(\oD)}\exp\left\{
\sum_{\bar e\in E(\oD)}J_{\bar e}\s_{\bar e}+
\sum_{\bar v\in V(\oD)}h_{\bar v}\s_{\bar v}\right\}.
\label{ihp2}
\end{multline}
The right side of \eqref{ihp2} is the partition 
function of the discrete Ising model on
the graph $G(\oD)$, with pair couplings $J_{\bar e}$ and external fields
$h_{\bar v}$.  We shall apply the random-current expansion of \cite{abf} to this
model.  

For convenience of exposition, we introduce the extended graph
\begin{align}
\wtilde G(\oD)&=(\wtilde V(\oD),\wtilde E(\oD))\label{ihp3}\\
&:=
\bigl(V(\oD)\cup\{\Gh\},E(\oD)\cup\{\bar v \Gh: \bar v \in V(\oD)\}\bigr)
\nonumber
\end{align}
where $\Gh$ is the
ghost-site.  We call members of $E(\oD)$ \emph{lattice-bonds},
and those of $\wtilde E(\oD)\sm E(\oD)$ \emph{ghost-bonds}.  
Let $\Psi(\oD)$ be the random multigraph with vertex set $\wtilde V(\oD)$ and
with each edge of $\wtilde E(\oD)$ replaced by a random number of parallel edges,
these numbers being independent and having the Poisson distribution, with
parameter $J_{\bar e}$ for lattice-bonds $\bar e$, and parameter $h_{\bar v}$
for ghost-bonds $\bar v\Gh$.  

Let $\{\partial\Psi(\oD)=A\}$ denote the event
that, for each $\bar v\in V(\oD)$, the total degree of $\bar v$ in $\Psi(\oD)$
\emph{plus} the number of elements of $A$ inside $\bar v$ (when regarded as an
interval) is even.  
(There is $\mu_\d$-probability $0$ that $A$ contains some endpoint of some $V(D)$,
and thus we may overlook this possibility.) 
Applying the discrete random-current expansion, and
in particular~\cite[eqn (9.24)]{grimmett_rcm}, we obtain by \eqref{o11} that
\begin{equation}
\sum_{\s\in\S(\oD)}\exp\left\{
\sum_{\bar e\in E(\oD)}J_{\bar e}\s_{\bar e}+
\sum_{\bar v\in V(\oD)}h_{\bar v}\s_{\bar v}\right\}=
c 2^{|V(\oD)|}P_D(\partial\Psi(\oD)=\es),
\end{equation}
where $P_D$ is the law of the edge-counts, and 
\begin{equation}
c=\exp\left\{\int_F \l(e)\,de + \int_K\g(x)\,dx\right\}.
\label{o21}
\end{equation}

By the same argument applied to the numerator in \eqref{st_Ising_eq}
(adapted to the measure on $\L$, see the remark after \eqref{o12}),
\begin{equation}
\label{rcr_step1_eq}
\el\s_A\er=
\frac{E(2^{|V(\oD)|}\one\{\partial\Psi(\oD)=A\})}
{E(2^{|V(\oD)|}\one\{\partial\Psi(\oD)=\es\})},
\end{equation}
where the expectation is with respect to $\mu_\d \times P_D$. 
The claim of the theorem will follow by an 
appropriate manipulation of \eqref{rcr_step1_eq}. 

Here is another way to sample $\Psi(\oD)$, which allows us to couple it with
the random colouring $\psi^A$.  Let $B\subseteq F$ and $G\subseteq K$
be finite sets sampled from $\mu_\l$ and $\mu_\g$ respectively.  
The number of points of $G$ lying
in the interval $\bar v=v\times J^v_k$ has the Poisson distribution with
parameter $h_{\bar v}$, and similarly the number of elements of $B$ lying
in $\bar e=e\times J^e_{k,l}\in E(\oD)$ has the Poisson distribution with
parameter $J_{\bar e}$.  Thus, for given $\oD$, the multigraph $\Psi(B,G,\oD)$,
obtained by replacing an edge of $\wtilde E(\oD)$ by parallel edges
equal in number to the corresponding number
of points from $B$ or $G$, respectively, has the same law as $\Psi(\oD)$.  
Using the
\emph{same} sets $B$, $G$ we may form the random colouring~$\psi^A$.  

The numerator of~\eqref{rcr_step1_eq} satisfies
\begin{align}
&E(2^{|V(\oD)|}\one\{\partial\Psi(\oD)=A\})\label{ihp5}\\
&\hskip1cm=\iint d\mu_\l(B)\, d\mu_\g(G)\,\int d\mu_\d(\oD)\, 
2^{|V(\oD)|}\one\{\partial\Psi(B,G,\oD)=A\}\nonumber\\
&\hskip1cm= \mu_\d(2^{|V(D)|}) 
\iint d\mu_\l(B)\, d\mu_\g(G)\,\wtilde\mu(\partial\Psi(B,G,\oD)=A),
\nonumber
\end{align}
where $\wtilde\mu$ is the probability measure on $\cF$ satisfying
\begin{equation}
\frac{d\wtilde\mu}{d\mu_\d}(D) \propto 2^{|V(\oD)|}.
\label{ihp10}
\end{equation}
Therefore, by \eqref{rcr_step1_eq},
\begin{equation}
\el\s_A\er=
\frac{\wtilde P(\partial\Psi(B,G,\oD)=A)}
{\wtilde P(\partial\Psi(B,G,\oD)=\es)},
\label{ihp9}
\end{equation}
where $\wtilde P$ denotes the probability under $\mu_\l\times\mu_\g\times\wtilde \mu$.
We claim that
\begin{equation}
\wtilde \mu(\pd\Psi(B,G,\oD)=A) = s M_{B,G}(\pd\psi^A(B,G)),
\label{ihp13}
\end{equation}
for all $B$, $G$, where $s$ is a constant,
and the expectation $M_{B,G}$ is over the uniform 
measure on the set of valid colourings. 
Claim \eqref{ihp8} follows from this,
and the remainder of the proof is to show \eqref{ihp13}.
The constants $s$, $s_j$ are permitted in the following to depend only on
$\L$, $\d$.  

Here is a special case:
\begin{equation}
\wtilde\mu(\partial\Psi(B,G,\oD)=A)=0
\end{equation}
if and only if some interval $\ol{I^v_k}$ contains an odd number of switching
points, if and only if $\psi^A(B,G) =\#$ and $\partial\psi^A(B,G)=0$.  
Thus \eqref{ihp13} holds in this case.

Another special case arises 
when $K_v=[0,\b)$ for all $v\in V$, that is, the `free boundary' case.  As remarked earlier,
there is a unique valid colouring $\psi^A=\psi^A(B,G)$.
Moreover, $|V(\oD)|=|\oD|+|V|$, whence from standard properties of Poisson
processes, $\wtilde\mu=\mu_{2\d}$.  It may be seen after some thought 
(possibly with the aid of a diagram) that, for given $B$, $G$,
the events $\{\partial\Psi(B,G,\oD)=A\}$ and $\{\oD\cap\odd(\psi^A)=\es\}$
differ by an event of $\mu_{2\d}$-probability $0$. Therefore,
\begin{align}
\wtilde\mu(\partial\Psi(B,G,\oD)=A)&=
\mu_{2\d}( \oD\cap\odd(\psi^A)=\es )\label{ihp7}\\
&=\exp\{-2\d(\odd(\psi^A))\}\nonumber\\
&= s_1\exp\{2\d(\ev(\psi^A))\}=s_1 \partial\psi^A,
\nonumber
\end{align}
with $s_1=e^{-2\d(K)}$. In this special case, \eqref{ihp13} holds.

For the general case, we first note some properties of
$\wtilde\mu$.  By the above, we may assume that
$B$, $G$ are such that $\wtilde\mu(\partial\Psi(B,G,\oD)=A)>0$,
which is to say that each $\ol{I_k^v}$ 
contains an even number of switching points. 
Let $W = \{v\in V: K_v=\SS\}$ and,  for
$v\in V$, let $\oD_v=D\cap(v\times K_v)$%
\nomenclature[D]{$\oD_v$}{Deaths in $K_v$}
 and $\od(v)=|\oD_v|$%
\nomenclature[d]{$\od(v)$}{Number of deaths in $K_v$}.
By \eqref{ihp10}, 
\begin{align*}
\frac{d\wtilde\mu}{d\mu_\d}(D) \propto 2^{|V(\oD)|}&=
\prod_{w\in W}2^{1\vee \od(w)}\prod_{v\in V\setminus W}2^{m(v)+\od(v)}\\
&\propto 2^{|\oD|} \prod_{w\in W} 2^{\one\{\od(w)=0\}},
\end{align*}
where $a\vee b = \max\{a,b\}$, and we recall the number $m(v)$ of intervals $I^v_k$ that constitute
$K_v$.  Therefore,
\begin{equation}
\frac{d\wtilde\mu}{d\mu_{2\d}}(D) \propto \prod_{w\in W} 2^{\one\{\od(w)=0\}}.
\end{equation}
Three facts follow.
\begin{letlist}
\item The sets $\oD_v$, $v\in V$ are independent under $\wtilde\mu$.
\item  For $v\in V\setminus W$, the law of $\oD_v$ under $\wtilde \mu$ is $\mu_{2\d}$.  
\item For $w\in W$, the law $\mu_w$ of 
$\oD_w$ is that of $\mu_{2\d}$ skewed by the Radon--Nikodym factor 
$2^{\one\{\od(w)=0\}}$, which is to say that
\begin{align}
\mu_w(\oD_w \in H) 
&= \frac1{\a_w}\Bigl[2\mu_{2\d}(\oD_w\in H,\,\od(w)=0) \label{ihp15}\\
&\hskip3cm +
\mu_{2\d}(\oD_w\in H,\,\od(w)\ge 1)\bigr],
\nonumber
\end{align}
for appropriate sets $H$, where 
$$
\a_w=\mu_{2\d}(\od(w)=0)+1.
$$
\end{letlist}

Recall the set $S=A\cup G\cup V(B)$ of switching points.  By (a) above, 
\begin{align}
\wtilde\mu(\partial\Psi(B,G,\oD)=A)&=
\wtilde\mu(\forall v,k:\, |S\cap \ol{J^v_k}|\mbox{ is even})\label{ihp11}\\
&= \prod_{v\in V} \wtilde\mu(\forall k:\, |S\cap \ol{J^v_k}|\mbox{ is even}).
\nonumber
\end{align}
We claim that
\begin{equation}
\wtilde\mu(\forall k:\, |S\cap \ol{J^v_k}|\mbox{ is even})= 
s_2(v)M_{B,G}\Bigl(\exp\bigl\{2\d\bigl(\ev(\psi^A)\cap(v\times K_v)\bigr)\bigr\}\Bigr),
\label{ihp12}
\end{equation}
where $M_{B,G}$ is as before. Recall that $M_{B,G}$ is a product measure.
Once \eqref{ihp12} is proved, \eqref{ihp13} follows by \eqref{def-wt} and \eqref{ihp11}.

For $v\in V\sm W$, the restriction of $\psi^A$ to $v\times K_v$
is determined given $B$
and $G$, whence by (b) above, and the remark prior to \eqref{ihp7},
\begin{align}
\wtilde\mu(\forall k:\, |S\cap \ol{J^v_k}|\mbox{ is even})&=
\mu_{2\d}(\forall k:\, |S\cap \ol{J^v_k}|\mbox{ is even})
\label{ihp14}\\
&= \exp\bigl\{-2\d\bigl(\odd(\psi^A)\cap(v\times K_v)\bigr)\bigr\}.
\nonumber
\end{align}
Equation \eqref{ihp12} follows with $s_2(v) =\exp\{-2\d(v\times K_v)\}$. 

For $w\in W$, by \eqref{ihp15},
\begin{align*}
&\wtilde\mu(\forall k:\, |S\cap J^w_k|\mbox{ is even})\\
&\hskip1cm =\frac1{\a_w}\Bigl[2\mu_{2\d}(\oD_w=\es)
+\mu_{2\d}(\oD_w\ne\es,\,\forall k:\, |S\cap J^w_k|\mbox{ is even})\Bigr]\\
&\hskip1cm=\frac1{\a_w}\Bigl[\mu_{2\d}(\oD_w=\es)+
\mu_{2\d}(\forall k:\, |S\cap J^w_k|\mbox{ is even})\Bigr].
\end{align*}
Let $\psi=\psi^A(B,G)$ be a valid colouring with $\psi(w,0) = \even$.  
The colouring $\ol\psi$, obtained from $\psi$ by flipping
all colours on $w \times K_w$, is valid also. We take into account the periodic boundary condition,
to obtain this time that 
\begin{align*}
&\mu_{2\d}(\forall k:\,|S\cap \ol{J^w_k}|\mbox{ is even})\\
&\quad =\mu_{2\d}\bigl(\{\oD_w\cap\odd(\psi)=\es\}\cup
\{\oD_w\cap\ev(\psi)=\es\}\bigr)\\
&\quad =\mu_{2\d}(\oD_w\cap\odd(\psi)=\es)+
\mu_{2\d}(\oD_w\cap\ev(\psi)=\es)
-\mu_{2\d}(\oD_w=\es),
\end{align*}
whence
\begin{align}
&\a_w \wtilde\mu(\forall k:\,|S\cap \ol{J^w_k}|\mbox{ is even})\\
&\hskip1cm =\mu_{2\d}(\oD_w\cap\odd(\psi)=\es)
+\mu_{2\d}(\oD_w\cap\ev(\psi)=\es)\nonumber\\
&\hskip1cm =2M_{B,G}\Bigl(\exp\bigl\{-2\d\bigl(\odd(\psi^A)\cap(w\times K_w)\bigr)\bigr\}\Bigr),
\nonumber
\end{align}
since $\odd(\psi^A) = \odd(\psi)$ with 
$M_{B,G}$-probability $\frac12$, and
equals $\ev(\psi)$ otherwise.
This proves \eqref{ihp12} with $s_2(w) = 2\exp\{-2\d(w\times K_w)\}/\a_w$. 
\end{proof}

By keeping track of the constants in the above proof, we arrive at the
following statement, which will be useful later.

\begin{lemma}\label{Z'}
The partition function $Z'=Z_K'$ of \eqref{o12} satisfies
$$
Z' = 2^N e^{\l(F)+\g(K)-\d(K)} E(\pd\psi^\es),
$$
where  $N=\sum_{v\in V}m(v)$
is the total number of intervals comprising~$K$.
\end{lemma}
We denote $Z_K=E(\partial\psi^\es)$,%
\nomenclature[Z]{$Z_K$}{$E(\partial\psi^\es)$}
 which is thus a constant multiple of~$Z'$.

\subsection{The backbone}

The concept of the backbone is key to the analysis of \cite{abf}, and
its definition there has a certain complexity. The corresponding
definition is rather easier in the current setting, because of the fact that
bridges, deaths, and sources have (almost surely) no common point. 

\begin{figure}[tbp]
\includegraphics{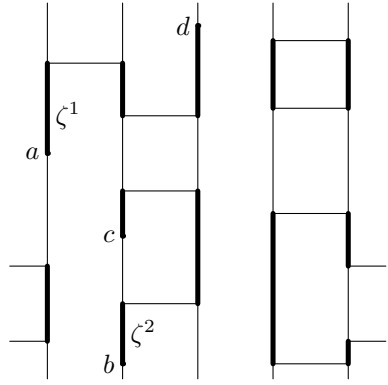}
\caption[The backbone]
{A valid colouring configuration $\psi$ with sources $A=\{a,b,c,d\}$, 
and its backbone  $\xi=\zeta^1\circ\zeta^2$.  Note that, in this 
illustration, bridges protruding from the sides `wrap around', and that 
there are no ghost-bonds.}
\label{backbone_fig}
\end{figure}

We construct a total order on $K$ by: first ordering the
vertices of $L$, and then using the natural order on $[0,\b)$.
Let $A\subseteq \ol K$, $B\subseteq F$ and $G\subseteq K$
be finite.  Let $\psi$ be a valid colouring.  We will
define a sequence of directed odd paths called the \emph{backbone} and denoted
$\xi=\xi(\psi)$%
\nomenclature[x]{$\xi(\psi)$}{Backbone}.
Suppose $A=(a_1,a_2,\dotsc,a_n)$ in the above ordering.  
Starting at $a_1$, follow the odd interval (in $\psi$) until you
reach an element of $S=A\cup G \cup V(B)$.  If the first
such point thus encountered is the endpoint of a
bridge, cross it, and continue along the odd interval;  continue likewise
until you first reach a point $t_1\in A\cup G$, at which point you stop.
Note, by the validity of $\psi$, that $a_1\ne t_1$.
The odd path thus traversed is denoted $\zeta^1$;%
\nomenclature[z]{$\zeta^k$}{Part of a backbone}
 we take $\zeta^1$ to be closed 
(when viewed as a subset of $\ZZ^d\times\RR$).  Repeat
the same procedure with $A$ replaced by $A\setminus\{a_1,t_1\}$, and iterate
until no sources remain.  The
resulting (unordered) set of paths $\xi=(\zeta^1,\dotsc,\zeta^k)$ 
is called the \emph{backbone} of
$\psi$.  The backbone will also be denoted at times as
$\xi=\zeta^1\circ\dotsb\circ\zeta^k$.  We define $\xi(\#)=\es$.
Note that, apart from the backbone, the remaining odd segments of
$\psi$ form disjoint self-avoiding cycles (or  `eddies').  
Unlike the discrete setting of \cite{abf}, 
there is a (a.s.) unique way of specifying
the backbone from knowledge of $A$, $B$, $G$ and the valid colouring $\psi$.
See Figure~\ref{backbone_fig}.

The backbone contains all the sources $A$ as endpoints, and the configuration outside $\xi$
may be any sourceless configuration.  Moreover, since $\xi$ is entirely odd,
it does not contribute to the weight $\partial\psi$ in \eqref{def-wt}.  
It follows, using properties of Poisson processes,  that the
conditional expectation $E(\partial\psi^A\mid\xi)$ equals the expected weight
of any sourceless colouring of $K\sm\xi$, which is to say that,
with $\xi:= \xi(\psi^A)$,
\begin{equation}\label{backb_cond_eq}
E(\partial\psi^A\mid\xi)=E_{K\setminus\xi}(\partial\psi^\es)
= Z_{K\setminus\xi}.
\end{equation}
Cf.\ \eqref{o12} and \eqref{ihp8}, and recall Remark~\ref{rem-as}.
We abbreviate $Z_K$ to $Z$, and recall
from Lemma \ref{Z'} that the $Z_R$ differ from the partition
functions $Z_R'$ by certain multiplicative constants.

Let $\Xi$ be the set of all possible backbones as $A$, $B$, and $G$ vary, 
regarded as sequences of
directed paths in $K$; these paths may, if required, be ordered by their
starting points.   
For $A \subseteq\ol K$ and $\nu\in \Xi$, we write
$A\sim\nu$ if there exist $B$ and $G$ such
that $M_{B,G}(\xi(\psi^A)=\nu)>0$.
We define the \emph{weight}
$\wt^A(\nu)$%
\nomenclature[w]{$\wt^A(\xi)$}{Weight of backbone}
 by 
\begin{equation}
\wt^A(\nu) = \wt^A_K(\nu):=
\begin{cases} \dfrac{Z_{K\sm \nu}}{Z} &\text{if } A \sim \nu,\\
0 &\text{otherwise}.
\end{cases}
\label{ihp16}
\end{equation}
By \eqref{backb_cond_eq} and Theorem \ref{rcr_thm},
with $\xi=\xi(\psi^A)$,
\begin{equation}\label{backbone_rep_eq}
E(\wt^A(\xi))=\frac{E(E(\pd\psi^A\mid\xi))}{Z} = 
\frac{E(\pd\psi^A)}{E(\pd\psi^\es)} =\el\s_A\er.
\end{equation}

For $\nu^1,\nu^2\in \Xi$ with $\nu^1\cap \nu^2=\es$
(when viewed as subsets of $K$), 
we write $\nu^1\circ \nu^2$ for the element of $\Xi$
comprising the union of $\nu^1$ and~$\nu^2$.

Let $\nu =\zeta^1\circ\dotsb\circ\zeta^k\in\Xi$ where $k \ge 1$.  If $\zeta^i$
has starting point $a_i$ 
and endpoint $b_i$, we write $\zeta^i:a_i\rightarrow b_i$, and also
$\nu:a_1\rightarrow b_1,\dotsc,a_k\rightarrow b_k$.  
If $b_i \in G$, we write $\zeta^i: a_i \rightarrow \Gh$.
There is a natural way to `cut' $\nu$ at points $x$ lying on
$\zeta^i$, say, where $x\ne a_i, b_i$: let $\bar\nu^1=\bar\nu^1(\nu,x)=
\zeta^1\circ\cdots\circ\zeta^{i-1}\circ \zeta^i_{\le x}$ and 
$\bar\nu^2=\bar\nu^2(\nu,x)=\zeta^i_{\ge x}\circ\zeta^{i+1}\circ\dots\circ\zeta^k$,
where $\zeta^i_{\le x}$ (\resp, $\zeta^i_{\ge x})$ is the closed sub-path of
$\zeta^i$ from $a_i$ to $x$ (\resp, $x$ to $b_i$).
We express this decomposition as $\nu=\bar\nu^1\circ\bar\nu^2$ where, this
time, each $\bar\nu^i$ may comprise a number of disjoint paths.
The notation $\ol\nu$ will be used only in a situation where there has been a
cut. 

We note two special cases. If $A=\{a\}$, then necessarily
$\xi(\psi^A):a\rightarrow \Gh$, so
\begin{equation}
\el\s_a\er=E\bigl(\wt^a(\xi)\cdot\one\{\xi:a\rightarrow \Gh\}\bigr).
\label{special1}
\end{equation}
If $A=\{a,b\}$ where $a<b$ in the ordering of $K$, then 
\begin{equation}
\el\s_a\s_b\er=E\bigl(\wt^{ab}(\xi)\cdot\one\{\xi:a\rightarrow b\}\bigr)
+E\bigl(\wt^{ab}(\xi)\cdot\one\{\xi:a\rightarrow \Gh,\,b\rightarrow \Gh\}\bigr).
\label{special2}
\end{equation}
The last term equals $0$ when $\g\equiv0$.

Finally, here is a lemma for computing the weight of $\nu$ in terms of its
constituent parts. 
The claim of the lemma is, as usual, valid only `almost surely'.

\begin{lemma}\label{backb2}
(a) Let $\nu^1,\nu^2 \in \Xi$ be disjoint, and $\nu=\nu^1\circ\nu^2$, $A\sim\nu$.  Writing $A^i=A\cap\nu^i$,
we have that 
\begin{equation}
\wt^A(\nu)=\wt^{A^1}(\nu^1)\wt^{A^2}_{K\setminus\nu^1}(\nu^2).
\end{equation}
(b) Let $\nu = \ol\nu^1\circ \ol\nu^2$ be a cut of the backbone $\nu$ at
the point $x$, and $A \sim \nu$. Then
\begin{equation}
\wt^A(\nu)=\wt^{B^1}(\ol\nu^1)\wt^{B^2}_{K\setminus\ol\nu^1}(\ol\nu^2).
\end{equation}
where $B^i=A^i\cup\{x\}$.
\end{lemma}

\begin{proof}
By \eqref{ihp16}, the first claim is equivalent to
\begin{equation}
\frac{Z_{K\setminus\nu}}{Z}\one\{A\sim\nu\}=
\frac{Z_{K\setminus\nu^1}}{Z}\one\{A^1\sim\nu^1\}
\frac{Z_{K\setminus(\nu^1\cup\nu^2)}}{Z_{K\setminus\nu^1}}\one\{A^2\sim\nu^2\}.
\end{equation}
The right side vanishes if and only if the left side vanishes. When
both sides are non-zero, their equality follows from the fact that
$Z_{K\setminus\nu}=Z_{K\setminus(\nu^1\cup\nu^2)}$. The second claim
follows similarly, on adding $x$ to the set of sources.
\end{proof}

\section{The switching lemma}\label{sw_sec}

We state and prove next the
principal tool in the random-parity representation, namely
the so-called `switching lemma'. In brief, this
allows us to take two independent colourings, with different sources, and to
`switch' the sources from one to the other in a measure-preserving way.
In so doing,  
the backbone will generally change.  In order to
preserve the measure, the \emph{connectivities}
inherent in the backbone must be retained. We begin by defining
two notions of connectivity in colourings. We work throughout this
section in the general set-up of Section~\ref{ssec-col}.

\subsection{Connectivity and switching}\label{ssec-switching}

Let $B\se F$, $G\se K$ be finite sets, let 
$A \subseteq\ol K$ be a finite set of sources, 
and write $\psi^A=\psi^A(B,G)$ for the colouring
given in the last section.  In what follows we think of the ghost-bonds
as bridges to the ghost-site~$\Gh$. 

Let 
$x,y\in K^\Gh:= K \cup\{\Gh\}$. A \emph{path} from $x$ to $y$ 
in the configuration $(B,G)$ is a 
self-avoiding path with endpoints $x$, $y$,
traversing intervals of $K^\Gh$, and possibly bridges in $B$ and/or ghost-bonds 
joining $G$ to $\Gh$.  
Similarly, a \emph{cycle} is a self-avoiding cycle in the
above graph. A \emph{route} is a path or a cycle. 
A route containing no ghost-bonds is called a
\emph{lattice-route}.  A route is called \emph{odd} (in the colouring $\psi^A$)
if $\psi^A$, when restricted to the route, takes only the value `odd'. 
The failed colouring $\psi^A=\#$ is deemed to contain no odd
paths.  

Let $B_1,B_2\se F$, $G_1,G_2\se K$, and
let $\psi_1^A=\psi_1^A(B_1,G_1)$ and $\psi_2^B=\psi_2^B(B_2,G_2)$
be the associated colourings.
Let $\D$%
\nomenclature[D]{$\D$}{Process of cuts}
 be an auxiliary Poisson process on $K$, with
intensity function $4\d(\cdot)$, that is independent of all other random
variables so far. We call points of
$\D$ \emph{cuts}.  A route of $(B_1\cup B_2, G_1\cup G_2)$
is said to be \emph{open} in the triple
$(\psi_1^A,\psi_2^B,\D)$ if it includes no 
sub-interval of $\ev(\psi_1^A)\cap\ev(\psi_2^B)$ containing
one or more  elements of $\D$.
In other words, the cuts break paths, but only when they
fall in intervals labelled `even' in \emph{both} colourings.  See Figure~\ref{connectivity_fig}.
In particular, if there is an odd path $\pi$ from $x$
to $y$ in $\psi_1^A$, then $\pi$ constitutes an open path in
$(\psi_1^A,\psi_2^B,\D)$ irrespective of $\psi_2^B$ and $\D$.  We let
\begin{equation}
\{x\lra y\mbox{ in }\psi_1^A,\psi_2^B,\D\}
\end{equation}
be the event that there exists an open path from $x$ to $y$ in 
$(\psi_1^A,\psi_2^B,\D)$. We may abbreviate this to
$\{x\lra y\}$ when there is no ambiguity. 

\begin{figure}[tbp]
\includegraphics{thesis.12}
\hspace{1.3cm} 
\includegraphics{thesis.13} 
\hspace{1.3cm} 
\includegraphics{thesis.14} 
\caption[Connectivity in the random-parity representation]
{Connectivity in pairs of colourings.
\emph{Left}:  $\psi_1^{ac}$.  Middle:  $\psi_2^\es$.
\emph{Right}:  the triple $\psi_1^{ac},\psi_2^\es,\D$. 
Crosses are elements of $\D$ and grey lines are
where either $\psi_1^{ac}$ or $\psi_2^\es$ is odd.  In
$(\psi_1^{ac},\psi_2^\es,\D)$ the following connectivities hold:
$a\nlra b$, $a\lra c$, $a\lra d$,
$b\nlra c$, $b\nlra d$,
$c\lra d$.  The dotted line marks $\pi$, one of the open paths from
$a$ to $c$.} 
\label{connectivity_fig}
\end{figure}

There is an analogy between open paths in the above construction
and the notion of connectivity in the
random-current representation of the discrete Ising model. Points
labelled `odd' or
`even' above may be considered as collections of infinitesimal parallel edges,
being odd or even in number, respectively.  If a point is `even', the
corresponding
number of edges may be $2,4,6,\dotsc$ \emph{or} it may be 0;  in the `union'
of $\psi_1^A$ and $\psi_2^B$, connectivity is broken at a point if and only if
both the corresponding numbers equal 0.  It turns out that the correct law for
the set of such points is that of~$\D$.

Here is some notation.  For any finite sequence $(a,b,c,\dots)$ of elements in $K$, 
the string  
$abc\dotsc$ will denote the subset of elements
that appear an odd number of times in the sequence.  
If $A\subseteq \ol K$ 
is a finite set with odd
cardinality, then for any pair $(B,G)$ for which there exists a valid 
colouring $\psi^A(B,G)$, the
number of ghost-bonds must be odd.  Thinking of these as bridges to $\Gh$, $\Gh$ may thus 
be viewed as an element of $A$, and we make the following remark.

\begin{remark}\label{gGremark}
For $A \subseteq \ol K$ with $|A|$ odd,
we shall use the expressions $\psi^A$ and $\psi^{A\cup\{\Gh\}}$ interchangeably.
\end{remark}
 
We call a function $F$, acting on $(\psi_1^A,\psi_2^B,\D)$, a
\emph{connectivity function} if it depends only on the connectivity
properties using open paths of $(\psi_1^A,\psi_2^B,\D)$, that is, the value
of $F$ depends only on the set $\{(x,y)\in (K^\Gh)^2: x \lra y\}$.
In the following, $E$ denotes expectation with respect to
$d\mu_\l d\mu_\g dM_{B,G} dP$, where $P$ is the law of~$\D$.

\begin{theorem}[Switching lemma]\label{sl}
Let $F$ be a connectivity function and $A,B\subseteq \ol K$ finite sets.
For $x,y\in K^\Gh$,
\begin{align}
\label{sw_eq_1}
&E\bigl(\partial\psi_1^A\partial\psi_2^B\cdot F(\psi_1^A,\psi_2^B,\D)
\cdot \one\{x\lra y\mbox{ in }\psi_1^A,\psi_2^B,\D\}\bigr)\\
&\hskip1cm =E\Big(\partial\psi_1^{A\sd xy}\partial\psi_2^{B\sd xy}
\cdot F(\psi_1^{A\sd xy},\psi_2^{B\sd xy},\D)\cdot \nonumber\\ 
&\hskip4cm \cdot \one\{x\lra y\mbox{ in }\psi_1^{A\sd xy},\psi_2^{B\sd xy},\D\}\Big).
\nonumber
\end{align}
In particular,
\begin{equation}\label{sw_eq_2}
E(\partial\psi_1^{xy}\partial\psi_2^B)
=E\bigl(\partial\psi_1^\es\partial\psi_2^{B\sd xy}
\cdot\one\{x\lra y\mbox{ in }
\psi_1^\es,\psi_2^{B\sd xy},\D\}\bigr).
\end{equation}
\end{theorem}

\begin{proof}
Equation \eqref{sw_eq_2} follows from \eqref{sw_eq_1} with $A=\{x,y\}$
and $F\equiv 1$, and so it suffices to prove \eqref{sw_eq_1}.
This is trivial if $x=y$, and we assume henceforth that $x\neq y$.
Recall that $W=\{v\in V:K_v=\SS\}$ and $|W|=r$.  

We prove \eqref{sw_eq_1} first for the special case when $F\equiv 1$, that is,
\begin{multline}\label{sw_eq_3}
E\bigl(\partial\psi_1^A\partial\psi_2^B
\cdot \one\{x\lra y\mbox{ in }\psi_1^A,\psi_2^B,\D\}\bigr)\\
=E\bigl(\partial\psi_1^{A\sd xy}\partial\psi_2^{B\sd xy}\cdot
\one\{x\lra y\mbox{ in }\psi_1^{A\sd xy},\psi_2^{B\sd xy},\D\}\bigr),
\end{multline}
and this will follow by conditioning on the pair $Q=(B_1\cup B_2,G_1\cup G_2)$.
 
Let $Q$ be given. Conditional on 
$Q$, the law of $(\psi_1^A,\psi_2^B)$ is given as follows.
First, we
allocate each bridge and each ghost-bond to either $\psi_1^A$ or
$\psi_2^B$ with equal probability (independently of one another).
If $W \ne \es$, then we must also allocate (uniform) random colours
to the points $(w,0)$, $w \in W$, for each of $\psi_1^A$, $\psi_2^B$.
If $(w,0)$ is itself a source, we work instead with $(w,0+)$.
(Recall that the pair $(B',G')$ may be reconstructed from
knowledge of a valid colouring $\psi^{A'}(B',G')$.)
There are $2^{|Q|+2r}$ possible outcomes of the above choices, and each is
equally likely.

The process $\D$ is independent of
all random variables used above.  Therefore, the conditional expectation,
given $Q$, of the random variable on the left side of \eqref{sw_eq_3} equals  
\begin{equation}\label{sw_cond_eq}
\frac{1}{2^{|Q|+2r}}\sum _{\Qab}\partial Q_1\partial Q_2\,
P(x\lra y\mbox{ in }Q_1,Q_2,\D),
\end{equation}
where the sum is over the set $\Qab=\Qab(Q)$ of all possible pairs $(Q_1,Q_2)$ of
values of $(\psi_1^A,\psi_2^B)$. 
The measure $P$ is that of~$\D$. 

We shall define an invertible (and therefore measure-preserving) map 
from $\Qab$ to $\Qabxy$. Let $\pi$ be a path of $Q$ with endpoints $x$ and $y$
(if such a path $\pi$ exists), and let $f_\pi:\Qab\to\Qabxy$
be given as follows.
Let $(Q_1,Q_2)\in\Qab$, say $Q_1=Q_1^A(B_1,G_1)$
and $Q_2=Q_2^B(B_2,G_2)$ where $Q=(B_1\cup B_2, G_1\cup G_2)$.
For $i=1,2$, let $B_i'$ (\resp, $G_i'$) be the set of bridges (\resp, ghost-bonds)
in $Q$ lying in exactly one
of $B_i$, $\pi$ (\resp, $G_i$, $\pi$). Otherwise expressed, $(B_i',G_i')$ is obtained
from $(B_i,G_i)$ by adding the bridges/ghost-bonds of $\pi$ `modulo 2'.
Note that $(B_1'\cup B_2', G_1' \cup G_2') = Q$.

If $W = \es$, we let $R_1=R_1^{A \sd xy}$ (\resp, $R_2^{B\sd xy}$) be the
unique valid colouring of $(B_1',G_1')$ with sources $A\sd xy$
(\resp, $(B_2',G_2')$ with sources $B\sd xy$), so
$R_1=\psi^{A\sd xy}(B_1',G_1')$, and similarly for $R_2$.
When $W\ne \es$ and $i=1,2$, we choose the colours of the $(w,0)$, $w\in W$, in $R_i$ in such a way
that $R_i \equiv Q_i$ on $K\sm\pi$. 

It is easily seen that the map $f_\pi:(Q_1,Q_2) \mapsto (R_1,R_2)$
is invertible, indeed its inverse is given by the same mechanism.
See Figure~\ref{connectivity_after_switch_fig}.

\begin{figure}[tbp]
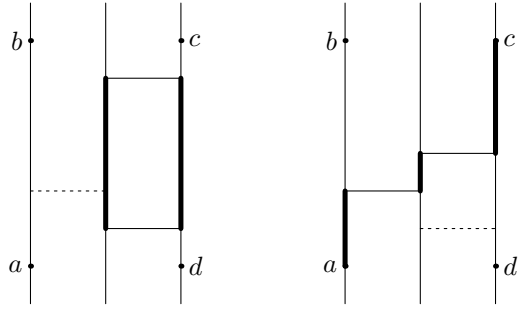

\includegraphics{thesis.15}
\hspace{1.3cm} 
\includegraphics{thesis.16} 
\caption[Switching]
{Switched configurations.  Taking $Q_1^{ac}$, $Q_2^\es$ and
$\pi$ to be $\psi_1^{ac}$, $\psi_2^\es$ and $\pi$ of
Figure~\ref{connectivity_fig}, respectively, this figure illustrates the
`switched' configurations $R_1^\es$ and $R_2^{ac}$ (left and
right, respectively).}
\label{connectivity_after_switch_fig}
\end{figure}

By \eqref{def-wt}, 
\begin{equation}
\partial Q_1\partial Q_2=
\exp\bigl\{2\d(\ev(Q_1))+2\d(\ev(Q_2))\bigr\}.
\label{o14}
\end{equation}
Now,
\begin{align}\label{o15}
\d(\ev(Q_i))
&=\d(\ev(Q_i)\cap\pi)+\d(\ev(Q_i)\setminus\pi)\\
&=\d(\ev(Q_i)\cap\pi)+\d(\ev(R_i)\setminus\pi), 
\nonumber
\end{align}
and
\begin{align*}
&\d(\ev(Q_1)\cap\pi)+\d(\ev(Q_2)\cap\pi) - 2\d\bigl(\ev(Q_1)\cap\ev(Q_2)\cap\pi\bigr) \\
&\hskip1cm =\d\bigl(\ev(Q_1)\cap\odd(Q_2)\cap\pi\bigr)+\d\bigl(\odd(Q_1)\cap\ev(Q_2)\cap\pi\bigr)\\
&\hskip1cm=\d\bigl(\odd(R_1)\cap\ev(R_2)\cap\pi\bigr)+\d\bigl(\ev(R_1)\cap\odd(R_2)\cap\pi\bigr)\\
&\hskip1cm =\d(\ev(R_1)\cap\pi)+\d(\ev(R_2)\cap\pi) -2\d\bigl(\ev(R_1)\cap\ev(R_2)\cap\pi\bigr),
\end{align*}
whence, by \eqref{o14}--\eqref{o15},
\begin{align}\label{switched_weights_eq}
\partial Q_1\partial Q_2=
\partial R_1\partial R_2&
\exp\bigl\{-4\d\bigl(\ev(R_1)\cap\ev(R_2)\cap\pi\bigr)\bigr\}\\ 
&\times \exp\bigl\{4\d\bigl(\ev(Q_1)\cap\ev(Q_2)\cap\pi\bigr)\bigr\}.
\nonumber
\end{align}

The next step is to choose a suitable path $\pi$.
Consider the final term in \eqref{sw_cond_eq}, namely 
\begin{equation}
P(x\lra y\mbox{ in }Q_1,Q_2,\D).
\end{equation}
There are finitely many paths in $Q$ from $x$ to $y$, let these paths be
$\pi_1,\pi_2,\dotsc,\pi_n$. Let 
$\cO_k=\cO_k(Q_1,Q_2,\D)$ be the event that
$\pi_k$ is the earliest such path that is open in $(Q_1,Q_2,\D)$.  Then 
\begin{align}\label{path_prob_eq}
&\hskip-1cm P(x\lra y\mbox{ in }Q_1,Q_2,\D)\\
&=\sum_{k=1}^n P(\cO_k)\nonumber\\
&=\sum_{k=1}^n P\bigl(\D\cap[\ev(Q_1)\cap\ev(Q_2)\cap \pi_k] = \es\bigr)
P(\wtilde\cO_k)\nonumber\\
&=\sum_{k=1}^n 
\exp\bigl\{-4\d\bigl(\ev(Q_1)\cap\ev(Q_2)\cap\pi_k\bigr)\bigr\}
P(\wtilde\cO_k),
\nonumber
\end{align}
where $\wtilde\cO_k = \wtilde\cO_k(Q_1,Q_2,\D)$ is the 
event  that each of  
$\pi_1,\dotsc,\pi_{k-1}$ is rendered non-open in $(Q_1,Q_2,\D)$
through the presence of elements of $\D$ lying in $K\sm\pi_k$.
In the second line of \eqref{path_prob_eq}, we have used the independence
of $\D \cap \pi_k$ and $\D \cap(K\sm \pi_k)$. 

Let $(R_1^k,R_2^k) =f_{\pi_k}(Q_1,Q_2)$.  
Since $R_i^k \equiv Q_i$ on $K \sm \pi_k$, we have that
$\wtilde\cO_k(Q_1,Q_2,\D) = \wtilde\cO_k(R_1^k,R_2^k,\D)$. 
By
\eqref{switched_weights_eq} and~\eqref{path_prob_eq}, the
summand in \eqref{sw_cond_eq} equals
\begin{align*}
& \sum_{k=1}^n \partial Q_1\partial Q_2
\exp\bigl\{-4\d\bigl(\ev(Q_1)\cap\ev(Q_2)\cap\pi_k\bigr)\bigr\}P(\wtilde\cO_k)\\
&\hskip1cm=\sum_{k=1}^n \partial R_1^k\partial R_2^k
\exp\bigl\{-4\d\bigl(\ev(R_1^k)\cap\ev(R_2^k)\cap\pi_k\bigr)\bigr\}
P(\wtilde\cO_k)\\
&\hskip1cm=\sum_{k=1}^n \partial R_1^k
\partial R_2^k \,P(\cO_k(R_1^k,R_2^k,\D)).
\end{align*}

Summing the above over $\Qab$, and remembering that each $f_{\pi_k}$
is
a bijection between $\Qab$ and $\Qabxy$, \eqref{sw_cond_eq} becomes
\begin{align*}
\frac1{2^{|Q|+2r}} \sum_{k=1}^n\, 
& \sum_{(R_1,R_2)\in\Qabxy}\partial R_1\partial R_2\,
P(\cO_k(R_1,R_2,\D))\\
&= \frac1{2^{|Q|+2r}} \sum_{\Qabxy}\partial R_1\partial R_2\,
P(x\lra y\mbox{ in }R_1,R_2,\D).
\end{align*}
By the argument leading to \eqref{sw_cond_eq}, this equals the right
side of \eqref{sw_eq_3}, and the claim is proved
when $F\equiv 1$.  

Consider now the case of general connectivity functions $F$ in \eqref{sw_eq_1}.
In~\eqref{sw_cond_eq}, the factor 
$P(x\lra y\mbox{ in }Q_1,Q_2,\D)$ is replaced by 
$$
P\bigl(F(Q_1,Q_2,\D)\cdot \one\{x\lra y\mbox{ in }Q_1,Q_2,\D\}\bigr),
$$
where $P$ is expectation with respect to $\D$.
In the calculation~\eqref{path_prob_eq}, we use the fact that
$$
P(F\cdot \one_{\cO_k})=P(F\mid\cO_k)P(\cO_k)
$$
and we deal with the factor $P(\cO_k)$ as before.  The result follows on noting
that, for each $k$,
$$
P\bigl(F(Q_1,Q_2,\D)\bigmid \cO_k(Q_1,Q_2,\D)\bigr)=
P\bigl(F(R_1^k,R_2^k,\D)\bigmid \cO_k(R_1^k,R_2^k,\D)\bigr).
$$
This holds because: (i) the configurations $(Q_1,Q_2,\D)$ and $(R_1^k,R_2^k,\D)$ are
identical off $\pi_k$, and (ii) in each, all points along $\pi_k$ are
connected. Thus the connectivities are identical in the two configurations.
\end{proof}

\subsection{Applications of switching}\label{sw_appl_sec}

In this section are presented a number of inequalities and identities
proved using the random-parity representation and 
the switching lemma.  With some exceptions 
(most notably~\eqref{dd_bound_eq}) the
proofs are adaptations of 
the proofs for the discrete Ising 
model that may be found in \cite{abf,grimmett_rcm}.

For $R\subseteq K$ a finite union of intervals, let
\begin{equation*}
\wtilde R:=\{(uv,t)\in F: \mbox{either } (u,t)\in R
\mbox{ or }(v,t)\in R\mbox{ or both}\}.
\end{equation*}
Recall that $W=W(K)=\{v\in V: K_v=\SS\}$, and $N=N(K)$ is the 
total number of intervals constituting~$K$.

\begin{lemma}\label{rw_mon_lem}
Let $R\subseteq K$ be finite union of intervals, 
and let $\nu\in\Xi$ be such that $\nu\cap
R=\es$.  If $A \subseteq \ol{K\sm R}$ is finite and  $A\sim\nu$, then
\begin{equation}
\wt^A(\nu)\leq 2^{r(\nu)-r'(\nu)}\wt^A_{K\setminus R}(\nu),
\end{equation}
where 
\begin{align*}
r(\nu) &= r(\nu,K) := |\{w\in W: \nu\cap (w\times K_w) \ne \es\}|,\\
r'(\nu) &= r(\nu,K\sm R).
\end{align*}
\end{lemma}%
\nomenclature[r]{$r(\nu)$}{The number of intersection points with $W$}

\begin{proof}
By \eqref{ihp16} and  Lemma~\ref{Z'},
\begin{align}\label{r0}
\wt^A(\nu)&=\frac{Z_{K\setminus\nu}}{Z_K}\\
&=2^{N(K)-N(K\sm\nu)}
e^{\l(\wtilde\nu)+\g(\nu)-\d(\nu)}\frac{Z'_{K\sm\nu}}{Z'_K}.
\nonumber
\end{align}
We claim that
\begin{equation}\label{r05}
\frac{Z'_{K\sm\nu}}{Z'_K}\leq\frac{Z'_{K\sm(R\cup\nu)}}{Z'_{K\sm R}},
\end{equation}
and the proof of this follows.

Recall the formula~\eqref{o12} for $Z'_K$ in terms of an integral over the
Poisson process $D$.  The set $D$ is the union of
independent Poisson processes $D'$ and $D''$, restricted respectively to $K\sm\nu$ and $\nu$.
We write $P'$ (\resp, $P''$) for the probability measure (and, on occasion,
expectation operator) governing $D'$ (\resp, $D''$).
Let $\S(D')$ denote the set of spin configurations on $K\sm\nu$
that are permitted by $D'$.  By \eqref{o12},
\begin{equation}
\label{Z'_split_eq}
Z'_K= P'\left( \sum_{\s'\in\S(D')}Z_\nu'(\s') 
\exp\left\{\int_{F\sm\wtilde\nu}\l(e)\s'_e\,de+\int_{K\sm\nu}\g(x)\s'_x\,dx\right\}\right),
\end{equation}
where
$$
Z_\nu'(\s') = P''\left(\sum_{\s''\in\wtilde\S(D'')}
\exp\left\{\int_{\wtilde\nu}\l(e)\s_e\,de+\int_{\nu}\g(x)\s_x\,dx\right\}\cdot \one_C(\s')\right)
$$
is the partition function on $\nu$ with boundary condition $\s'$,
and where $\s$, $\wtilde\S(D'')$, and $C=C(\s',D'')$ are given as follows.

The set $D''$ divides $\nu$, in the usual way, into a collection 
$V_\nu(D'')$ of intervals.  From the set of endpoints of
such intervals, we distinguish the subset $\cE$ that: (i) lie in $K$, and (ii) 
are endpoints of some interval of $K\sm \nu$. For $x\in\cE$,
let $\s'_x =\lim_{y\to x} \s'_y$, where the limit is taken over $y\in K\sm\nu$.
Let $\wtilde V_\nu(D'')$ be the subset of $V_\nu(D'')$ containing
those intervals with no endpoint in $\cE$,
and let $\wtilde\S(D'') =\{-1,+1\}^{\wtilde V_\nu(D'')}$.

Let $\s'\in \S(D')$, and 
let $\cI$ be the set of maximal sub-intervals $I$ of $\nu$
having both endpoints in $\cE$, and such that $I \cap D''=\es$.
Let $C=C(D'')$ be the set of $\s'\in\S(D')$ such that, for all $I\in\cI$,
the endpoints of $I$ have equal spins under $\s'$.
Note that
\begin{equation}
\label{o30}
\one_C(\s') = \prod_{I\in\cI} \tfrac12(\s'_{x(I)}\s'_{y(I)} + 1),
\end{equation}
where $x(I)$, $y(I)$ denote the endpoints of~$I$.

Let $\s''\in \wtilde\S(D'')$. The conjunction $\s$ of $\s'$ and $\s''$ is defined
except on sub-intervals of $\nu$ lying in $V_\nu(D'')\sm \wtilde V_\nu(D'')$. 
On any such sub-interval with exactly one endpoint $x$ in $\cE$, we set
$\s\equiv \s'_x$. On the event $C$, an interval of $\nu$
with both endpoints $x(I)$, $y(I)$ in $\cE$
receives the spin $\s\equiv \s_{x(I)} = \s_{y(I)}$.
Thus, $\s\in \S(D'\cup D'')$ is well defined for $\s'\in C$.

By \eqref{Z'_split_eq}, 
\begin{equation*}
\frac{Z'_K}{Z'_{K\sm\nu}}=\el Z'_\nu(\s')\er_{K\sm\nu}.
\end{equation*}
Taking the expectation $\el\cdot\er_{K\sm\nu}$ inside the integral, the last expression becomes
\begin{equation*}
P''\left(\sum_{\s''\in\wtilde\S(D'')}\left\el
\exp\left\{\int_{\wtilde\nu}\l(e)\s_e\,de\right\}
\exp\left\{\int_{\nu}\g(x)\s_x\,dx\right\}
\cdot \one_C(\s')\right\er_{K\sm\nu}\right)
\end{equation*}
The inner expectation may be expressed as a sum over $k,l\geq 0$ 
(with non-negative coefficients) of iterated integrals of the form
\begin{equation}
\frac1{k!}\,\frac1{l!}\,\iint\limits_{\wtilde\nu^k\times\nu^l}\l(\mathbf e)\g(\mathbf x)
\el\s_{e_1}\cdots\s_{e_k}\s_{x_1}\cdots\s_{x_l}\cdot \one_C\er_{K\sm\nu}
\,d\mathbf e \,d\mathbf x,
\label{o32}
\end{equation}
where we have written $\mathbf e=(e_1,\dotsc,e_k)$, and $\l(\mathbf e)$
for $\l(e_1)\dotsb\l(e_k)$ (and similarly for $\mathbf x$). 
We may write
\begin{equation*}
\el\s_{e_1}\cdots\s_{e_k}\s_{x_1}\cdots\s_{x_l}\cdot \one_C\er_{K\sm\nu}
=\el\s'_S\s''_T\cdot \one_C\er_{K\sm\nu}=\s''_T\el\s_S'\cdot \one_C\er_{K\sm\nu},
\end{equation*}
for sets $S\subseteq \ol{K\sm\nu}$, $T\subseteq \nu$
 determined by $e_1,\dotsc,e_k,x_1,\dotsc,x_l$ and $D''$ only.
We now bring the sum over $\s''$ inside the integral of \eqref{o32}.  For $T\neq\es$,
\begin{equation*}
\sum_{\s''\in\wtilde\S(D'')}\s''_T\el\s_S'\cdot \one_C\er_{K\sm\nu}=0,
\end{equation*}
so any non-zero term is of the form
\begin{equation}\label{r1}
\el\s_S'\cdot \one_C\er_{K\sm\nu}.
\end{equation}

By \eqref{o30}, \eqref{r1} 
may be expressed in the form 
\begin{equation}
\sum_{i=1}^s2^{-a_i}\el\s'_{S_i}\er_{K\sm\nu}
\label{o35}
\end{equation}
for appropriate sets $S_i$ and integers $a_i$.  By Lemma~\ref{cor_mon_lem},
\begin{equation*}
\el\s'_{S_i}\er_{K\sm\nu}\geq\el\s'_{S_i}\er_{K\sm(R\cup\nu)}.
\end{equation*}
On working backwards, we obtain \eqref{r05}.

By \eqref{r0}--\eqref{r05},
\begin{equation*}
\wt^A(\nu)\leq 2^U\wt^A_{K\setminus R}(\nu),
\end{equation*}
where 
\begin{align*}
U&=\bigl[N(K)-N(K\sm\nu)\bigr]- \bigl[N(K\sm R)-N(K\sm (R\cup\nu) )\bigr]\\
&=r(\nu)-r'(\nu)
\end{align*}
as required.  
\end{proof}

For distinct $x,y,z\in K^\Gh$, let
\begin{align*}
\el\s_x;\s_y;\s_z\er &:=
\el\s_{xyz}\er -\el\s_{x}\er\el\s_{yz}\er\\
&\hskip1.5cm -\el\s_{y}\er\el\s_{xz}\er
-\el\s_{z}\er\el\s_{xy}\er
+2\el\s_{x}\er\el\s_{y}\er\el\s_{z}\er.
\end{align*}

\begin{lemma}[\ghs\ inequality]\label{ghs_lem}
For distinct $x,y,z\in K^\Gh$, 
\begin{equation}\label{ghs_1_eq}
\el\s_x;\s_y;\s_z\er\leq 0.
\end{equation}
Moreover, $\el\s_x\er$ is concave in $\g$ in the sense that, for
bounded, measurable functions $\g_1,\g_2: K\to\RRp$ satisfying
$\g_1\le\g_2$, and
$\theta\in[0,1]$, 
\begin{equation}
\theta\el\s_x\er_{\g_1}+(1-\theta)\el\s_x\er_{\g_2}\leq
\el\s_x\er_{\theta\g_1+(1-\theta)\g_2}.
\end{equation}
\end{lemma}

\begin{proof}
The proof of this follows very closely the corresponding proof for the
classical Ising model~\cite{ghs}. We include it here because it allows us
to develop the technique of `conditioning on clusters', which will be 
useful later.

We prove~\eqref{ghs_1_eq} via the following more general result.
Let $(B_i,G_i)$, $i=1,2,3$, be independent sets of bridges/ghost-bonds,
and write $\psi_i$, $i=1,2,3$, for corresponding colourings (with sources
to be specified through their superscripts).
We claim that, for any four points $w, x, y, z\in K^\Gh$,
\begin{equation}
\label{ghs_2_eq}
\begin{split}
&E\bigl(\partial\psi_1^\es\partial\psi_2^\es \partial\psi_3^{wxyz}\bigr)-
E\bigl(\partial\psi_1^\es\partial\psi_2^{wz}\partial\psi_3^{xy}\bigr)
\\
&\quad\leq E(\partial\psi_1^\es \partial\psi_2^{wx} \partial\psi_3^{yz})
+E(\partial\psi_1^\es\partial\psi_2^{wy}\partial\psi_3^{xz})
-2E(\partial\psi_1^{wx}\partial\psi_2^{wy}\partial\psi_3^{wz}).
\end{split}
\end{equation}
Inequality \eqref{ghs_1_eq} follows by Theorem \ref{rcr_thm} on letting
$w=\Gh$. 

The left side of \eqref{ghs_2_eq} is
\begin{align*}
&E(\partial\psi_1^\es)\bigl[
E(\partial\psi_2^\es\partial\psi_3^{wxyz})-
E(\partial\psi_2^{wz}\partial\psi_3^{xy})\bigr]\\
&\hskip3cm =
Z\, E\bigl(\partial\psi_2^\es\partial\psi_3^{wxyz}
\cdot\one\{w\nlra z\}\bigr),
\end{align*}
by the switching lemma \ref{sl}.  When $\partial\psi_3^{wxyz}$  is
non-zero, parity constraints imply that at least one of  
$\{w\lra x\}\cap \{y\lra z\}$ and $\{w\lra y\}\cap 
\{x\lra z\}$ occurs, but that, in the presence of the indicator function 
they cannot both occur.  Therefore,
\begin{align}
\label{ghs_pf_1_eq}
&E(\partial\psi_2^\es\partial\psi_3^{wxyz}
\cdot\one\{w\nlra z\})\\
&\hskip1cm
=E\bigl(\partial\psi_2^\es\partial\psi_3^{wxyz}
\cdot\one\{w\nlra z\}
\cdot\one\{w\lra x\}\bigr)\nonumber\\
&\hskip3cm +
E\bigl(\partial\psi_2^\es\partial\psi_3^{wxyz}
\cdot\one\{w\nlra z\}
\cdot\one\{w\lra y\}\bigr).
\nonumber
\end{align}
Consider the first term. By the switching lemma,
\begin{equation}
E\bigl(\partial\psi_2^\es\partial\psi_3^{wxyz}
\cdot\one\{w\nlra z\}
\cdot\one\{w\lra x\}\bigr)=
E\bigl(\partial\psi_2^{wx}\partial\psi_3^{yz}
\cdot\one\{w\nlra z\}\bigr).
\label{o17}
\end{equation}

We next `condition on a cluster'.  Let 
$C_z=C_z(\psi_2^{wx},\psi_3^{yz},\D)$ be the set of all points of
$K$ that are connected by open paths to $z$.  Conditional on 
$C_z$, define new independent colourings
$\mu_2^\es$, $\mu_3^{yz}$ on the domain $M=C_z$.  Similarly,
let $\nu_2^{wx}$, $\nu_3^\es$ be independent colourings on the domain
$N=K\setminus C_z$, that are also independent of the $\mu_i$.  It is
not hard to see that, if $w\nlra z$ in $(\psi_2^{wx},\psi_3^{yz},\D)$, then, conditional on $C_z$,
the law of $\psi_2^{wx}$ equals 
that of the superposition of $\mu_2^\es$ and $\nu_2^{wx}$;
similarly the conditional law of $\psi_3^{yz}$ is the same as that of the
superposition of $\mu_3^{yz}$ and $\nu_3^\es$. 
Therefore, almost surely on the event $\{w \nlra z\}$,
\begin{align}
E(\partial\psi_2^{wx}\partial\psi_3^{yz}\mid C_z)&=
E'(\pd\mu_2^\es)E'(\pd\nu_2^{wx})E'(\pd\mu_3^{yz})E'(\pd\nu_3^\es)\label{o19}\\
&=\el\s_{wx}\er_N E'(\pd\mu_2^\es)E'(\pd\nu_2^\es)
E'(\pd\mu_3^{yz})E'(\pd\nu_3^\es)\nonumber\\
&\leq
\el\s_{wx}\er_KE(\partial\psi_2^\es\partial\psi_3^{yz}\mid C_z),
\nonumber\end{align}
where $E'$ denotes expectation conditional on $C_z$,
and we have used Lemma \ref{cor_mon_lem}. 
Returning to \eqref{ghs_pf_1_eq}--\eqref{o17},
\begin{align*}
&E\bigl(\partial\psi_2^\es\partial\psi_3^{wxyz}
\cdot\one\{w\nlra z\}
\cdot\one\{w\lra x\}\bigr)\\
&\hskip2cm \leq
\el\s_{wx}\er E(\partial\psi_2^\es\partial\psi_3^{yz}
\cdot\one\{w\nlra z\}).
\end{align*}
The other term in \eqref{ghs_pf_1_eq} satisfies the same inequality
with $x$ and $y$ interchanged.
Inequality \eqref{ghs_2_eq} follows on applying the switching lemma to the right sides
of these two last inequalities, and adding them.

The concavity of $\el\s_x\er$ follows from the fact that, if 
\begin{equation}
T=\sum_{k=1}^n a_k\one_{A_k} 
\end{equation}
is a step function on $K$ with $a_k\ge 0$ for all $k$,
and $\g(\cdot)=\g_1(\cdot)+\a T(\cdot)$, then 
\begin{equation}
\frac{\partial^2}{\partial \a^2}\el\s_x\er
=\sum_{k,l=1}^na_ka_l\iint_{A_k\times A_l} dy\,dz\, \el\s_x;\s_y;\s_z\er\leq 0.
\end{equation}
Thus, the claim holds whenever $\g_2-\g_1$ is a step function.  The general
claim follows by approximating $\g_2-\g_1$ by step functions, and applying the
dominated convergence theorem. 
\end{proof}

For the next lemma we assume for simplicity that $\g\equiv 0$ (although
similar results can easily be proved for $\g\not\equiv 0$).  We let $\bar\d\in\RR$ be
an upper bound for $\d$, thus $\d(x)\leq\bar\d<\oo$ for all $x\in K$.
Let $a,b\in K$ be two distinct points.  A closed set $T\subseteq K$
is said to \emph{separate} $a$ from 
$b$ if every lattice path from $a$ to $b$ (whatever the set of bridges) intersects $T$.  
Moreover, if $\eps>0$
and $T$ separates $a$ from $b$, we say that $T$ is an \emph{$\eps$-fat
separating set} if every point in $T$ lies in a closed sub-interval of $T$ of length at
least~$\eps$.

\begin{lemma}[Simon inequality]\label{simon_lem}
Let $\g\equiv 0$.
If $\eps>0$ and $T$ is an $\eps$-fat separating set for $a,b\in K$,
\begin{equation}
\el\s_a\s_b\er\leq\frac{1}{\eps}\exp(8\eps\bar\d)
\int_T\el\s_a\s_x\er\el\s_x\s_b\er\, dx.
\end{equation}
\end{lemma}

\begin{proof}
By Theorems \ref{rcr_thm} and \ref{sl},
\begin{equation}
\el\s_a\s_x\er\el\s_x\s_b\er=
\frac{1}{Z^2}E(\partial\psi_1^\es\partial\psi_2^{ab}
\cdot\one\{a\lra x\}),
\end{equation}
and, by Fubini's theorem,
\begin{equation}
\int_T\el\s_a\s_x\er\el\s_x\s_b\er\;dx=
\frac{1}{Z^2}E(\partial\psi_1^\es\partial\psi_2^{ab}
\cdot|\what T|),
\end{equation}
where $\what T=\{x\in T:a\lra x\}$ and $|\cdot|$ denotes Lebesgue
measure.  Since $\g\equiv 0$, the backbone $\xi = \xi(\psi_2^{ab})$ 
consists of a single (lattice-) path from $a$ to $b$
passing through $T$.  Let $U$ denote the set of points in $K$
that are separated from $b$ by $T$, and let $X$ be the point at which $\xi$
exits $U$ for the first time.  Since $T$ is
assumed closed, $X\in T$.
See Figure \ref{simon_fig}.  

\begin{figure}[tbp]
\includegraphics{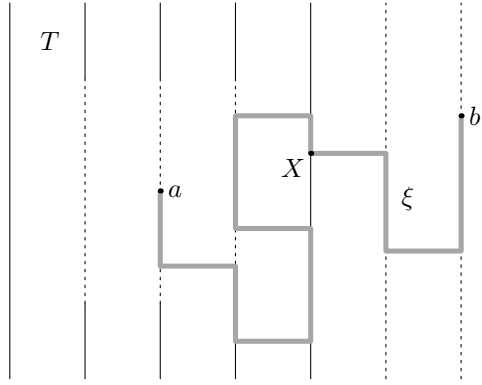}
\caption[The Simon inequality]
{The Simon inequality.  The separating set $T$ is drawn with solid
  black lines, and the backbone $\xi$ with a grey line.}
\label{simon_fig}
\end{figure}

For  
$x\in T$, let $A_x$ be the event that there is no element of $\D$ within
the interval of length $2\eps$ centered at $x$.  Thus,
$P(A_x)=\exp(-8\eps\bar\d)$.  On the
event $A_X$, we have that $|\what T|\geq\eps$, whence
\begin{align}\label{simon_eq1_eq}
E(\partial\psi_1^\es\partial\psi_2^{ab}
\cdot|\what T|)&\geq
E(\partial\psi_1^\es\partial\psi_2^{ab}
\cdot|\what T|\cdot\one\{A_X\})\\
&\geq\eps
E(\partial\psi_1^\es\partial\psi_2^{ab}
\cdot\one\{A_X\}).
\nonumber
\end{align}
Conditional on $X$, the event $A_X$ is independent of $\psi_1^\es$ and
$\psi_2^{ab}$, so that
\begin{equation}
E(\partial\psi_1^\es\partial\psi_2^{ab}
\cdot|\what T|)\geq\eps\exp(-8\eps\bar\d)
E(\partial\psi_1^\es\partial\psi_2^{ab}),
\label{o18}
\end{equation}
and the proof is complete.
\end{proof}

Just as for the classical Ising model, only a small amount of extra work is
required to deduce the following improvement of Lemma~\ref{simon_lem}.

\begin{lemma}[Lieb inequality]\label{lieb_lem}
Under the assumptions of Lemma \ref{simon_lem},
\begin{equation}
\el\s_a\s_b\er\leq\frac{1}{\eps}\exp(8\eps\bar\d)
\int_T\el\s_a\s_x\er_{\overline{T}}\,\el\s_x\s_b\er\;dx,
\end{equation}
where $\el\cdot\er_{\ol T}$ denotes expectation with respect to the
measure restricted to~$\ol T$.
\end{lemma}

\begin{proof}
Let $x \in T$, let $\ol\psi_1^{ax}$ denote a colouring on the restricted region $U$, and let $\psi_2^{xb}$ denote a colouring on the
full region $K$ as before.  We claim that
\begin{equation}\label{lieb_sw_eq}
E(\partial\ol\psi_1^{ax}\partial\psi_2^{xb})=
E\bigl(\partial\ol\psi_1^\es\partial\psi_2^{ab}
\cdot \one\{a\lra x\mbox{ in }\overline{T}\}\bigr).
\end{equation}
The use of the letter $E$ is an abuse of notation,
since the $\ol\psi$ are colourings of $U$ only.

Equation \eqref{lieb_sw_eq} may be established using a slight variation in the proof of the
switching lemma.  We follow the proof of that lemma, first
conditioning on the set $Q$ of all bridges and ghost-bonds in the two
colourings taken together, and then allocating them to the colourings $Q_1$ and $Q_2$,
uniformly at random.  We then order the paths $\pi$ of $Q$ from $a$ to
$x$, and add the earliest open path to both $Q_1$ and $Q_2$ `modulo 2'.  There
are two differences here:  firstly, any element of $Q$ that is not
contained in $U$ will be allocated to $Q_2$, 
and secondly, we only consider paths $\pi$ that lie 
inside $U$.  Subject to these two changes, we follow the argument of
the switching lemma to arrive at~\eqref{lieb_sw_eq}.  

Integrating \eqref{lieb_sw_eq} over $x \in T$,
\begin{equation}
\int_T\el\s_a\s_x\er_{\overline{T}}\,\el\s_x\s_b\er\;dx=
\frac{1}{Z_{\ol T}Z}E(\partial\ol\psi_1^\es\partial\psi_2^{ab}
\cdot|\what T|),
\end{equation}
where this time 
$\what T=\{x\in T:a\lra x\mbox{ in }U\}$.  
The proof is completed as in \eqref{simon_eq1_eq}--\eqref{o18}.
\end{proof}

For the next lemma we specialize to the situation that is the main focus of
this chapter, namely the following. 
Similar results are valid for other lattices and for
summable translation-invariant interactions.

\begin{assumption}\label{periodic_assump}\hspace{1cm}
\begin{itemize}
\item The graph $L=[-n,n]^d \subseteq \ZZ^d$ where $d \ge 1$,
with periodic boundary condition.
\item The parameters $\l$, $\d$, $\g$ are non-negative constants.
\item The set $K_v=\SS$ for every $v\in V$.
\end{itemize}
\end{assumption}
Under the periodic boundary condition, two vertices of $L$ are 
joined by an edge whenever there exists $i\in\{1,2,\dots,d\}$
such that their $i$-coordinates differ by exactly~$2n$.  

Under Assumption~\ref{periodic_assump}, 
the model is invariant under automorphisms of $L$ and, furthermore,
the quantity
$\el\s_x\er$ does not depend on the choice of $x$.  
Let $0$ denote some fixed but arbitrary point of $K$, and let
$M=M(\l,\d,\g)=\el\s_0\er$ denote the common value of the $\el\s_x\er$.  

For $x,y\in K$, we write $x\sim y$ if $x=(u,t)$ and
$y=(v,t)$ for some $t \ge 0$ and $u,v$ adjacent in $L$.  We write 
$\{x\lrao z y\}$ for the complement of the event that there exists
an open path from $x$ to $y$ not containing $z$.
Thus, $x\lrao z y$ if: either $x \nlra y$, or $x \lra y$ and every
open path from $x$ to $y$ passes through~$z$.

\begin{theorem}\label{three_ineq_lem}
Under Assumption~\ref{periodic_assump}, the following hold. 
\begin{align}\label{dg_bound_eq}
\frac{\partial M}{\partial\g}&=\frac{1}{Z^2}\int_K dx\;
E\bigl(\partial\psi_1^{0x}\partial\psi_2^\es
\cdot \one\{0\nlra \Gh\}\bigr)
\leq \frac{M}{\g}.\\
\label{dl_bound_eq}
\frac{\partial M}{\partial\l}&=\frac{1}{2Z^2}\int_K dx
\sum_{y\sim x} E\bigl(\partial\psi_1^{0xy\Gh}\partial\psi_2^\es
\cdot \one\{0\nlra \Gh\}\bigr)
\leq 2dM\frac{\partial M}{\partial\g}.\\
\label{dd_bound_eq}
-\frac{\partial M}{\partial\d}&=
\frac{2}{Z^2}\int_K dx \: E\bigl(\partial\psi_1^{0\Gh}\partial\psi_2^\es
\cdot \one\{0\overset{x}{\lra} \Gh\}\bigr)
\leq \frac{2M}{1-M^2}\frac{\partial M}{\partial\g}.
\end{align}
\end{theorem}

\begin{proof}
With the exception of~\eqref{dd_bound_eq}, the proofs mimic those 
of \cite{abf} for the classical Ising model.
For the equality in~\eqref{dg_bound_eq}, note that
\[
\frac{\partial M}{\partial \g}=\int_K  \el \s_0;\s_x\er \,dx.
\]
Now
\[
\el\s_0;\s_x\er=\el \s_0\s_x\er - \el \s_0\er\el \s_x\er=
\frac{1}{Z^2}(E(\partial\psi_1^{0x}\partial\psi_2^\varnothing)-
E(\partial\psi_1^{0\G}\partial\psi_2^{x\G}))
\] 
and the difference $E(\partial\psi_1^{0x}\partial\psi_2^\varnothing)-
E(\partial\psi_1^{0\G}\partial\psi_2^{x\G})$ on the right hand side
equals
\[
E(\partial\psi_1^{0x}\partial\psi_2^\varnothing)-
E(\partial\psi_1^{0x}\partial\psi_2^\varnothing\cdot
\one\{0\leftrightarrow \G\})
=E(\partial\psi_1^{0x}\partial\psi_2^\varnothing
\cdot \one\{0\not\leftrightarrow \G\})
\]
by the switching lemma.
For the inequality in~\eqref{dg_bound_eq}, the
concavity of $M$ in $\g$ means that for all $\g_2\geq\g_1>0$,
\begin{equation}
\frac{\partial M}{\partial\g}\leq
\frac{M(\l,\d,\g_2)-M(\l,\d,\g_1)}{\g_2-\g_1}.
\end{equation}
Letting $\g_1\rightarrow0$ and using the continuity of $M$ and the fact
that $M(\l,\d,0)=0$ for all $\l,\d>0$, the result follows. 

Similarly, for the equality in~\eqref{dl_bound_eq} we note that
\[
\frac{\partial M}{\partial \l}=\int_F  \el \s_0;\s_e\er\,de
=\frac{1}{2}\int_K dx\sum_{y\sim x}
(\el \s_0\s_x\s_y\er - \el \s_0\er\el \s_x\s_y\er).
\]
Again
\[\begin{split}
\el \s_0\s_x\s_y\er - \el \s_0\er\el \s_x\s_y\er&=
\frac{1}{Z^2}(E(\partial\psi_1^{0xy\G}\partial\psi_2^\varnothing)-
E(\partial\psi_1^{0\G}\partial\psi_2^{xy}))\\
&=E(\partial\psi_1^{0xy\G}\partial\psi_2^\varnothing
\cdot \one\{0\not\leftrightarrow \G\})
\end{split}\]
by the switching lemma.  For the inequality,
\begin{equation}
\begin{split}
\frac{\partial M}{\partial \l} &=
\frac{1}{2}\int_K dx\sum_{y\sim x} 
\big(\el\s_0\s_x\s_y\er-\el\s_0\er\el\s_x\s_y\er\big)\\
&\leq\frac{1}{2}\int_K dx\sum_{y\sim x} 
\big(\el\s_x\er\el\s_0\s_y\er+\el\s_y\er\el\s_0\s_x\er
-2\el\s_0\er\el\s_x\er\el\s_y\er\big)\\
&=\int_K dx\:\el\s_0;\s_x\er\sum_{y\sim x} 
\el\s_y\er\\
&=2dM\int_K dx\:\el\s_0;\s_x\er=
2dM\frac{\partial M}{\partial \g},
\end{split}
\end{equation}
where we have used the \ghs-inequality and translation invariance.

Here is the proof of \eqref{dd_bound_eq}. Let $|\cdot|$ denote Lebesgue
measure.  By differentiating 
\begin{equation}
M=\frac{E(\partial\psi^{0\Gh})}{E(\partial\psi^\es)}=
\frac{E(\exp(2\d|\ev(\psi^{0\Gh})|))}{E(\exp(2\d|\ev(\psi^\es)|))},
\end{equation}
with respect to $\d$, we obtain that
\begin{align}
\frac{\partial M}{\partial \d}&=
\frac{2}{Z^2}E\bigl(\partial\psi_1^{0\Gh}\partial\psi_2^\es\cdot
\bigl[|\ev(\psi_1^{0\Gh})|-|\ev(\psi_2^\es)|\bigr]\bigr)\label{ihp17}\\
&=\frac{2}{Z^2}\int dx\,
E\bigl(\partial\psi_1^{0\Gh}\partial\psi_2^\es\cdot
\bigl[\one\{x\in\odd(\psi_2^\es)\}-\one\{x\in\odd(\psi_1^{0\Gh})\}\bigr]\bigr).
\nonumber
\end{align}

Consider the integrand in \eqref{ihp17}.
Since $\psi_2^\es$ has no sources, all odd routes in $\psi_2^\es$ are
necessarily cycles.  If $x\in\odd(\psi_2^\es)$, then $x$ lies in an odd
cycle.  We shall assume that $x$ is not the endpoint of a bridge, since this
event has probability 0.  It follows that, on the event $\{0\lra \Gh\}$,
there exists an open path from $0$ to $\Gh$ that avoids $x$
(since any path can be re-routed around the odd cycle of $\psi_2^\es$ containing
$x$). 
Therefore, the event $\{0\lrao{x} \Gh\}$ does not occur, and hence
\begin{align} 
&E\bigl(\partial\psi_1^{0\Gh}\partial\psi_2^\es\cdot
\one\{x\in\odd(\psi_2^\es)\}\bigr)\label{dd_pf_1_eq}\\
&\hskip2cm =E\bigl(\partial\psi_1^{0\Gh}\partial\psi_2^\es\cdot
\one\{x\in\odd(\psi_2^\es)\}\cdot 
\one\{0\lrao{x} \Gh\}^\tc\bigr).
\nonumber
\end{align}
  
We note next that, if $\partial\psi_1^{0\Gh}\ne 0$ and
$0\lrao x\Gh$, then necessarily
$x\in\odd(\psi_1^{0\Gh})$. Hence,
\begin{align}\label{dd_pf_2_eq}
&E\bigl(\partial\psi_1^{0\Gh}\partial\psi_2^\es\cdot
\one\{x\in\odd(\psi_1^{0\Gh})\}\bigr)\\
&\hskip1cm =E\bigl(\partial\psi_1^{0\Gh}\partial\psi_2^\es\cdot
\one\{x\in\odd(\psi_1^{0\Gh})\}\cdot \one\{0\overset{x}{\lra} \Gh\}^\tc\bigr)\nonumber\\
&\hskip5cm +E\bigl(\partial\psi_1^{0\Gh}\partial\psi_2^\es\cdot 
\one\{0\overset{x}{\lra} \Gh\}\bigr).
\nonumber
\end{align}
We wish to switch the sources $0\Gh$ from $\psi_1$ to $\psi_2$ in the right
side of \eqref{dd_pf_2_eq}.
For this we need to adapt some details of the proof of the 
switching lemma to this situation.  
The first step in the proof of that lemma
was to condition on the union $Q$ of the bridges and
ghost-bonds of the two colourings;  then, the paths from $0$ to $\Gh$ in $Q$ were listed
in a fixed \emph{but arbitrary} order. 
We are free to choose this
ordering in such a way that paths not containing $x$ have precedence, and
we assume henceforth that the ordering is thus chosen.  The next step is
to find the earliest open path $\pi$, and `add $\pi$ modulo 2' to both
$\psi_1^{0\Gh}$ and $\psi_2^\es$.  On the event 
$\{0\overset{x}{\lra} \Gh\}^\tc$, this earliest path $\pi$ does not
contain $x$, by our choice of ordering.  Hence, in the new colouring
$\psi_1^\es$,  $x$ continues to lie in an `odd' interval (recall
that, outside $\pi$, the colourings are unchanged by the switching procedure).
Therefore,
\begin{align}
&E\bigl(\partial\psi_1^{0\Gh}\partial\psi_2^\es\cdot
\one\{x\in\odd(\psi_1^{0\Gh})\}\cdot \one\{0\overset{x}{\lra} \Gh\}^\tc\bigr)\\
&\hskip 2cm =E\bigl(\partial\psi_1^\es\partial\psi_2^{0\Gh}\cdot
\one\{x\in\odd(\psi_1^\es)\}\cdot \one\{0\overset{x}{\lra} \Gh\}^\tc\bigr).
\nonumber
\end{align}
Relabelling, putting the last expression into~\eqref{dd_pf_2_eq}, and
subtracting~\eqref{dd_pf_2_eq} from~\eqref{dd_pf_1_eq}, we obtain
\begin{equation}\label{dd_pf_3_eq}
\frac{\partial M}{\partial \d}=
-\frac{2}{Z^2}\int dx\:
E\bigl(\partial\psi_1^{0\Gh}\partial\psi_2^\es\cdot 
\one\{0\overset{x}{\lra} \Gh\}\bigr)
\end{equation}
as required.

Turning to the inequality, let  $C^x_z$ denote the set of
points that can be reached from $z$ along open paths 
\emph{not containing $x$}.  When conditioning 
$E\bigl(\partial\psi_1^{0\Gh}\partial\psi_2^\es\cdot 
\one\{0\overset{x}{\lra} \Gh\}\bigr)$ on $C^x_0$ as in the proof of the
\ghs\ inequality, we find that $\psi_1^{0\Gh}$ is a combination of two independent
colourings, one inside $C^x_0$ with sources $0x$, and one outside $C^x_0$
with sources $x\Gh$.  As in \eqref{o19}, using Lemma \ref{cor_mon_lem} as there,
\begin{align}
\label{dd_pf_4_eq}
E\bigl(\partial\psi_1^{0\Gh}\partial\psi_2^\es
\cdot \one\{0\overset{x}{\lra} \Gh\}\bigr)&=
E\bigl(\partial\psi_1^{0x}\partial\psi_2^\es\el\s_x\er_{K\setminus C^x_0}
\cdot \one\{0\overset{x}{\lra} \Gh\}\bigr)\\ 
&\leq M\cdot E\bigl(\partial\psi_1^{0x}\partial\psi_2^\es
\cdot \one\{0\overset{x}{\lra} \Gh\}\bigr).
\nonumber
\end{align}
We split the expectation on the right side according to whether or not
$x\lra \Gh$.  Clearly,
\begin{equation}\label{dd_pf_5_eq}
E\bigl(\partial\psi_1^{0x}\partial\psi_2^\es
\cdot \one\{0\overset{x}{\lra} \Gh\}
\cdot \one\{x\nlra \Gh\}\bigr)\leq
E\bigl(\partial\psi_1^{0x}\partial\psi_2^\es
\cdot \one\{x\nlra \Gh\}\bigr).
\end{equation}
By the switching lemma \ref{sl}, the other term satisfies
\begin{equation}
E\bigl(\partial\psi_1^{0x}\partial\psi_2^\es
\cdot \one\{0\overset{x}{\lra} \Gh\}
\cdot \one\{x\lra \Gh\}\bigr)=
E\bigl(\partial\psi_1^{0\Gh}\partial\psi_2^{x\Gh}
\cdot \one\{0\overset{x}{\lra} \Gh\}\bigr).
\end{equation}
We again condition on a cluster, this time $C^x_\Gh$, to obtain as in
\eqref{dd_pf_4_eq} that
\begin{equation}\label{dd_pf_6_eq}
E\bigl(\partial\psi_1^{0\Gh}\partial\psi_2^{x\Gh}
\cdot \one\{0\overset{x}{\lra} \Gh\}\bigr)\leq 
M\cdot E\bigl(\partial\psi_1^{0\Gh}\partial\psi_2^\es
\cdot \one\{0\overset{x}{\lra} \Gh\}\bigr).
\end{equation}
Combining \eqref{dd_pf_4_eq}, \eqref{dd_pf_5_eq}, \eqref{dd_pf_6_eq} with \eqref{dd_pf_3_eq}, 
we obtain by \eqref{dg_bound_eq} that 
\begin{equation}
-\frac{\partial M}{\partial \d}\leq 2M\frac{\partial M}{\partial \g}+
M^2\Big(-\frac{\partial M}{\partial \d}\Big),
\end{equation}
as required.
\end{proof}

\section{Proof of the main differential inequality}\label{pf_sec}

In this section we will prove Theorem \ref{main_pdi_thm},
the differential inequality which, in
combination with the inequalities of the previous section, will
yield information about the critical behaviour of the space--time Ising model.  The
proof proceeds roughly as follows.  In the
random-parity representation of $M=\el\s_0\er$, there is a backbone from $0$
to $\Gh$ (that is, to some point $g \in G$).  
We introduce two new  sourceless configurations;
depending on how the backbone interacts with these configurations, the
switching lemma allows a decomposition into a combination of
other configurations which, via Theorem \ref{three_ineq_lem}, may be 
transformed into derivatives of the magnetization.

Throughout this section we work under Assumption~\ref{periodic_assump},
  that is, \emph{we work with a translation-invariant model on 
a cube in the $d$-dimensional
  lattice}, while noting that our conclusions are valid
for more general interactions with similar symmetries.
The arguments in this section borrow heavily from~\cite{abf}.
As in Theorem \ref{three_ineq_lem}, the main novelty in the proof
concerns connectivity in the `vertical' direction 
(the term $R_v$ in \eqref{ihp19}--\eqref{ihp20}
below). 

\begin{proof}[Proof of Theorem~\ref{main_pdi_thm}]
By Theorem \ref{rcr_thm},
\begin{equation}
M=\frac{1}{Z}E(\partial\psi_1^{0\Gh})
=\frac{1}{Z^3}E(\partial\psi_1^{0\Gh}\partial\psi_2^\es
\partial\psi_3^\es).
\end{equation}
We shall consider the backbone
$\xi=\xi(\psi_1^{0\Gh})$ and the open cluster $C_\Gh$ of $\Gh$ in 
$(\psi_2^\es,\psi_3^\es,\D)$.
All connectivities will refer
to the triple $(\psi_2^\es,\psi_3^\es,\D)$.  Note that $\xi$ consists of
a single path with endpoints $0$ and $\Gh$.
There are four possibilities, illustrated in Figure~\ref{decomp_fig}, for the way in which $\xi$,
viewed as a directed path from $0$ to $\Gh$, interacts with $C_\Gh$:  
\begin{romlist}
\item $\xi\cap C_\Gh$ is empty, 
\item $0 \in \xi\cap C_\Gh$, 
\item $0 \notin \xi\cap C_\Gh$, and $\xi$ first meets $C_\Gh$ immediately
  after a bridge, 
\item $0 \notin \xi\cap C_\Gh$, and $\xi$ first meets $C_\Gh$ at a cut, which
necessarily belongs to $\ev(\psi_2^\es)\cap\ev(\psi_3^\es)$.  
\end{romlist}

Thus,
\begin{equation}
M=T+R_0+R_h+R_v,
\label{ihp19}
\end{equation}
where
\begin{equation}
\label{ihp20}
\begin{split}
T&=\frac{1}{Z^3}E\bigl(\partial\psi_1^{0\Gh}\partial\psi_2^\es
\partial\psi_3^\es\cdot \one\{\xi\cap C_\Gh=\es\}\bigr),\\
R_0&=\frac{1}{Z^3}E\bigl(\partial\psi_1^{0\Gh}\partial\psi_2^\es
\partial\psi_3^\es\cdot \one\{0\lra \Gh\}\bigr),\\
R_h&=\frac{1}{Z^3}E\bigl(\partial\psi_1^{0\Gh}\partial\psi_2^\es
\partial\psi_3^\es\cdot 
\one\{\mbox{first point on $\xi\cap C_\Gh$ is a bridge of $\xi$}\}\bigr),\\
R_v&=\frac{1}{Z^3}E\bigl(\partial\psi_1^{0\Gh}\partial\psi_2^\es
\partial\psi_3^\es\cdot 
\one\{\mbox{first point on $\xi\cap C_\Gh$ is a cut}\}\bigr).
\end{split}
\end{equation}
We will bound each of these terms separately.

\begin{figure}[tbp]
\includegraphics{thesis.19}
\hspace{1cm}
\includegraphics{thesis.20}
\hspace{1cm}
\includegraphics{thesis.21}
\hspace{1cm}
\includegraphics{thesis.22}
\caption[The four cases in Theorem~\ref{main_pdi_thm}]
{Illustrations of the four possibilities for $\xi\cap C_\Gh$.  
Ghost-bonds in $\psi^{0\Gh}$ are labelled $g$. The backbone $\xi$ is drawn as a
  solid black line, and $C_\Gh$ as a grey rectangle.}
\label{decomp_fig}
\end{figure}

By the switching lemma,
\begin{align}
R_0&=\frac{1}{Z^3}E\bigl(\partial\psi_1^{0\Gh}\partial\psi_2^\es
\partial\psi_3^\es\cdot \one\{0\lra \Gh\}\bigr)\label{ihp21}\\
&=\frac{1}{Z^3}E\bigl(\partial\psi_1^{0\Gh}\partial\psi_2^{0\Gh}
\partial\psi_3^{0\Gh}\bigr)=M^3.
\nonumber
\end{align}

Next, we bound $T$. The letter $\xi$ will always denote the backbone
of the first colouring $\psi_1$, with corresponding sources. 
Let $X$ denote the location of the ghost-bond that ends
$\xi$.  By conditioning on $X$, 
\begin{equation}
\begin{split}
T&=\frac{1}{Z^3}\int P(X\in dx)\,E\bigl(\partial\psi_1^{0\Gh}\partial\psi_2^\es
\partial\psi_3^\es\cdot \one\{\xi\cap C_\Gh=\es\}\bigmid X=x\bigr)\\
&\leq \frac{\g}{Z^3}\int dx\, E\bigl(\partial\psi_1^{0x}\partial\psi_2^\es
\partial\psi_3^\es\cdot \one\{\xi\cap C_\Gh=\es\}\bigr).
\end{split}
\label{o51}
\end{equation}
We study the last expectation by conditioning on $C_\Gh$ and bringing one of
the factors $1/Z$ inside.  By \eqref{backb_cond_eq}--\eqref{ihp16} and conditional
expectation, 
\begin{align}
&\frac{1}{Z}E\bigl(\partial\psi_1^{0x}\cdot \one\{\xi\cap C_\Gh=\es\}\bigmid C_\Gh\bigr)\label{o23}\\
&\hskip3cm =E\Bigl( Z^{-1}E(\pd\psi_1^{0x} \mid \xi, C_\Gh)\one\{\xi\cap C_\Gh=\es\} \Bigmid C_\Gh\Bigr)
\nonumber\\
&\hskip3cm =E\bigl(\wt^{0x}(\xi)\cdot \one\{\xi\cap C_\Gh=\es\}\bigmid C_\Gh\bigr).
\nonumber
\end{align}
By Lemma \ref{rw_mon_lem},
\begin{equation}
\wt^{0x}(\xi) \le 2^{r(\xi)-r'(\xi)}
\wt_{K\setminus C_\Gh}^{0x}(\xi)\quad\mbox{on}\quad \{\xi\cap C_\Gh = \es\},
\label{o22}
\end{equation}
where
\begin{equation*}
r(\xi) = r(\xi,K),\qquad 
r'(\xi) = r(\xi,K\sm C_\Gh).
\end{equation*}
Using~\eqref{special2} and~\eqref{backbone_rep_eq}, we have
\begin{align}
&E\bigl(\wt^{0x}(\xi)
\cdot \one\{\xi\cap C_\Gh=\es\}\bigmid C_\Gh\bigr)\label{o24} \\
&\hskip3cm \leq E\bigl( 2^{r(\xi)-r'(\xi)}\wt_{K\setminus C_\Gh}^{0x} (\xi)\cdot
\one\{\xi\cap C_\Gh=\es\}\bigmid C_\Gh\bigr)\nonumber \\
&\hskip3cm  \le\el\s_0\s_x\er_{K\setminus C_\Gh}.
\nonumber\end{align}

The last step merits explanation. 
Recall that $\xi=\xi(\psi_1^{0x})$, and assume $\xi\cap C_\Gh=\es$.
Apart from the randomization that takes place when $\psi_1^{0x}$ is one of
several valid colourings,  the law of $\xi$, $P(\xi\in d\nu)$,
 is a function of the positions of
bridges and ghost-bonds along $\nu$ only, that is, the existence of
bridges where needed, and the non-existence of ghost-bonds along $\nu$.
By \eqref{o22} and Lemma~\ref{rw_mon_lem},
with $\Xi_{K\sm C} := \{\nu\in \Xi: \nu\cap C = \es\}$ and $P$
the law of $\xi$,
\begin{align*}
&E\bigl(\wt^{0x}(\xi)\cdot \one\{\xi\cap C_\Gh=\es\}\bigmid C_\Gh\bigr)\\
&\hskip3cm =\int_{\Xi_{K\sm C_\Gh}} w^{0x}(\nu)\, P(d\nu)\\
&\hskip3cm \le \int_{\Xi_{K\sm C_\Gh}}2^{r(\nu)-r'(\nu)} w^{0x}_{K\sm C_\Gh}(\nu) 
\left(\tfrac12\right)^{r(\nu)} \mu(d\nu)  
\end{align*}
for some measure $\mu$, where the factor $(\frac12)^{r(\nu)}$ arises from the
possible existence of more than one valid colouring.  
Now, $\mu$ is a measure on paths which by the remark above depends only locally
on $\nu$, in the sense that $\mu(d\nu)$ depends only on the bridge- and
ghost-bond configurations along $\nu$.  In particular, the same measure $\mu$
governs also the law of the backbone in the \emph{smaller} region 
$K\setminus C_\Gh$.  More explicitly, by \eqref{backbone_rep_eq} with
$P_{K\sm C_\Gh}$ the law of the backbone of the colouring
$\psi_{K\sm  C_\Gh}^{0x}$ defined on $K\sm C_\Gh$, we have
\begin{align*}
\el\s_0\s_x\er_{K\setminus C_\Gh}
&=\int_{\Xi_{K\sm C_\Gh}} w^{0x}_{K\sm C_\Gh}(\nu)\, P_{K\sm C_\Gh}(d\nu)\\
&= \int_{\Xi_{K\sm C_\Gh}} w^{0x}_{K\sm C_\Gh}(\nu) 
\left(\tfrac12\right)^{r'(\nu)}\, \mu(d\nu).
\end{align*}
Thus \eqref{o24} follows. 

Therefore, by \eqref{o51}--\eqref{o24},
\begin{align}
T&\leq \frac{\g}{Z^2}\int dx\: E\bigl(\partial\psi_2^\es
\partial\psi_3^\es\el\s_0\s_x\er_{K\setminus C_\Gh}\cdot
\one\{0\nlra \Gh\}\bigr)\\
&=\g\int dx\: \frac{1}{Z^2}E\bigl(\partial\psi_2^{0x}
\partial\psi_3^\es\cdot \one\{0\nlra \Gh\}\bigr)\nonumber\\
&=\g\frac{\partial M}{\partial \g},
\nonumber
\end{align}
by `conditioning on the cluster' $C_\Gh$
and Theorem~\ref{three_ineq_lem}.

Next, we bound $R_h$.  Suppose that the bridge bringing $\xi$
into 
$C_\Gh$ has endpoints $X$ and $Y$, where we take $X$ to be the endpoint not in
$C_\Gh$.  When the bridge $XY$ is removed, the backbone $\xi$ consists of two
paths:  $\zeta^1:0\rightarrow X$ and $\zeta^2:Y\rightarrow \Gh$.  Therefore,
\begin{align*}
R_h&= \frac{1}{Z^3}\int P(X\in dx)\,E\bigl(\partial\psi_1^{0\Gh}
\partial\psi_2^\es\partial\psi_3^\es\bigmid X=x\bigr)\\
&\le\frac{\l}{Z^3}\int dx \, \sum_{y\sim x}E\bigl(\pd\psi_1^{0xy\Gh}
\partial\psi_2^\es\partial\psi_3^\es\cdot\one\{0\nlra\Gh,\,y\lra\Gh\}\cdot \one\{J_\xi\}\bigr),
\nonumber
\end{align*}
where $\xi=\xi(\psi_1^{0xy\Gh})$ and
$$
J_\xi=\bigl\{\xi=\zeta^1\circ\zeta^2,
\,\zeta^1:0\rightarrow x,\,\zeta^2:y\rightarrow \Gh,\,
\zeta^1\cap C_\Gh=\es\bigr\}.
$$
As in \eqref{o23},
\begin{equation}
R_h
\leq \frac{\l}{Z^2}\int dx\:\sum_{y\sim x}
E\bigl(\partial\psi_2^\es\partial\psi_3^\es
\cdot \one\{0\nlra \Gh,\, y\lra \Gh\}\cdot
\wt^{0xy\Gh}(\xi)\cdot\one\{J_\xi\}\bigr).
\label{ihp23}
\end{equation}

By Lemmas \ref{backb2}(a) and \ref{rw_mon_lem}, on the event $J_\xi$,
\begin{align*}
\wt^{0xy\Gh}(\xi) &= \wt^{0x}(\zeta^1) \wt^{y\Gh}_{K\sm \zeta^1}(\zeta^2)\\
&\le 2^{r-r'}\wt^{0x}_{K\sm C_\Gh}(\zeta^1)\wt^{y\Gh}_{K\sm \zeta^1}(\zeta^2),
\end{align*}
where $r = r(\zeta^1,K)$ and $r'= r(\zeta^1, K\sm C_\Gh)$.
By Lemma~\ref{cor_mon_lem} and the reasoning after~\eqref{o24},
\begin{align*}
E\bigl(\wt^{0xy\Gh}(\xi) \cdot \one\{J_\xi\}\bigmid \zeta^1, C_\Gh\bigr) 
&\leq 2^{r-r'}
\wt_{K\setminus C_\Gh}^{0x}(\zeta^1)\cdot \el\s_y\er_{K\sm \zeta^1}\\
&\leq M\cdot 2^{r-r'} \wt_{K\setminus C_\Gh}^{0x}(\zeta^1),
\end{align*}
so that, similarly, 
\begin{equation}
E\bigl(\wt^{0xy\Gh}(\xi)\cdot \one\{J_\xi\}\bigmid C_\Gh\bigr) \le 
M \cdot \el\s_0\s_x\er_{K\setminus C_\Gh}.
\label{o25}
\end{equation}
We substitute into the summand in \eqref{ihp23}, using the switching
lemma, conditioning on the cluster $C_\Gh$, and the bound 
$\el\s_y\er_{C_\Gh}\leq M$, to obtain the upper bound
\begin{align}
&M\cdot E\bigl(\partial\psi_2^\es\partial\psi_3^\es\cdot 
\one\{0\nlra \Gh,\, y\lra \Gh\}\cdot
\el\s_0\s_x\er_{K\setminus C_\Gh}\bigr)\\
&\hskip2.5cm =M\cdot E\bigl(\partial\psi_2^{y\Gh}\partial\psi_3^{y\Gh}
\cdot \one\{0\nlra \Gh\}\cdot \el\s_0\s_x\er_{K\setminus C_\Gh}\bigr)\nonumber\\
&\hskip2.5cm =M\cdot E\bigl(\partial\psi_2^{0xy\Gh}\partial\psi_3^\es
\el\s_y\er_{C_\Gh} \cdot \one\{0\nlra \Gh\}\bigr)\nonumber\\
&\hskip2.5cm \leq M^2\cdot E\bigl(\partial\psi_2^{0xy\Gh}\partial\psi_3^\es
\cdot \one\{0\nlra \Gh\}\bigr).
\nonumber
\end{align}
Hence, by \eqref{dl_bound_eq},
\begin{align*}
R_h &\leq \l M^2\frac{1}{Z^2}\int dx\,\sum_{y\sim x}
E\bigl(\partial\psi_2^{0xy\Gh}\partial\psi_3^\es
\one\{0\nlra \Gh\}\bigr)\\
&=2\l M^2 \frac{\partial M}{\partial \l}.
\end{align*}

Finally, we bound $R_v$.  Let 
$X\in\D\cap\ev(\psi_2^\es)\cap\ev(\psi_3^\es)$ be the first
point of $\xi$ in $C_\Gh$.  In a manner similar to that used for $R_h$ 
at \eqref{ihp23} above,
and by cutting the backbone $\xi$ at the point $x$,
\begin{equation}
R_v\le\frac{1}{Z^2}\int P(X\in dx)\,
E\bigl(\partial\psi_2^\es\partial\psi_3^\es
\cdot \one\{0\nlra \Gh,\,x\lra \Gh\}\cdot
\wt^{0\Gh}(\xi)\cdot \one\{J_\xi\}\bigr),
\label{ihp22}
\end{equation}
where 
$$
J_\xi=  1\bigl\{\xi=\ol\zeta^1\circ\ol\zeta^2,\, \ol\zeta^1:0\rightarrow x,\,\ol\zeta^2:x\rightarrow \Gh,
\, \zeta^1\cap C_\Gh=\es\bigr\}.
$$
As in \eqref{o25},
\begin{align*}
E(\wt^{0\Gh}(\xi)\cdot \one\{J_\xi\} \mid C_\Gh)
&= E\bigl(E(\wt^{0\Gh}(\xi)\cdot \one\{J_\xi\}\mid \ol\zeta^1,C_\Gh)\bigmid C_\Gh\bigr)\\
&\leq E\bigl(\el\s_0\s_x\er_{K\setminus C_\Gh}\cdot \el\s_x\er_{K\setminus\zeta^1}\bigmid C_\Gh\bigr)\\
&\leq \el\s_0\s_x\er_{K\setminus C_\Gh}\cdot M.
\end{align*}
By \eqref{ihp22} therefore,
\begin{equation*}
R_v\leq M\frac{1}{Z^2}\int P(X\in dx)\,
E\bigl(\partial\psi_2^\es\partial\psi_3^\es
\cdot \one\{0\nlra \Gh,\,x\lra \Gh\}
\el\s_0\s_x\er_{K\setminus C_\Gh}\bigr).
\end{equation*}
By removing the cut at $x$, the origin $0$ becomes connected to $\Gh$, but only
via $x$.  Thus,
\begin{equation*}
R_v\leq 4\d M\frac{1}{Z^2}\int dx\:
E\bigl(\partial\psi_2^\es\partial\psi_3^\es
\cdot \one\{0\overset{x}{\lra} \Gh,\,x\lra \Gh\}
\el\s_0\s_x\er_{K\setminus C^x_\Gh}\bigr),
\end{equation*}
where $C^x_\Gh$ is the set of points reached from $\Gh$ along open
paths not containing $x$.  By the switching lemma, and conditioning
twice on the cluster $C_\Gh^x$,
\begin{align*}
R_v&\leq4\d M\frac{1}{Z^2}\int dx\:
E\bigl(\partial\psi_2^{x\Gh}\partial\psi_3^{x\Gh}
\cdot \one\{0\overset{x}{\lra} \Gh\}
\el\s_0\s_x\er_{K\setminus C^x_\Gh}\bigr)\\
&=4\d M\frac{1}{Z^2}\int dx\, E\bigl(\partial\psi_2^{0\Gh}\partial\psi_3^{x\Gh}
\cdot \one\{0\overset{x}{\lra} \Gh\}\bigr)\\
&=4\d M\frac{1}{Z^2}\int dx\, E\bigl(\partial\psi_2^{0\Gh}\partial\psi_3^\es
\cdot \one\{0\overset{x}{\lra} \Gh\}\el\s_x\er_{C^x_\Gh}\bigr)\\
&\leq 4\d M^2\frac{1}{Z^2}\int dx\,
E\bigl(\partial\psi_2^{0\Gh}\partial\psi_3^\es
\cdot \one\{0\overset{x}{\lra} \Gh\}\bigr)\\
&=-2\d M^2 \frac{\partial M}{\partial \d},
\end{align*}
by \eqref{dd_bound_eq}, as required.
\end{proof}

\section{Consequences of the inequalities}\label{cons_sec}

In this section we formulate the principal results
of this chapter, and show how the differential inequalities of
Theorems \ref{main_pdi_thm} and \ref{three_ineq_lem} 
may be used to prove them.  We will rely in this section on the
results in Section~\ref{inf_potts_sec}, and
we work under Assumption~\ref{periodic_assump}, unless
otherwise stated.  It is sometimes inconvenient to use periodic boundary
conditions, and we revert to the free condition where necessary.

We shall consider the infinite-volume limit as $L \uparrow \ZZ^d$;
the ground state is obtained by letting $\b \to\oo$ also.
Let $n$ be a positive integer, and set $L_n = [-n,n]^d$ with
periodic boundary condition.  
Let $\L_n^\b := [-n,n]^d \times[-\frac12 \b, \frac12\b]$.
The symbol $\b$ will appear as superscript in the following; 
the superscript $\oo$ is to be interpreted
as the ground state. Let $0=(0,0)$ and
$$
M^\b_{n}(\l,\d,\g) =\el\s_0\er_{L_{n}}^\b=\el\s_0\er_{\L^\b_{n}}
$$  
be the magnetization in $\L_{n}^\b$, noting that $M_n^\b\equiv 0$ when $\g=0$.

We have from the results in Section~\ref{press_sec} that
the limits
\begin{equation}\label{m_lim_eq}
M^\b := \lim_{n\to\oo} M^\b_n,\quad 
M^\oo := \lim_{n,\b\to\oo} M^\b_n,
\end{equation}
exist for all $\g\in\RR$ (where, in the second limit,
$\b=\b_n$ is comparable to $n$ in the sense that
Assumption~\ref{vanhove} holds).
Note that $M^\b(\l,\d,0)=0$ for $\b\in(0,\oo]$.
Recall that  we set $\d=1$, $\rho=\l/\d$,
and write
$$
M^\b(\rho,\g)= M^\b(\rho,1,\g),\qquad \b\in(0,\oo],
$$
with a similar notation
for other functions.  

Recall the following facts.
From Theorem~\ref{ising_summary_thm}  there is a unique
infinite-volume state $\el\cdot\er^\b$ at every $\g>0$.
Letting  $\el\cdot \er^\b_+$ be
the limiting state as $\g \downarrow 0$,
there is a unique state at $(\rho,0)$ if and only if
\[
M^\b_+(0):=\el\s_0\er^\b_+=0.
\]
From~\eqref{plus_ising} the state $\el\cdot\er^\b_+$ may
alternatively be obtained as the infinite volume limit of the $+$ boundary
states taken with $\g=0$.   The critical value
\begin{equation}
\bc^\b:=\inf\{\rho>0:M^\b_+(\rho)>0\},
\label{crit_val_defs_eq}
\end{equation}
see also \eqref{o1} and \eqref{critvals}. We shall have need later for
the infinite-volume limit $\el\cdot\er^{\free,\b}$,
as $n\to\oo$, with \emph{free} boundary condition in the
$\ZZ^d$ direction (or in both directions, if $\b\rightarrow\oo$).
This limit exists by Theorem~\ref{potts_lim_thm}.
Note from Theorem~\ref{ising_summary_thm} that 
\begin{equation}
\el\cdot\er^{\free,\b}_{\g=0} = \el\cdot\er^\b_{\g=0}=\el\cdot \er^\b_+ 
\quad \text{if} \quad M^\b_+(\rho)=0.
\label{o40}
\end{equation}
The superscript `f' shall always indicate the free boundary condition.

For $\b\in(0,\oo]$,
let $\phi_\rho^{b,\b}$, $b\in\{\free,\wired\}$, be the $q=2$ \rc\ measures 
of Theorem~\ref{conv_lem}, with $\g=0$.
By Theorem~\ref{correlation}, 
these measures are non-decreasing in $\rho$, and, as we saw
in~\eqref{seb5},
\begin{equation}
\label{mel61}
\fr^{\wired,\b} \le \frr^{\free,\b}, \qquad\text{when }0\leq\rho<\rho'.
\end{equation}
As in Remark~\ref{inf_vol_rk}, for $\b\in(0,\oo]$,
\begin{equation}
\fr^{\wired,\b}(x\lra y) = \el\s_x\s_y\er^\b_+, 
\quad \fr^{\wired,\b}(0\lra\oo) =M_+(\rho).
\label{mel60}
\end{equation}
By \eqref{mel60}, the
{\fkg} inequality (Theorem~\ref{rc_fkg}), 
and the uniqueness of the unbounded cluster 
(Theorem~\ref{inf_clust_uniq}),
\begin{equation}
\el\s_x \s_y\er^\b_+ \ge \fr^{\wired,\b}(x\lra\oo)
\fr^{\wired,\b}(y\lra\oo)= M^\b_+(\rho)^2.
\label{o50}
\end{equation}

Let $\b\in(0,\oo)$. Using the concavity of $M^\b$ implied
by Lemma~\ref{ghs_lem}, as well as the properties of convex functions 
in Proposition~\ref{conv_prop},
the derivative $\pd M^\b/\pd\g$ exists for all 
$\g\in\cC\subseteq(0,\infty)$, where $\cC$ is a set 
whose complement has measure zero.
When $\g\in\cC$,
\begin{equation}\label{chi_lim_eq}
\chi^\b_n(\rho,\g):=\frac{\partial M^\b_n}{\partial\g}\rightarrow
\chi^\b(\rho,\g):=\frac{\partial M^\b}{\partial\g}< \infty.
\end{equation}%
\nomenclature[c]{$\chi$}{Magnetic susceptibility}%
The corresponding conclusion holds also as $n,\b\to\oo$.  
Furthermore, by the \ghs-inequality, Lemma~\ref{ghs_lem},
$\chi^\b$ is decreasing in $\g\in\cC$, which implies
that the limits
$$
\chi^\b_+(\rho) := \lim_{\g\downarrow0} \chi^\b(\rho,\g), \qquad \b\in(0,\oo].
$$
exist when taken along sequences in~$\cC$. 

The limit
\begin{align}
\chi^{\rf,\b}(\rho,0) &:= 
\lim_{n\to\oo}\left(\left.\frac{\pd M^{\rf,\b}_n}{\pd \gamma}\right|_{\g=0}\right)
\label{mel65}\\
&=\lim_{n\to\oo}\int_{\L_n^\b} \el\s_0\s_x\er_{n,\g=0}^{\rf,\b}\,dx
=\int\el\s_0\s_x\er^{\rf,\b}_{\g=0}\,dx
\nonumber
\end{align} 
exists by monotone convergence, see Lemma \ref{cor_mon_lem}.
Let
\begin{equation}
\bs^\b :=\inf\{\rho>0:\chi^{\rf,\b}(\rho,0)=\infty\},\qquad \b\in(0,\oo].
\label{o27}
\end{equation}
We shall see in Theorem \ref{0mass} that $\chi^{\rf,\b}(\bs^\b,0)=\oo$.

It will be useful later to note that
\begin{equation}\label{i2}
\chi_+^\b(\rho)\geq\chi^{\free,\b}(\rho,0)\quad\text{whenever }
M^\b_+(\rho)=0,\qquad \b\in(0,\oo].
\end{equation}
To see this, let $\g\in\cC$ and first note from Fatou's lemma
that
\begin{equation}
\chi^\b(\rho,\g)\geq\int \el\s_0;\s_x\er^\b_\g \, dx,
\end{equation}
where we have written $\el\cdot\er^\b_\g$ for the unique state at
$\g$.  Hence, using also the monotone convergence theorem and 
the \ghs-inequality,
\begin{equation}
\chi^\b_+(\rho)=\lim_{\substack{\g\downarrow 0\\\g\in\cC}}
\chi^\b(\rho,\g)\geq\lim_{\substack{\g\downarrow 0\\\g\in\cC}}
\int\el\s_0;\s_x\er^\b_\g\, dx=\int\el\s_0;\s_x\er^\b_+\, dx.
\end{equation}
When $M_+(0)=0$ there is a unique state at $\g=0$, so that 
$\el\s_0;\s_x\er^\b_+=\el\s_0\s_x\er^{\rf,\b}_{\g=0}$
which by~\eqref{mel65} gives~\eqref{i2}.
It will follow in particular from Theorem~\ref{0mass}
that $\chi^\b_+(\rho^\b_\rs)=\oo$.  Of course, similar arguments are
valid for the limit $n,\b\rightarrow\oo$.

By~\eqref{mel65} and Lemma~\ref{cor_mon_lem} we have that
$\chi^{\free,\b}(\rho,0)$ is increasing in $\rho$.  We claim that
\begin{equation}
\bs^\b \le \bc^\b;
\label{mel64}
\end{equation}
it will follow that there is a unique equilibrium state when $\g=0$ 
and $\rho<\bs^\b$.
First note that, by~\eqref{mel61} and~\eqref{mel60}, if
$\rho<\rho'<\bs^\b$ then
\begin{equation}
M_+(\rho)=\phi^{\wired,\b}_\rho(0\leftrightarrow\oo)
\leq\phi^{\free,\b}_{\rho'}(0\leftrightarrow\oo),
\end{equation}
so it suffices to show that
$\phi^{\free,\b}_{\rho}(0\leftrightarrow\oo)=0$ if $\rho<\bs^\b$.
To see this, note that if 
$\phi^{\free,\b}_{\rho}(0\leftrightarrow\oo)>0$ then certainly 
\begin{equation}
\chi^{\rf,\b}(\rho,0)=
\int_{\ZZ^d\times[-\frac12\b,\frac12\b]}
\el\s_0\s_x\er^{\free,\b} \;dx=
\phi^\free_\rho(|C_0|)=\infty,
\end{equation}
where $C_0$ denotes the cluster at the origin, and $|\cdot|$ denotes
Lebesgue measure.

For $x\in\ZZ^d\times\RR$, let $\|x\|$ denote the supremum norm of~$x$.

\begin{theorem}\label{exp_decay_cor}
Let $\b\in(0,\oo]$ and $\rho <\bs^\b$.  There exists $\a=\a(\rho)>0$ such that
\begin{equation}
\el\s_0\s_x\er^\b \leq e^{-\a\|x\|},\qquad x\in\ZZ^d\times\RR.
\end{equation}
\end{theorem}
\begin{proof}
Fix $\b\in(0,\oo)$ and $\g=0$, and let $\rho<\bs^\b$, so that \eqref{o40} applies.  By the uniqueness of the equilibrium state, we have that
\begin{equation}
\chi^{\rf,\b}(\rho,0)=
\int_{\ZZ^d\times[-\frac12\b,\frac12\b]}
\el\s_0\s_x\er^{\b} \;dx=
\sum_{k\geq1}\int_{C_k^\b}\el\s_0\s_x\er^{\b} \;dx,
\end{equation}
where $C_k^\b:=\L^\b_k\setminus\L^\b_{k-1}$.
Since $\rho<\bs^\b$, the last 
summation converges, whence, for sufficiently large $k$,
\begin{equation}\label{exp_cond_eq}
\int_{C_k^\b}\el\s_0\s_x\er^{\b} \, dx<e^{-8}.
\end{equation}
The result follows from the the Simon inequality,
Lemma~\ref{simon_lem}, with the 1-fat separating sets $C_k^\b$ 
using standard arguments (see
\cite[Corollary~9.38]{grimmett_rcm} for more details on the method).
A similar argument holds when $\b=\oo$.
\end{proof}

Let $\b\in(0,\oo]$, $\g=0$ and define the \emph{mass}
\begin{equation}
m^\b(\rho):=\liminf_{\|x\|\rightarrow\infty}
\left(-\frac{1}{\|x\|}\log\el\s_0\s_x\er^\b_\rho\right)
\end{equation}
By Theorem \ref{exp_decay_cor} and \eqref{o50},
\begin{equation}
m^\b(\rho) 
\begin{cases} >0 &\text{if  }\rho<\bs^\b,\\
=0 &\text{if } \rho>\bc^\b.
\end{cases}
\label{mel66}
\end{equation}

\begin{theorem}\label{0mass}
Except when $d=1$ and $\b<\oo$, $m^\b(\bs^\b)=0$ 
and $\chi^{\rf,\b}(\bs^\b,0)=\oo$. 
\end{theorem}

\begin{proof}
Let $d \ge 2$, $\g=0$, and fix $\b\in(0,\oo)$.
We use the Lieb inequality,
Lemma~\ref{lieb_lem}, and the argument of~\cite{lieb80,simon80},
see also \cite[Corollary 9.46]{grimmett_rcm}.
It is necessary and sufficient for $m^\b(\rho)>0$ that
\begin{equation}\label{exp_cond_eq_2}
\int_{C_n^\b}\el\s_0\s_x\er^{\rf,\b}_{n,\rho}\, dx<e^{-8}\quad\mbox{for some }n.
\end{equation}
Necessity holds because the integrand is no greater than $\el\s_0\s_x\er^\b$.
Sufficiency follows from Lemma \ref{lieb_lem}, as in the proof of Theorem \ref{exp_decay_cor}.

By \eqref{st_Ising_eq},
\begin{align*}
\frac{\partial}{\partial\rho}\el\s_0\s_x\er^{\rf,\b}_{n,\rho}
&=
\frac{1}{2} \int_{\L_n^\b}dy\,\sum_{z\sim y}\el\s_0\s_x;\s_y\s_z\er^{\rf,\b}_{n,\rho}\\
&\leq d\b(2n+1)^d.
\end{align*}
Therefore, if $\rho'>\rho$,
\begin{equation}
\int_{C_n^\b}\el\s_0\s_x\er^{\rf,\b}_{n,\rho'}\,dx\leq d[\b(2n+1)^d]^2(\rho'-\rho)
+\int_{C_n^\b}\el\s_0\s_x\er^{\rf,\b}_{n,\rho}\,dx.
\end{equation}
Hence, if \eqref{exp_cond_eq_2} holds for some $\rho$, then it holds for $\rho'$
when $\rho'-\rho>0$ is sufficiently
small.  

Suppose $m^\b(\bs^\b)>0$. Then $m^\b(\rho')>0$ for some $\rho'>\bs^\b$,
which contradicts $\chi^{\rf,\b}(\rho',0)=\oo$, and the first claim of the theorem follows.
A similar argument holds when $d=1$ and $\b=\oo$.
The second claim follows similarly: if $\chi^{\rf,\b}(\bs^\b,0)<\oo$, then
\eqref{exp_cond_eq_2} holds with $\rho=\bs^\b$, whence $m^\b(\rho')>0$
and $\chi^{\rf,\b}(\rho',0)<\oo$
for some $\rho'>\bs^\b$, a contradiction. (See also \cite{aizenman_tree-graph}.)
\end{proof}

We are now ready to state the main results.
We will adapt the arguments of~\cite[Lemmas~4.1, 5.1]{ab} 
(see also \cite{abf,grimmett_perc})
to prove the following. 
\begin{theorem}\label{ab_thm}
There are constants $c_1$, $c_2>0$ such that, for $\b\in(0,\oo]$, 
\begin{align}\label{ab_1_eq}
M^\b(\bs,\g)&\geq c_1\g^{1/3},\\
\label{ab_2_eq}
M^\b_+(\rho,0)&\geq c_2(\rho-\bs^\b)^{1/2},
\end{align}
for small, positive $\g$ and $\rho-\bs^\b$, \resp. 
\end{theorem}

This is vacuous when $d=1$ and $\b<\oo$; see \eqref{critvals}.
The exponents in the above inequalities are presumably sharp
in the corresponding mean-field model 
(see \cite{abf,af} and Remark \ref{remark_mf}).
It is standard that a number of important results follow from
Theorem~\ref{ab_thm}, of which we state the following here.
\begin{theorem}\label{eq_cor}
For $d \ge 1$ and $\b\in(0,\oo]$, we have that $\bc^\b=\bs^\b$.
\end{theorem}
\begin{proof}
Except when $d=1$ and $\b<\oo$, this is immediate 
from \eqref{mel64} and \eqref{ab_2_eq}.
In the remaining case, $\bc^\b=\bs^\b=\oo$.
\end{proof}

\begin{proof}[Proof of Theorem~\ref{ab_thm}]
We will describe the case when $\b<\oo$ is fixed;  the ground
state case is proved by a similar method.  The argument is based
on~\cite{ab}.

We start by proving~\eqref{ab_1_eq}.  If $M^\b_+(\rho_\rs,0)>0$ there
is nothing to prove, so we assume that $M^\b_+(\rho_\rs,0)=0$.
The inequalities of Theorems~\ref{three_ineq_lem} and~\ref{main_pdi_thm}
may be combined to obtain
\begin{equation}\label{i1}
M^\b_n\leq (M^\b_n)^3+\chi^\b_n\cdot
\left(\g+4d\l (M^\b_n)^3+4\d \frac{(M^\b_n)^3}{1-(M^\b_n)^2}\right).
\end{equation}
Set $\d=1$ and $\rho=\rho^\b_\rs$, and write
$f_n(\g)=2M^\b_n(\rho^\b_\rs,\g)$.  Recall that 
the sequence $f_n(\g)$ converges 
as $n\rightarrow\oo$ to some $f(\g)$
for all $\g\geq0$,  and that the derivatives $f_n'=2\chi^\b_n$ 
converge for $\g\in\cC$ to some $g(\g)$ which is decreasing in $\g$. 
Moreover, from the discussion around~\eqref{i2} and the assumption
that $M^\b_+(\rho_\rs,0)=0$ it follows that
\begin{equation}
\lim_{\substack{\g\downarrow 0\\\g\in\cC}} g(\g)=\oo.
\end{equation}

Multiplying through
by $1-(M^\b_n)^2$ and discarding non-positive terms on the right
hand side, we may deduce from~\eqref{i1} that the functions $f_n$
satisfy the inequality
\begin{equation}
f_n(\g)\leq\g\cdot f'_n(\g)+a\cdot f'_n(\g)f_n(\g)^3
+f_n(\g)^3,\qquad \g\geq0,
\end{equation}
where $a>0$ is an appropriate constant depending on $\l$ and $d$ only. 
For $\g>0$ we may rewrite this as
\begin{equation}
\frac{1}{f'_n(\g)}\frac{d}{d\g}\Big[\frac{\g}{f_n(\g)}\Big]\leq
f'_n(\g)\Big(a+\frac{1}{f'_n(\g)}\Big).
\end{equation}
Letting $\g>\eps>0$ and integrating from $\eps$ to $\g$ it follows
that
\begin{equation}
\frac{\g}{f_n(\g)}-\frac{\eps}{f_n(\eps)}\leq
\int_\eps^\g f'_n(x)f_n(x)\Big(a+\frac{1}{f'_n(x)}\Big)\,dx.
\end{equation}
Using~\eqref{dg_bound_eq} of Theorem~\ref{three_ineq_lem}, 
it follows on letting $\eps\downarrow0$ that
\begin{equation}\label{i3}
\frac{\g}{f_n(\g)}-\frac{1}{f'_n(0)}\leq
\int_0^\g f'_n(x)f_n(x)\Big(a+\frac{1}{f'_n(x)}\Big)\,dx.
\end{equation}

Now suppose that $\g>0$ lies in $\cC$.  If $\g$ is sufficiently
small then $g(\g)\geq1.1$, and for such a $\g$ fixed we have for
sufficiently large $n$ that $f'_n(\g)\geq1$.  Since $f'_n$ is
decreasing in $\g$ we may deduce from~\eqref{i3} that
\begin{equation}
\frac{\g}{f_n(\g)}-\frac{1}{f'_n(0)}\leq(a+1)
\int_0^\g f'_n(x)f_n(x)\,dx=\frac{a+1}{2}f_n(\g)^2
\end{equation}
Letting $n\rightarrow\oo$ it follows that
$$
\frac{\g}{f(\g)}\leq\frac{a+1}{2}f(\g)^2
$$
as required.

Let us now turn to~\eqref{ab_2_eq}.  Note first that if $\rho=\l/\d$ then
\begin{equation}
\frac{\partial M^\b_n}{\partial\l}=\frac{1}{\d}
\frac{\partial M^\b_n}{\partial\rho}\quad\text{and}\quad
\frac{\partial M^\b_n}{\partial\d}=-\frac{\l}{\d^2}
\frac{\partial M^\b_n}{\partial\rho},
\end{equation}
so that the inequality of Theorem~\ref{main_pdi_thm} may be
rewritten as
\begin{equation}
M^\b_n\leq\g\frac{\partial M^\b_n}{\partial\g}+(M^\b_n)^3
+2\rho (M^\b_n)^2\frac{\partial M^\b_n}{\partial\rho}.
\end{equation}
This may in turn be rewritten as
\begin{equation}
\frac{\partial}{\partial\g}(\log M^\b_n)+
\frac{1}{\g}\frac{\partial}{\partial\rho}
(\rho (M^\b_n)^2-\rho)\geq0.
\end{equation}
We wish to integrate this over the rectangle
$[\rho^\b_\rs,\rho']\times[\g_0,\g_1]$ for $\rho'>\rho^\b_\rs$
and $\g_1>\g_0>0$.  Since $M^\b_n$ is increasing in $\rho$ and
in $\g$ we deduce, after discarding a term
$-\rho_\rs^\b M^\b_n(\rho^\b_\rs,\g)^2$, that
\begin{equation}\label{i4}
(\rho'-\rho^\b_\rs)\log\Big(
\frac{M^\b_n(\rho',\g_1)}{M^\b_n(\rho^\b_\rs,\g_0)}\Big)+
(\rho' M^\b_n(\rho',\g_1)^2-\rho'+\rho^\b_\rs)
\log\frac{\g_1}{\g_0}\geq0.
\end{equation}
We may let $n\rightarrow\oo$ in~\eqref{i4}, to deduce that the
same inequality is valid with $M^\b_n$ replaced by $M^\b$.  
It follows from~\eqref{ab_1_eq} that
\begin{equation}
\liminf_{\g_0\downarrow 0}
\frac{\log\Big(
\frac{M^\b_n(\rho',\g_1)}{M^\b_n(\rho^\b_\rs,\g_0)}\Big)}{\log(\g_1/\g_0)}
\leq \frac{1}{3}.
\end{equation}
It follows that
\begin{equation}
\frac{1}{3}(\rho'-\rho^\b_\rs)+\rho' M^\b(\rho',\g_1)
-(\rho'-\rho^\b_\rs)\geq0,
\end{equation}
which on letting $\g_1\downarrow0$ gives the result.
\end{proof}

\begin{remark}\label{remark_mf}
Let $\b\in(0,\oo]$. Except when $d=1$ and $\b<\oo$,
one may conjecture the existence of exponents $a=a(d,\b)$, $b=b(d,\b)$ 
such that
\begin{alignat}{2}
M^\b_+(\rho)&=(\rho-\bc^\b)^{(1+\o(1))a }
\qquad&&\mbox{as }\rho\downarrow\bc^\b,\\
M^\b(\bc^\b,\g)&=\g^{(1+\o(1))/b}\qquad&&\mbox{as }\g\downarrow 0.
\end{alignat}
Theorem~\ref{ab_thm} would then imply that $a\leq \frac12$ and $b \geq 3$.
In~\cite[Theorem~3.2]{chayes_ioffe_curie-weiss} 
it is proved for the ground-state 
quantum Curie--Weiss, or mean-field, model that the corresponding
$a =\frac12$.  It may be conjectured that the values $a=\frac12$ and
$b=3$ are attained for the space--time 
Ising model on $\ZZ^d\times [-\frac12 \b,\frac12\b]$
for $d$ sufficiently large, as proved for the classical Ising model
in \cite{af}.  See also Section~\ref{af_sec}.
\end{remark}

Finally, a note about \eqref{o5}. The \rc\ measure corresponding to
the quantum Ising model is \emph{periodic} in
both $\ZZ^d$ and $\b$ directions, and this complicates 
the infinite-volume limit.
Since the periodic \rc\ measure dominates the 
free \rc\ measure, for $\b\in(0,\oo)$,
as in \eqref{mel61} and \eqref{o50},
\begin{alignat*}{2}
\liminf_{n\to\oo} \tau^\b_{L_n}(u,v) &\ge 
\el\s_{(u,0)}\s_{(v,0)}\er_{+,\rho'}^\b \qquad&&\text{for } \rho'<\rho\\
&\to M_+(\rho-)^2 \qquad &&\text{as } \rho'\uparrow \rho,
\end{alignat*}
and a similar argument holds in the ground state also.

\chapter{Applications and extensions}
\label{appl_ch}

\begin{quote}
{\it Summary.}  First we prove that the critical ratio for the 
ground state quantum
Ising model on $\ZZ$ is $\rho^\oo_\crit=2$;  we then extend this
result to more complicated (`star-like') graphs.
Next we discuss the possible applications of `reflection positivity'
to strengthen the results of Chapter~\ref{qim_ch} when $d\geq 3$,
and conclude with a discussion of versions of the random-parity 
representation of the Potts model.  
\end{quote}

\section{In one dimension}\label{sec_1d}

The quantum Ising model on $\ZZ$
has been thoroughly studied in the mathematical physics literature.
It is an example of what is called an `exactly solvable model':
using transfer matrices and related techniques, the critical ratio
and other important quantities have been computed, see for 
example~\cite{pfeuty70} or~\cite{sachdev99} and references therein.
In this section we prove by graphical methods that the critical
value coincides with the self-dual value of Section~\ref{duality_sec}.
The graphical method is valuable in that it extends to
more complicated geometries, as in the next section.
In the light of \eqref{critvals}, we shall study only
the ground state, and we shall suppress the superscript
$\oo$. 
\begin{theorem}\label{crit_val_cor}
Let $\LL=\ZZ$.  Then $\rho_\crit=2$, and the transition is
of second order in that $M_+(2)=0$.
\end{theorem}

We mention an application of this theorem. 
In an account \cite{GOS} of so-called `entanglement' in the 
quantum Ising model on the subset $[-m,m]$ of $\ZZ$,
it was shown that the reduced density matrix $\nu_{m}^L$
of the block $[-L,L]$ satisfies
$$
\|\nu_m^L - \nu_n^L\| \le \min\{2, C L^{\a}e^{-cm}\},\qquad 2\le m<n<\oo,
$$ 
where $C$ and $\a$ are constants depending on $\rho=\l/\d$,
and $c=c(\rho)>0$ whenever $\rho < 1$. Using Theorems \ref{exp_decay_cor}
and \ref{crit_val_cor}, we have that $c(\rho)>0$ for $\rho<2$. 

\begin{proof}
We adapt the well-known methods~\cite[Chapter~6]{grimmett_rcm}
for the discrete random-cluster model.
Write $\phi_\rho^\free$ and  $\phi_\rho^\wired$ for the free and
wired $q=2$ \rc\ measures, respectively.  By Theorem~\ref{potts_lim_thm}
and Remark~\ref{inf_vol_rk}, and the  
representation~\eqref{plus_ising} of the state $\el\cdot\er_+$,
we have that
\begin{equation}
\el\s_x\s_y\er_+=\fr^\wired(x\lra y),
\qquad\el\s_x\er_+=\fr^\wired(x\lra\infty).
\label{o26}
\end{equation}
Recall from Theorem~\ref{duality_thm} that the measures 
$\phi_\rho^\free$ and $\phi_{4/\rho}^\wired$ are mutually
dual.  By Zhang's argument, Theorem~\ref{zhang_thm}, we
know of the self-dual point $\rho=2$ that
\begin{equation}
\label{o52}
\phi_2^\free(0\lra \oo) = 0
\end{equation}
and hence that $\rho_\crit\geq 2$.

We show next that $\rho_\crit\le 2$, following the method developed for
percolation to be found in \cite{grimmett_perc,grimmett_rcm}. 
Suppose that $\bc>2$.  Consider the `lozenge' $D_n$ of side length $n$, as 
illustrated in Figure~\ref{zhang_fig}
on p.~\pageref{zhang_fig}, and its `dual' $D_n^\dual$.
Let $A_n$ denote the event that there is an open path from the bottom
left to the top right of $D_n$ in $\om$, and let $A_n^\dual$ be the 
`dual' event that there is in $\om_\dual$ an open path from the
top left to the bottom right of $D_n^\dual$.  The events $A_n$ and
$A_n^\dual$ are complementary, so we have by duality and symmetry
under reflection that
\begin{equation}
1=\phi_2^\free(A_n)+\phi_2^\free(A_n^\dual)=
\phi_2^\free(A_n)+\phi_2^\wired(A_n)\leq 2\phi_2^\wired(A_n).
\end{equation}
However, if $2<\rho_\crit$ then we have by~\eqref{o26} and 
Theorem~\ref{exp_decay_cor} that $\phi_2^\wired(A_n)$
decays to zero in the manner of $C n^2 e^{-\a n}$ as $n\to\oo$, 
a contradiction.

We show that $M_+(2)=0$ by adapting a simple argument 
developed by Werner in~\cite{WW08} for the classical Ising model on $\ZZ^2$.
Certain geometrical details are omitted.
Let $\pi^\free$ be the Ising state
obtained with free boundary condition, as
in Theorem~\ref{potts_lim_thm}.   Recall that $\pi^\free$
may be obtained from the random-cluster measures 
$\phi_2^\free$ by assigning to the clusters
spin $\pm 1$ independently at random, with probability $1/2$ each.
By Lemma~\ref{ergod_lem}, $\pi^\free$ is ergodic.

The binary relations $\lrao \pm$ are defined as follows.
A \emph{path} of $\ZZ\times\RR$ is a path of $\RR^2$ that:
traverses a finite number of line-segments of $\ZZ\times\RR$, 
and is permitted to connect them by passing between any two points of
the form $(u,t)$, $(u\pm 1,t)$. For $x,y\in\ZZ\times\RR$, we write 
$x\lrao + y$ (\resp, $x\lrao - y$) if there exists a 
path with endpoints $x$, $y$ all of
whose elements are labelled $+1$ (\resp, $-1$).
(In particular, for any $x$ we have that $x\lrao + x$
and $x\lrao - x$.)
Let $N^+$ (\resp, $N^-$) be the number of unbounded $+$ (\resp, $-$)
Ising clusters with connectivity relation $\lrao +$ (\resp, $\lrao -$).
By the Burton--Keane argument, as in Theorem~\ref{inf_clust_uniq},
one may show that either
$\pi^\free(N^+=1) = 1$ or $\pi^\free(N^+=0)=1$. The former would entail
also that $\pi^\free(N^-=1)=1$, by the $\pm$ symmetry in the coupling with
the random-cluster measure.  With an application of Zhang's argument
as in Theorem~\ref{zhang_thm}, however, one can show that this is 
impossible.  Therefore, 
\begin{equation}
\pi^\free(N^\pm = 0) = 1.
\label{o31} 
\end{equation}

Recall that $\el\cdot\er^+=\pi^\wired$.
There is a standard argument for deriving $\pi^\free=\el\cdot\er^+$ 
from \eqref{o31}, of which the idea is roughly as follows.
(See
\cite{ACCN} or \cite[Thm 5.33]{grimmett_rcm} for examples of similar arguments
applied to the \rc\ model.) 
Let $\L_n=[-n,n]^2 \subseteq \ZZ\times \RR$, and let $m<n$.  We call
a set $S\subseteq\L_n$ a \emph{separating set} if any path from $\L_m$
to $\partial\L_n$ contains an element of $S$.
We adopt the harmless
convention that, for any spin-configuration $\s$, 
the subset of $\L_n$ labelled $+1$ is closed,
compare Remark~\ref{rem-as}.
By \eqref{o31}, for given $m$, and for $\eps>0$ and large $n$,
the event $A_{m,n}= \{\L_m \lrao -\pd\L_n\}^\tc$ satisfies 
$\pi^\free(A_{m,n}) > 1-\eps$.
On $A_{m,n}$, there is a separating set labelled entirely $+$;  
let us call any
such separating set a $+$-separating set.  Let $U$ denote 
the set of all points in $\L_n$
which are $-$-connected to $\partial\L_n$ (note that
this includes $\partial\L_n$ itself).  Write $S=S(\s)$ for 
$\partial (\L_n\sm U)$.  Note that
$S\subseteq \L_n\sm\L_m$ is a $+$-separating set.  
See Figure~\ref{2nd_order_fig}.
\begin{figure}[hbt]
\centering
\includegraphics{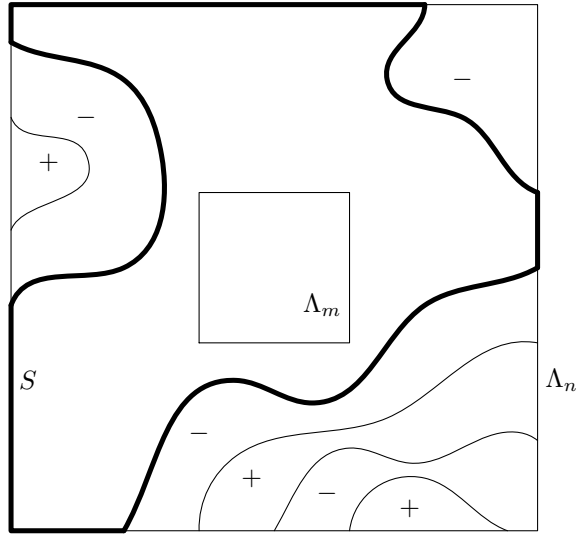}
\caption[First order phase transition when $d=1$]
{Sketch of an Ising configuration $\s$, 
with the set $S(\s)$ drawn bold; $S$ is
a $+$-separating set.}\label{2nd_order_fig} 
\end{figure}

For any closed separating set $S_1$, define $\hat S_1$ to be the union
of $S_1$ with the unbounded component of $(\ZZ\times\RR)\sm S_1$.  
Also let $\tilde S_1$ be the set of points in $\L_n$ that 
are separated from $\partial\L_n$ by $S_1$.  The event
$\{S(\s)=S_1\}$ is $\cG_{\hat S_1}$-measurable, i.e. it depends only
on the restriction of $\s$ to $\hat S_1$.  
Let $B\se\L_m$ be a finite set, and recall the notation
$\nu'_B$ at~\eqref{mrk1}.
By the \dlr-property of Lemma~\ref{potts_cond} 
(the natural extension of which holds also for
infinite-volume measures)  we deduce that
\begin{equation}
\pi^\free(\nu'_B\mid A_{m,n},S)=\pi^\wired_{\tilde S}(\nu'_B\mid A_{m,n}).
\end{equation}
Let $n\rightarrow\oo$ to deduce,
using also the \fkg-inequality
of Lemma~\ref{ising_fkg}, that
\begin{equation}
\pi^\free(\nu'_B\mid S)=\pi^\wired_{\tilde S}(\nu'_B)\geq\pi^\wired(\nu'_B).
\end{equation}
By integrating, and letting $m\rightarrow\oo$, 
we obtain that $\pi^\free(\nu'_B)\geq\pi^\wired(\nu'_B)$ 
for all finite sets $B\se\ZZ\times\RR$.
Since the reverse inequality 
$\pi^\free(\nu'_B)\le \pi^\wired(\nu'_B)$ always holds
(by  Lemma~\ref{potts_cond} and Lemma~\ref{ising_fkg} again), 
we deduce that $\pi^\free=\pi^\wired$ as claimed.

One way to conclude that $M_+(2)=0$ is to use
the \rc\ representation again.
By \eqref{o52} and the above, 
$$
\phi_2^\free(0\lra \oo) = \phi_2^\wired(0\lra\oo)=0,
$$
whence $M_+(2) = \phi_2^\wired(0\lra \oo) = 0$.
\end{proof}

\section{On star-like graphs}\label{starlike_sec}

We now extend Theorem~\ref{crit_val_cor} of the previous section, to
show that the critical ratio $\rho_\crit(2)=2$ for a larger class of
graphs than just $\ZZ$.  This section forms the contents of the 
article~\cite{bjo0}.

\begin{figure}[hbt]
\centering
\includegraphics{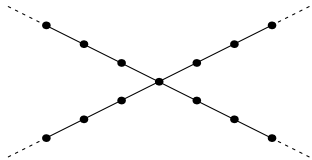}
\caption[The star graph]
{The star graph has a central vertex of degree $k\geq3$ and $k$
  infinite arms, on which each vertex has degree 2.  In this illustration,
  $k=4$.}\label{star_fig} 
\end{figure}

The class of graphs for which we prove that the critical
ratio is $2$ includes for example the 
\emph{star graph}, which is the junction 
of several copies of $\ZZ$ at a single point.  See Figure~\ref{star_fig}.  
It also includes many other planar graphs (see 
Definition~\ref{starlike_def}).
The result for the star is perhaps not unexpected, 
since the star is only `locally'
different from $\ZZ$:  if you go far enough out on one of the
`arms' then the star `looks like' $\ZZ$.  However, as pointed out
before, the quantum Ising model on the star, unlike on $\ZZ$,
is not exactly solvable, and graphical methods are the only known
way to prove this result.

The Ising model on the star-graph has recently arisen in
the study of boundary effects in the two-dimensional classical Ising model,
see for example~\cite{trombettoni,martino_etal}.  Similar geometries have
also arisen in different problems in quantum theory, such as transport
properties of quantum wire systems,
see~\cite{chamon03,hou08,lal02}.

\begin{figure}[hbt]
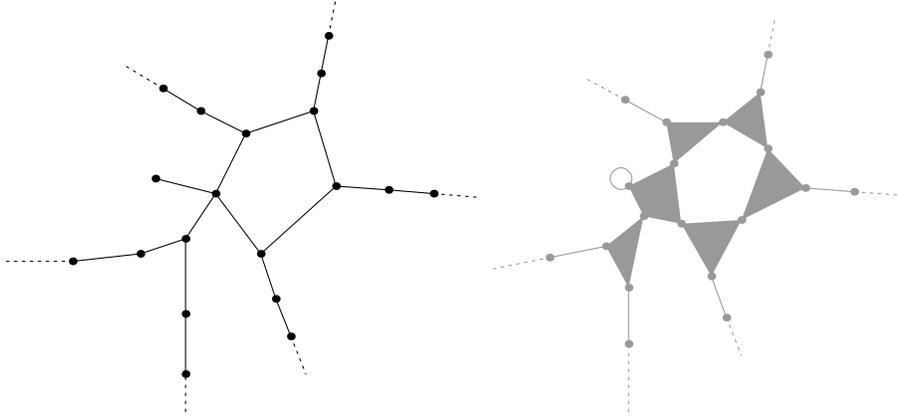

\centering
\includegraphics{thesis.24}
\includegraphics{thesis.25}
\caption[A star-like graph and its line-hypergraph]
{A star-like graph $\LL$ (left) and its 
  line-hypergraph $\HH$ (right).
  Any vertex of degree $\geq3$ in $\LL$ is associated with a ``polygonal''
  (hyper)edge in $\HH$.}\label{graph_fig}
\end{figure}

Throughout this section we consider the ground-state only, that
is to say we let $\b=\oo$;  reference to $\b$ will be suppressed.
We also let $\l,\d>0$ be constant and $\g=0$.
Let $\LL=(\VV,\EE)$ be a 
fixed \emph{star-like graph}:
\begin{definition}\label{starlike_def}
A star-like graph is a countably infinite connected planar graph, in which all
vertices have finite degree and only finitely many vertices have
degree larger than two.
\end{definition}
Such a graph is illustrated in Figure~\ref{graph_fig};  note that
the star graph of Figure~\ref{star_fig} 
is an example in which exactly one vertex
has degree at least three.

The following is the main result of this section. 
\begin{theorem}\label{starlike_thm}
Let $\LL$ be any star-like graph.  Then the critical ratio of the 
ground state quantum Ising model on $\LL$ is $\rho_\crit(2)=2$. 
\end{theorem}
Simpler arguments than those presented here can be used
to establish the analogous result when $q=1$, namely that $\rho_\crit(1)=1$.
Also, the same arguments can be used to calculate the critical probability 
of the discrete graphs $\LL\times\ZZ$ when $q=1,2$.  
As in the case $\LL=\ZZ$, an essential ingredient of the proof is the
exponential decay of correlations below~$\rho_\crit$.

Recall that a \emph{hypergraph} is a set $\WW$ 
together with a collection $\BB$ of 
subsets of $\WW$, called \emph{edges}
(or hyperedges).  A graph is a hypergraph in which all
edges contain two elements.  In our analysis we will use  a
suitably defined hypergraph `dual' of $\LL$.
To be precise, let $\HH=(\WW,\BB)$%
\nomenclature[H]{$\HH$}{Hypergraph}
 be the 
\emph{line-hypergraph} of $\LL$, given by letting $\WW=\EE$%
\nomenclature[W]{$\WW$}{Vertices of $\HH$}
 and
letting the set $\{e_1,\dotsc,e_n\}\subseteq \EE=\WW$ be an edge 
(that is, an element of $\BB$)%
\nomenclature[B]{$\BB$}{Edge set of $\HH$}
 if and only if $e_1,\dotsc,e_n$ are all the edges 
of $\LL$ adjacent to some particular vertex of $\LL$.  Note 
that only finitely many edges of $\HH$ have size larger than two,
since $\LL$ is star-like.

Fix an arbitrary 
planar embedding of $\LL$ into $\RR^2$;  we will typically
identify $\LL$ with its embedding.  
We let $\cO$ denote an arbitrary but fixed vertex of 
$\LL$ which has degree at least two;  we think of $\cO$ as the
`origin'.  There is a
natural planar embedding of $\HH$ defined via the embedding $\LL$, in which an
edge of size more than two is represented as a polygon.  See
Figure~\ref{graph_fig}.  In this section we will use the symbol
$\XX$%
\nomenclature[X]{$\XX$}{Product $\LL\times\RR$ for $\LL$ star-like}
 in place of $\LLam$ for $\LL\times\RR$, and will identify $\XX$
with the corresponding subset of $\RR^3$.  Similarly, we write
$\YY=\HH\times\RR$%
\nomenclature[Y]{$\YY$}{Dual of $\XX$}
 for the `dual' of $\XX$, also thought of as
a subset of $\RR^3$.  
We will often identify $\om=(B,D)\in\Om$ with its embedding, 
$\om\equiv(\XX\setminus D)\cup B$.  We
let $\L_n$ be the simple region corresponding to $\b=n$ and 
$L$ the subgraph of $\LL$ induced by the vertices at graph distance
at most $n$ from $\cO$, see~\eqref{def-oldL}.  
Note that $\L_n\uparrow\XX$.
In this section we let uppercase
$\Phi^b_n$%
\nomenclature[F]{$\Phi^b$}{Random-cluster measure on $\XX$}
 denote the random-cluster measure on $\L_n$
with parameters $\l,\d>0$, $\g=0$, $q=2$ and boundary
condition $b\in\{0,1\}$, where, as in Section~\ref{duality_sec},
we let $0$ and $1$ denote the free and wired boundary conditions,
respectively.

Given any configuration $\om\in\Om$, one may as in the case $\LL=\ZZ$
associate with it a \emph{dual} configuration on $\YY$ by placing a death
wherever 
$\om$ has a bridge, and a (hyper)bridge wherever $\om$ has a death. 
Recall Figure~\ref{duality_fig} on p.~\pageref{duality_fig}.
More precisely, we let $\Om_\dual$ be the set of pairs of locally finite
subsets of $\BB\times\RR$ and $\WW\times\RR$, 
and for each $\om=(B,D)\in\Om$ we
define its dual to be $\om_\dual:=(D,B)$.  As before, we may identify
$\om_\dual$ with its 
embedding in $\YY$, noting that some bridges may be embedded as polygons.  We
let $\Psi^b_n$%
\nomenclature[P]{$\Psi^b$}{Dual of $\Phi^{1-b}$}
 and $\Psi^b$ denote the laws of $\om_\dual$ under
$\Phi_n^{1-b}$ and $\Phi^{1-b}$ respectively.  

We will frequently be comparing the
random-cluster  measures on $\XX$ and $\YY$ with the random-cluster
measures on $\ZZ\times\RR$;  the latter may be regarded as a subset 
of both $\XX$ and $\YY$ (in a sense made more precise below).
We will reserve the lower-case
symbols $\phi^b_n,\phi^b$ for the random-cluster measures 
on $\ZZ\times\RR$ with the same parameters as $\Phi^b_n$
(where $\phi^b_n$ lives on the simple region given by
$\b=n$ and $L=[-n,n]$).
We will write $\psi_n^{1-b},\psi^{1-b}$ for the dual measures of 
$\phi^b_n,\phi^b$ on $\ZZ\times\RR$; thus 
by Theorem~\ref{duality_thm}, the measures $\psi_n^{1-b},\psi^{1-b}$
are random cluster measures with parameters
$q'=q$, $\l'=q\d$ and $\d'=\l/q$, and boundary condition $1-b$.

Here is a brief outline of the proof of Theorem~\ref{starlike_thm}.  
First we make
the straightforward observation that $\rho_\crit(2)\leq 2$.
Next, we use exponential decay to establish the
existence of certain infinite paths in the dual 
model on $\YY$ when $\l/\d<2$.  
Finally, we show how to put these paths together to
form  `blocking circuits' in $\YY$, which prevent the existence of infinite
paths in $\XX$ when $\l/\d<2$.  Parts of the argument are inspired
by~\cite{gkr}.  

\begin{lemma}\label{upperbound_lem}  
For $\LL$ any star-like graph, $\rho_\crit(2)\leq 2$.
\end{lemma}
\begin{proof}
Since $\LL$ is star-like, it
contains an isomorphic copy of $\ZZ$ as a subgraph.  Let 
$Z$ be such a subgraph;  we may assume that $\cO\in Z$.  We may
identify $\phi^b_n,\phi^b$ with the random-cluster measures on $Z\times\RR$.
For each $n\geq 1$, let $C_n$ be the event that no two points in
$\L_n\cap(Z\times\RR)$ are connected by a path which leaves $Z\times\RR$.
Each $C_n$ is a decreasing event.  It follows from  the
\dlr-property, Lemma~\ref{del_contr}, 
that $\Phi^b_n(\cdot\mid C_n)=\phi^b_n(\cdot)$.  If $A$ is an increasing 
local event defined on $Z\times\RR$, this means that
\begin{equation}
\phi^b_n(A)=\Phi^b_n(A\mid C_n)\leq\Phi^b_n(A),
\end{equation}
i.e. $\phi^b_n\leq\Phi^b_n$ for all $n$.  Letting $n\rightarrow\infty$ it
follows that $\phi^b\leq\Phi^b$.  If $\l/\d>2$ then 
$\phi^b((\cO,0)\leftrightarrow\infty)>0$ and it follows that also 
\begin{equation}
\Phi^b((\cO,0)\leftrightarrow\infty)>0,
\end{equation}
which is to say that $\rho_\crit(2)\leq2$.
\end{proof}

\subsection{Infinite paths in the half-plane}

Let us now establish some facts about the random-cluster model on
the `half-plane'
$\ZZ_+\times\RR$ which will be useful later.  
Our notation is as follows:  for $n\geq 1$, let
\begin{equation}
\begin{split}
S_n&=\{(a,t)\in\ZZ\times\RR: -n\leq a\leq n, |t|\leq n\}\\
S_n(m,s)&=S_n+(m,s)=\{(a+m,t+s)\in\ZZ\times\RR: (a,t)\in S_n\}.
\end{split}
\end{equation}%
\nomenclature[S]{$S_n$}{Region in $\ZZ\times\RR$}%
For brevity write $T_n=S_n(n,0)$.%
\nomenclature[T]{$T_n$}{$S_n(n,0)$}
 For $b\in\{0,1\}$ and $\D$ one of 
$S_n,T_n$, we let $\phi^b_\D$ denote the
$q=2$ random-cluster measure on the simple region in $\XX$
with $K=\D$ with boundary condition $b$ and
parameters $\l,\d$.  Note that
\begin{equation}
\phi^b=\lim_{n\rightarrow\infty}\phi^b_{S_n},\qquad
\psi^b=\lim_{n\rightarrow\infty}\psi^b_{S_n}.
\end{equation}
We will also be using the limits
\begin{equation}
\phi^\swired=\lim_{n\rightarrow\infty} \phi^1_{T_n},\qquad
\psi^\sfree=\lim_{n\rightarrow\infty} \psi^0_{T_n},
\end{equation}%
\nomenclature[s]{$\swired$}{`Side wired' boundary condition}%
\nomenclature[s]{$\sfree$}{`Side free' boundary condition}%
which exist by similar arguments to Theorem~\ref{conv_lem}.
(The notation `$\swired$' and `$\sfree$' is short for 
`side wired' and `side free', respectively.)
These are measures on configurations $\om$ on $\ZZ_+\times\RR$;  
standard arguments let us deduce all the properties 
of $\phi^\swired$ and $\psi^\sfree$ that we need.  In
particular $\psi^\sfree$ and $\phi^\swired$ 
are mutually dual (with the obvious
interpretation of duality) and they enjoy the positive
association property of Theorem~\ref{rc_fkg} and 
the finite energy property of Lemma~\ref{fin_en}.

Let $W$ be the `wedge'
\begin{equation}
W=\{(a,t)\in\ZZ_+\times\RR: 0\leq t\leq a/2+1\},
\end{equation}
and write $0$ for the origin $(0,0)$.
\begin{lemma}\label{wedge_lem}  Let $\l/\d<2$.  Then
\begin{equation}
\psi^\sfree(0\leftrightarrow\infty\mbox{ in } W)>0.
\end{equation}
\end{lemma}

Here is some intuition behind the proof of Lemma~\ref{wedge_lem}.  
The claim is well-known with $\psi^0$ 
in place of $\psi^\sfree$, by standard
arguments using duality and exponential decay.  However, $\psi^\sfree$ is
stochastically smaller than $\psi^0$, so we cannot deduce the result
immediately.  Instead we pass to the dual $\phi^\swired$ and 
establish directly a
lack of blocking paths.  The problem is the presence of the infinite `wired
side';  we get the required fast decay of two-point functions by using the
following result.

\begin{proposition}\label{exp_prop}
Let $\l/\d<2$.  There is $\a>0$ such that for all~$n$,
\begin{equation}
\phi^1_{S_n}(0\leftrightarrow \partial S_n)\leq e^{-\a n}.
\end{equation}
\end{proposition}

In words, correlations decay exponentially under finite volume
measures if they do so under infinite volume measures.
Results of this type for the classical Ising and random-cluster 
models appear in many places.  In~\cite{campanino_ioffe_velenik_08} 
and~\cite{cerf_messikh_06}
it is proved for general $q\geq1$ random-cluster models in two
dimensions, and more general results about the two-dimensional
case appear in~\cite{alexander04}.  
A proof of general results of this type
for the classical Ising model in any
dimension appears in~\cite{higuchi93_ii}.  Below we adapt the argument
in~\cite{higuchi93_ii} to the current setting, with the difference that
we shorten the proof by using the Lieb inequality, Lemma~\ref{lieb_lem}, 
in place
of the \ghs-inequality;  use of the Lieb-inequality was suggested by Grimmett
(personal communication).  Note that the same argument works 
on $\ZZ^d$ for any $d\geq1$.
\begin{proof}
Let $\hat S_n\supseteq S_n$ denote the `tall' box 
\begin{equation}
\hat S_n=\{(a,t)\in\ZZ\times\RR: -n\leq a\leq n, |t|\leq n+1\}.
\end{equation}
We will use a random-cluster measure on $\hat S_n$ 
which has non-constant $\l,\d$, and nonzero $\g$.
The particular intensities we use are these.  Fix $n$, and fix
$m\geq0$, which we think of as large.  Let
$\l(\cdot)$, $\d(\cdot)$ and $\g_m(\cdot)$ be given by 
\begin{equation}
\begin{split}
\d(a,t)&=\left\{
\begin{array}{ll}
\d, & \mbox{if } (a,t)\in S_n \\
0, & \mbox{otherwise},
\end{array}
\right.\\
\l(a+1/2,t)&=\left\{
\begin{array}{ll}
\l, & \mbox{if } (a,t)\in S_n \mbox{ and } (a+1,t)\in S_n\\
0, & \mbox{otherwise},
\end{array}
\right.\\
\g_m(a,t)&=\left\{
\begin{array}{ll}
\l, & \mbox{if exactly one of } 
      (a,t) \mbox{ and } (a+1,t) \mbox{ is in } S_n\\ 
m, & \mbox{if } (a,t)\in \hat S_n\setminus S_n\\
0, & \mbox{otherwise}.
\end{array}
\right.
\end{split}
\end{equation}
In words, the intensities are as usual `inside' $S_n$ and in particular
there is no external field in the interior;  on the left and right sides of
$S_n$, the external field simulates the wired boundary condition;  and on top
and bottom, the external field simulates an approximate wired boundary (as
$m\rightarrow\infty$). 
We let $\tilde \phi^b_{m,n}$
denote the random-cluster measure on $\hat S_n$ with intensities
$\l(\cdot),\d(\cdot),\g_m(\cdot)$ and boundary condition $b\in\{0,1\}$. 
Note that $\tilde\phi^0_{m,n}$ and $\phi^0_{S_n}$ 
agree on events defined on
$S_n$, for any~$m$.  

Let $X$ denote $\hat S_n\setminus S_n$ together with the left and right
sides of $S_n$.
By the Lieb inequality, Lemma~\ref{lieb_lem}, we have that
\begin{equation}
\begin{split}
\tilde\phi^1_{m,n}(0\leftrightarrow \G)&\leq
e^{8\d}\int_X dx\;\tilde\phi^0_{m,n}(0\leftrightarrow x)
\tilde\phi^1_{m,n}(x\leftrightarrow \G) \\
&\leq e^{8\d}\int_X dx\;\tilde\phi^0_{m,n}(0\leftrightarrow x),
\end{split}
\end{equation}
since (with these intensities)  $X$
separates $0$ from $\G$.  Therefore, by stochastic domination by the 
infinite-volume measure,
\begin{equation}\label{o76}
\tilde\phi^1_{m,n}(0\leftrightarrow \G)\leq
e^{8\d}\int_X dx\;\phi^0(0\leftrightarrow x).
\end{equation}
All the points $x\in X$ are at distance at least $n$ from the origin.  
By
exponential decay in the infinite volume, Theorem~\ref{exp_decay_cor}, 
it follows from~\eqref{o76} that
there is an absolute constant $\tilde\a>0$ such that  
\begin{equation}
\tilde\phi^1_{m,n}(0\leftrightarrow \G)\leq e^{8\d}|X|e^{-\tilde\a n}=
e^{8\d}(8n+2)e^{-\tilde\a n}.
\end{equation}
Now let $C$ be the event that all of $\hat S_n\setminus S_n$ belongs to
the connected component of $\G$, which is to say that all points on   
$\hat S_n\setminus S_n$ are linked to $\G$.  Then by the \dlr-property
of random-cluster measures the conditional measure
$\tilde\phi^1_{m,n}(\cdot\mid C)$ agrees 
with $\phi^1_{S_n}(\cdot)$ on events
defined on $S_n$.  Therefore
\begin{equation}
\begin{split}
\phi^1_{S_n}(0\leftrightarrow \partial S_n)&=
\tilde\phi^1_{m,n}(0\leftrightarrow \partial S_n\mid C)
=\tilde\phi^1_{m,n}(0\leftrightarrow \G\mid C)\\
&\leq\frac{\tilde\phi^1_{m,n}(0\leftrightarrow \G)}{\tilde\phi^1_{m,n}(C)}
\leq\frac{e^{8\d}}{\tilde\phi^1_{m,n}(C)}\cdot (8n+2)e^{-\tilde\a n}.
\end{split}
\end{equation}
Since $\tilde\phi^1_{m,n}(C)\rightarrow 1$ as 
$m\rightarrow\infty$ we conclude
that 
\begin{equation}
\phi^1_{S_n}(0\leftrightarrow \partial S_n)
\leq e^{8\d}(8n+2)e^{-\tilde\a n}.
\end{equation}
Since each $\phi^1_{S_n}(0\leftrightarrow \partial S_n)<1$ it is a simple
matter to tidy this up to get the result claimed.
\end{proof}

\begin{proof}[Proof of Lemma~\ref{wedge_lem}]
Let $T=\{(a,a/2+1): a\in\ZZ_+\}$ be the `top' of the wedge $W$.  
We claim that
\begin{equation}
\sum_{n\geq 1}\phi^\swired((n,0)\leftrightarrow T \mbox{ in } W)<\infty.
\end{equation}
Once this is proved, it follows from the Borel--Cantelli lemma that with
probability one under $\phi^\swired$, at most
finitely many of the points $(n,0)$ are connected to $T$ inside $W$.  Hence
under the dual measure $\psi^\sfree$ there is an infinite path inside $W$ with
probability one, and by the \dlr- and positive association properties it
follows that 
\begin{equation}
\psi^\sfree(0\leftrightarrow\infty\mbox{ in } W)>0,
\end{equation}
as required.

To prove the claim we note that, if $n$ is larger than some constant, then
the event `$(n,0)\leftrightarrow T \mbox{ in } W$' implies the event
`$(n,0)\leftrightarrow \partial S_{n/3}(n,0)$'.  The latter event, being
increasing, is more likely under the measure $\phi^1_{S_{n/3}(n,0)}$ than
under $\phi^\swired$.  But by Proposition~\ref{exp_prop},
\begin{equation}
\phi^1_{S_{n/3}(n,0)}((n,0)\leftrightarrow \partial S_{n/3}(n,0))
=\phi^1_{S_{n/3}}(0\leftrightarrow \partial S_{n/3})\leq e^{-\a n/3},
\end{equation}
which is clearly summable.
\end{proof}

The next lemma uses a variant of standard blocking arguments.
\begin{lemma}\label{circuits_lem}
Let $\l/\d<2$.  There exists $\eps>0$ such that for each~$n$,
\begin{equation}
\psi^\sfree((0,2n+1)\leftrightarrow (0,-2n-1)\mbox{ off } T_n)\geq \eps.
\end{equation}
\end{lemma}
\begin{proof}
Let $L_n=\{(a,n):a\geq0)\}$ be the horizontal line at height $n$, and let
$\eps>0$ be such that 
$\psi^\sfree(0\leftrightarrow\infty\mbox{ in } W)\geq\sqrt\eps$.
We claim that 
\begin{equation}
\psi^\sfree((0,-2n-1)\leftrightarrow L_{2n+1}\mbox{ off } T_n)\geq \sqrt\eps.
\end{equation}
Clearly $\psi^\sfree$ is invariant under reflection in the $x$-axis and
under vertical translation, see Lemma~\ref{invar_lem}.
Thus once the claim is proved we get that
\begin{equation}
\begin{split}
\psi^\sfree((0,2n+1)&\leftrightarrow (0,-2n-1)\mbox{ off } T_n)\\
&\geq\psi^\sfree((0,-2n-1)\leftrightarrow L_{2n+1}\mbox{ off } T_n\\
&\qquad\qquad\mbox{ and } (0,2n+1)\leftrightarrow L_{-2n-1}\mbox{ off } T_n)\\
&\geq (\sqrt\eps)^2,
\end{split}
\end{equation}
as required.  See Figure~\ref{paths_fig}.
\begin{figure}[hbt]
\centering
\includegraphics{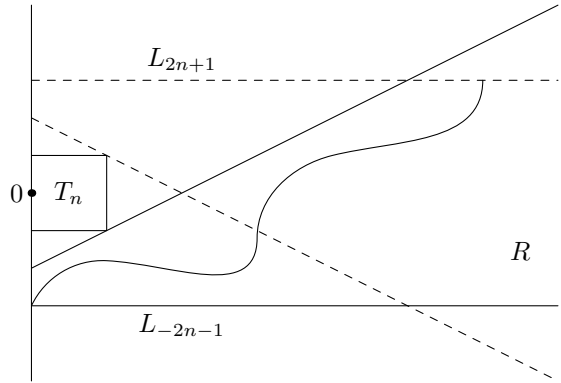}
\caption[Constructing circuits in the half-plane]
{Construction of a `half-circuit' in $\ZZ_+\times\RR$.  With
  probability one, any infinite path in the lower wedge must reach
  the line $L_{2n+1}$, and similarly for any infinite path in the upside-down
  wedge.  Any pair of such paths starting on the horizontal axis must
  cross.}\label{paths_fig}  
\end{figure}

The claim follows if we prove that
\begin{equation}\label{strip_eqn}
\psi^\sfree(0\leftrightarrow\infty\mbox{ in } R)=0,
\end{equation}
where $R$ is the strip
\begin{equation}
R=\{(a,t):a\geq 0,-2n-1\leq t\leq 2n+1\}.
\end{equation}
However,~\eqref{strip_eqn} follows from the \dlr-property,
Lemma~\ref{del_contr}, the
stochastic domination of Theorem~\ref{comp_perc},
and the Borel--Cantelli lemma;  these combine to show that the event
`no bridges between $\{k\}\times[-2n-1,2n+1]$ and 
$\{k+1\}\times[-2n-1,2n+1]$' must happen for infinitely many $k$ with
$\psi^\sfree$-probability one.  In more detail:  we have 
that $\psi^\sfree\leq\mu$, where $\mu$ is the percolation measure with
parameters $\l,\d$;  under $\mu$ the events above are independent, so
\begin{equation}
\psi^\sfree(0\leftrightarrow\infty\mbox{ in } R)\leq
\mu(0\leftrightarrow\infty\mbox{ in } R)=0.
\end{equation}
\end{proof}

\subsection{Proof of Theorem~\ref{starlike_thm}}

We may assume that $\LL\neq\ZZ$, since the case $\LL=\ZZ$ is known.
Let $\l/\d<2$, and recall that $\LL$ consists of finitely many infinite
`arms', where each vertex has degree two, together with a `central'
collection of other vertices.  On each of the arms, let us fix one
arbitrary vertex (of degree two) and call it an \emph{exit point}.  Let $U$
denote the set of exit points of~$\LL$.  

Given an exit point $u\in U$, call its
two neighbours $v$ and $w$;  we may assume that they are labelled so that only
$v$ is connected to the origin $\cO$ by a path not including $u$.  
If the edge $uv$ were
removed from $\LL$, the resulting graph would consist of two components, where
we denote by $J_u$ the component containing $w$.
Let $\hat \Phi^b_n,\hat \Phi^b$ denote the marginals of
$\Phi^b_n,\Phi^b$ on $X_u:=J_u\times\RR$;  similarly let
$\hat\Psi^b_n,\hat\Psi^b$ 
denote the marginals of the dual measures.  Of course $X_u$ is isomorphic to
the half-plane graph considered in the previous subsection.  
By positive association and the \dlr-property 
of random-cluster measures, $\hat\Phi_n^0\leq\phi^1_{T_n(u)}$, so letting
$n\rightarrow\infty$ also $\hat\Phi^0\leq\phi^\swired$.  
Passing to the dual, it
follows that $\hat\Psi^1\geq\psi^\sfree$.  The (primal) edge $uv$ is a
\emph{vertex} in the line-hypergraph;  denoting it still by $uv$ we therefore
have by Lemma~\ref{circuits_lem} that
there is an $\eps>0$ such that for all  $n$,
\begin{equation}
\Psi^1((uv,-2n-1)\leftrightarrow(uv,2n+1)
\mbox{ off }T_n(u)\mbox{ in }X_u)\geq\eps.
\end{equation}
Here $T_n(u)$ denotes the copy of the box $T_n$ contained in $X_u$.
Letting $A$ denote the intersection of the events above over all exit points
$u$, and letting $A_1=A_1(n)$ be the dual event
$A_1=\{\om_\dual:\om\in A\}$,  it follows from positive association that
$\Phi^0(A_1)\geq\eps^k$, where $k=|U|$  is the number of exit points.  Note
that $A_1$ is a decreasing event in the primal model.  

\begin{figure}[htb]
\centering
\includegraphics{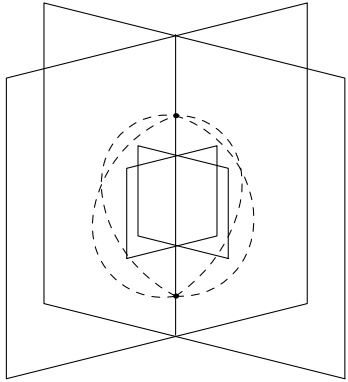}
\caption[Blocking in star-like graphs]
{The dashed lines indicate dual paths that block any primal connection
  from the interior to $\infty$.  Note that this figure illustrates only the
  simplest case when $\LL$ is a junction of lines at a single
  point.}\label{blocking_fig} 
\end{figure}

On $A_1$, no point in $T_n(u)$ can reach $\infty$ without passing the line 
$\{u\}\times[-2n-1,2n+1]$, since there is a dual blocking path in $X_u$.
Let $I$ denote the (finite) subgraph of $\LL$ spanned by the complement of
all the $J_u$ for $u\in U$, and let $A_2=A_2(n)$ denote the event that for all
vertices  $v\in I$, the intervals $\{v\}\times[2n+1,2n+2]$ and 
$\{v\}\times[-2n-1,-2n-2]$
all contain at least one death and the endpoints of no bridges (in the primal
model).  There is $\eta>0$ independent of $n$ such
that $\Phi^0(A_2)\geq\eta$. 
So by positive association $\Phi^0(A_1\cap A_2)\geq\eta\eps^k>0$.  
We have that $A_1\cap A_2\se A_3$, where $A_3$ is the event
that no point inside the union of 
$I\times[-n,n]$ with $\cup_{u\in U}T_n(u)$ lies in an unbounded 
connected component.  
See Figure~\ref{blocking_fig}. 
Taking the intersection of the $A_3=A_3(n)$ over all $n$, it follows
that  
\begin{equation}
\Phi^0(\mbox{there is no unbounded connected component})\geq\eta\eps^k.
\end{equation}
The event that there is no unbounded connected component is a tail event.
By tail-triviality, Proposition~\ref{tail_triv}, 
it follows that whenever $\l/\d<2$ then
\begin{equation}
\Phi^0(0\not\leftrightarrow\infty)=1.
\end{equation}
In other words, $\rho_\crit(2)\geq 2$.  Combined with the opposite bound in
Lemma~\ref{upperbound_lem}, this gives the result.
\qed

One may ask if, as in the case $\LL=\ZZ^d$, 
the phase transition on star-like
graphs is of second order, and if there is exponential decay of 
correlations below the critical point.  We do not know how to prove such
results:  Zhang's argument (Theorem~\ref{zhang_thm}) fails on
star-like graphs, and so do the arguments for
Theorem~\ref{main_pdi_thm}, due to the lack of symmetry.

\section{Reflection positivity}\label{af_sec}

The theory of reflection positivity was first developed 
in~\cite{frohlich_lieb78,froehlich_etal78,froehlich_etal80},
originally as a way to prove the existence of discontinuous phase
transitions in a wide range of models in statistical physics.
A model which is reflection positive
(see definitions below) will satisfy what are called 
`Gaussian domination bounds' and `chessboard estimates'.  The latter
will not be touched upon here, see the review~\cite{biskup_rp} and
references therein.  One may think of the Gaussian domination
bounds, and the related `infrared bound', as a way of bounding
certain quantities in the model by corresponding quantities in
another, simpler, model, namely what is called the `Gaussian free
field'.  Very roughly, existence of a phase transition in the Gaussian 
free field therefore implies existence of a phase transition in your
reflection positive model.

In~\cite{af}, it was shown that Gaussian domination bounds could also
be used in another way for the Ising model.  By relating the bounds
to quantities that appear naturally in the random-current representation
of the Ising model, Aizenman and Fern\'andez were able to establish
that the behaviour of the classical Ising model on $\ZZ^d$ resembles that of
the `mean field' Ising model when $d$ is large, in fact already when
$d\geq4$.  In this section we will state more precisely the sense
in which `large $d$ resembles mean field', and give a very brief sketch 
of the arguments involved.  We will also indicate how one might
extend the results of~\cite{af} to the quantum Ising model;  this
is currently work in progress.

In this section we will only be considering the case when
$\LL=\ZZ^d$ for some $d\geq1$, and $L=(V,E)=[-n,n]^d$ for some $n$,
with periodic boundary (see Assumption~\ref{periodic_assump}).  
For $j=1,\dotsc,d$, write $e_j$ for the element of $V$ whose $j$th
coordinate is $1$ and whose other coordinates are zero.
For $\s\in\{-1,+1\}^V$,
we write its classical Ising weight in this section as
\begin{equation}\label{q1}
\exp\Big(\b\sum_{xy\in E}J_{xy}\s_x\s_y+\g\sum_{x\in V}\s_x\Big),
\end{equation}
where $\b,\g,J_e\geq0$.  We assume that the model is translation invariant
in that $J_{xy}\equiv J_{y-x}$, where, for $z\in V$, 
$J_z\geq0$ and $J_z=0$ unless $z=e_j$ for some $j$.
We also assume that $J_{e_j}=J_{-e_j}$ for all $j=1,\dotsc,d$.

The classical Ising model displays a phase-transition in $\b$ when
$\g=0$, at the critical value $\b_\crit$.  As in the quantum Ising
model (Theorem~\ref{ab_thm}), the infinite-volume magnetization 
$M=M(\b,\g)$ satisfies the inequalities
\begin{align}\label{q0}
M&\geq c_2(\b-\b_\crit)^{1/2}, 
&\text{for $\g=0$ and $\b\downarrow\b_\crit$}, \\
M&\geq c_1\g^{1/3}, &\text{for $\b=\b_\crit$ and $\g\downarrow 0$},
\nonumber
\end{align}
for some constants $c_1,c_2$ (this was first proved in~\cite{abf}).
As mentioned in Remark~\ref{remark_mf}, it is conjectured that the
limits
\begin{equation}\label{q-1}
a=\lim_{\b\downarrow\b_\crit}\frac{\log M(\b,0)}{\log(\b-\b_\crit)},\qquad
\frac{1}{b}=\lim_{\g\downarrow 0}\frac{\log M(\b_\crit,\g)}{\log\g}
\end{equation}
exist.  Using the random-current representation coupled with 
results from reflection positivity, \cite{af} shows that these
limits do indeed exist when $d\geq 4$, and that~\eqref{q0}
is sharp in that $a=1/2$ and $b=3$.  The values $a=1/2$
and $b=3$ are called the `mean field' values because they are
known to be the correct critical exponents for the Ising 
model on the complete graph (this result is `well-known', but
see~\cite{fisher64,fisher67} for reviews).
Intuitively, complete graphs are infinite-dimensional, so the
higher $d$ is the closer one may expect the behaviour to be to
that on the complete graph.  The results of~\cite{af} confirm
this, and show that the `critical dimension' is at
most $d=4$.  Their method is roughly as follows.

For $j=1,\dotsc,d$ we let $P_i=\{x=(x_1,\dotsc,x_d)\in V: x_j=0\}$,
and we let $P_j^+=\{x\in V:x_j>0\}$ and 
$P_j^-=\{x\in V:x_j<0\}$.  The symbol $\theta_i$ will denote 
reflection in $P_i$, thus
$\theta_j(x_1,\dotsc,x_j,\dotsc,x_d)=(x_1,\dotsc,-x_j,\dotsc,x_d)$.
Write $\cF_{P_j^+}$ and $\cF_{P^-_j}$ for the $\s$-algebras
of events defined on $P_j^+$ and $P_j^-$, respectively.  

Although
we will be using the concept of reflection positivity only
for the Ising measure~\eqref{q1}, the definition makes sense in greater
generality, as follows.  
Let $S\subseteq\RR$ be a compact set, and endow
$S^V$ with the product $\s$-algebra.  Fix $j\in\{1,\dotsc,d\}$, and let
$\psi$ denote a probability measure on $S^V$ which is invariant
under $\theta_j$.  For $s=(s_x:x\in V)\in S^V$, write 
$\theta_j(s)=(s_{\theta_j(x)}:x\in V)$, and for $f:S^V\rightarrow\RR$
define $\theta_j f(s)=f(\theta_j (s))$.
\begin{definition}
The probability measure $\psi$ is \emph{reflection positive}
with respect to $\theta_j$ if for all $\cF_{P_j^+}$-measurable
$f:S^V\rightarrow\RR$, we have that
\[
\psi(f\cdot\theta_j f)\geq 0.
\]
\end{definition}

\begin{lemma}\label{rp_lem}$\;$
\begin{itemize}
\item Any product measure on $S^V$ invariant under $\theta_j$
is reflection positive with respect to $\theta_j$,
\item The Ising measure~\eqref{q1} is reflection positive
with respect to all the $\theta_j$.
\end{itemize}
\end{lemma}
For a proof of this standard fact, see for example~\cite{biskup_rp}.
It follows from Lemma~\ref{rp_lem} that the Ising model satisfies the
following `Gaussian domination' bounds.  For $p\in[-\pi,\pi]^d$, let
\begin{equation}
G(p):=\sum_{x\in V}\el\s_o\s_x\er_{\g=0}e^{ip\cdot x}
\end{equation}
be the Fourier transform of $\el\s_0\s_x\er_{\g=0}$, where
$i=\sqrt{-1}$ and $p \cdot x$ denotes the  usual dot product.
Due to our symmetry assumptions we see that the complex
conjugate $\ol{G(p)}=G(-p)=G(p)$ so that $G(p)\in \RR$.
Also define
\begin{equation}
E(p):=\frac{1}{2}\sum_{x\in V}(1-e^{ip\cdot x})J_x;
\end{equation}
similarly we see that $E(p)\in\RR$.
\begin{proposition}[Gaussian domination]\label{gd_prop}
\[
G(p)\leq \frac{1}{2\b E(p)}.
\]
\end{proposition}
Before we describe how this relates to the random-current
representation, we note that a simple calculation shows that
$E(p)\geq c\sum_{j=1}^dp_j^2$, which at least
gives some indication of why Gaussian domination may be
particularly useful for large~$d$.

The link to the random-current representation is roughly as
follows.  Define the \emph{bubble diagram}
\begin{equation}
B_0=\sum_{x\in V}\el\s_0\s_x\er_{\g=0}^2.
\end{equation}
Recall that $M=\el\s_0\er$ and that we write
$\chi=\partial M/\partial\g$.  We saw in Section~\ref{sw_appl_sec}
that random-current arguments imply the \ghs-inequality,
namely that $\partial \chi/\partial\g\leq 0$.  
In~\cite{af}, elaborations of such arguments (for the discrete
model) show that in fact 
\begin{equation}\label{q2}
\frac{\partial\chi}{\partial\g}\leq
-\frac{|1-\tanh(\g)B_0/M|_+^2}{96B_0(1+2\b B_0)^2}
\tanh(\g)\chi^4,
\end{equation}
where $|x|_+=x\vee 0$.  The bubble diagram appears here as it becomes
necessary to consider the existence of two independent
currents between sites $0$ and $x$.  Inequality~\eqref{q2} is an improvement
on the \ghs-inequality if $B_0$ is finite;  thus the first
task is to obtain bounds on $B_0$.  Such bounds are provided
primarily by Gaussian domination.
The link is provided via Parseval's identity:
\begin{equation}
B_0=\frac{1}{(2\pi)^d}\int_{[-\pi,\pi]^d} G(p)^2\,dp.
\end{equation}
By careful use of Gaussian domination and other bounds, one may
establish bounds on $B_0$ for $\b$ close to the critical value
$\b_\crit$.  More precisely, one may show that there are constants
$0<c_1,c_2<\oo$ such that
\begin{align*}
B_0&\leq c_1,  &\text{if } d>4,\\
B_0&\leq c_2 |\log(\b_\crit-\b)|,  &\text{if } d=4,
\end{align*}
as $\b\uparrow\b_c$.  Careful manipulation and integration
of~\eqref{q2} then gives that there are constants
$c_1',c_2',c_1'',c_2''$ such that the infinite-volume
magnetization $M$ satisfies the following.
First, as $\b\downarrow\b_\crit$ for $\g=0$,
\begin{align*}
M&\leq c_1'(\b-\b_\crit)^{1/2},  &\text{if } d>4,\\
M&\leq c_2'(\b-\b_\crit)^{1/2} |\log(\b-\b_\crit)|^{3/2},  &\text{if } d=4,
\end{align*}
and second, for $\b=\b_\crit$ and $\g\downarrow0$,
\begin{align*}
M&\leq c_1''\g^{1/3},  &\text{if } d>4,\\
M&\leq c_2'' \g^{1/3}|\log \g|,  &\text{if } d=4.
\end{align*}
These are the complementary bounds to~\eqref{q0} needed to
show that the limits~\eqref{q-1} exist and take the values
$a=1/2$ and $b=3$.

There are two main steps to extending the results of~\cite{af}
to the quantum (or space--time) Ising model:  first, to establish
reflection positivity and the related Gaussian domination bound, and
second, to verify that the random-parity representation can produce
an inequality of the form~\eqref{q2}.  There is essentially only
one known way of showing that a measure is reflection positive, which
is to show that it has a density against a product measure
which is of a prescribed form~\cite[Lemma~4.4]{biskup_rp}.  
Preliminary calculations suggest that
this method works also for the space--time Ising model.  Although
the random-current manipulations in~\cite{af} leading up 
to~\eqref{q2} are considerably more delicate than those presented
in Chapter~\ref{qim_ch} of this work and involve some new ideas
such as `dilution', preliminary calculations again
suggest that it should be possible to extend them as required.

\section{Random currents in the Potts model}

The main results of this work have relied on the
random-parity representation for the space--time Ising model.
It is natural to ask if there is a similar representation for 
the $q\geq 3$ Potts model.  Here we will discuss this question,
to start with in the context of the \emph{classical} (discrete)
Potts model on a finite graph  $L=(V,E)$.  
For simplicity we will assume free boundary condition and zero
external field;  it is easy to adapt the results here to positive
fields.  

It is shown in~\cite[Chapter~9]{grimmett_rcm}
(see also~\cite{essam_tsallis,magalhaes_essam})
that the $q$-state Potts model with $q\geq3$ possesses a
\emph{flow representation}, which is 
akin to the random-current representation, in that the
two-point correlation function may be written as the ratio of two
expected values.  This representation is as follows.  

Let the integer $q\geq 2$ be fixed.
For $\ul n=(n_e:e\in E)$ a vector of non-negative integers,
define the graph $L_{\ul n}=(V,E_{\ul n})$ 
by replacing each edge $e$ of $L$ by $n_e$ parallel edges.  
If $P=(P_e:e\in E)$ is a collection of finite sets with
$|P_e|=n_e$, we identify $L_P$ with $L_{\ul n}$, and interpret
$P_e$ as the set of edges replacing $e$.
We assign to the elements of $E_{\ul n}$ arbitrary directions
and write $\vec{e}$ for directed elements of $E_{\ul n}$;  if
$\vec{e}$ is adjacent to a vertex $x\in V$ and is directed
into $x$ we write $\vec{e}\mapsto x$, and if $\vec{e}$ is 
directed out of $x$ we write $\vec{e}\mapsfrom x$.
We say that a function $f:E_{\ul n}\rightarrow\{1,\dotsc,q-1\}$
is a (nonzero)
\emph{mod $q$ flow on $L_{\ul n}$} 
(or $q$-flow for short) if for all $x\in V$ we
have that
\begin{equation}
\sum_{\substack{\vec{e}\in E_{\ul n}: \\e\mapsfrom x}} f(\vec{e})-
\sum_{\substack{\vec{e}\in E_{\ul n}: \\e\mapsto x}} f(\vec{e})
\equiv 0 \quad\text{(mod $q$)}.
\end{equation}
Let $C(L_{\ul n};q)$ denote the number of mod $q$ flows on 
$L_{\ul n}$ (this is called the flow polynomial of $L_{\ul n}$).
It is easy to see that this number does not depend on the directions
chosen on the edges  (if the direction of an edge $\vec{e}$
is reversed we can replace $f(\vec{e})$ by $q-f(\vec{e})$).

For each $e\in E$, let $\b'_e\geq0$, and recall that the Potts
weight of an element $\nu\in\{1,\dotsc,q\}^V=\cN$ is 
\begin{equation}\label{p1}
\exp\Big(\sum_{e=xy\in E}\b'_e\d_{\nu_x,\nu_y}\Big),
\end{equation}
so that the partition function is
\begin{equation}\label{p2}
Z=\sum_{\nu\in\cN} \exp\Big(\sum_{e=xy\in E}\b'_e\d_{\nu_x,\nu_y}\Big).
\end{equation}
Let $\b_e=\b_e'/q$ and let the collection
$P=(P_e:e\in E)$ of finite sets be given by letting the $|P_e|$ be 
independent Poisson random variables, each with
parameter $\b_e$.  Write $\PP_\b$ for the probability measure
governing the $P_e$ and $\EE_\b$ for the corresponding expectation
operator.  

The flow representation of $Z$ is
\begin{equation}\label{fl_z}
Z=\exp\Big(2\sum_{e\in E}\b_e\Big)q^{|V|}\EE_\b[C(L_P;q)].
\end{equation}
In fact, more is true.  For $x,y\in V$,
let $L_{\ul n}^{xy}=(V,E_{\ul n}\cup\{xy\})$ 
denote the graph $L_{\ul n}$ with an edge added from $x$ to $y$.
Write $\el\cdot\er$ for the expected value under the
$q$-state Potts measure defined by~\eqref{p1}--\eqref{p2}.  Then 
for any $x,y\in V$ we have that
\begin{equation}\label{fl_2p}
q\el\one\{\nu_x=\nu_y\}\er-1=
\frac{\EE_\b[C(L_P^{xy};q)]}{\EE_\b[C(L_P;q)]}.
\end{equation}

Here is a simple observation that changes the expected value
in~\eqref{fl_z} into a probability.  For $\ul n\in\ZZ_+^E$,
let $F_{q}(\ul n)$ denote the set
of functions $f:V\rightarrow\{1,\dotsc,q-1\}$.  Then
\begin{equation}\label{fl_prob}
\begin{split}
\EE_\b[C(L_P;q)]&=\sum_{\ul n\in\ZZ_+^E}\prod_{e\in E}
\frac{\b_e^{n_e}}{n_e!}e^{-\b_e}\sum_{f\in F_{q}(\ul n)}
\one\{f\text{ is $q$-flow}\}\\
&=\exp\Big((q-2)\sum_{e\in E}\b_e\Big)
\sum_{\ul n\in\ZZ_+^E}\prod_{e\in E}
\frac{((q-1)\b_e)^{n_e}}{n_e!}e^{-(q-1)\b_e}\\
&\qquad\cdot\frac{1}{(q-1)^{\sum_{e\in E}n_e}}
\sum_{f\in F_{q}(\ul n)}\one\{f\text{ is $q$-flow}\}\\
&=\exp\Big((q-2)\sum_{e\in E}\b_e\Big)
\PP(\psi\text{ is $q$-flow on $L_{P'}$}),
\end{split}
\end{equation}
where, under $\PP$, the collection $P'=(P'_e:e\in E)$ is given by letting
the $|P'_e|$ be independent Poisson random variables with parameters
$(q-1)\b_e$ respectively, and $\psi$ is, given $P'$, a uniformly
chosen element of $F_q(P')$.  (As before, arbitrary directions are
assigned to the elements of $E_{P'}$, but the probability that $\psi$
is a $q$-flow does not depend on the choice of directions.)

We now show that a similar representation to~\eqref{fl_prob}
holds for the
two-point correlation functions~\eqref{fl_2p}, and indeed for
more general correlation functions.
As in  Section~\ref{gks_sec} we will use the variables
\[
\s_x=\exp\Big(\frac{2\pi i\nu_x}{q}\Big), \qquad\nu_x=1,\dotsc,q.
\]
We write $Q\subseteq\CC$ for the set of $q$th roots of unity, and 
$\S=Q^V$.  For $\ul r\in\ZZ^V$ and $\s\in\S$ we let 
\[
\s^{\ul r}=\prod_{x\in V}\s_x^{r_x}.
\]
Note that it is equivalent  to regard $r_x$ as an
element of $\ZZ/(q\ZZ)$, the integers modulo $q$.
Let $\PP$, $P'$ and $\psi$ be as in~\eqref{fl_prob},
and write $\{\psi\equiv 0\}$ for the event that
$\psi$ is a $q$-flow.  More generally, write
$\{\psi+\ul r\equiv 0\}$ for the event that for each
$x\in V$, 
\begin{equation}
\sum_{\substack{\vec{e}\in E_{P'}: \\e\mapsfrom x}} 
\psi(\vec{e})-
\sum_{\substack{\vec{e}\in E_{P'}: \\e\mapsto x}} 
\psi(\vec{e})
\equiv -r_x \quad\text{(mod $q$)}.
\end{equation}
(Recall that we have assigned arbitrary directions to the elements
of $E_{P'}$.)  
\begin{theorem}\label{potts_flow_thm}
In the discrete Potts model with zero field and coupling 
constants $\b'_e$,
\[
\el\s^{\ul r}\er=
\frac{\PP(\psi+\ul r\equiv 0)}{\PP(\psi\equiv 0)}.
\]
\end{theorem}
Before proving this, note that if
$\s\in\cN$ and $x,y\in V$, then $\tau_{xy}:=\s_x\s_y^{-1}$
has the property that $\tau_{xy}=1$
if and only if $\s_x=\s_y$, and in fact
\[
\frac{1}{q}\sum_{r=0}^{q-1}\tau_{xy}^r=\d_{\s_x,\s_y}.
\]
Thus the partition function~\eqref{p2} may be written
\begin{equation}\label{p3}
\begin{split}
Z&=\sum_{\s\in\S}\exp\Big(\sum_{e=xy\in E}\b'_e\d_{\s_x,\s_y}\Big)\\
&=\sum_{\s\in\S}\exp\Big(\frac{1}{2}\sum_{x,y\in V}
\b_{xy}\sum_{r=1}^{q-1}\tau_{xy}^r\Big)
\cdot\exp\Big(\sum_{e\in E}\b_e\Big),
\end{split}
\end{equation}
where the first sum inside the exponential is over all
ordered pairs $x,y\in V$, and we set $\b_{xy}=\b_e$ if $e\in E$
is an edge between $x$ and $y$, and $\b_{xy}=0$ otherwise.
Note finally that $\tau_{xy}\neq\tau_{yx}$ in general.

\begin{proof}
We perform a calculation on the factor
\[
\sum_{\s\in\S}\exp\Big(\frac{1}{2}\sum_{x,y\in V}
\b_{xy}\sum_{r=1}^{q-1}\tau_{xy}^r\Big)
\]
which appears on the right-hand-side of~\eqref{p3};  this 
will only re-prove the relation~\eqref{fl_prob}, but
it will be clear that a simple extension of the calculation 
will give the  result.

Let us write $\tilde\b_{xy}=\b_{xy}/2$.  We have that
\begin{equation}\label{rcr1}
\begin{split}
\sum_{\s\in\S}\exp\Big(\frac{1}{2}\sum_{x,y\in V}
\b_{xy}\sum_{r=1}^{q-1}\tau_{xy}^r\Big)
&=\sum_{\s\in\S}\prod_{x,y\in V}\prod_{r=1}^{q-1}
\sum_{m\geq 0}\frac{1}{m!}(\tilde\b_{xy}\tau_{xy}^r)^m\\
&=\sum_{\s\in\S}\sum_{\ul m} w(\ul m)
\prod_{x,y\in V}\prod_{r=1}^{q-1}(\tau_{xy}^r)^{m_{x,y,r}},
\end{split}
\end{equation}
where the vector $\ul m=(m_{x,y,r}:x,y\in V,r=1,\dotsc,q-1)$ 
consists of non-negative integers and
\[
w(\ul m)=\prod_{x,y\in V}\prod_{r=1}^{q-1}
\frac{\tilde\b_{xy}^{m_{x,y,r}}}{m_{x,y,r}!}
\]
is an un-normalized Poisson weight on $\ul m$.
Reordering~\eqref{rcr1} we obtain
\begin{equation}\label{rcr2}
\sum_{\s\in\S}\exp\Big(\frac{1}{2}\sum_{x,y\in V}
\b_{xy}\sum_{r=1}^{q-1}\tau_{xy}^r\Big)=
\sum_{\ul m}w(\ul m)\sum_{\s\in\S}\prod_{x,y\in V}
\tau_{xy}^{M_{xy}}
\end{equation}
where 
\[
M_{xy}=\sum_{r=1}^{q-1} r\cdot m_{x,y,r}.
\]

We may interpret $m_{x,y,r}$ as a random number of edges, each
of which is directed
from $x$ to $y$ and receives flow value $r$.  
Then $M_{xy}$ is the total flow from $x$ to $y$.  Up to the constant 
multiple $\exp\big((q-1)\sum_e\b_e\big)$, 
the quantity~\eqref{rcr2} equals the expected value of the quantity
\begin{equation}\label{rcr3}
\sum_{\s\in\S}\prod_{x,y\in V}
\tau_{xy}^{M_{xy}}
\end{equation}
when the $m_{x,y,r}$ have the Poisson distribution with parameter
$\tilde\b_{xy}$ and are chosen independently.

The quantity~\eqref{rcr3} simplifies, as follows.  
Let $a\in V$ be fixed, and let 
$L_a=(V_a,E_a)$ denote $L$ with $a$ removed.  Then
\begin{equation}
\begin{split}
\sum_{\s\in\S}\prod_{x,y\in V}\tau_{xy}^{M_{xy}}&=
\sum_{\s\in\S}\Big(\prod_{b\sim a}\tau_{ab}^{M_{ab}}
\tau_{ba}^{M_{ba}}\Big)
\prod_{x,y\in V_a}\tau_{xy}^{M_{xy}}\\
&=\sum_{\s\in\S}\Big(\prod_{b\sim a}\s_a^{M_{ab}-M_{ba}}
\s_b^{M_{ba}-M_{ab}}\Big)
\prod_{x,y\in V_a}\tau_{xy}^{M_{xy}}.
\end{split}
\end{equation}
Write $M_a=\sum_{b\sim a}(M_{ab}-M_{ba})$.  We may now take out the
factor
\begin{equation}\label{p4}
\sum_{\s_a\in Q}\s_a^{M_a}=q\cdot \one_{\{M_a\equiv 0\text{ (mod $q$)}\}}.
\end{equation}
Proceeding as above with the remaining vertices of $L$ we obtain
that
\begin{equation}
\sum_{\s\in\S}\prod_{x,y\in V}\tau_{xy}^{M_{xy}}=
q^{|V|}\cdot \one\{M_a\equiv 0\text{ (mod $q$) for all }a\in V\}.
\end{equation}
Thus
\begin{equation}
Z=q^{|V|}\exp\Big(q\sum_{e\in E}\b_e\Big)
\Pr(M_a\equiv 0\;\forall a\in V)
\end{equation}

It remains to show that the distribution of $M$ coincides with
that of $\psi$.  This is easy:  given $P'$, do the following.
First, assign for all $e\in E$ each of the $|P'_e|$ edges 
replacing $e$ a direction uniformly a random;  the number of
edges directed from $x$ to $y$ then has the Poisson distribution
with parameter $(q-1)\b_e/2$.  Next, assign each directed edge
a value $1,\dotsc,q-1$ uniformly at random;  the number of 
edges directed from $x$ to $y$ with value $r$ then has the Poisson
distribution with parameter $\tilde\b_e$.  The 
corresponding element of $F_q(P')$ is uniformly chosen
given the edge numbers and directions, and since the probability
of obtaining a $q$-flow does not depend on the choice of directions,
we are done.

To obtain the full result in the theorem, repeat the above steps
with the numerator of $\el\s^{\ul r}\er$.  The quantity $M_a$
in~\eqref{p4} must then be replaced by $M_a+r_a$, but the rest
of the calculation is as before.  It follows that
\begin{equation}\label{p5}
\el\s^{\ul r}\er=
\frac{q^{|V|}\exp\Big(q\sum_{e\in E}\b_e\Big)\PP(\psi+\ul r\equiv 0)}
{q^{|V|}\exp\Big(q\sum_{e\in E}\b_e\Big)\PP(\psi\equiv 0)}
=\frac{\PP(\psi+\ul r\equiv 0)}{\PP(\psi\equiv 0)}
\end{equation}
\end{proof}

It is straightforward to extend Theorem~\ref{potts_flow_thm} 
to an analogous representation for the space--time model, and we
sketch this here.  First, by conditioning on the set $D$, 
one obtains (as in~\eqref{ihp3}) a discrete graph $G(D)=(V(D),E(D))$.
By applying the formulas in the numerator and denominator
of~\eqref{p5} on the graph $G(D)$, one obtains a representation 
of the form~\eqref{rcr_step1_eq}.  One may then repeat the procedure
in the proof of Theorem~\ref{rcr_thm} to obtain a formula
in terms of weighted labellings;  these labellings are defined
as follows. 

Let $\L=(K,F)$ and $\b$ be as in Chapter~\ref{qim_ch}.
Fix an arbitrary ordering of the vertices $V$ of $L$.
Let $B\subseteq F$ be a Poisson process with rate
$(q-1)\l$.  We assign directions to the elements
of $B$ by letting a bridge between $(u,t)\in K$ 
and $(v,t)\in K$ be directed from $u$ to $v$ if $u$ comes before
$v$ in the ordering of $V$.  We then assign to each element of $B$
a weight from $\{1,\dotsc,q-1\}$ uniformly at random, these choices
being independent.

Let $A\subseteq K^\circ$ be a finite set (which lies in
the interior of $K$ only for convenience of exposition).
Let $\ul r=(r_x:x\in A)$ be a vector of integers, indexed
by $A$,  and let $S\se K$ denote the union of $A$ with the set
of endpoints of bridges in $B$.  
Given the above, a labelling $\psi^{\ul r}$ is a map
$K\rightarrow\ZZ/(q\ZZ)$, which is constrained to
be `valid' in that:
\begin{enumerate}
\item on each subinterval of each $K_v$,
the label is constant between elements of $S$,
\item as we move along a subinterval of $K_v$ 
($v\in V$) in the increasing $\b$
direction, the label changes at elements of $S$;  if the label
is $t$ before reaching $x\in S$, then the label just after $x$ is
\begin{itemize}
\item $t+r$ if $x$ is the endpoint of a bridge directed into $x$
and which has weight $r$,
\item $t-r$ if $x$ is the endpoint of a bridge directed out of $x$
with weight $r$,
\item $t-r_x$ if $x\in A$,
\end{itemize}
\item as one moves towards an endpoint of an interval $I^v_k\neq\SS$
(in either direction) the label converges to $0$.
\end{enumerate}
As for the random-parity representation of the space--time Ising model,
these conditions do not uniquely define $\psi^{\ul r}$ if there is 
a $v\in V$ such that $K_v=\SS$.  If this is the case, the label
at $0$ is chosen uniformly at random for each such $v$, these choices
being independent.

A valid labelling is given the weight 
\[
\partial\psi^{\ul r}:=\exp(q\d(\cL_0(\psi^{\ul r}))),
\]
where $\cL_0(\psi^{\ul r})$ is the set labelled $0$ in $\psi^{\ul r}$.
In the following, $\ul r=0$ denotes the 
vector which takes the value $0$ at all $x\in A$;
we let $E(\cdot)$ denote the expectation over $B$
as well as the weights assigned to
the elements of $B$, and the randomization which takes place when
there are several valid labellings. 
\begin{theorem}\label{st_potts_flow_thm} 
In the space--time Potts model,
\[
\el\s^{\ul r}\er=\frac{E(\partial\psi^{\ul r})}
{E(\partial\psi^0)}.
\]
\end{theorem}

The usefulness of Theorems~\ref{potts_flow_thm} 
and~\ref{st_potts_flow_thm} 
when $q\geq3$ is questionable.  Mod $q$ flows with $q\geq3$
are considerable more complicated than mod $2$ flows, and there
does not seem to be a useful switching lemma (along the lines
of Theorem~\ref{sl} or its discrete version~\cite{abf}) 
for general~$q$.

\backmatter

\appendix

\chapter{The Skorokhod metric and tightness}
\label{skor_app}

In this appendix we define carefully the Skorokhod metric on
$\Om$ and show that the sequence $\phi^b_n$ of random-cluster
measures in Section~\ref{rc_wl_sec} is tight, 
proving Lemma~\ref{tight_lem}.  We will rely partly on the
notation and results in~\cite[Chapter~3]{ethier_kurtz};  see 
also~\cite[Appendix~1]{lindvall}.

A function $f:\RR\rightarrow\RR$ is called 
\emph{c\`adl\`ag}
if it is right-continuous and has left limits.  We let
$\cD^0_\ZZ(\RR)$ denote the set of increasing 
c\`adl\`ag step functions
on $\RR$ with values in $\ZZ$, and which take the value
$0$ at $0$.  It is straightforward to modify the definitions 
and results of~\cite[Chapter~3]{ethier_kurtz},
which concern c\`adl\`ag functions on $[0,\oo)$ with values in
some metric space $E$, to apply to the set $\cD^0_\ZZ(\RR)$.
Specifically, we define the Skorokhod metric on $\cD^0_\ZZ(\RR)$
as follows.   Let $U$ denote the set of
strictly increasing bijections $u:\RR\rightarrow\RR$
which are Lipschitz continuous and for which the quantity
\begin{equation}
\a(u):=\sup_{t>s}\;\log\Big| \frac{u(t)-u(s)}{t-s}\Big|
\end{equation}
is finite.  For $a,b\in\ZZ$ let $r(a,b)=\d_{a,b}$, and note
that $r$ is a metric on $\ZZ$.  The Skorokhod metric on
$\cD^0_\ZZ(\RR)$ is by definition given by
\begin{equation}
d'(f,g)=\inf_{u\in U}
\Big[\a(u)\wedge \int_{-y}^y e^{-|y|}d'(f,g,u,y)\,dy\Big],
\end{equation}
where 
\begin{equation}
d'(f,g,u,y)=\sup_{t\in\RR}
r(f((t\wedge y)\vee -y),g((u(t)\wedge y)\vee -y)).
\end{equation}
It may be checked, as in~\cite[pp.~117]{ethier_kurtz}, that
$d'$ is indeed a metric, and that the metric space $(\cD^0_\ZZ(\RR),d')$
is complete and separable.  

Recall that we are given a countable graph $\LL=(\VV,\EE)$.
Let $\TT$ denote the countable set 
\[
\TT=(\VV\times\{\mathrm{d}\})\cup(\VV\times\{\mathrm{g}\})\cup\EE,
\]
and let $\upsilon:\TT\rightarrow\{1,2,\dotsc\}$ denote an arbitrary
bijection.  Then we formally define the set $\Om$ to be the product 
space $\Om=\cD^0_\ZZ(\RR)^\TT$.  For $\om\in\Om$ and $x\in\TT$,
the restriction $\om_x$ of $\om$ to $x\times\RR$ (not to be
confused with the $\om_x$ of Section~\ref{rcm_si}) is to be interpreted 
as: the process of deaths on $x\times\RR$ if $x\in\VV\times\{\mathrm{d}\}$,
or the process of ghost-bonds on $x\times\RR$ if
$x\in\VV\times\{\mathrm{g}\}$,
or the process of bridges on $x\times\RR$ if $x\in\EE$.  In this section
we do \emph{not} overlook events of probability zero, that is 
Remark~\ref{rem-as} does not apply.
\begin{definition}
We define the Skorokhod metric $d$ on $\Om$ by
\[
d(\om,\om')=\sum_{x\in\TT}e^{-\upsilon(x)}d'(\om_x,\om_x').
\]
\end{definition}
Note that the sum is absolutely convergent since $d'$ is bounded, and 
in fact also $d$ is bounded.  It is straightforward to check that $d$ is
indeed a metric on $\Om$, and (using the dominated convergence theorem)
that $(\Om,d)$ is a complete metric space.  
It is also separable, hence Polish.  
The $\s$-algebra $\cF$ on $\Om$ generated by
$d$ agrees with that generated by all the coordinate functions
$\pi_{x,t}:\om\mapsto\om_x(t)$ for $x\in\TT$ and $t\in\RR$, 
see~\cite[Proposition~3.7.1]{ethier_kurtz}.  The fact that all
finite tuples of such coordinate functions forms a convergence
determining class (a fact used in Theorem~\ref{conv_lem})
follows as in~\cite[Theorem~3.7.8]{ethier_kurtz}. 

In order to establish tightness of the sequence $\phi^{b_n}_n$ we must
find compact sets in $\Om$.  Since $(\Om,d)$ is a metric space,
compactness is equivalent to sequential compactness.  If for each
$x\in\TT$, the set $A_x$ is (sequentially)
compact in $(\cD^0_\ZZ(\RR),d')$, 
then by a straightforward diagonal argument the set 
$A=\bigotimes_{x\in\TT}A_x$ is a compact subset of $(\Om,d)$.

\begin{proof}[Proof of Lemma~\ref{tight_lem}]
As a witness for the tightness of $\{\phi^{b_n}_n:n\geq 1\}$ we  will
use the product $A$ of the following compact sets $A_x$.
For each $x\in\TT$, let $\xi_x:[0,\oo)\rightarrow(0,\oo)$ be a strictly
positive function, to be specified later.  Let $A_x$ be the
set of $\om\in\Om$ such that for all $t>0$, all jumps of $\om_x$
in the interval $[-t,t]$ are separated from each other 
by at least $\xi_x(t)$.  It follows from the characterization 
in~\cite[Theorem~3.6.3]{ethier_kurtz} that $A_x$ is compact
(alternatively, it is not hard to deduce the sequential compactness
of $A_x$ using a diagonal argument).

It remains to show that we can choose the functions $\xi_x$
so as to get a uniform lower bound on $\phi^{b_n}_n(A)$ which is
arbitrarily close to 1.  We can use stochastic domination,
Corollary~\ref{comp_perc}, to reduce this to checking the tightness
of \emph{a single percolation measure}, as follows.  If 
$x\in\VV\times\{\mathrm{d}\}$ 
then the event $A_x$ is increasing, otherwise it is
decreasing.  Thus $A=\bigcap_{x\in\TT}A_x=A^+\cap A^-$ 
where
\[
A^+=\bigcap_{x\in\VV\times\{\mathrm{d}\}}A_x\qquad\text{and}\qquad
A^-=\bigcap_{x\in(\VV\times\{\mathrm{g}\})\cup\EE}A_x
\]
are increasing and decreasing events, respectively.  
We have that
\begin{equation}\label{ek1}
\phi^{b_n}_n(A)\geq\phi^{b_n}_n(A^+)+\phi^{b_n}_n(A^-)-1.
\end{equation}
The events $A^+,A^-$ are not local events, but by writing
them as decreasing limits of local events it is easy to
justify the following application of Corollary~\ref{comp_perc} 
to~\eqref{ek1}.  For suitable choices of the parameters
$\l_i,\d_i,\g_i$ ($i=1,2$) which are multiples of the original parameters
$\l,\d,\g$ we have that
\begin{equation}\label{ek2}
\phi^{b_n}_n(A)\geq\mu_{\l_1,\d_1,\g_1}(A^+)+\mu_{\l_2,\d_2,\g_2}(A^-)-1.
\end{equation}
Clearly, any lower bound on the right-hand-side 
of~\eqref{ek2} is a uniform lower bound on the $\phi^{b_n}_n(A)$.

Let us focus on $A^+$, since $A^-$ is similar.  Suppose we can,
for any $\eps>0$, choose $\xi_x$ so that 
\[
\mu_{\l_1,\d_1,\g_1}(A_x)\geq e^{-\eps/\upsilon(x)^2}.
\]
Then, since the $A_x$ are independent under $\mu_{\l_1,\d_1,\g_1}$,
we will have that
\[
\mu_{\l_1,\d_1,\g_1}(A^+)\geq\exp\Big(-\eps\frac{\pi^2}{6}\Big),
\]
which is enough.  The event $A_x$ concerns only the process
$D$ of deaths on $x\times\RR$.   We may replace $\d_1$ by a constant
upper bound.  By adjusting parameters it follows that we are done
if we prove the following:  for any $\eps>0$ we have that
\begin{equation}\label{ek3}
P(N\in A_x)\geq 1-\eps,
\end{equation}
where $P$ is the measure governing the Poisson process $N$ of rate
1 on $\RR$.  The proof of~\eqref{ek3} is a straightforward exercise
on Poisson processes, but we include it for completeness.

For $I\subseteq\RR$ and $a\in\RR$ we  write $aI=\{at:t\in I\}$.
Define $I_1^+=I_1^-=[-1,1]$ and for $k\geq2$ let
$I^+_k$ be the closed interval of length $1/k$ with left endpoint
$1+1/2+1/3+\dotsb+1/(k-1)$;  let $I^-_k=-I^+_k$.  Since the series
$\sum\tfrac{1}{k}$ diverges, the $I^\pm_k$ ($k\geq1$) cover $\RR$.
Next let $J^+_k$ ($k\geq1$) be the closed interval whose left and right
endpoints are at the midpoints of $I^+_k$ and $I^+_{k+1}$ respectively;
let $J^-_k=-J^+_k$.  Note that 
$|J^\pm_k|=(|I^\pm_k|+|I^\pm_{k+1}|)/2\geq \tfrac{1}{k+1}$.
Let $\eps>0$ and let $A'$ be the event that each $\eps I^\pm_k$ and 
each $\eps J^\pm_k$ ($k\geq1$) contains at most one element of $N$.

We claim that $A'\subseteq A_x$ for $\xi_x(t)=\eps e^{-t/\eps}/4$.
Suppose $A'$ happens and $s\in N$.  We may assume $s\in\eps I^+_k$ with
$k\geq2$ (the other cases are similar).  Then $s$ also lies in either
$\eps J^+_{k-1}$ or $\eps J^+_k$.  
Hence the closest possible other point of
$N$ is a distance at least $\tfrac{\eps}{2(k+1)}$ from $s$.  Let
$t>0$ and suppose $s\in N\cap[0,t]$.  Let $k$ be maximal with 
$I^+_k\cap[0,t]\neq\varnothing$.  Then
\[
t\geq\eps\sum_{i=1}^{k-1}\frac{1}{i}\geq\eps \log k,
\]
and the closest point to $s$ in $N$ is a distance
at least 
\[
\frac{\eps}{2(k+1)}\geq\frac{\eps}{2(e^{t/\eps}+1)}\geq
\frac{\eps}{4}e^{-t/\eps}.
\]
Similarly if $s<0$.  Hence $A'\subseteq A_x$ as claimed.

It is well-known that there is an absolute constant $C$ such 
that for $\eta>0$ small and $I$ a fixed interval
of length at most $\eta$, we have that $P(|N\cap I|\geq 2)\leq C\eta^2$.
Clearly we have 
\begin{equation}
\begin{split}
P(N\in A')&\geq 1-2\sum_{k\geq 2} P(|N\cap \eps I^+_k|\geq 2)-\\
&\qquad-2\sum_{k\geq 1} P(|N\cap \eps J^+_k|\geq 2) 
 - P(|N\cap \eps I_1|\geq 2)\\
&\geq 1-\eps^2 C\cdot 2\pi^2/3.
\end{split}
\end{equation}
This proves the result.
\end{proof}

\chapter{Proof of Proposition~\ref{cond_meas_rcm}}
\label{pfs_app}

\begin{proof}[Proof of Lemma~\ref{cond_meas_rcm}]
This is essentially straightforward, but notationally intricate.
We write $(\eta,\om,\tau)_{\L,\D}$ for the
configuration which equals $\eta$ inside the smallest set $\L$, equals $\om$ in
the intermediate region $\D\setminus\L$, and equals $\tau$ outside $\D$.  
For readability, let us write $k_\L(\cdot;\tau)$ in place of
$k_\L^\tau(\cdot)$ in what follows. 

Let $A'\in\cF_{\D\setminus\L}$.  Then
\begin{multline}\label{abc0}
  \phi^\tau_\D(\one_{A'}(\cdot)\phi^{(\cdot,\tau)_\D}_\L(A))=
  \iint\one_{A'}(\om)\one_A((\eta,\om,\tau)_{\L,\D})
  \:d\phi^{(\om,\tau)_\D}_\L(\eta)d\phi^\tau_\D(\om)
  \\=\iint\one_{A\cap {A'}}((\eta,\om,\tau)_{\L,\D})
  \frac{q^{k_\L(\eta;(\om,\tau)_\D)}}{Z^{(\om,\tau)_\D}_\L}
  \frac{q^{k_\D(\om;\tau)}}{Z^{\tau}_\D}\:d\mu(\eta)d\mu(\om).
\end{multline}
Note that if $(\a,\b)_\L\in\Om$ then
\begin{equation}\label{cpt_count_eq}
  k_\D((\a,\b)_\L;\tau)=k_\L(\alpha;(\beta,\tau)_\D)+\tilde k_\D(\beta;\tau),
\end{equation}
where $\tilde k_\D$ counts the number of components 
in $\D$ which do not intersect $\L$.  Let $\a,\b$ be independent
with law $\mu$;  then $\om$ has the law of $(\a,\b)_\L$.
Use~\eqref{cpt_count_eq} on each power of $q$ in~\eqref{abc0}
to see that
\begin{align*}
\phi^\tau_\D\big(&\one_{A'}(\cdot)\phi^{(\cdot,\tau)}_\L(A)\big)\\
&=\iiint\one_{A\cap {A'}}((\eta,\beta,\tau)_{\L,\D})
\frac{q^{k_\L(\alpha;(\beta,\tau)_\D)}q^{k_\D((\eta,\beta)_\L;\tau)}}
     {Z_\L^{(\beta,\tau)_\D}Z_\D^\tau}
\:d\mu(\eta)d\mu(\alpha)d\mu(\beta)\\ 
&=\int\one_{A\cap {A'}}((\om',\tau)_\D)
\frac{q^{k_\D(\om';\tau)}}{Z^\tau_\D}
\bigg(\int\frac{q^{k_\D(\alpha;(\om',\tau))}}{Z^{(\om',\tau)}_\L}
\:d\mu(\alpha)\bigg)\:d\mu(\om')\\
&=\phi_\D^\tau(A\cap {A'}),
\end{align*}
where $\om'=(\eta,\beta)_\L$.  This proves the claim.  
\end{proof}

\bibliographystyle{amsplain}
\bibliography{thesis}

\end{document}